\documentclass[a4paper, 12pt]{report}
\usepackage
[
        a4paper,
        left=2.5cm,
        right=1.0cm,
        top=2.0cm,
        bottom=2.0cm,
]
{geometry}

\usepackage[cp1251]{inputenc}
\usepackage[english]{babel}
\usepackage{amsthm}
\usepackage{cite}
\usepackage{hyphenat}
\usepackage{amsmath,amsfonts,amssymb,amscd,amsthm,xspace,latexsym}
\usepackage{tikz-cd}
\usepackage{amsmath}
\usepackage{amsfonts}
\usepackage{amssymb}
\usepackage[colorlinks = true, urlcolor = black, linkcolor = black, citecolor = black]{hyperref}
\usepackage[]{authblk}
\usepackage{epsfig,graphicx}
\usepackage{indentfirst}
\usepackage{setspace}
\usepackage{fancyhdr}

\usepackage{titlesec}
\titleformat{\chapter}[display]
{\normalfont\Huge\bfseries\centering}{\chaptertitlename\ \thechapter}{10pt}{\Huge}

\titlespacing*{\chapter}
  {0pt}{2\baselineskip}{3\baselineskip}

\newcommand\Proof{\noindent{\it Proof.\/ }}
\newcommand\Endproof{\hfill$\Box$\medbreak}

\DeclareMathOperator{\rank}{rank}
\DeclareMathOperator{\Frac}{Frac}
\DeclareMathOperator{\Trdeg}{Trdeg}
\DeclareMathOperator{\tr}{tr}
\DeclareMathOperator{\Spec}{Spec}
\DeclareMathOperator{\Aut}{Aut}
\DeclareMathOperator{\NAut}{NAut}
\DeclareMathOperator{\End}{End}
\DeclareMathOperator{\Char}{char}
\DeclareMathOperator{\Ch}{char}
\DeclareMathOperator{\Det}{det}

\DeclareMathOperator{\Ker}{Ker}

\DeclareMathOperator{\Tr}{tr}
\DeclareMathOperator{\Fr}{Fr}

\DeclareMathOperator{\Deg}{deg}

\DeclareMathOperator{\Sp}{Sp}
\DeclareMathOperator{\Id}{Id}
\DeclareMathOperator{\Ht}{ht}

\DeclareMathOperator{\J}{J}
\DeclareMathOperator{\SL}{SL}
\DeclareMathOperator{\GL}{GL}
\DeclareMathOperator{\diag}{diag}

\DeclareMathOperator{\gr}{gr}
\DeclareMathOperator{\Ind}{Ind}
\DeclareMathOperator{\Pic}{Pic}

\DeclareMathOperator{\Sympl}{Sympl}
\DeclareMathOperator{\supp}{supp}
\newcommand{\TAut}{\operatorname{TAut}}

\def\JC{\operatorname{JC}}

\def\ch{\operatorname{Char}}
\def\Ker{\operatorname{Ker}}

\def\Id{\operatorname{Id}}
\def\Ad{\operatorname{Ad}}
\def\End{\operatorname{End}}
\def\Nilp{\operatorname{Solv}}
\def\Aut{\operatorname{Aut}}
\def\PI{\operatorname{PI}}
\def\Var{\operatorname{Var}}
\def\goth{\mathfrak}
\newcommand\C{{\mathbb C}}

\newcommand\A{{\mathbb A}}
\renewcommand\k{{\bf k}}
\newcommand\spec{{\mathsf{Spec}\,}}

\makeatletter
\newsavebox{\@brx}
\newcommand{\llangle}[1][]{\savebox{\@brx}{\(\m@th{#1\langle}\)}%
  \mathopen{\copy\@brx\kern-0.5\wd\@brx\usebox{\@brx}}}
\newcommand{\rrangle}[1][]{\savebox{\@brx}{\(\m@th{#1\rangle}\)}%
  \mathclose{\copy\@brx\kern-0.5\wd\@brx\usebox{\@brx}}}
\makeatother

\newtheorem{thm}{Theorem}[section]
\newtheorem{theorem}[thm]{Theorem}
\newtheorem{lem}[thm]{Lemma}
\newtheorem{lemma}[thm]{Lemma}

\newtheorem{prop}[thm]{Proposition}
\newtheorem{proposition}[thm]{Proposition}
\newtheorem{cor}[thm]{Corollary}
\newtheorem{corollary}[thm]{Corollary}
\newtheorem{conj}[thm]{Conjecture}
\newtheorem{Def}[thm]{Definition}
\newtheorem{definition}[thm]{Definition}
\newtheorem{remark}[thm]{Remark}
\newtheorem{rem}[thm]{Remark}

\newtheorem{problem}[thm]{Problem}
\newtheorem{example}[thm]{Example}

\newtheorem{defn}{Definition}[thm]
\newtheorem{statement}{Proposition}[thm]
\begin{document}
\pagestyle{empty}

\begin{titlepage}

\vspace{3.0cm}

\center
{\bfseries \LARGE  Polynomial Automorphisms, Quantization and Jacobian Conjecture Related Problems}\\[1.0cm]

\vspace{1.5cm}

\begin{flushleft}
Alexei Belov-Kanel (SZU)\\
Andrey Elishev (MIPT)\\
Farrokh Razavinia(MIPT)\\
Jie-Tai Yu (SZU)\\
Wenchao Zhang (BIU)
\end{flushleft}

College of Mathematics and Statistics, Shenzhen University, Shenzhen, 518061, China

Laboratory of Advanced Combinatorics and Network Applications, Moscow Institute of Physics and
Technology, Dolgoprudny, Moscow Region, 141700, Russia
College of Mathematics and Statistics, Shenzhen University, Shenzhen, 518061, China

Mathematics Department, Bar-Ilan University, Ramat-Gan, 52900, Israel

\end{titlepage}

\pagestyle{fancy} \pagenumbering{arabic} \setcounter{page}{2}

\fancyhead{}  
\rhead{\thepage}  
\lhead{}  

\pagestyle{fancy}  

\newpage

\fancypagestyle{plain}{}

\tableofcontents \fontsize{12}{12pt}\selectfont
\renewcommand{\to}{\mapsto}

\sloppy \setstretch{1.3}

\chapter*{Preface}
\addcontentsline{toc}{chapter}{Preface}

The purpose of this review is the collection and systematization of results concerning the quantization approach to the
Jacobian conjecture.

The Jacobian conjecture of O.-H. Keller remains, as of the writing of this text, an open and apparently unassailable
problem. Various possible approaches to the Jacobian conjecture have been explored, resulting in accumulation of a
substantial bibliography, while the development of vast parts of modern algebra and algebraic geometry were in part
stimulated by a search for an adequate framework in which the Jacobian conjecture could be investigated. This has
engendered a situation of simultaneous existence of circumstantial evidence in favor and against the positivity of this
conjecture.

One of the more established plausible approaches to the Jacobian problem concerns the study of infinite-dimensional
algebraic semigroups of polynomial endomorphisms and groups of automorphisms of associative algebras, as well as
mappings between those. The foundation for this approach was laid by I.R. Shafarevich. During the last several
decades, the theory was developed and vastly enriched by the works of Anick, Artamonov, Bass, Bergman, Dicks,
Dixmier, Lewin, Makar-Limanov, Czerniakiewicz, Shestakov, Umirbaev, Bia\l{}ynicki-Birula, Asanuma,
Kambayashi, Wright, and many others. In particular, the results of Anick, Makar-Limanov, Shestakov and
Umirbaev established a connection between the Jacobian conjecture for the commutative polynomial algebra and its
associative analogues on the one hand with combinatorial and geometric properties (stable tameness, approximation)
of the spaces of polynomial automorphisms on the other.

More recently, the stable equivalence between the Jacobian conjecture and a conjecture of Dixmier on the
endomorphisms of the Weyl algebra has been discovered by Kanel-Belov and Kontsevich and, independently, by
Tsuchimoto. The cornerstone of this rather surprising feature is a certain mapping (sometimes referred to as the
anti-quantization map) from the semigroup of Weyl algebra endomorphisms (a quantum object) to the semigroup of
endomorphisms of the corresponding Poisson algebra (the appropriate classical object). In view of that, it seems
reasonable to think there are insights to be gained by studying quantization of spaces of polynomial mappings and
properties of the corresponding quantization morphisms.

In this direction, one of the larger milestones is given by a series of conjectures of Kontsevich concerning
equivalences between polynomial symplectomorphisms, holonomic modules over algebras of differential operators,
and automorphisms of such algebras. Another rather non-trivial side of the quantization program rests upon the
interaction with universal algebra.

In this review we present some of our progress regarding quantization, Kontsevich conjecture, as well as recall some
of our recent results on the geometry of $\Ind$-scheme automorphisms, approximation by tame automorphisms
together with its symplectic version, and torus actions on free associative algebras. We also provide a review of
work of Kanel-Belov, Bokut, Rowen and Yu, which sought to connect the Jacobian problem with various problems
in universal algebra, as conceived by the brilliant late mathematician A.V. Yagzhev.

We have benefitted greatly from extensive and fruitful discussions with E. Aljadeff, I.V. Arzhantsev, V.A.
Artamonov, E.B. Vinberg, A.E. Guterman, V.L. Dolnikov, I.Yu Zhdanovskii, A.B. Zheglov, D. Kazhdan, R.N.
Karasev, I.V. Karzhemanov, V.O. Manturov, A.A. Mikhalev, S.Yu. Orevkov, E.B. Plotkin, B.I. Plotkin, A.M.
Raigorodskii, E. Rips, A.L. Semenov, N.A. Vavilov, G.B. Shabat, U. Vishne and G.I. Sharygin. It is a pleasant task
to express our utmost gratitude to our esteemed colleagues.

This work is supported by the Russian Science Foundation grant No. 17-11-01377.

\chapter{Introduction}
\label{Chapter1} \lhead{Chapter 1. \emph{Introduction}}

\section{Quantization and algebra problems}
This section provides the overview of the Jacobian conjecture together with motivation for the theory of
$\Ind$-schemes and quantization, as well as some necessary preliminaries on the proof of Bergman's centralizer
theorem. Throughout this paper, all rings are associative with multiplicative identity.

\subsection{Free algebras}
A free algebra is a noncommutative analogue of a polynomial ring since its elements may be described as
"polynomials" with non-commuting variables, while the free commutative algebra is the polynomial algebra. Let us
first give the definition of a free monoid, which is needed in our definition of free algebras \cite{rowen2008}.

\begin{definition}
	Let $X=\{x_i:i\in I\}$. A \textit{word} is a string with elements in $X$. A \textit{free associative monoid} on a set $X$, namely $X^*$, is the set of words in $X$, including the empty product to represent 1. The multiplication on $X^*$ is given by the juxtaposition of words.
\end{definition}

Next we can naturally give a definition of the free associative algebra respect to a generating set over a commutative
ring.

\begin{definition}
	Let $C$ be a commutative ring with multiplicative identity. A \textit{free associative $C$-algebra} $C\langle X\rangle$ with respect to a generating set $X=\{x_i:i\in I\}$ is the free $C$-module with base $X^*$.
\end{definition}

\begin{remark}
	This $C$-module becomes a $C$-algebra by defining a multiplication as follows: the product of two basis elements is the concatenation of the corresponding words and the product of two arbitrary $C$-module elements are thus uniquely determined. Note that $C\langle X\rangle:=\bigoplus_{w\in X^*}Cw$ and the elements of $C\langle X\rangle$ are called \textit{noncommutative polynomials} over $C$ generated by $X$.
\end{remark}

By the same token, we can also define the free associative algebra respect to a generating set $X=\{x_i:i\in I\}$ over
an arbitrary field $k$, namely $k\langle X\rangle$.

\begin{remark}
	If $C$ is an integral domain, then the product of leading monomials of two noncommutative polynomials $f$ and $g$ in $C\langle X\rangle$ is the leading monomial of $fg$. It follows that $C\langle X\rangle$ is a domain (but still noncommutative) as well.
\end{remark}

We finally remark that we will only discuss free associative $k$-algebras with respect to a generating set
$X=\{x_1,\dots,x_s\}$ (for $s\geq2$) over a field $k$ instead of a commutative ring throughout this review.

\subsection{Matrix representations of algebras}\label{matrep}

Let $A$ be a $k$-algebra, and let $K$ be a field extension of $k$. We talk about finite dimensional representations
of $A$ in this work, so when we mention a representation, we mean it is a finite-dimensional representation.

\begin{definition}
	An $n$-dimensional \textit{matrix representation} over $K$ is a $k$-homomor-\\phism $\rho: A\to M_n(K)$ to the matrix algebra over $K$.
\end{definition}

\begin{remark}
	Two representations $\rho',\rho$ are equivalent if they are conjugate, namely $\rho'=\tau\rho\tau^{-1}$ for some invertible matrix $\tau\in M_n(K)$.
\end{remark}

The representation is \textit{irreducible} if the images of $A$ generates the matrix algebra as $K$-algebra, or if the
map $A\otimes K\to M_n(K)$ is surjective. Usually, we study the case when $K=k$. With this assumption, we call
a representation $\rho:A\to M_n(k)$ is irreducible if and only if it is surjective.

\subsection{Algebra of generic matrices}
In order to use the concept of generic matrices, we need to first introduce the matrix representation of free
associative algebra \cite{Artin1999}. We have introduced matrix representations of any $k$-algebras in Section
\ref{matrep}. A matrix representation of the free associative ring $k\langle X\rangle=k\langle x_1,\dots,x_s\rangle$
over $k$ generated by a finite set $X=\{x_1,\dots,x_s\}$ of $s$ ($s\geq2$) indeterminates is given by assigning
arbitrary matrices as images of the variables. In itself, this is not very interesting. However, when one asks for
equivalence classes of irreducible representations, the other is directed to an interesting problem in invariant theory.
We will discuss this topic in the following.
\medskip

\begin{definition}
Let $n$ be a positive integer, and let $\{x_{ij}^{(\nu)}|1\leq i,j\leq n,\nu\in \mathbb{N}\}$ be independent
commuting indeterminates over $k$. Then $$X_{\nu}:=(x_{ij}^{(\nu)})\in M_n(k[x_{ij}^{(\nu)}])$$ is called an
$n\times n$ \textit{generic matrix} over $k$, and the $k$-subalgebra of $M_n(k[x_{ij}^{(\nu)}])$ generated by
the $X_{\nu}$ is called the \textit{algebra of generic matrices} and will be denoted by $k\langle
X_1,\dots,X_s\rangle$ or simply $k\{X\}$.
\end{definition}
The algebra of generic matrices is a basic object in the study of the polynomial identities and invariants of $n\times n$
matrices.

\medskip

There is a canonical homomorphism

\begin{equation}\label{eqngeneric}
\pi: k\langle x_1,\dots,x_s\rangle\to k\langle X_1,\dots,X_s\rangle
\end{equation}
from the free associative ring on variables $x_1,\dots,x_s$ to this ring.

If $u_1,\dots,u_s$ are $n\times n$ matrices with entries in a commutative $k$-algebra $R$, then we can substitute
$u_j$ for $X_j$ , and thereby obtain a homomorphism $$k\langle X_1,\dots, X_s\rangle \to M_n(R).$$

There is an important property of the homomorphism $\pi$: an element $f$ of the free associative algebra is in the
kernel of the map $\pi$, if and only if it vanishes identically on $M_n(R)$ for every commutative $k$-algebra $R$,
and this is true if and only if $f$ vanishes identically on $M_n(k)$. In addition, the (irreducible) matrix
representations of the free ring $A$ of dimension $\leq n$ correspond bijectively to the (irreducible) matrix
representations of the ring of generic matrices. This result is a core tool in our proof.

Let $u_1,\dots,u_N$ ($N=n^2$) be a basis for the matrix algebra $M_n(\bar{K})$, and let $z_1,\dots,z_N$ be
indeterminates. Then the entries of the matrix $Z=\sum z_ju_j$ are all algebraically independent. Moreover, we have
the famous Amitsur's theorem \cite{Artin1999} as follows.

\begin{theorem}[Amitsur]\label{amitthm}
The algebra $k\langle X_1,\dots,X_s\rangle$ of generic matrices is a \\domain.
\end{theorem}

\begin{proof} cf. \cite{Artin1999} Theorem V.10.4 or \cite{Zhang-mas} Theorem 3.2.
\end{proof}

\subsection{The Amitsur-Levitzki theorem}
For the free associative algebra $A=k\langle X\rangle$, the \textit{commutator} of two elements in $A$ is defined
by $[x,y]=xy-yx$. The commutator has analogues for more variables, called \textit{generalized commutators}
\cite{Artin1999} of elements $x_1,\dots,x_n$ of $A$,
\begin{equation}
S_n(x_1,x_2,\dots,x_n):=\sum (-1)^{\sigma}x_{\sigma1}\cdots x_{\sigma n},
\end{equation}
where $\sigma$ runs over the groups of all permutations. It is clear that $S_2(x,y)=[x,y]$. Note that the generalized
commutators are multilinear and alternating polynomials in the variables. Moreover, a general multilinear polynomial
in $n$ variables has the form $p(x_1,\dots,x_n)=\sum c_{\sigma}x_{\sigma1}\cdots x_{\sigma n},$ where the
coefficients $c_{\sigma}$ are elements of $k$.

There is an important and powerful result \cite{AmLev1} which is first proved by A. S. Amitsur and J. Levitzki in
1950.

\begin{theorem}[Amitsur-Levitzki]\label{amit-lev}
	Let $R$ be a commutative ring, and let $r$ be an integer. \\
	(i) If $r\geq 2n$, then $S_r(a_1,\dots,a_r)=0$ for every set $a_1,\dots,a_r$ of $n\times n$ matrices with entries in $R$. \\
	(ii) let $p(x_1,\dots,x_r)$ be nonzero multilinear polynomial. If $r<2n$, then there exist $n\times n$ matrices $a_1,\dots,a_r$ such that $p(a_1,\dots,a_r)\neq 0$. In particular, $S_r(x_1,\dots,x_r)$ is not identically zero.
\end{theorem}

The identity $S_{2n}\equiv0$ is called the \textit{standard identity} of $n\times n$ matrices. Note that
$S_2\equiv0$ is the commutative law, which holds for any $1\times 1$ matrices but not for any $n\times n$ matrices
if $n>1$.

\begin{remark}
	The Amitsur-Levitzki theorem is quite important \cite{Artin1999}. Suppose we study a representation $A\to M_n(k)$ of a $k$-algebra $A$. Let $I\subset A$ be the ideal generated by all substitutions of elements of $A$ into $S_{2n}$, and let $\tilde{A}=A/I$. The Amitsur-Levitzki theorem tells us that $S_{2n}=0$ is true in $M_d(k)$ if $d\leq n$ whereas it is not true if $d>n$. Killing $I$ has the effect of keeping the representations of dimensions $\leq n$, and cutting out all irreducible representations of higher dimension.
\end{remark}

The original proof of the Amitsur-Levitzki theorem by Amitsur and Levitzki is a direct proof, which is quite involved.
Rosset (1976) gives a short proof \cite{rosset1976} using the exterior algebra of a vector space of dimension $2n$.
This proof can be also found in \cite{Zhang-mas} Theorem 1.7. Then we obtain the following proposition.

\begin{proposition}\label{propCH}
Let $k\{X\}$ be the algebra of generic matrices.
\begin{enumerate}
\item[a)] Every minimal polynomial of $A\in k\{X\}$ is irreducible. In particular, $A$ is diagonalizable.
\item[b)] Eigenvalues of $A\in k\{X\}$ are roots of irreducible minimal polynomial of $A$, and every eigenvalue
    appears same amount of times.
\item[c)] The characteristic polynomial of $A$ is a power of minimal polynomial of $A$.
\end{enumerate}
\end{proposition}

There is an important open problem well-known in the community.

\begin{problem}
Whether for big enough $n$, every non-scalar element in the algebra of generic matrices has a minimal polynomial
which always coincides with its characteristic polynomial.
\end{problem}

This is an important open problem. For small $n$, Galois group of extension quotient field of center of algebra of
generic matrices might not be symmetry. But it still unknown for big enough $n$.

\medskip

From above Proposition \ref{propCH} c), for $n=p$,  a big enough prime, we can obtain following corollary.

\begin{corollary}\label{mincha}
Let $k\{X\}$ be the algebra of generic matrices of a big enough prime order $n:=p$. Assume $A$ is a non-scalar
element in $k\{X\}$, then the minimal polynomial of $A$ coincides with its characteristic polynomial.

\end{corollary}

\begin{proof}
Let $m(A)$ and $c(A)$ be the minimal polynomial and the characteristic polynomial of $A$ respectively. Note that
$\deg{c(A)}=n$, and $c(A)=(m(A))^k$. Because $A$ is not scalar, hence $\deg{m(A)}>1$. Since $n$ is a prime,
$k$ divides $n$. Hence $k=1$.
\end{proof}

Let us here remind the following fact:

\begin{proposition}\label{ProPcomm}
Every matrix with same eigenvectors as a matrix $A$ commutes with $A$.
\end{proposition}
\begin{proof}
	Let $A, B\in M_{n\times n}(k)$ and having $n$ eigenvectors meaning that the have $n$ linearly independent eigenvectors. And since $A$ and $B$ are $n\times n$ matrices with $n$ eigenvectors, then they both are diagonalizable and hence $A=Q^{-1}D_AQ$ and $B=P^{-1}D_BP,$ for $Q$ and $P$ are matrices whose columns are eigenvectors of $A$ and $B$ associated with the eigenvalues listed in the diagonal matrices $D_A$ and $D_B$ respectively. But $A$ and $B$ according to the hypothesis, have same eigenvectors and hence $P=Q=:S$. And hence $A=S^{-1}D_AS$ and $B=S^{-1}D_BS$ and so $AB=S^{-1}D_ASS^{-1}D_BS = S^{-1}D_AD_BS$ and in a same way we have $BA= S^{-1} D_BD_AS$ and since $D_A$ and $D_B$ are diagonal matrices, then commute and hence so do $A$ and $B$.
\end{proof}

\begin{proposition}\label{propeigen}
If $n$ is prime, $A$ is a non-scalar element of the algebra of generic matrices, then all eigenvalues of $A$ are
pairwise different.
\end{proposition}

\begin{proof}
It follows from Proposition \ref{propCH} and Corollary \ref{mincha} directly.
\end{proof}

Proposition \ref{propeigen} implies following results.

\begin{proposition}
The set of generic matrices commuting with $A$ are diagonalizable with $A$ simultaneously in the same
eigenvectors basis as $A$.
\end{proposition}

If $A$ is a non-scalar matrix, then we have following.

\begin{proposition}
$A$ is a non-scalar element of the algebra of generic matrices $k\{X\}$, then every eigenvalue of $A$ is
transcendental over $k$.
\end{proposition}

\subsection{Deformation quantization}


\subsubsection{Literature review}
In general, the ``quantization problem" can be stated as following. Given a classical physical model (Hamiltonian
system,  Lagrange system on a Riemannian manifold etc.), quantization amounts to replacing the observable functions
with operators acting on a Hilbert space, such that they satisfy some specific quantization conditions. In quantum
mechanics, this quantization condition is called the \textit{canonical commutation relation}, which is the fundamental
relation between canonical conjugate quantities. For example, the commutation relation between different
components of position and momentum can be expressed as $[\hat P_i,\hat Q_j] = i\hbar \delta_{ij}$, where $i$ is
the imaginary unit and $\delta_{ij}$ is the Kronecker delta. M. Hermann Weyl studied the Heisenberg uncertainty
principle in quantum mechanics by considering the operator ring generated by $P$ and $Q$. For any $2n$
dimensional linear space $V$, the Kronecker delta can be realized as a symplectic form $\omega $ such that
$u\otimes v-v\otimes u=\omega(u,v)$ defines a Weyl algebra $W(V)$ over $V$. In this sense, classical mechanics
corresponds to symmetric algebra, while the Weyl algebra is the "quantization" of symmetric algebra.

\medskip


In 1940s, J. E. Moyal \cite{moyal1949} conducted a more in-depth study of the Weyl quantization. Unlike Weyl,
the object he was interested in is not operators, but the classical function space: Weyl ignores the Poisson structure
of the classical function space. Instead of building a Hilbert space from a Poisson manifold and associating an
algebra of operators to it, He was only concerned with the algebra. He used the star product and Moyal bracket to
define a Poisson algebraic structure named Moyal algebra over the classical function space. Through the
investigation of the Moyal algebra, F. Bayen \cite{bayen1978} et al. raised that the quantum algebra can be
regarded as the deformation of the classical algebra if we think of $\hbar$ as the deformation parameter. In
particular, they proved that for the classical Poisson algebraic structure on the symmetric algebra over $\mathbb
R^{2n}$, the Moyal algebra is the only possible deformation in the sense of normative equivalence. That is, quantum
mechanics is the only possible ``deformation" of classical mechanics.

\medskip

We use the Poisson bracket to ``deform" the ordinary commutative product of observables in classical mechanics,
elements of our function algebra, and obtain a noncommutative product suitable for quantum mechanics. In order to
make deformation, we ask that the Moyal product is not only an asymptotic expansion, but also a real analytical
expansion. There is no a prior guarantee for this. From the Darboux theorem, the local Poisson algebra structure on
the symplectic manifold can always be deformed into the Moyal algebra. We only need to extend this local
deformation to the entire manifold after equipping a flat symplectic connection. However, for a typical Poisson
manifold, the situation is much more complicated.

\medskip

In mid 1970s, the existence of star-products for symplectic manifolds whose third cohomology group is trivial was
proved, but this restriction turned out to be merely technical. In the early 1980s the existence of star-products for
larger and larger classes of symplectic manifolds was proved, and finally it was shown that any symplectic manifold
can be ``quantized". A further generalization was achieved with \cite{fedosov1994simple} where Fedosov proved
that the results about the canonical star-product on an arbitrary symplectic manifold can be used to prove that all
regular Poisson manifolds can be quantized. However, in physics we sometimes require manifolds which have a
degenerate Poisson bracket and so are not sympletic. Therefore all the results mentioned above provided only a
partial answer to the problem of quantization.
\medskip

In 1993-1994 M. Kontsevich proposed the \textit{Formality Conjecture} which would imply the desired result. If
the Formality Conjecture could be proved, this would infer that any finite-dimensional Poisson manifold can be
canonically quantized in the sense of deformation quantization. The Formality Conjecture is proved by Kontsevich in
\cite{kontsevich2003}. Kontsevich then derived an explicit quantization formula which gives a formal definition of
the Moyal product for any Poisson manifold. However, it is not clear whether it gives the only possible deformation
quantization in the sense of canonical equivalence.

\medskip

Another direction in which research in deformation quantization has developed is strict deformation quantization in
which the parameter  is no longer a formal parameter, but a real one. In a way, the deformed algebras
$A[\![\hbar]\!]$ are identified with the original algebra $A$.

\subsubsection{Definitions and basic results}

\begin{definition}[Ring of formal power series]
Let $R$ be a commutative ring with identity, $R[\![X]\!]$ is said to be a ring of formal power series in the variable
$X$ over $R$ if and only if any element of $R[\![X]\!]$ is of form $\sum_{i\in\mathbb{N}}a_i X^i$ and satisfying
following:
\begin{align}
    &\sum_{i\in\mathbb{N}}a_i X^i+\sum_{i\in\mathbb{N}}b_i X^i=\sum_{i\in\mathbb{N}}(a_i+b_i) X^i\\
    &\sum_{i\in\mathbb{N}}a_i X^i\times\sum_{i\in\mathbb{N}}b_i X^i=\sum_{{n\in\mathbb{N}}}\left(\sum_{k=0}^n a_kb_{n-k}\right)X^n\label{cauchyprod}
\end{align}
\end{definition}

The above product \ref{cauchyprod} of coefficients is called the Cauchy product of the two sequences of
coefficients, and is a sort of discrete convolution. Note that the zero element and the multiplicative identity of the ring
of formal power series are the same as ring $R$'s.

\begin{remark}
The series $A=\sum_{n\in\mathbb{N}}a_n X^n\in R[\![X]\!]$ is invertible if and only if its constant coefficient
$a_0$ is invertible in $R$. The inverse series of an invertible series $A$ is $B=\sum_{n\in\mathbb{N}}b_n X^n\in
R[\![X]\!]$ with:
\begin{align*}
    &b_0=\frac{1}{a_0}\\
    &b_n=-\frac{1}{a_0}\sum_{i=1}^{n}a_ib_{n-i},\quad n\geq1
\end{align*}
An important example is the geometric series
$$(1-X)^{-1}=\sum_{n=0}^\infty X^n.$$
If $R$ is a field, then a series is invertible if and only if the constant term is non-zero.
\end{remark}

\begin{definition}
A \textit{Lie algebra} is a vector space with a skew-symmetric bilinear operation $(f, g)\to [f, g]$ satisfying the
Jacobi identity $$[[f, g], h] + [[g, h], f] + [[h, f], g] = 0$$
\end{definition}
\begin{definition}
A \textit{Poisson algebra} is a vector space equipped with a commutative associative algebra structure $(f,g)\to fg$
and a Lie algebra structure $(f,g)\to \{f,g\}$  satisfying the \textit{Leibniz rule}
$$\{fg,h\}=f\{g,h\}+\{f,h\}g$$
\end{definition}

\begin{definition}
A \textit{Poisson manifold} is a manifold $M$ whose function space $C^{\infty}(M)$ is a Poisson algebra with the
pointwise multiplication as commutative product.
\end{definition}

Let $k$ be an arbitrary field, and $A$ be a unitary $k$-algebra. Denote by $k[\![\hbar]\!]$ the \textit{ring of formal
power series} in an indeterminate $\hbar$, and by $A[\![\hbar]\!]$ the $k[\![\hbar]\!]$-module of formal power
series with coefficients in $A$.

\begin{definition}
A \textit{formal deformation} or \textit{star product} of the algebra $A$ is an associative, $\hbar$-adic continuous,
$k[\![\hbar]\!]$ bilinear product
\begin{equation*}
    \star: A[\![\hbar]\!]\times A[\![\hbar]\!]\to A[\![\hbar]\!]
\end{equation*}
satisfying the following rule on $A$:
\begin{equation}
    f\star g=\sum_{n=0}^\infty B_n(f,g)\hbar^n=fg+\sum_{n=1}^\infty B_n(f,g)\hbar^n \quad \forall f,g\in A
\end{equation}
where $B_n:A\times A\to A$ are bilinear operators.
\end{definition}

\begin{remark}
We usually want the bilinear operators $B_n$ to be bidifferential operators, i.e. bilinear maps which are differential
operators with respect to each argument.
\end{remark}

\begin{remark}The formal deformation extends $k[\![\hbar]\!]$-linearity in $A[\![\hbar]\!]$ with respect to: $$\left(\sum_{k=0}^{\infty}f_k\hbar^k\right)\star \left(\sum_{m=0}^{\infty}g_m\hbar^m\right)=\sum_{n=0}^{\infty}\left(\sum_{m+k+r=n}^{\infty}B_r(f_k,g_m)\right)\hbar^n.$$
\end{remark}

\medskip

There is a natural gauge group acting on star-products. This group consists of automorphisms of $A[\![\hbar]\!]$
considered as an $k[\![\hbar]\!]$-module, of the following form:
\begin{align*}
    &f\mapsto f+\sum_{n=0}^{\infty} D_n(f)\hbar^n, \quad \forall f\in A\subset A[\![\hbar]\!] \\
    &\sum_{n=0}^{\infty}f_n\hbar^n\mapsto \sum_{n=0}^{\infty}f_n\hbar^n+\sum_{n=0}^{\infty}\sum_{m=1}^{\infty}D_m(f_m)\hbar^{n+m},\quad \forall \sum_{n=0}^{\infty}f_n\hbar^n\in A[\![\hbar]\!]
\end{align*}
where $D_i: A\to A$ are differential linear operators.

\begin{definition}
$D(\hbar)$ as defined above is called a \textit{gauge transformation} in $A$. The set of such $D(\hbar)$ is naturally
a group.
\end{definition}

If $D(\hbar)=1+\sum_{m=1}^{\infty}D_m\hbar^m$ is such an automorphism, then it defines an equivalence and
acts on the set of star products as
\begin{equation}
    \star\mapsto\star', f(\hbar)\star'g(\hbar):=D(\hbar)(D(\hbar)^{-1}(f(\hbar))\star D(\hbar)^{-1}(g(\hbar))), \forall f(\hbar),g(\hbar)\in A[\![\hbar]\!]
\end{equation}



Each associative formal deformation $\star$ of the multiplication of $A$ admits a unit element $1_{\star}$.
Moreover, such an associative formal deformation $\star$ is always equivalent to another formal deformation
$\star'$ with $1_{\star'}=1_A$, where $1_A$ is the unit element of $A$. We are interested in star products up to
gauge equivalence.

\medskip

The following lemma gives a Poisson structure for an associative formal deformation of the multiplication of an
associative and commutative $k$-algebra $A$.

\begin{lemma}[\cite{keller2003}, lemma 1.1]
Let $\star$ be an associative formal deformation of the multiplication of an associative and commutative $k$-algebra
$A$. For $f,g\in A$, put $$\{f,g\}:=B_1(f,g)-B_1(g,f).$$ Then the map $\{,\}$ is a Poisson bracket on $A$, i.e., a
$k$-linear map such that the bracket is a Lie bracket and satisfies Leibniz rule. In addition, the bracket is dependent
only on the equivalence class of $\star$.
\end{lemma}

\begin{proof}
For simplicity, we write $[f,g]_{\star}$ the commutator of the star product $f\star g-g\star f$ for short. The map
\begin{equation}\label{starLie}
    (f,g)\mapsto \frac{1}{\hbar}[f,g]_\star
\end{equation}

clearly defines a Lie bracket on $A[\![\hbar]\!]$. The bracket $\{,\}$ equals the reduction modulo $\hbar$ of this
Lie bracket, i.e. it satisfies

\begin{equation}
\frac{1}{\hbar}[f,g]_\star\equiv \{f,g\}\mod \hbar A[\![\hbar]\!].
\end{equation}

We may write it in the another form as follows
\begin{equation}
    \{f,g\}:= \left.\frac{[f,g]_\star}{\hbar}\right|_{\hbar=0}=B_1(f,g)-B_1(g,f)
\end{equation}
Therefore the bracket $\{,\}$ is still a Lie bracket, and it also satisfies the Leibniz rule because the Lie bracket
defined in \ref{starLie} obeys the rule by associativity of the star product.

Suppose $D(\hbar)$ is an automorphism which yields equivalence of $\star$ and $\star'$, then we have
$$B_1(f,g)+D_1(fg)=B_1'(f,g)+D_1(f)g+fD_1(g)$$ for all $f,g\in A$. Thus the difference $B_1(f,g)-B_1'(f,g)$ is symmetric in $f,g$ and does not contribute to $\{,\}$.
\end{proof}

One can also decompose the operator $B_1$ into the sum of the symmetric part and of the anti-symmetric part:
$$B_1=B_1^+ +B_1^-, B_1^{+}(f,g)=B_1^{+}(g,f), B_1^{-}(f,g)=-B_1^{-}(g,f).$$

Then gauge automorphisms affect only the symmetric part of $B_1$, i.e. $B_1^-=(B_1')^{-}$. The symmetric part
is killed by a gauge automorphism. In this notation, we infer that $$\{f,g\}=B_1(f,g)-B_1(g,f)=2 B_1^{-}(f,g).$$

Thus, gauge equivalence classes of star products modulo $\hbar^2 A[\![\hbar]\!]$ are classified by Poisson
structures. However, it is not clear whether there exists a star product for a given Poisson structure. Moreover, we
may ask whether there exists a preferred choice of an equivalence class of star products. As we mentioned before,
Maxim Kontsevich \cite{kontsevich2003} showed that there is a canonical construction of an equivalence class of
star products for any Poisson manifold.

\subsubsection{Formal deformation quantization}

In this section, we may assume that $A$ is the algebra of smooth functions on a Poisson manifold $M$.

\begin{definition}
A \textit{deformation quantization} of a Poisson manifold $M$ is a star product on $A$ such that $2 B_1^{-} =
\{,\}$.
\end{definition}

We will not reproduce Kontsevich's proof here. His proof that we will not deal with is in terms of the cohomology of
the Hochschild complex. From the following theorem given by M. Kontsevich \cite{kontsevich2003}, there is a
surjection from the equivalence classes of formal deformations of $A$ onto Poisson brackets on $A$.

\begin{theorem}[Kontsevich \cite{kontsevich2003}]
Let $M$ be a smooth manifold and $A=C^{\infty}(M)$. Then there is a natural isomorphism between equivalence
classes of deformations of the null Poisson structure on $M$ and equivalence classes of smooth deformations of the
associative algebra $A$.

In particular, any Poisson bracket on $M$ comes from a canonically defined (modulo equivalence) star product.

\end{theorem}


Moreover, Kontsevich constructs a section of map, and his construction is canonical up to equivalence for general
manifolds $M$. A later result shows that in addition to the existence of a canonical way of quantization, we can
define a universal infinite-dimensional manifold parametrizing quantizations.
\medskip

The simplest example of a deformation quantization is the Moyal product for the Poisson structure on
$\mathbb{R}^n$. This is the first known example of a non-trivial deformation of the Poisson bracket.

\begin{example}
Let $M=\mathbb{R}^n$ and consider a Poisson structure with constant coefficients
$$\alpha=\sum_{i,j}\alpha^{ij}\partial_i\wedge\partial_j, \alpha^{ij}=-\alpha^{ji}\in \mathbb{R}$$

where $\partial_i=\partial/\partial x^i$ is the partial derivative in the derivation of coordinate $x^i$, $i=1,2,\dots,n$.
In such a case, we could have
$$\{f,g\}=\sum_{i,j}\alpha^{ij}\partial_i(f)\partial_j(g).$$
The Moyal $\star$-product is then given by exponentiating this Poisson operator

\begin{align*}
    f \star g  &= e^{\hbar \alpha}(f,g)\\
    & = fg+\hbar \sum_{i,j}\alpha^{ij}\partial_i(f)\partial_j(g)+\frac{h^2}{2}\sum_{i,j,k,l}\alpha^{ij}\alpha^{kl}\partial_i\partial_k(f)\partial_j\partial_l(g)+\dots\\
    & = \sum_{n=0}^\infty \frac{\hbar^n}{n!}\sum_{i_1,\dots,i_n;j_1,\dots,j_n}\prod_{k=1}^n\alpha^{i_kj_k}\prod_{k=1}^n\partial_{i_k}(f)\prod_{k=1}^n\partial_{j_k}(g).
\end{align*}
The Moyal product is a deformation of $(M, \alpha)$ but this formula is only valid when $\alpha$ has constant
coefficients.
\end{example}

In particular,

\begin{example}
Let $M=\mathbb{R}^2$. Consider the Poisson bracket given by
$$\{f,g\}=\mu\circ\left(\frac{\partial}{\partial x_1}\wedge \frac{\partial}{\partial x_2}\right)(f\otimes g)=\frac{\partial f}{\partial x_1}\frac{\partial g}{\partial x_2}-\frac{\partial f}{\partial x_2}\frac{\partial g}{\partial x_1}$$

where $\mu$ is the multiplication of functions on $M$. Then Kontsevich's construction yields the associative formal
deformation given by
$$f\star g=\sum_{n=0}^{\infty}\frac{\partial^n f}{\partial x_1^n}\frac{\partial^n g}{\partial x_2^n}\frac{\hbar^n}{n!}.$$
\end{example}

\medskip

The explicit construction of Kontsevich's formal quantization uses combinatorics, such as quivers. We close this
section here since we do not need to construct an explicit formula of deformation quantization in our proof.


\subsection{Algebraically closed skew field}
The role of algebraically closed fields in commutative algebra is well known. There are some parallel generalizations
of the concept of an algebraically closed skew field to non-commutative skew fields have proved useful for settling
various questions in ring theory. However, there are various definitions. The diversity of definitions of algebraically
closed skew fields is based on different choices of some particular characteristic of a commutative algebraically
closed field. A most natural generalization is in the sense of solvability of arbitrary equations which was brought in
sight by Bokut \cite{bokut62alg,bokut62lie,bokut66}. In \cite{bokut66}, in particular, Bokut raises a question
whether algebraically closed skew fields exist or not. The affirmative answer to the question is given by L.
Makar-Limanov \cite{makar85skew}. His result is one of the fundamental contributions to the theory of
non-commutative algebraically closed skew fields. In \cite{cohn92skew}, P. M. Cohn outlined a wide research
program for skew fields that are algebraically closed in the various senses. Note that not every associative algebras
can be embedded into an algebraically closed one, in the sense of solvability of arbitrary equations. For example, the
``Metro-Equation" $ax-xa=1$ (cf. \cite{cohn74progress}) is never solvable in any extension of a quaternionic skew
field. In \cite{kole2000}, P. S. Kolesnikov re-prove the Makar-Limanov theorem on the existence of an
algebraically closed skew field in the sense of there being a solution for any generalized polynomial equation. He
employs a simpler argument for proving that the skew field constructed is algebraically closed.

\medskip

\subsubsection{Existence of algebraically closed skew field}

We construct a non-commutative skew field $A$ satisfying the following (cf. \cite{kole2000}):

\begin{definition}\label{algskewfied}
$A$ with center $F$ is said to be \emph{algebraically closed} if, for any $S(x)\in A * F [x]\setminus A$, there
exists an element $a\in A$ such that $S(a)=0$; here, $*$ stands for a free product.
\end{definition}

It is easy to see that if $A$ is a field, that is, $A = F$, then Definition \ref{algskewfied} checks with the usual
definition of an algebraically closed field.

Let $F$ be an algebraically closed field of characteristic $0$ and $G$ be a commutative group generated by the
elements $$p_1^{\lambda_1},q_{1}^{\mu_1},p_2^{\lambda_2},q_{2}^{\mu_2},\ldots,$$ where
$\lambda_i,\mu_i\in\mathbb{Q}$, and $p_1,q_i$ are symbols in some countable alphabet. The group is isomorphic
to a direct sum of countably many additive groups $\mathbb{Q}$ of rational numbers. Then we define the
lexicographic order on $G$ by setting $p_1\ll q_1\ll p_2\ll \cdots < 1$, where $a\ll b$ means that $a^n <b$ for all
$n>0$. Correspondingly, $p_1^{-1}\gg q_1^{-1}\gg p_2^{-1}\gg \cdots > 1$. Put $G_n=\langle
p_n^{\lambda_n},q_{n}^{\mu_n},p_{n+1}^{\lambda_{n+1}},\ldots,\rangle$ and $G_{(m)}=\langle
p_1^{\lambda_1},q_1^{\mu_1},\dots,q_{m}^{\lambda_m}\rangle$. Obvious, $G_n$ is isomorphic to $G$.

\medskip

Given $G$ and $F$, we construct a set $\mathcal{A}$ of Maltsev-Neumann series. Elements $a\in\mathcal{A}$
has the form $$a=\sum_{g\in H_a} a(g)g, \;\; H_a\in G\;\; \text{is well ordered}, a(g)\in F\setminus \{0\},$$ the set
$H_a$ is denoted by $\mathrm{supp} a$. Choose a subset $A$ of $\mathcal{A}$ so that
$$A=\{a\in\mathcal{A}|\mathrm{supp}a \subset G_{(n(a)}\}.$$ Accordingly, put $A_n=\{a\in A| \mathrm{supp}a
\subset G_{n}\}$ and $A_{(n)}=\{a\in A| \mathrm{supp}a \subset G_{(n)}\}$. The set A constructed is exactly
the universe of the desired algebraic system. For the series on $A_n$, we define ordinary addition and multiplication,
and also derivations $\left(\frac{\partial}{\partial p_n},\frac{\partial}{\partial q_n}\right)$. Derivatives of the
elements $g\in G_n$ w.r.t. $p_1$ and $q_1$ are elements of $A_n$. There are several formula related to those
derivations which is omitted here. Following \cite{makar85skew,kole2000}, multiplication $*$ on A is defined thus:
$$a,b  \in A, a*b=\sum_{i\geq0} \frac{1}{i!}\frac{\partial ^i a}{\partial q_1^i}\frac{\partial ^i b}{\partial p_1^i}.$$

The $*$ is well defined and associative (cf. \cite{kole2000}). Then the system $\langle A, +, *, ||\rangle$ is an
associative algebra with valuation. That this is a skew field follows from the fact that $a*x = 1$ has a solution in
$A$. Moreover, $A$ does not satisfy any generalized polynomial identity, i.e. for every non-trivial generalized
polynomial $S(x)\in A * F [x]\setminus A$, there exists an element $a\in A$ such that $S(a)\neq 0$ (cf.
\cite{kole2000} Lemma 1.3.).

We bring the following notion that generalizes the concept of an homogeneous polynomial in $A * F [x]$.

\begin{definition}
An \emph{homogeneous operator} over $A_n$ is $$S_n(x)=\sum_{\imath,\jmath}
f_{\imath,\jmath}x^{(i_1,j_1)}\ldots x^{(i_k,j_k)},$$ where $\imath=(i_1\ldots i_k), \jmath=(j_1\ldots j_k),
f_{\imath,\jmath}\in A_n$, and $x$ is a common element in $A_n$, if the following conditions hold:
\smallskip

(1) there exists an $m$ such that $f_{\imath,\jmath}\in A_{(m)}$ for all $\imath,\jmath$;

(2) for ant $g\in G_n$, $x\in A_n$, the following inequality holds only for finitely many summands in $S_n(x)$:
$$|f_{\imath,\jmath}x^{(i_1,j_1)}\dots x^{(i_k,j_k)}|\leq g;$$

(3) all summands have the same degree over $x$, denoted deg $S_n(x) (= k)$.
\end{definition}

In \cite{kole2000}, Kolesnikov solves $|S(x)|=g$ and then $S_{1}(x)=f_1$. In his proof, there is a modification of
Makar-Limanov's original proof which is expedient for it compensates for this loss by instilling much more simplicity
in the argument for algebraic closedness.

In \cite{kole2001}, Kolesnikov shows that every polynomial equation containing more than one homogeneous
component over such a skew field has a non-zero solution necessarily. Precisely, he obtains following proposition:

\begin{proposition}[cf. \cite{kole2001} Theorem 1.]
Let $S_i(x), i=1,\dots,n$ be homogeneous operators over $A$, where $n\geq 1$, and $T(x)$ be a homogeneous
operator such that $\Deg S_i < \Deg T$. Then the equation $\sum_{i} S(x) = T (x)$ has a solution $x \in A, x \neq
0$.
\end{proposition}

\subsubsection{Algebraically closed skew field in the sense of matrices}
Another conception of algebraic closedness is associated with the notion of singular eigenvalues of matrices. The
definitions are given in Cohn \cite{cohn92skew}.

Let $D$ be a skew field with center $k$. Denote by $M_n(D)$ a ring of all $n \times n$-matrices over $D$. A
matrix $A\in M_n(D)$ is said to be \emph{singular} if there exists a non-zero column $u \in D^n$ such that $Au =
0$. A square matrix is singular if and only if it is not invertible. The property of being singular for a matrix is
preserved under left or right multiplication by an invertible one, in particular, under elementary transformations of
columns with coefficients from the skew field on the right, and of rows - on the left.

An element $\lambda\in D$ is called a \emph{singular eigenvalue} of $A$ if $A-\lambda I$ is a singular matrix. It is
worth mentioning that singular eigenvalues of matrices are not always preserved under similarity transformations, but
central eigenvalues are invariant in this sense.

The following definition of algebraically closed skew field is due to P. Cohn \cite{cohn92skew}.
\begin{definition}
A skew field $D$ is said to be \emph{algebraically closed} in the sense of Cohn (written $AC$) if every square
matrix over $D$ has a singular eigenvalue in that skew field. $D$ is said to be \emph{fully algebraically closed}
(written $FAC$) if every matrix $A \in M_n(D)$, which is not similar to a triangular matrix over the center of $D$,
has a non-zero singular eigenvalue in $D$.
\end{definition}

The definition of $FAC$ skew field is equivalent to the following:

\begin{definition}
A skew field $D$ is fully algebraically closed if every matrix $A\in M_n(D)$ which is not nilpotent has a non-zero
singular eigenvalue in $D$.
\end{definition}

Consequently, if $A$ is similar to a triangular matrix over the center of $D$ then either it is nilpotent or has a
non-zero eigenvalue. Conversely, if A is nilpotent then it is similar to its canonical form containing only $1$ on a
secondary diagonal.

\begin{definition}
We say that $D$ is an $AC_n$ (resp., $FAC_n$) skew field if every (non-nilpotent) matrix $A \in M_m(D), m \leq
n$, has a non-zero eigenvalue in $D$.
\end{definition}

\begin{proposition}[cf. \cite{kole2001} Theorem 2.]
Let $D$ be an $FAC_n$ skew field and $a_i,b_i,c\in D, i=1,\dots,n$. Then the equation
$$L_n(x)=\sum_{i=1}^{n}a_i x b_i=c$$ has a solution $x=x_0\in D$ if $L_n(x)\equiv0$ for all $x\in D$.
\end{proposition}

\section{Automorphisms of polynomial algebras and Kontsevich Conjecture}

One of the main objects of study in the theory of polynomial mappings are given by $\Ind$-schemes, the points of
which are automorphisms of various algebras with polynomial identities. The latter are normally algebras of
commutative polynomials $ \mathbb{K} [x_1, \ldots, x_n] $ of $ n $ variables, the free algebras $ \mathbb{K}
\langle x_1, \ldots, x_n \rangle $ of $ n $ generators, some selected quotients, as well as algebras with additional
structure -- the prominent example is the polynomial algebra equipped with the Poisson bracket. This area of
research is rooted in the widely famous Jacobian Conjecture. Thanks to the relatively recent progress of A.
Belov-Kanel and M. Kontsevich \cite{K-BK1, K-BK2} and Y. Tsuchimoto \cite{Tsu2, Tsu1}), as well as in
connection with earlier studies, a significant place in the scientific program regarding the Jacobian Conjecture has
come to be occupied by questions related to the quantization of classical algebras.

Studying the geometry and topology of $\Ind$-schemes of automorphisms, development of approximation theory of
symplectomorphisms by tame symplectomorphisms, as well as the construction of a correspondence between plane
algebraic curves and holonomic modules (over the corresponding the Weyl algebra) are the basis of the approach to
solving the Conjecture of A. Belov-Kanel and M. Kontsevich on automorphisms of the Weyl algebra, built by A.
Elishev, A. Kanel-Belov and J.-T. Yu in articles \cite{K-BE4, K-BE2} (cf. also \cite{Eli-phd}).

\subsection{Jacobian Conjecture}

One of the most well-known unresolved problems in the theory of polynomials in several variables is the so-called
the \emph{Jacobian Conjecture}, formulated in 1939 by O.-H. Keller \cite{keller1939}. Let $ \mathbb{K} $ be
the main field, and for a fixed positive integer $ n $ are given $ n $ polynomials
$$ f_1 (x_1, \ldots, x_n), \ldots, f_n (x_1, \ldots, x_n) $$
of$ n $ variables $ x_1, \ldots, x_n $. Any such system of polynomials defines a unique image
\emph{endomorphism} of the algebra $ \mathbb{K} [x_1, \ldots, x_n] $
$$
F: \mathbb{K} [x_1, \ldots, x_n] \rightarrow \mathbb{K} [x_1, \ldots, x_n]
$$

$$
F \leftrightarrow (F (x_1), \ldots, F (x_n)) \equiv (f_1 (x_1, \ldots, x_n), \ldots, f_n (x_1, \ldots, x_n),
$$
the  $\mathbb{K}$-endomorphism $ F $ of polynomial algebra is determined by its action on the set of generators.
Let $ J (F) $ denote \emph{Jacobian} (the determinant of the Jacobi matrix) of the map $ F $:
\begin{equation*}
J (F) = \Det \begin{bmatrix}
\frac{\partial f_1}{\partial x_1} & \cdots & \frac{\partial f_1}{\partial x_n} \\
\vdots & \ddots & \vdots \\
\frac{\partial f_n}{\partial x_1} & \cdots & \frac{\partial f_n}{\partial x_n}
\end{bmatrix}
\end{equation*}
The Jacobian Conjecture is as follows.

\begin{conj} [The Jacobian Conjecture, $ JC_n $] \label{jc}
Let the characteristic of the base field $ \mathbb{K} $ be equal to zero. Then, if the Jacobian $ J (F) $ of the
endomorphism $ F $ is equal to a nonzero constant (that is, it belongs to the set $ \mathbb{K} ^{\times} $), then $
F $ is an automorphism.
\end{conj}

An elementary exercise is to verify the statement that automorphisms of polynomial algebra always have a nonzero
Jacobian constant. The Conjecture \ref{jc} is thus partially inverse statement of this property. It is also easy to see
that if a polynomial  endomorphism $F$ is invertible, then the inverse will also be a polynomial endomorphism.

\medskip

The Jacobian conjecture is trivial for $ n = 1 $. On the other hand, when the field $ \mathbb{K} $ has positive
characteristic, the Jacobian Conjecture formulated as Conjecture \ref{jc} is incorrect even in the case of $ n = 1 $.
Indeed, if $ \Char \mathbb{K} = p $ and $ n = 1 $, we can take $ \varphi (x) = x - x ^ p $; the Jacobian of such a
mapping is equal to unity, but it is irreversible.

\medskip

Despite the apparent simplicity of wording and context, the Jacobian Conjecture is one of the most difficult open
questions of modern algebraic geometry. This problem has become the subject of numerous studies and has greatly
contributed to the development of related fields of algebra, algebraic geometry and mathematical physics, which are
also of independent interest.

\medskip

The literature on the Jacobian Conjecture, its analogues and related problems is extensive. A detailed discussion of
the results established in the context of the Jacobian Conjecture is beyond the scope of this work; Below we give a
brief overview of some results directly related to the Jacobian Conjecture (i.e., for the algebra of polynomials in
commuting variables). Among studies of issues similar to the Jacobian Conjecture in associative algebra, it is worth
noting the work of Dicks \cite{Di} and Dicks and Levin \cite{DiLev} on an analogue of the Jacobian Conjecture for
free associative algebras, the proof by U.U. Umirbaev \cite{umirbaev1995ext} of an analogue of the Jacobian
Conjecture for the free metabelian algebra, as well as the deep and extremely non-trivial work of A.V. Yagzhev
\cite{Yag1, Yag2, Yag3, Yag4} (see also \cite{BBRY}).

\subsection{Some results related to the Jacobian Conjecture}

While the general case of the Jacobian Conjecture (or even the Jacobian Conjecture on the plane) remains, at the
time of writing this text, an open problem, various partial results are known. We mention but a few of them.

\medskip

S.-S. Wang \cite{wang1980jacobian} established the Jacobian conjecture for the case of endomorphisms defined
by polynomials of degree $ 2 $. Also, H. Bass, E.H. Connell, and D. Wright \cite{BCW} showed that the general
case of the Jacobian Conjecture would follow from the special case of the Jacobian Conjecture for the so-called
endomorphisms of homogeneous cubic type, which are defined as mappings of the form
$$
(x_1, \ldots, x_n) \mapsto (x_1 + H_1, \ldots, x_n + H_n)
$$
where the polynomials $ H_k $ are homogeneous of degree $ 3 $.

\medskip

Moreover, L.M. Dru{\.z}kowski \cite{Druz1} proved that the previous hypothesis can be weakened, by
considering as $ H_k $ only polynomials that are cubes of linear homogeneous polynomials.

\medskip

In the works of M. de Bondt and A. van den Essen \cite{dBvdE1, dBvdE2}, as well as in the work of
Dru{\.z}kowski \cite{Druz2}, it was shown that the Jacobian Conjecture was enough to be established for
endomorphisms of homogeneous cubic type with a symmetric Jacobi matrix.

\medskip

Suppose, as before, that the polynomial endomorphism $ F $ is given by the set of images of the generators:
$$
F \leftrightarrow (F (x_1), \ldots, F (x_n)) \equiv (F_1, \ldots, F_n).
$$
Then $ F $ is invertible if and only if the algebras
$$
\mathbb{K} [x_1, \ldots, x_n] \; \; \text{and} \; \; \mathbb{K} [F_1, \ldots, F_n]
$$
are isomorphic to each other. Keller's original paper \cite{keller1939} considered a rational analogue of the
presented criterion, i.e. case of isomorphism of function fields
$$
\mathbb{K} (x_1, \ldots, x_n) \; \; \text{and} \; \; \mathbb{K} (F_1, \ldots, F_n)
$$
and the invertibility following from the existence of an isomorphism is established by L.A. Campbell \cite{Campbell}.
A generalization of Keller's original result to the case when $ \mathbb{K} (x_1, \ldots,x_n) $ is a Galois extension
of the field $ \mathbb{K} (F_1, \ldots, F_n) $ (see also the works of M. Razar \cite{Razar} and D. Wright
\cite{Wright} generalizing the result mentioned).

\medskip

In addition, some efforts were aimed at testing the fulfillment of the Jacobian Conjecture for all endomorphisms
defined by polynomials of degree not higher than some fixed number. Moh \cite{Moh1, Moh2} performed a similar
test for polynomials of two variables of degree not exceeding $100$.

\medskip

Despite the existence of the results described above (as well as some other similar theorems), the general case of the
Jacobian Conjecture remains not only open, but, apparently, at the moment unassailable.

\medskip

On the other hand, there are situations in which mappings, by their geometric properties close to polynomial
endomorphisms, are nevertheless not invertible. S.Yu. Orevkov \cite{Orevkov} points to the following reformulation
of the Jacobian Conjecture, leading to a similar situation. Let $ l $ be an infinitely distant line in the complex
projective plane $ \mathbb{C} P ^ 2 $, $ U $ be its tubular neighborhood, $ f_1 $, $ f_2 $ are meromorphic
functions on $ U $, holomorphic on $ U \backslash l $ and defining a locally one-to-one mapping
$$
F: U \backslash l \rightarrow \mathbb{C} ^ 2.
$$
The Jacobian conjecture is equivalent to the statement about the injectivity of mappings of this kind. S.Yu. Orevkov
\cite{Orevkov} constructed the following example.
\begin{thm} [S.Yu. Orevkov, \cite{Orevkov}] \label{thmorevkov}
There is a smooth, non-compact complex analytic surface $ \tilde{X} $, on which there is a smooth curve $
\tilde{L}$, isomorphic to the projective line, with the self-intersection index $ + 1 $, and two functions $ f_1 $, $
f_2 $, meromorphic on $ \tilde{X} $ and holomorphic on $ \tilde{X} \backslash \tilde{L} $, such that the mapping
defined by
$$
F: \tilde{X} \backslash \tilde{L} \rightarrow \mathbb{C} ^ 2
$$
is locally one-to-one, but not injective.
\end{thm}
As noted in \cite{Orevkov}, if $ \tilde{U} $ is a tubular neighborhood of the curve $ \tilde{L} $, then the pairs $
(U, l) $ (as above) and $ (\tilde{U}, \tilde{L}) $ are diffeomorphic, which implies the existence of a smooth
immersion in a two-dimensional complex exterior space of a ball which is geometrically similar to a polynomial map
and non-invertible. Also, if the pairs $ (U, l) $ and $ (\tilde{U}, \tilde{L}) $ were biholomorphic to each other, then
from the example of Orevkov one would derive the existence of a counterexample to the Jacobian Conjecture. This
consideration allows one to conclude that suspicion in favor of the negativity of the Jacobian Conjecture are
generally warranted.
\medskip

In his classic work, D. Anick \cite{An} developed a theory of approximation of polynomial endomorphisms by tame
automorphisms (one of the main results of this paper is the proof of a symplectic analogue of Anick's main theorem).
In connection with D. Anick's theorem on approximation, the Jacobian Conjecture can be reduced to solving the
question of the invertibility of limits of sequences of tame automorphisms. Moreover, a symplectic analogue of
Anick's theorem gives a natural (albeit requiring non-trivial extensions) idea to solve the lifting problem of polynomial
symplectomorphisms to automorphisms of the Weyl algebra in order to prove the Kanel-Belov -- Kontsevich
Conjecture (often called the Kontsevich Conjecture), an overview of which we provide in the latter section of this
paper.

\medskip

The central question in the approach to the Jacobian Conjecture and to the Kontsevich Conjecture based on
approximation by sequences of tame automorphisms is the proof of polynomial nature of the resulting limit. While in
the case of the lifting of symplectomorphisms (Kontsevich Conjecture, \cite{K-BE4}) the proof of the correctness
of the construction seems to be possible (a significant role in it is played by the invertibility of the sequence limits,
which is obviously not the case for the Jacobian Conjecture), in the context of the Jacobian Conjecture there is no
clarity in the matter, and considerations following from S.Yu. Orevkov \cite{Orevkov}, indicate possible significant
obstacles.

\medskip

The Jacobian Conjecture is studied by the methods of covering groups in S.Yu. Orevkov \cite{Orevkov2,
Orevkov3} and A.G. Vitushkin \cite{Vit1, Vit2}.

\medskip

The Jacobian Conjecture is also the subject of highly non-trivial work of Vik S. Kulikov \cite{Kulik01, Kulik92}.

\medskip

A number of other difficult problems from the theory of polynomial automorphisms are closely connected with the
Jacobian Conjecture and with affine algebraic geometry. These problems are important in the general mathematical
context. For example, a special case of the classical Abyankar-Sataye Conjecture \cite{Za96, Popov} posits
isomorphisms of all embeddings of the complex affine line into three-dimensional space (in other words, it is a
conjecture about the possibilities of formally algebraic definition of the knot).

\subsection{Ind-schemes and varieties of automorphisms}

One of the essential areas of algebraic geometry, the development of which was motivated Jacobian Conjecture is
the theory of infinite-dimensional algebraic groups. The main reference is the seminal article of I.R. Shafarevich
\cite{Shafarevich}, in which he defined concepts that allowed one to study questions about some natural
infinite-dimensional groups -- for example, the group of automorphisms of an algebra of polynomials in several
variables -- using tools from algebraic geometry. In particular, Shafarevich defines \emph{infinite-dimensional
varieties} as inductive limits of directed systems of the form
$$
\lbrace X_i, f_{ij}, \; i, j \in I \rbrace
$$
where $ X_i $ are algebraic varieties (more generally, algebraic sets) over a field $ \mathbb{K} $, and the
morphisms $f_{ij} $ (defined for $ i \leq j $) are closed embeddings. The inductive limit of a system of topological
spaces carries a natural topology, and therefore the natural questions about connectivity and irreducibility arise,
which were also studied in \cite{Shafarevich}.

\medskip

Following generally accepted terminology, we will call the direct limit of systems of varieties and closed embeddings
an $\Ind$-\emph{variety}, and the corresponding limits of systems of schemes and morphisms of schemes an
$\Ind$-\emph{scheme}.

\medskip

The Jacobian Conjecture has the following elementary connection with $\Ind$-schemes. Since the algebra of
polynomials $ \mathbb{K} [x_1, \ldots, x_n] $ can be endowed with a natural $ \mathbb{Z} $ - grading in total
degree $ \Deg $, which is defined as the appropriate monoid homomorphism  by the requirement $ \Deg x_i = 1 $,
we can define the degree of endomorphism $ \varphi $: namely, if
$$
\varphi = (\varphi (x_1), \ldots, \varphi (x_n))
$$
defined by its action on algebra generators, then the degree $\Deg \varphi$ is the maximum value of $ \Deg $ on the
polynomials $ \varphi (x_1), \ldots, \varphi (x_n) $. It defines an increasing filtration
$$
\End ^{\leq N} \mathbb{K} [x_1, \ldots, x_n], \; N \geq 0
$$
on the set $ \End \mathbb{K} [x_1, \ldots, x_n] $ of endomorphisms of the polynomial algebra; points
$$ \End ^{\leq N} \mathbb{K} [x_1, \ldots, x_n] $$ are endomorphisms of degree at most $ N $. It is easy to see that the algebraic sets $$ \End ^{\leq N} \mathbb{K} [x_1, \ldots, x_n] $$ are isomorphic to affine spaces of appropriate dimension; the coordinates of the point $ \varphi $ are the coefficients of the polynomials $$ \varphi (x_1), \ldots, \varphi (x_n) $$, and for $$ \End \mathbb{K} [x_1, \ldots, x_n] $$ these coordinates are not connected by any relations.

The total degree filtration also enables endowing the sets of automorphisms with the Zariski topology as follows (see
also \cite{Shafarevich}): if $\varphi$ is a polynomial automorphism, then consider a set of polynomials
$$
(\varphi (x_1), \ldots, \varphi (x_1), \varphi ^{- 1} (x_1), \ldots, \varphi ^{- 1} (x_n))
$$
- the images of generators under the action of the automorphism and its inverse. The coefficients of these
polynomials serve as coordinates of $ \varphi $ as a point of some affine space.

Define the subsets

$$ \Aut ^{\leq N} \mathbb{K} [x_1, \ldots, x_n] = \lbrace \varphi \in \Aut \mathbb{K} [x_1, \ldots, x_n]: \Deg \varphi , \Deg \varphi ^{- 1} \leq N \rbrace $$
as sets of automorphisms such that all coefficients of polynomials in the presentation above for degrees greater than
$ n $ are zero.

The sets $ \Aut ^{\leq N} \mathbb{K} [x_1, \ldots, x_n] $ are algebraic sets; Indeed, the identities that define the
points $ \Aut ^{\leq N} $ are derived from the identity
$$
\varphi \circ \varphi ^{- 1} = \Id
$$
and, it is easy to see, are specified by polynomials.

\medskip

Now let $ \mathfrak{J} ^{\leq N} $ denote a subset of $$ \End ^{\leq N} \mathbb{K} [x_1, \ldots, x_n], $$
whose points are endomorphisms with a Jacobian equal to a nonzero constant.

Then Conjecture \ref{jc} can be clearly reformulated as follows
$$
\forall \varphi \in \mathfrak{J} ^{\leq N} \Rightarrow \varphi \in \Aut \mathbb{K} [x_1, \ldots, x_n], \; \; \forall N, \; \; \text{for} \; \; \Char \mathbb{K} = 0.
$$

\subsection{Conjectures of Dixmier and Kontsevich}

J. Dixmier \cite{Dix} in his seminal study of Weyl algebras found a connection between the Jacobian Conjecture and
the following Conjecture. Let $ W_ {n, \mathbb {K}} $ denote the \emph {$n$-th Weyl algebra} over the field $
\mathbb {K} $ defined as the quotient algebra of the free algebra $$ F_ {2n} = \mathbb {K} \langle a_1, \ldots,
a_n, b_1, \ldots, b_n \rangle $$ of $ 2n $ generators by the two-sided ideal $ I_W $, generated by polynomials
$$
a_ia_j-a_ja_i, \; \; b_ib_j-b_jb_i, \; \; b_ia_j-a_jb_i- \delta_ {ij} \; \; (1 \leq i, j \leq n)
$$
($ \delta_ {ij} $ is the Kronecker symbol). The Dixmier Conjecture states:
\begin {conj} [Dixmier Conjecture, $ DC_n $] \label{dc}
Let $ \Char \mathbb {K} = 0 $. Then $ \End W_ {n, \mathbb {K}} = \Aut W_ {n, \mathbb {K}} $.
\end {conj}
In other words,  the Dixmier conjecture asks whether every endomorphism of the Weyl algebra over a field of
characteristic zero is in fact an automorphism.

\medskip

The Dixmier conjecture for $ n $ variables, $ DC_n $, implies the Jacobian Conjecture $ JC_n $ for $ n $ variables
(see, for example, \cite{vdE}). Significant progress in recent years in the study of Conjecture \ref{jc} has been
achieved by A. Kanel-Belov (Belov) and M.L. Kontsevich \cite{K-BK2} - and independently by Y. Tsuchimoto
\cite{Tsu2} (also see \cite{Tsu1}) - in the form of the following theorem.

\begin{thm} [A.Ya. Kanel-Belov and M.L. Kontsevich \cite{K-BK2}, Y. Tsuchimoto \cite{Tsu2}] \label{dcjc}
$ JC_{2n} $ implies $ DC_n $.
\end{thm}

In particular, the Theorem \ref{dcjc} implies the stable equivalence of the Jacobian and Dixmier conjectures - i.e.
the equivalence of the conjectures $ JC_{\infty} $ and $ DC_{\infty} $, where $ JC_{\infty} $ denotes the
conjunction corresponding conjectures for all finite $ n $.

\medskip

Theorem \ref{dcjc} laid the foundation for the research into Jacobian Conjecture based on the study of the behavior
of varieties of endomorphisms and automorphisms of algebras under deformation quantization. The principal
reference in this direction is given by an article by A. Kanel-Belov and M.L. Kontsevich \cite{K-BK1}; in it several
conjectures concerning $\Ind$-varieties of automorphisms of the corresponding algebras are formulated . The main
Conjecture is called the Kontsevich Conjecture and is as follows.

\medskip

\begin{conj} [Kontsevich Conjecture, \cite{K-BK1}] \label{mainconj}
Let $ \mathbb {K} = \mathbb {C} $ be the field of complex numbers. The automorphism group $ \Aut W_ {n,
\mathbb {C}} $ of the $ n $ -th Weyl algebra over $ \mathbb {C} $ is isomorphic to the automorphism group $
\Aut P_ {n, \mathbb {C}} $ of the so-called $ n $ -th (commutative) Poisson algebra $ P_ {n, \mathbb {C}} $:
$$
\Aut W_ {n, \mathbb {C}} \simeq \Aut P_ {n, \mathbb {C}}
$$
\end{conj}
The algebra $ P_ {n, \mathbb {C}} $ is by definition the polynomial algebra
$$
\mathbb {C} [x_1, \ldots, x_n, p_1 \ldots, p_n]
$$
of $ 2n $ variables, equipped with the Poisson bracket - a bilinear operation $ \lbrace \;, \; \rbrace $, which is a Lie
bracket satisfying the Leibniz rule and acting on generators of the algebra in the following way:
$$
\lbrace x_i, x_j \rbrace = 0, \; \; \lbrace p_i, p_j \rbrace = 0, \; \; \lbrace p_i, x_j \rbrace = \delta_ {ij}.
$$
Endomorphisms of the algebra $ P_n $ are endomorphisms of the algebra of polynomials that preserve the Poisson
bracket (which we sometimes call the Poisson structure in the text). Elements of $\Aut P_ {n, \mathbb {C}} $ are
called polynomial symplectomorphisms; the choice of name is due to the existence of an (anti-) isomorphism
between the group $ \Aut P_ {n, \mathbb {C}} $ and the group of polynomial symplectomorphisms of the affine
space $ \mathbb {A} ^ {2n} $.

\medskip

Kontsevich conjecture is true for $ n = 1 $. The proof of this result is a direct description of automorphism groups $
\Aut P_ {1, \mathbb {C}} $ and $ \Aut W_ {1, \mathbb {C}} $, contained in the classical works of Jung
\cite{Jung}, Van der Kulk\cite{VdK}, Dixmier \cite{Dix} and Makar-Limanov \cite{ML1} (see also \cite{ML2}).
Namely, consider the following transformation groups: the group $ G_1 $ is a semidirect product
$$
\SL (2, \mathbb {C}) \rtimes \mathbb {C} ^ 2
$$
whose elements are called special affine transformations, and the group $ G_2 $ by definition consists of the
following, ``triangular" substitutions:
$$
(x, p) \mapsto (\lambda x + F (p), \lambda ^ {- 1} p), \; \; \lambda \in \mathbb {C} ^ {\times}, \; \; F \in \mathbb {C} [t].
$$

The automorphism group of the algebra $ P_ {1, \mathbb {C}} $ then \cite{Jung} is isomorphic to the quotient
group of the free product of the groups $ G_1 $ and $ G_2 $ by their intersection. J. Dixmier \cite{Dix} and, later,
L.G. Makar-Limanov \cite{ML1} showed that if in the description above one replaces the commuting Poisson
generators with their quantum (Weyl) analogues, one obtains a description of the group of automorphisms of the first
Weyl algebra $ W_{1,\mathbb {C}} $.

\medskip

\begin{remark} \label{rem2D}
The theorems of Jung, van der Kulk, Dixmier and Makar-Limanov also mean that all automorphisms of the
polynomial algebra of two variables and the first Weyl algebra $W_1$ are \emph{tame} (definition of the concept of
tame automorphism, which plays a significant role in this study, we provide in the sections below). Also,
Makar-Limanov \cite{ML2} and A. Czerniakiewicz \cite{Czer1,Czer2} proved that all automorphisms of the free
algebra $ \mathbb {K} \langle x, y \rangle $ are tame.

In view of these circumstances, the case of two variables is to be considered exceptional. However, the Jacobian
Conjecture is a hard open problem even in this case.
\end{remark}

\medskip

Recently A. Kanel-Belov, together with A. Elishev and J.-T. Yu, have suggested a proof of the general case of
Kontsevich conjecture (\cite{K-BE4,K-BE2}). An independent proof of a closely related result (based on a study
of the properties of holonomic $ \mathcal {D} $ - modules) was proposed by C. Dodd \cite{Dodd}.

\medskip

In contrast to the Jacobian Conjecture, which is an extremely difficult problem, in the study of the Kontsevich
Conjecture there are several conceivable approaches. First of all, in the program article \cite{K-BK1} Kanel-Belov
and Kontsevich have formulated several generalizations of the Conjecture \ref{mainconj}. In \cite{K-BK2} and
\cite{Tsu2}, which is devoted to the proof of Theorem \ref{dcjc}, the construction of homomorphisms
$$
\phi: \Aut W_ {n, \mathbb{C}} \rightarrow \Aut P_ {n, \mathbb{C}}
$$
and
$$
\phi: \End W_ {n, \mathbb{C}} \rightarrow \End P_ {n, \mathbb{C}},
$$
involved in the construction, from a counterexample to $ DC_n $, of an irreversible endomorphism with a single
Jacobian, has been presented. A straightforward strengthening of Conjecture \ref{mainconj} is the statement that the
homomorphism $ \phi $ realizes the isomorphism of the Kontsevich Conjecture. Also - namely, in Chapter 8 of the
article \cite{K-BK1}, an approach to resolving the problem of lifting of polynomial symplectomorphisms to
automorphisms of the Weyl algebra (i.e., constructing a homomorphism inverse to $\phi$) was discussed.
Conjecture 5 of the paper \cite{K-BK1}, along with Conjecture 6, which is a weaker form of Conjecture
\ref{mainconj}, make up the essential contents of the construction proposed in \cite{K-BK1}. To solve the problem
of lifting of symplectomorphisms in the sense of these conjectures, it is necessary to study the properties of $
\mathcal{D} $ - modules --  (left) modules over the Weyl algebra. The work of Dodd \cite{Dodd} is based on this
approach.

\medskip

\subsection{Approximation by tame automorphisms}

Tsuchimoto \cite{Tsu1, Tsu2}, and independently Kanel-Belov and Kontsevich \cite{K-BK2}, found a deep
connection between the Jacobian conjecture and a celebrated conjecture of Dixmier \cite{Dix} on endomorphisms
of the Weyl algebra, which is stated as in Conjecture \ref{dc}.

The correspondence between the two open problems, in the case of algebraically closed $K$, is based on the
existence of a composition-preserving map

\begin{equation*}
\End W_n(\mathbb{K})\rightarrow \End \mathbb{K}[x_1,\ldots,x_{2n}]
\end{equation*}

which is a homomorphism for the corresponding automorphism groups. Furthermore, the mappings that belong to
the image of this homomorphism preserve the canonical symplectic form on $\mathbb{A}^{2n}_{\mathbb{K}}$.
In accordance with this, Kontsevich and Kanel-Belov \cite{K-BK1} formulated several conjectures on
correspondence between automorphisms of the Weyl algebra $W_n$ and the Poisson algebra $P_n$ (which is the
polynomial algebra $\mathbb{K}[x_1,\ldots,x_{2n}]$ endowed with the canonical Poisson bracket) in
characteristic zero. In particular, there is a

\begin{conj}
The automorphism groups of the $n$-th Weyl algebra and the polynomial algebra in $2n$ variables with Poisson
structure over the rational numbers are isomorphic:
\begin{equation*}
\Aut W_n(\mathbb{Q})\simeq \Aut P_n(\mathbb{Q})
\end{equation*}
\end{conj}
Relatively little is known about the case $\mathbb{K}=\mathbb{Q}$, and the proof techniques developed in
\cite{K-BK1} rely heavily on model-theoretic objects such as infinite prime numbers (in the sense of non-standard
analysis); that in turn requires the base field $\mathbb{K}$ to be of characteristic zero and algebraically closed
(effectively $\mathbb{C}$ by the Lefschetz principle). However, even the seemingly easier analogue of the above
conjecture, the case $\mathbb{K}=\mathbb{C}$, is known (and positive) only for $n=1$.

\medskip

In the case $n=1$, the affirmative answer to the Kontsevich conjecture, as well as positivity of several isomorphism
statements for algebras of similar nature, relies on the fact that all automorphisms of the algebras in question are tame
(see definition below). Groups of tame automorphisms are rather interesting objects. Anick \cite{An} has proved
that the group of tame automorphisms of $\mathbb{K}[x_1,\ldots,x_N]$ is dense (in power series topology) in the
subspace of all endomorphisms with non-zero constant Jacobian. This fundamental result enables one to reformulate
the Jacobian conjecture as a statement on invertibility of limits of tame automorphism sequences.

\medskip

Another interesting problem is to ask whether all automorphisms of a given algebra are tame \cite{Jung, VdK,
Czer1, Czer2, SU2}. For instance, it is the case \cite{ML2, ML3} for $\mathbb{K}[x,y]$, the free associative
algebra $\mathbb{K}\langle x,y\rangle$ and the free Poisson algebra $\mathbb{K}\lbrace x,y\rbrace$. It is also the
case for free Lie algebras (a result of P. M. Cohn). On the other hand, tameness is no longer the case for
$\mathbb{K}[x,y,z]$ (the wild automorphism example is provided by the well-known Nagata automorphism, cf.
\cite{Shes2}).

\medskip

Anick's approximation theorem was established for polynomial automorphisms in 1983. We obtain the
approximation theorems for polynomial symplectomorphisms and Weyl algebra automorphisms. These new cases
are established after more than 30 years. The focus of this paper is \textbf{the problem of lifting of
symplectomorphisms}:

\medskip

\begin{quote}
    Can an arbitrary symplectomorphism in dimension $2n$ be lifted to an automorphism of the $n$-th Weyl
    algebra in characteristic zero?
\end{quote}

\medskip

The lifting problem is the milestone in the Kontsevich conjecture. The use of tame approximation is advantageous
due to the fact that tame symplectomorphisms correspond to Weyl algebra automorphisms: in fact \cite{K-BK1},
the tame automorphism subgroups are isomorphic when $\mathbb{K}=\mathbb{C}$.

\medskip

The problems formulated above, as well as other statements of similar flavor, outline behavior of algebra-geometric
objects when subject to quantization. Conversely, quantization (and anti-quantization in the sense of Tsuchimoto)
provides a new perspective for the study of various properties of classical objects; many of such properties are of
distinctly K-theoretic nature. The lifting problem is a subject of a thorough study of Artamonov
\cite{art1978,art1985,art1991,art1998}, one of the main results of which is the proof of an analogue of the
Serre-Quillen-Suslin theorem for metabelian algebras. The possibility of lifting of (commutative) polynomial
automorphisms to automorphisms of metabelian algebra is a well-known result of Umirbaev, cf.
\cite{umirbaev1995ext}; the metabelian lifting property was instrumental in Umirbaev's resolution of the Anick's
conjecture (which says that a specific automorphism of the free algebra $\mathbb{K}\langle x,y,z\rangle$, $\Char
\mathbb{K}=0$ is wild). Related to that also is a series of well-known papers \cite{Shes1, Shes2, SU2}.

\medskip

In this note we establish the approximation property for polynomial symplectomorphisms and comment on the lifting
problem of polynomial symplectomorphisms and Weyl algebra automorphisms. In particular, the main results
discussed here are as follows.

\medskip

\begin{theorem}\label{approx-main-1}
Let $\varphi=(\varphi(x_1),\;\ldots,\;\varphi(x_N))$ be an automorphism of the polynomial algebra
$\mathbb{K}[x_1,\ldots,x_N]$ over a field $\mathbb{K}$ of characteristic zero, such that its Jacobian
\begin{equation*}
\J(\varphi)=\Det \left[\frac{\partial \varphi(x_i)}{\partial x_j}\right]
\end{equation*}
is equal to $1$. Then there exists a sequence $\lbrace \psi_k\rbrace\subset \TAut \mathbb{K}[x_1,\ldots,x_N]$ of
tame automorphisms converging to $\varphi$ in formal power series topology.
\end{theorem}

\medskip

D. Anick \cite{An} proved this tame approximation theorem for polynomial automorphisms. In this paper we get the
approximation theorems for polynomial symplectomorphisms and Weyl algebra automorphisms.

\medskip
\begin{thm}\label{approx-main-2}
Let $\sigma=(\sigma(x_1),\;\ldots,\;\sigma(x_n),\;\sigma(p_1),\;\ldots,\;\sigma(p_n))$ be a symplectomorphism of
$\mathbb{K}[x_1,\ldots,x_n,p_1,\ldots,p_n]$ with unit Jacobian. Then there exists a sequence $\lbrace
\tau_k\rbrace\subset \TAut P_n(\mathbb{K})$ of tame symplectomorphisms converging to $\sigma$ in formal
power series topology.
\end{thm}


\medskip
\begin{theorem}\label{approx-main-3}
Let $\mathbb{K}=\mathbb{C}$ and let $\sigma:P_n(\mathbb{C})\rightarrow P_n(\mathbb{C})$ be a
symplectomorphism over complex numbers. Then there exists a sequence
\begin{equation*}
\psi_1,\;\psi_2,\;\ldots,\;\psi_k,\;\ldots
\end{equation*}
of tame automorphisms of the $n$-th Weyl algebra $W_n(\mathbb{C})$, such that their images $\sigma_k$ in
$\Aut P_n(\mathbb{C})$ converge to $\sigma$.
\end{theorem}

\medskip

The last proposition is of main concern to us. As we shall see, sequences of tame symplectomorphisms lifted to
automorphisms of Weyl algebra (either by means of the isomorphism of \cite{K-BK1}, or explicitly through
deformation quantization $P_n(\mathbb{C})\rightarrow P_n(\mathbb{C})[[\hbar]]$) are such that their limits may
be thought of as power series in Weyl algebra generators. If we could establish that those power series were actually
polynomials, then the Dixmier conjecture would imply the Kontsevich's conjecture (with $\mathbb{Q}$ replaced by
$\mathbb{C}$). Conversely, approximation by tame automorphisms provides a possible means to attack the
Dixmier conjecture (and, correspondingly, the Jacobian conjecture).

\medskip

\subsection {Holonomic \texorpdfstring{$\mathcal{D}$}{D}-modules, Lagrangian subvarieties and the work of Dodd}

The following general conjecture holds (\cite{K-BK1}, see also \cite{Kon}).
\begin{conj} \label{genconj}
Let $ X $ be a smooth variety. There is a one-to-one correspondence between (irreducible) holonomic $
\mathcal{D} (X) $ - modules and Lagrangian subvarieties $ T ^ * X $ of the corresponding dimension.
\end{conj}

Kontsevich \cite{Kon} introduces the general definition of the holonomic $ \mathcal {D} $ - module as follows. Let
$ X $ be a smooth affine algebraic variety of dimension $ n $ over the field $ \mathbb {K} $. Consider the $
\mathbb {K} $ - algebra $ \mathcal {D} (X) $ of differential operators - the algebra of operators, acting on the ring
$ \mathcal {O} (X) $ generated by functions and $ \mathbb {K} $ - derivations:
$$
f \mapsto gf, \; \; f \mapsto \xi (f), \; \; g \in \mathcal {O} (X), \; \; \xi \in \Gamma (X, T_ {X / \Spec \mathbb {K}} )
$$
The natural filtration is defined on the algebra
$$
\mathcal {D} (X) = \cup_ {k \geq 0} \mathcal {D} _ {\leq k} (X)
$$
with respect to the order of operators, the associated graded algebra is naturally isomorphic to the algebra of
functions on the cotangent bundle $ T ^ * X $. Let $ M $ be a finitely generated module over $ \mathcal {D} (X) $,
and $ V $ be a finite-dimensional subspace of elements generating $ M $. It induces a filtration
$$
M _ {\leq k} = \mathcal {D} _ {\leq k} (X) V \subset M, \; \; k \geq 0,
$$
such that the associated graded module $ \gr (M) $ is finitely generated over $ \mathcal {O} (T ^ * X) $. It is
known (this result belongs to O. Gabber, see \cite{Kon}) that its support
$$
\supp (\gr (M)) \subset T ^ * X
$$
is a coisotropic variety; in particular, the dimension of any of its irreducible components is not less than $ n $. The
support is independent of the choice of the subspace $ V $ (and is denoted in the original article \cite{Kon} via $
\supp (M) $).

A finitely generated module $ M $ is called \emph {holonomic} if, by definition, the dimension of its support is $ n $.

\medskip

Conjecture \ref{genconj} (which can also be called the Kontsevich Conjecture) generalizes the Conjecture
\ref{mainconj}, as well as Conjectures 5 and 6 of the article \cite{K-BK1} in the context of the lifting of
symplectomorphisms. Namely, with any symplectomorphism one may naturally associate a Lagrangian subvariety
(namely, its graph). On the other hand, holonomic $ \mathcal {D} $ - modules correspond to autoequivalences of
Weyl algebra, from which in principle (taking into account Conjecture 5 of article \cite{K-BK1}) one can get a
correspondence with automorphisms.

\medskip

In connection with these circumstances, the necessity to study the holonomic $ \mathcal {D} $ - modules has
naturally presented itself. The problems of lifting of polynomial symplectomorphisms in the case of low dimensions -
namely, for $ n = 1 $, which corresponds to the well-known case of Kontsevich's conjecture, has become the prime
candidate for testing these new deep insights. Some progress in this direction has been achieved in the paper
\cite{K-BE} by Kanel-Belov and  Elishev. The general case of arbitrary dimension was investigated by Kontsevich
in the main article. \cite{Kon} (see also \cite{Bit}); significant results on the Conjecture \ref{genconj} were obtained
(according to our understanding) by Dodd \cite{Dodd}.

\medskip

Namely, Dodd devised the proof of the following result.
\begin{thm} [C. Dodd, \cite{Dodd}] \label {dodd1}
Let $ X $ be a smooth variety over $ \mathbb {C} $, $ L \subset T ^ * X $ be a Lagrangian subvariety of the
cotangent space. Suppose that:

1. The projection $ \pi: L \rightarrow X $ is a dominant mapping.

2. The first singular homology group $ H_1 ^ {sing} (L, \mathbb {Z}) $ is trivial.

3. There exists a smooth projective compactification $ \bar {L} $ of the variety $ L $ with trivial $ (0,2) $ - Hodge
cohomology.

Then there exists a unique irreducible holonomic $ \mathcal {D} (X) $ - module $ M $ with constant arithmetic
support \footnote {For the definition of arithmetic support, see \cite{Kon}.}, equal to $ L $, with  multiplicity $ 1 $.
\end{thm}

This theorem partially resolves the problem of finding sufficient conditions for the correspondence between
holonomic modules and Lagrangian varieties as formulated in the Conjecture \ref{genconj}. Dodd also notes that in
the case when $ X = \mathbb {A} ^ n $ is an affine space, condition 2 of Theorem \ref{dodd1} can be dropped, in
connection with which there is a corollary:
\begin{cor}[C. Dodd, \cite{Dodd}] \label {dodd2}
Let $ L \subset T ^ * \mathbb {A} ^ m $ be a smooth Lagrangian subvariety satisfying conditions 2 and 3 of
Theorem \ref{dodd1}. Then there exists a unique irreducible holonomic $ \mathcal {D} (\mathbb {A} ^ m)
$-module  $ M $ whose arithmetic support is $ L $, with multiplicity $ 1 $.
\end{cor}

This result is closely related to the construction studied in \cite{K-BE}.

As Dodd notes, Theorem \ref{dodd1} and Corollary \ref{dodd2} allow us to give a description of the Picard group
$ \Pic (W_ {n, \mathbb {C}}) $ of the Weyl algebra. Recall that the Picard group of an associative algebra is
defined as a group of classes (modulo isomorphism) of invertible bimodules over a given algebra, with a group
operation given by the tensor product of modules.

Consider polynomial symplectomorphisms of the variety $ T ^ * \mathbb {A} ^ m $. It is easy to show that the
graph of any symplectomorphism $ \varphi $ is a Lagrangian subvariety of $ L ^ {\varphi} $ in $ T ^ * \mathbb {A}
^ {2m} $, isomorphic to $ \mathbb {A} ^ {2m} $ and, therefore, satisfies cohomological conditions of Theorem
\ref{dodd1}. Applying Corollary \ref{dodd2}, we obtain (uniquely identified) $ \mathcal {D} (\mathbb {A} ^
{2m}) \simeq \mathcal {D} (\mathbb {A} ^ m) \otimes \mathcal {D} (\mathbb {A} ^ m) ^ {op} $ - module $ M ^
{\varphi} $ corresponding to $ L ^ {\varphi} $. One can check \cite{Dodd} that the inverse symplectomorphism $
\varphi ^ {- 1} $ in such a construction corresponds to inverse bimodule.

From these considerations, Dodd obtains the following result.
\begin{thm}[C. Dodd, \cite{Dodd}] \label {dodd3}
There is an isomorphism of groups (over $ \mathbb{C}$)
$$
\Pic (\mathcal {D} (\mathbb {A} ^ m)) \simeq \Sympl (T ^ * \mathbb {A} ^ m),
$$
where $ \Sympl (T ^ * \mathbb {A} ^ m) $ denotes the group of polynomial symplectomorphisms  (this group is a
geometric analogue of the group $ \Aut P_ {m, \mathbb {C}} $).
\end{thm}

In the case $ m = 1 $, it is known (Dixmier, \cite{Dix}) that $$ \Pic (\mathcal {D} (\mathbb {A} ^ 1)) = \Aut
(\mathcal {D} (\mathbb {A} ^ 1)), $$ and the algebra $ \mathcal {D} (\mathbb {A} ^ 1) $ is isomorphic to the first
Weyl algebra $ W_ {1, \mathbb {C}} $. This means that we are in the situation of Conjecture \ref {mainconj} for $
m = 1 $.

\medskip

\subsection {Tame automorphisms and the Kontsevich Conjecture}

Dodd's constructions are deep in content and, apparently, prove the Kontsevich Conjecture on the correspondence
between Lagrangian varieties and holonomic modules (more precisely, its essential part). On the other hand, starting
from Theorem \ref{dodd3} we cannot immediately arrive at the general case of Conjecture \ref{mainconj} - proof
of Conjecture 1 of article \cite{K-BK1} requires a solution to the lifting problem of symplectomorphisms to
automorphisms of the corresponding Weyl algebra.

One of the main results of the paper \cite{K-BK1} was the proof of the following homomorphism properties
$$
\phi: \Aut W_ {n, \mathbb {C}} \rightarrow \Aut P_ {n, \mathbb {C}}
$$
constructed in \cite{K-BK1} and \cite{Tsu2}. First, let $ \varphi $ be an automorphism of the polynomial algebra $
\mathbb {K} [x_1, \ldots, x_n] $. We call $ \varphi $ \emph {elementary} if it has the form
\begin {equation*}
\varphi = (x_1, \ldots, \; x_ {k-1}, \; ax_k + f (x_1, \ldots, x_ {k-1}, \; x_ {k + 1}, \; \ldots, \; x_n), \; x_ {k + 1},
\; \ldots, \; x_n).
\end {equation*}
In particular, automorphisms given by linear substitutions  of generators are elementary. Denote by $\TAut \mathbb
{K} [x_1, \ldots, x_n] $ the subgroup generated by all elementary automorphisms. Elements of this subgroup are
called \emph {tame automorphisms} of the algebra of polynomials; non-tame automorphisms  are called \emph {wild
automorphisms}.

\medskip

Tame automorphisms of the algebra $ P_ {n, \mathbb {K}} $ are, by definition, compositions of those tame
elementary automorphisms which preserve the Poisson bracket. Tame automorphisms of the Weyl algebra are
defined $ W_ {n, \mathbb {K}} $ similarly.

\medskip

The following theorem is proved in \cite{K-BK1}.
\begin{thm} [A. Kanel-Belov and M.L. Kontsevich, \cite{K-BK1}] \label {thmtame}
The homomorphism $$ \phi: \Aut W_ {n, \mathbb {C}} \rightarrow \Aut P_ {n, \mathbb {C}} $$ restricts to
isomorphism
$$
\phi_ {| \TAut}: \TAut W_ {n, \mathbb {C}} \rightarrow \TAut P_ {n, \mathbb {C}}
$$
between subgroups of tame automorphisms.
\end{thm}

In particular, due to the tame nature of automorphism groups of Weyl and Poisson algebras for $ n = 1 $, the
homomorphism $ \phi $ gives an isomorphism of the Kontsevich Conjecture between $ \Aut W_ {1, \mathbb {C}}
$ and $ \Aut P_ {1, \mathbb {C}} $.

\medskip

It is not known whether all automorphisms of the Poisson and Weyl algebras are tame for $ n> 1$, or even stably
tame (an automorphism is called \emph {stably tame} if it becomes tame after adding dummy variables and
extending the action on them by means of the identity automorphism). For the algebra of polynomials in three
variables, the Nagata automorphism
$$
(x, y, z) \mapsto (x-2 (xz + y ^ 2) y- (xz + y ^ 2) ^ 2z, y + (xz + y ^ 2) z, z)
$$
is wild (the famous result due to  I.P. Shestakov and U.U. Umirbaev, \cite{Shes2}).

\medskip

Nevertheless, tame automorphisms turn out to play a significant role in the context of the Kontsevich Conjecture and
Jacobian Conjecture, due to the following reason. D. Anick \cite{An} showed that the set of tame automorphisms of
the algebra of polynomials $ \mathbb {K} [x_1, \ldots, x_n] $ ($ n \geq 2 $) is dense in the topology of formal
power series in the space $\mathfrak{J} $ of polynomial endomorphisms with nonzero constant Jacobian. In
particular, for any automorphism of a polynomial algebra there exists a sequence of tame automorphisms converging
to it in this topology - in other words, Anick's theorem implies the existence of \emph{approximations} of
automorphisms, or \emph{approximations} by tame automorphisms (and in general, endomorphisms with nonzero
constant Jacobian). In view of the theorem of Anick, the Jacobian Conjecture can be formulated as a problem of
invertibility of limits of sequences of tame automorphisms (this is discussed in the conclusion of the article \cite{An}).
This formulation of the Jacobian Conjecture can be directly generalized to the case of a field of arbitrary
characteristic,  see more below as well as in \cite{KBYu}.

\medskip

Anick's results, together with Theorem \ref {thmtame}, suggest the idea of solving the lifting problem of polynomial
symplectomorphisms to automorphisms of the Weyl algebra, alternative to that proposed in \cite{K-BK1}
constructs. Namely, if there is a symplectic analogue of Anick's theorem - that is, if there is an approximation of
polynomial symplectomorphisms by tame symplectomorphisms, then, taking a sequence of tame
symplectomorphisms converging to a given point, we can take the sequence their pre-images under the isomorphism
$\phi_{| \TAut} $ and try to prove that its limit exists and is an automorphism of the Weyl algebra. The symplectic
analogue of Anick's theorem  was proved in \cite{KGE}. The application of approximation theory to the lifting
problem constitutes the main idea of the proof of Conjecture \ref{mainconj} in \cite{K-BE4}.

\medskip

However, the direct application of the main result of \cite{KGE} to the solution of the lifting problem does not
achieve the desired result, since the homomorphism $\phi$ does not preserve the topology of formal power series
(due to commutation relations in the Weyl algebra). In connection with this circumstance, the naive approximation
approach needs some modification. It turns out that such a modification is possible (see \cite{K-BE4}). The nature
of this modification is significant. It is connected with the geometric properties of $\Ind$-schemes of automorphisms
of the corresponding algebras. \emph {Therefore, the study of the geometry of $\Ind$-schemes of automorphisms is
justified in the framework of  Kontsevich Conjecture.}

\subsection {Homomorphism \texorpdfstring{$\phi$}{} and Ind-schemes}

It is necessary to mention one more circumstance justifying the study of $\Ind$-schemes in the context of Kontsevich
Conjecture. In the construction of the homomorphism

$$
\phi: \Aut W_ {n, \mathbb {C}} \rightarrow \Aut P_ {n, \mathbb {C}}
$$
as in the papers \cite{K-BK1} and \cite{Tsu2}, tools from model theory and non-standard analysis are used. In
particular, the fixed nonprincipal ultrafilter $ \mathcal{U} $ on the index set as well as a fixed infinitely large prime
number $ [p] $ are involved in the construction of $ \phi $.

\medskip

Briefly, the construction of the homomorphism $ \phi $ proceeds as follows. First of all, we make the following
observations (cf.  \cite{K-BK1, K-BK2, Tsu1, Tsu2}).

1. Over the field $ \mathbb {F}_p $ of positive characteristic $ p $ the Weyl algebra $ W_ {n, \mathbb {F} _p} $
has a large center given by a subalgebra
$$
\mathbb {F} _p [x_1 ^ p, \ldots, x_n ^ p, d_1 ^ p, \ldots, d_n ^ p]
$$
($ x_i $, $ d_j $ are generators of the Weyl algebra; thus, the center is the polynomial algebra).

2. Any endomorphism $ \varphi $ of the algebra $ W_ {n, \mathbb {F} _p} $ restricts to an endomorphism $
\varphi ^ c $ center (i.e. maps the center into itself).

3. On the center $ \mathbb {F} _p [x_1 ^ p, \ldots, x_n ^ p, d_1 ^ p, \ldots, d_n ^ p] $ there is a natural Poisson
bracket which descends from the ambient algebra. Namely, we take in $ \mathbb {F} _p $ the prime subfield $
\mathbb {Z} _p $  and define the projection
$$
\pi: W_ {n, \mathbb {Z}} \rightarrow W_ {n, \mathbb {Z} _p}
$$
from the Weyl algebras over the ring of integers. Let $ a, b $ belong to the center of $ W_ {n, \mathbb {Z} _p} $.
Take arbitrary inverse images $ a ', b' $ of elements $ a, b $ under reduction. Consider the expression
$$
\lbrace a, b \rbrace = \pi \left (- \frac {[a ', b']} {p} \right).
$$
It can be shown that it correctly defines the Poisson bracket on the center of the algebra $ W_ {n, \mathbb {Z} _p}
$. It can be extended linearly to the entire algebra $ W_ {n, \mathbb {F} _p} $. Since the Poisson bracket is
induced by the commutator, for any endomorphism $\varphi$ its restriction to the center$\varphi ^ c$ preserves the
Poisson bracket. Thus, in characteristic $ p $ there is a homomorphism from the group of automorphisms of the
Weyl algebra to the group of polynomial symplectomorphisms. (Strictly speaking, the result of the mapping $ \varphi
\mapsto \varphi ^ c $ has to be untwisted via the action of the inverse Frobenius automorphism; it is done in order to
get rid of the $p$-th powers of the coefficients in the image).

The case of characteristic zero reduces to positive characteristic as follows. It is well known (see \cite{LLS} and
\cite{FMS}) that the set of minimal prime ideals of the ring
$$
A = \prod_ {I} A_i
$$
where $ A_i $ are domains and $ I $ is an arbitrary nonempty set of indices (finite or infinite), is in one-to-one
correspondence with the set of all ultrafilters on the set of subsets of $ I $. Namely, if $ \mathcal {U} $ is an
ultrafilter, then we put
$$
(\mathcal {U}) = \lbrace a \in A \; | \; \theta (a) \in \mathcal {U} \rbrace,
$$
Where
$$
\theta (a) = \lbrace i \in I \; | \; a_i = 0 \rbrace.
$$
Then $ (\mathcal {U}) $ is a minimal prime ideal in $ A $, and all minimal prime ideals are of this form. Furthermore,
if $ \mathbb {P} ^ {\mathbb {N}} $ is the set of sequences of prime numbers, and $ \mathcal {U} $ is a
nonprincipal ultrafilter on subsets of natural numbers (playing here the role of the index set $ I $), then it determines
the equivalence relation on prime number sequences: two sequences are equivalent if the set of indices on which the
entries of the two sequences coincide lies in $ \mathcal {U} $. Consider the set of equivalence classes:
$$
{}^* \mathbb {P} \equiv \mathbb {P} ^ {\mathbb {N}} / \sim _ {\mathcal {U}}.
$$
Among points in that set there are points corresponding to prime numbers (classes of stationary sequences), while
there are also other classes. The latter are called \emph {infinitely large primes}. Such a definition corresponds to the
concept of a prime number in the ring of hyperinteger numbers (see \cite{robinson2016non} pp. 432-440).

\medskip

If we fix the nonprincipal ultrafilter $ \mathcal {U} $ (on $ \mathbb {N} $) and consider any infinitely large prime $
[p] $ represented by the sequence $ (p_m) $, then we can consider the ring
$$
\prod_ {m \in \mathbb {N}} \mathbb {F} _ {p_m}
$$
(direct product of a countable number of algebraically closed fields of positive characteristics). The Krull dimension
of the direct product (of any set) of fields is zero, therefore the minimal ideal $ (\mathcal {U}) $ will also be the
maximal, so that the quotient ring is a field. When $ [p] $ is infinitely large relative to $ \mathcal {U} $, this factor set
has the power of the continuum (the proof can be found, for example, in \cite{Tsu2} or in \cite{K-BK2}). It is also
easy to see that such a field is algebraically closed and of characteristic zero. Therefore, due to the well-known result
from field theory, this field will be isomorphic to the field of complex numbers.

Now, consider the ultraproduct (a direct product factorized with respect to the equivalence relation induced $
\mathcal {U} $) of Weyl algebras
$$
\left (\prod_ {m \in \mathbb {N}} W_ {n, \mathbb {F} _ {p_m}} \right) / (\mathcal {U}).
$$
The center of this algebra contains a (proper) subalgebra isomorphic to the algebra of polynomials
$$
\mathbb {C} [\xi_1, \ldots, \xi_ {2n}]
$$
(it is obtained due to the isomorphism of the base field and observation 1 above). Taking an arbitrary endomorphism
of algebra Weyl in characteristic zero, one can lift it to an automorphism of the ultraproduct and then apply the
results of observations 1 - 3 to the endomorphisms in the ultraproduct decomposition and restrict these
symplecto-endomorphisms onto the subalgebra $ \mathbb {C} [\xi_1, \ldots, \xi_ {2n}] $ (which is possible due to
the fact that the original endomorphism is of finite degree). Thus constructed mapping defines the Belov - Kontsevich
homomorphism $\phi$ (the construction given here follows Y. Tsuchimoto \cite{Tsu2}).

\medskip

The nature of infinitely large primes (as well as the nature of the ultrafilter $ \mathcal {U} $) is fundamentally
non-constructible. In view of that, the question of  independence of the homomorphism $ \phi = \phi _ {[p]} $ of the
choice of infinitely large prime number $ [p] $ naturally arises. This issue was resolved in the paper \cite{K-BK2}.
In the proof presented in \cite{K-BK2}, the properties of the so-called loop morphism
$$
\Phi: \Aut W_ {n, \mathbb {C}} \rightarrow \Aut W_ {n, \mathbb {C}},
$$
obtained from $ \phi _ {[p]} $ and $ \phi _ {[p ']} $ ($ [p] \neq [p'] $) under the assumption of the Conjecture
\ref{mainconj} (which in turn is justified in the article \cite{K-BE4}) are studied. The loop morphism provides an
example of an $\Ind$-automorphism of $\Ind$-schemes, and some of its essential properties (for example, local
unipotency) can be proved using the algebraic-geometric methods examined in the article \cite{KBYu}. In this
sense, the situation with the Kontsevich Conjecture is similar to the investigation of $\Ind$-schemes in \cite{KBYu},
however, it should be noted that the geometry of $\Ind$-schemes $ \Aut P_ {n, \mathbb {C}} $ (and $\Aut W_ {n,
\mathbb {C}}$) is more complicated than that of their analogues in the commutative and free associative cases.

\medskip

The study of the geometry of $\Ind$-schemes of automorphisms was conceived separately from the Jacobian
Conjecture and related conjectures in the classical works of B. I. Plotkin \cite{BIP1, BIP2}. In particular, Plotkin
studied the problem of the structure of sets
$$
\Aut \Aut \; \; \text {and} \; \; \Aut \End
$$
whose elements are automorphisms of automorphism groups (respectively, semigroups endomorphisms) of algebras
with polynomial identities. Significant results in this direction were obtained by A. Kanel-Belov, R. Lipyanski and A.
Berzins \cite{BelovLiapiansk2, KBLBerz, Berz}. A description of the sets
$$
\Aut \End \mathbb {K} [x_1, \ldots, x_n]
$$
and
$$
\Aut \End \mathbb {K} \langle x_1, \ldots, x_n \rangle
$$
of automorphisms of semigroups of endomorphisms of polynomial algebra and free algebra was obtained.

\medskip

Similar questions for automorphism groups can be resolved at the level of automorphisms of $\Ind$-schemes given
by
$$
\Aut _ {\Ind} \Aut
$$
for the polynomial algebra and the free associative algebra. This was done by A. Kanel-Belov, A. Elishev and J.-T.
Yu in the paper \cite{KBYu}. The fact that the study of the sets of $\Ind$-automorphisms utilizes approximation by
tame automorphisms along with techniques from algebraic geometry and topology (namely, the study of the
properties of curves in $ \Aut $), which find their application in the (modified) context of the lifting problem and
Kontsevich Conjecture. In this sense, the study of the geometry of $\Ind$-schemes is necessary to prove
Kontsevich Conjecture and - in the sense of the lifting problem - that study precedes it.

\medskip

Generally speaking, the geometry of automorphism groups of affine varieties, going back to I. R. Shafarevich and B.
I. Plotkin, is an actively developing field, regardless of progress in the work on the Jacobian Conjecture. Various
topics of independent interest are discussed in the review article by J.-P. Furter and H. Kraft \cite{FurKr}. The
study of questions similar to the results discussed in the present review was the subject of the works of T.
Kambayashi \cite{Kam96, Kam03}, H. Kraft and A. Regeta \cite{KrReg}, H. Kraft and I. Stampfli
\cite{KraftStampfli}, S. Kovalenko, A.Yu. Perepechko and M.G. Zaidenberg \cite{KPZ}, I.V. Arzhantsev, K.G.
Kuyumzhiyan and M.G. Zaidenberg \cite{AKZ}.

\medskip

Of particular interest, in the light of results obtained in \cite{KGE,K-BE4}, is the recent work by K. Urech and S.
Zimmermann \cite{UrZ}, which proves the following result: if an automorphism of the Cremona group of arbitrary
rank is also a topology homeomorphism in either Zariski or Euclidean topology, then it is an inner automorphism (up
to action of an automorphism induced by an automorphism of the base field). Also in \cite{UrZ} a similar result is
established, obtained by replacing the Cremona group with the group of polynomial automorphisms of an affine
space, thus the work of Urech and Zimmermann is a generalization of \cite{KBYu}.

\medskip

\subsection {Approximation and analysis of curves in Aut}

In the proof of Kontsevich conjecture in \cite{K-BE4}, namely in the proof of the correctness of the lifting of
augmented symplectomorphisms, a certain kind of analysis of curves of automorphisms is used. In essence, it is a
technique that allows, by considering curves in an infinite-dimensional variety $ \Aut $, one to control the degrees of
large-order terms in the automorphism and its image under a morphism $ \Phi $, as long as the initial automorphism
is close enough to the identity automorphism. Also, this procedure plays a significant role in the proof of the main
results of \cite{KBYu}.

\medskip

As we mentioned before, the naive approximation by tame symplectomorphisms does not achieve the resolution of
the lifting problem to the Weyl algebra. In connection with the search for a stronger approximation theory, it became
necessary to introduce the deformation (or augmentation, see \cite{K-BE4}) of algebras and, accordingly, of the
power series topology. In the augmented case, the singularity analysis procedure presented below (which we often
call the "singularity trick", eng. singularity trick) allows one to establish the correctness of the lifting procedure, which
in turn allows one to apply the theory of approximation by tame symplectomorphisms to the construction of the lifting
map acting on augmented symplectomorphisms.

\medskip

This idea amounts to a very useful procedure when working with infinite-dimensional varieties. It was first described
in \cite{KBYu} (Theorem 3.2, Lemmas 3.5, 3.6, and 3.7).

\medskip

For the case of the commutative algebra of polynomials $ \mathbb {K} [x_1, \ldots, x_n] $ the singularity trick has
the following form.

\medskip

Let $ L = L (t) $ be a curve whose points are linear automorphisms, i.e. a curve
$$ L \subset \Aut (K [x_1, \dots, x_n]), $$ the points of which are given by linear changes of variables. Suppose that as $ t $ tends to zero, the $i$-th eigenvalue of the matrix $ L (t) $ (corresponding to linear substitutions) also tends to zero as $ t ^ {k_i} $, $ k_i \in \mathbb { N} $. Such a curve always exists.

Now note that the orders $ \lbrace k_i, \; i = 1, \ldots n \rbrace $ of singularities of eigenvalues at zero are such that
for any pair $ (i, j) $, if $ k_i \neq k_j $, there exists a positive integer $ m $ what
$$
\text {either \; \;} k_im \leq k_j \; \; \text {or \; \;} k_jm \leq k_i.
$$
\begin{Def} \label{defsingtrickdeg}
The largest such $ m $ will be called the order of the curve $ L (t) $ at $ t = 0 $.
\end{Def}
Since all $ k_i $ are natural numbers, the order equals the integer part of $ \frac {k _ {\text {max}}} {k _ {\text
{min}}} $.

Let $ M \in \Aut_0 (K [x_1, \dots, x_n]) $ be a polynomial automorphism preserving the origin. The following
statement holds.

\begin{lem} \label{Lm2init} The curve $ L (t) ML (t) ^ {- 1} $ has no singularity at zero for any everywhere diagonalizable conjugating curve $ L (t) $ of order $ \leq N $ if and only if $ M \in \hat {H} _N $, where $ \hat {H} _N $ denotes the subgroup of automorphisms that are homothety \footnote{Detailed definitions are given in the first section of Chapter 2.} modulo $ N $ degrees of the expanding ideal $ (x_1, \ldots, x_n) $.
\end{lem}

From a topological point of view, this is a criterion for a point to belong to some sufficiently small neighborhood of
the identity automorphism. In particular, variations of this singularity trick are useful in proving the continuity of
mappings with sufficient regularity. The proof of Lemma \ref{Lm2init} is given in \cite{KBYu}.

\medskip

\subsection{Quantization of classical algebras}

As already noted, the approach to the Jacobian conjecture, using techniques from the theory of deformation
quantization - namely, the approach based on stable equivalence between the Jacobian conjecture and the Dixmier
conjecture as well as, to a somewhat lesser extent, the Kontsevich Conjecture - is currently one of the more
promising approaches to finding a possible solution to the Jacobian Conjecture. However, as in questions of the
geometric theory of $\Ind$-schemes and infinite-dimensional algebraic groups, the issues arising in connection with
the application of quantization methods, due to their nontriviality and depth, is a direction whose value may well be
comparable with the value of a possible solution to the original problem.

\medskip

Analogues of JC and DC for algebras of quantum polynomials are not obvious and often do not admit a naive
transfer of formulations (for example, E. Backelin \cite{Backelin} wrote about the $q$-quantum version of the
Dixmier conjecture). On the other hand, the well known theorem of Umirbaev \cite {U}, showing the validity of an
analogue of the Jacobian conjecture for free metabelian algebras, may be considered as an argument in favor of the
validity of the Jacobian Conjecture.

\medskip

Significant development of algebra and non-commutative geometry of quantum polynomials has been achieved by
V.A. Artamonov \cite{art1978,art1985,art1991,art1998,art2008}. In particular, he proved \cite{art1998} the
quantum-algebra analogue of the Serre conjecture (Quillen--Suslin theorem) -- the result which is extremely
non-trivial even in the commutative case.

\medskip

In connection with the Jacobian Conjecture, we mention the works of Dicks \cite{Di}, Dicks and Lewin
\cite{DiLev} as well as Yagzhev \cite{Yag1, Yag2, Yag3, Yag4}. In a sense, they can be interpreted as works
consistent with the point of view on the Jacobian problem as a problem related to quantization. Regarding the
practical benefits of studying relationships induced by quantization-type correspondences, there are known examples
of application of the elements of the quantization procedure to some (previously proven by other means) problems of
general algebra. An example \cite{zhang2017,Zhang-mas} is a new proof of Bergman's centralizer theorem
\ref{berg} of the free associative algebra, based on the deformation quantization procedure, which we discuss in this
work.

\medskip

\section{Torus actions on free associative algebras and the Bia\l{}ynicki-Birula theorem}

In the proof of the results concerning the geometry of $\Ind$-schemes automorphisms, we use the famous A.
Bia\l{}ynicki-Birula theorem \cite{BialBir1, BialBir2} on the linearizability of regular actions of a maximal torus on
an affine space merits consists in the following.

\medskip

Let $ \mathbb {K} $ be the base field, and let $ \mathbb{K}^{\times} = \mathbb{K}\backslash \{0 \} $ be the
multiplicative group of the field, considered as an algebraic $\mathbb{K}$ - group.


We call a $ n $ -dimensional algebraic $ \mathbb {K} $ - torus a group
$$
\mathbb {T} _n \simeq (\mathbb {K} ^ {\times}) ^ n
$$
(with obviously certain multiplication).

\begin{Def}
An $n$-dimensional algebraic $\mathbb{K}$-torus is a group
$$
\mathbb{T}_n\simeq (\mathbb{K}^{\times})^n
$$
(with obvious multiplication).
\end{Def}
Denote by $\mathbb{A}^n$ the affine space of dimension $n$ over $\mathbb{K}$.
\begin{Def}
A (left, geometric) torus action is a morphism
$$
\sigma: \mathbb{T}_n\times \mathbb{A}^n\rightarrow \mathbb{A}^n.
$$
that fulfills the usual axioms (identity and compatibility):
$$
\sigma(1,x)=x,\;\;\sigma(t_1,\sigma(t_2,x))=\sigma(t_1t_2,x).
$$

The action $\sigma$ is \textbf{effective} if for every $t\neq 1$ there is an element $x\in \mathbb{A}^n$ such that
$\sigma(t,x)\neq x$.
\end{Def}

In \cite{BialBir1}, Bia\l{}ynicki-Birula proved the following two theorems, for $\mathbb{K}$ algebraically closed.

\begin{thm}
Any regular action of $\mathbb{T}_n$ on $\mathbb{A}^n$ has a fixed point.
\end{thm}

\begin{thm}
Any effective and regular action of $\mathbb{T}_n$ on $\mathbb{A}^n$ is a representation in some coordinate
system.
\end{thm}

The notion of regular action means regularity in the sense of algebraic geometry (preservation of regular functions;
Bia\l{}ynicki-Birula also considered birational actions in \cite{BialBir1}). The last theorem states that any effective
regular action of the maximal  torus on an affine space is conjugate to a linear action (representation) - in other
words, such action \emph{admits linearization}.

\medskip
An algebraic group action on $\mathbb{A}^n$ is the same as the corresponding action by automorphisms on the
algebra
$$
\mathbb{K}[x_1,\ldots,x_n]
$$
of coordinate functions. In other words, it is a group homomorphism
$$
\sigma: \mathbb{T}_n\rightarrow \Aut \mathbb{K}[x_1,\ldots,x_n].
$$
An action is effective if and only if $\Ker\sigma = \lbrace 1\rbrace$.

The polynomial algebra is a quotient of the free associative algebra
$$
F_n = \mathbb{K}\langle z_1,\ldots,z_n\rangle
$$
by the commutator ideal $I$ (it is the two-sided ideal generated by all elements of the form $fg-gf$). The definition
of torus action on the free algebra is thus purely algebraic.

The following result has been established in \cite{TA1,TA2}.

\begin{thm} \label{BBfree_copy}
Suppose given an action $\sigma$ of the algebraic $n$-torus $\mathbb{T}_n$ on the free algebra $F_n$. If
$\sigma$ is effective, then it is linearizable.
\end{thm}

\medskip

The theory of  algebraic group actions on varieties is a substantial part of the study of $\Ind$-varieties. Among the
significant works on this subject, the reader is well advised to consult the papers of T. Kambayashi and P. Russell
\cite {KR}, M. Koras and P. Russell \cite {KoRu2}, T. Asanuma \cite {Asanuma}, G. Schwartz \cite {Sch} and
H. Bass \cite {Bass84}.


The group action linearity problem asks, generally speaking, whether any action of a given algebraic group on an
affine space is linear in some suitable coordinate system (or, in other words, whether for any such action there exists
an automorphism of the affine space such that it conjugates the action to a representation). This subject owes its
existence largely to the classical work of A. Bia\l{}ynicki-Birula \cite{BialBir1}, who considered regular (i.e. by
polynomial mappings) actions of the $n$-dimensional torus on the affine space $\mathbb{A}^n$  (over algebraically
closed ground field) and proved that any faithful action is conjugate to a representation (or, as we sometimes say,
linearizable). The result of Bia\l{}ynicki-Birula had motivated the study of various analogous instances, such as those
that deal with actions of tori of dimension smaller than that of the affine space, or, alternatively, linearity conjectures
that arise when the torus is replaced by a different sort of algebraic group. In particular, Bia\l{}ynicki-Birula himself
\cite{BialBir2} had proved that any effective action of $(n-1)$-dimensional torus on $\mathbb{A}^n$ is
linearizable, and for a while it was believed \cite{KR} that the same was true for arbitrary torus and affine space
dimensions. Eventually, however, the negation of this generalized linearity conjecture was established, with
counter-examples due to Asanuma \cite{Asanuma}.

\smallskip

More recently, the linearity of effective torus actions has become a stepping stone in the study of geometry of
automorphism groups. In the paper \cite{KBYu}, the following result was obtained.
\begin{thm}        \label{commautautthm}
 Let $\mathbb{K}$ be algebraically closed, and let $n\geq 3$. Then any $\Ind$-variety automorphism $\Phi$ of the $\Ind$-group
$\Aut(K[x_1,\dots,x_n])$ is inner.
\end{thm}
The notions of $\Ind$-variety (or $\Ind$-group in this context) and $\Ind$-variety morphism were introduced by
Shafarevich \cite{Shafarevich}: an $\Ind$-variety is the direct limit of a system whose morphisms are closed
embeddings. Automorphism groups of algebras with polynomial identities, such as the (commutative) polynomial
algebra and the free associative algebra, are archetypal examples; the corresponding direct systems of varieties
consist of sets $\Aut^{\leq N}$ of automorphisms of total degree less or equal to a fixed number, with the degree
induced by the grading. The morphisms are inclusion maps which are obviously closed embeddings.

\medskip

Theorem \ref{commautautthm} is proved by means of tame approximation (stemming from the main result of
\cite{An}), with the following Proposition, originally due to E. Rips, constituting one of the key results.
\begin{prop} \label{proprips}
Let $\mathbb{K}$ be algebraically closed and $n\ge 3$ as above, and suppose that $\Phi$ preserves the standard
maximal torus action on the commutative polynomial algebra\footnote{That is, the action of the $n$-dimensional
torus on the polynomial algebra $\mathbb{K}[x_1,\ldots,x_n]$, which is dual to the action on the affine space.}.
Then $\Phi$ preserves all tame automorphisms.
\end{prop}

The proof relies on the Bia\l{}ynicki-Birula theorem on the maximal torus action. In a similar fashion, the paper
\cite{KBYu} examines the $\Ind$-group $\Aut \mathbb{K}\langle x_1,\ldots, x_n\rangle $ of automorphisms of the
free associative algebra $\mathbb{K}\langle x_1,\ldots, x_n\rangle$ in $n$ variables, and establishes results
completely analogous to Theorem \ref{commautautthm} and Proposition \ref{proprips}. \footnote{The free
associative case was amenable to the above approach when $n>3$.} In accordance with that, the free associative
analogue of the Bia\l{}ynicki-Birula theorem was required.

\smallskip

Such an analogue is indeed valid, and we have established it in our notes \cite{TA1,TA2} on the subject. We will
provide the proof of this result in the sequel.

\smallskip

Given the existence of a free algebra version of the Bia\l{}ynicki-Birula theorem, one may inquire whether various
other instances of the linearity problem (such as the Bia\l{}ynicki-Birula theorem on the action of the
$(n-1)$-dimensional torus on $\mathbb{K}[x_1,\ldots,x_n]$) can be studied. As it turns out, direct adaptation of
proof techniques from the commutative realm is sometimes possible. There are certain limitations, however. For
instance, Bia\l{}ynicki-Birula's proof \cite{BialBir2} of linearity of $(n-1)$-dimensional torus actions uses
commutativity in an essential way. Nevertheless, a neat workaround of that hurdle can be performed when $n=2$,
as we show in this note. Also, a special class of torus actions (positive-root actions) turns out to be linearizable.
Finally, some of the proof techniques developed by Asanuma \cite{Asanuma} admit free associative analogues; this
will allow us to prove the existence of non-linearizable torus actions in positive characteristic, in complete analogy
with Asanuma's work.

\chapter{Quantization proof of Bergman's centralizer theorem} 
\label{Chapter2} \lhead{Chapter 2. \emph{Quantization proof of Bergman's centralizer theorem}}

We first give a brief summary of the background of two well-known centralizer theorems in the power series ring
and in the free associative algebra, i.e., Cohn's centralizer theorem and Bergman's centralizer theorem.

\section{Centralizer theorems}

This section is a relatively independent part of the paper, and only sketches proofs with classic tools, while the
following sections will focus on the new proof of Bergman's centralizer theorem.
\medskip

Throughout this section, $X$ is a finite set of noncommutative variables, and $k$ is a field. Let $X^*$ denote the
free monoid generated by $X$. Let $k\langle X\rangle$ and $k\langle X\rangle$ denote the $k$-algebra of formal
series and noncommutative polynomials (i.e., the free associative algebra over $k$) in $X$, respectively.
\medskip
Both elements of $k\llangle X\rrangle$ and $k\langle X\rangle$ have the form $a=\sum_{\omega\in
X^*}a_{\omega}\omega$, where $a_{\omega}\in k$ is the coefficient of the word $w$ in $a$, but they have
different details inside the above formula. An element of $k\langle X\rangle$ is only a finite sum of words, while there
are infinitely many terms of the sum for an element in $k\llangle X\rrangle$. The multiplication of elements in
$k\llangle X\rrangle$ is the concatenation of words and normal multiplication of coefficients. We can only combine
the coefficients which have the same corresponding words for addition.
\medskip
The \textit{length} $|\omega|$ of a word $\omega\in X^*$ is the number of letters inside $\omega$. Now we can
define the valuation
$$\nu: k\llangle X\rrangle\to \mathbb{N}\cup\{\infty\}$$ as follows: $\nu(0)=\infty$ and if $a=\sum_{\omega\in X^*}a_{\omega}\omega\neq 0$, then $\nu(a)=\min\{|\omega|:a_{\omega}\neq0\}$. Note that if $\omega$ is constant, then $\nu(\omega)=0$ and $\nu(ab)=\nu(a)+\nu(b)$ for all $a,b$ in $k\llangle X\rrangle$.

For the words valuation, there is an easy but quite useful lemma \cite{sha2013}.

\begin{lemma}[Levi's Lemma]
Let $\omega_1,\omega_2,\omega_3,\omega_4\in X^*$ be nonzero with $|\omega_2|\geq |\omega_4|$. If
$\omega_1\omega_2=\omega_3\omega_4$, then $\omega_2=\omega\omega_4$ for some $\omega\in X^*$.
\end{lemma}

The proof is trivial by backward induction on $|\omega_2|$ since $\omega_2$ has the same last letter as
$\omega_4$. Next lemma extends Levi's lemma to $k\llangle X\rrangle$, and we post the result as follows.

\begin{lemma}[\cite{lothaire1997}, Lemma 9.1.2]\label{lemloth}
Let $a,b,c,d\in k\llangle X\rrangle$ be nonzero. If $\nu(a)\geq\nu(c)$ and $ab=cd$, then $a=cq$ for some $q\in
k\llangle X\rrangle$.
\end{lemma}

\begin{proof}
We can fix a word $u$ which appears in $b$ and $|u|=\nu(b)$. Suppose $v$ is any nonzero word appearing in $d$,
then we have
\begin{equation}\label{levieqn}
    |v|\geq \nu(d)=\nu(a)+\nu(b)-\nu(c)\geq \nu(b)=|u|.
\end{equation}
Let $w$ be any word in $X^*$. The coefficient of $wu$ in $ab$ is $\sum_{rs=wu}a_r b_s$, where $a_r$ and
$b_s$ are the coefficients of the words $r,s$ which appear in $a,b$ respectively. Similarly, the coefficient of $wu$
in $cd$ is $\sum_{yz=wu}c_y d_z$. Since $ab=cd$, we have
\begin{equation}\label{leviequal}
    \sum_{rs=wu}a_r b_s=\sum_{yz=wu}c_y d_z.
\end{equation}
By the inequality \ref{levieqn}, we have $|z|\geq |u|$, and $|s|\geq |u|$ by the definition of $u$. Thus $rs=wu$ and
$yz=wu$ imply $s=s_1u$ and $z=z_1u$ for some $s_1,z_1\in X^*$, by Levi's lemma. Hence $rs_1=yz_1=w$ and
we can rewrite the formula \ref{leviequal} as
\begin{equation}\label{newleviequal}
    \sum_{rs_1=w}a_r b_{s_1 u}=\sum_{yz_1=w}c_y d_{z_1 u}.
\end{equation}
Let $b'=\sum_{s_1\in X^*}b_{s_1u}s_1$ and $d'=\sum_{z_1\in X^*}b_{z_1u}z_1$. Then the equation gives
$ab'=cd'$. The constant term of $b'$ is $b_u\neq 0$ and hence $b'$ is invertible in $k\llangle X\rrangle$. Hence if
we let $q=d'b'^{-1}$, then $a=cq$.
\end{proof}

\subsection{Cohn's centralizer theorem}

With the help of the preceding lemmas, we could post and prove this well-known centralizer theorem of $k$-algebra
of formal series by P. M. Cohn.

\begin{theorem}[Cohn's Centralizer Theorem, \cite{cohn1964}]\label{cohnthm}
If $a\in k\llangle X\rrangle$ is not a constant,  then the centralizer $C(a;k\llangle X\rrangle)\cong k[\![x]\!]$, where
$k[\![x]\!]$ is the algebra of formal power series in the variable $x$.
\end{theorem}

\begin{proof}
Let $C:=C(a;k\llangle X\rrangle)$. Let $a_0$ be the constant term of $a$, then it is clear that $C=C(a-a_0;k\llangle
X\rrangle)$. So we may assume that the constant term of $a$ is zero. Thus we have a nonempty set $A=\{c\in
C:\nu(c)>0\}$ because $a\in C$ and so there exists $b\in A$ such that $\nu(b)$ is minimal. An easy observation is
that $k[\![b]\!]\cong k[\![x]\!]$. Because suppose $\sum_{i\geq m}\beta_i b^i=0,\beta_i\in k, \beta_m\neq0$, then
we must have $\infty=\nu(\sum_{i\geq m}\beta_i b^i)=\nu(b^m)=m\nu(b)$, which is absurd. So we just need to
show that $C=k[\![b]\!]$. Assume that an element $c\in C$ is not constant. Our first claim is that there exist
$\beta_i\in k$ such that
\begin{equation}\label{val1}
    \nu(c-\sum_{i=0}^n\beta_ib^i)\geq(n+1)\nu(b).
\end{equation}
The proof is by induction on $n$. let $\beta_0$ be the constant term of $c$. Then $c-\beta_0\in A$ and thus
$\nu(c-\beta_0)\geq\nu(b)$, by the minimality of $b$. This proves the $n=0$ case for the inequality \ref{val1}.

Now we need the second claim to complete this induction proof. Our second claim is following: suppose that the
constant term of an element $a\in k\llangle X\rrangle$ is zero and $b,c\in C\setminus\{0\}$. If $\nu(c)\geq\nu(b)$,
then $c=bd$ for some $d\in C$. In fact, since the constant term of an element $a\in k\llangle X\rrangle$ is zero we
have $\nu(a)\geq1$. Thus for $n$ large enough, we have $\nu(a^n)=n\nu(a)\geq\nu(c)$. we also have $a^nc=ca^n$
because $c\in C$. Thus, by lemma \ref{lemloth}, $a^n=cq$ for some $q\in k\llangle X\rrangle$. Hence,
$cqb=a^nb=ba^n$ and since $\nu(c)\geq\nu(b)$, we have $c=bd$, for some $d\in k\llangle X\rrangle$, by lemma
\ref{lemloth}. Finally, $$bad=abd=ac=ca=bda,$$ which gives $ad=da$, i.e. $d\in C$.

Now let us continue to prove the first claim. Suppose we have found $\beta_0,\dots,\beta_n\in k$ such that
$\nu(c-\sum_{i=0}^n\beta_ib^i)\geq(n+1)\nu(b)$. Then since $(n+1)\nu(b)=\nu(b^{n+1})$, we have
$c-\sum_{i=0}^n\beta_ib^i=b^{n+1}d$ for some $d\in C$, by the second claim we proved above. If $d$ is a
constant, we are done because then $c\in k[b]\subset k[\![b]\!]$. Otherwise, let $\beta_{n+1}$ be the constant
term of $d$. Then $d-\beta_{n+1}\in A$ and hence $\nu(d-\beta_{n+1})>\nu(b)$ by the minimality of $b$.
Therefore, by the first claim, $d-\beta_{n+1}=bd'$ for some $d'\in C$. Hence
$$c-\sum_{i=0}^n\beta_i b^i=b^{n+1}d=b^{n+1}(bd'+\beta_{n+1})=b^{n+2}d'+\beta_{n+1}b^{n+1},$$
which gives $c-\sum_{i=0}^{n+1}\beta_i b^i=b^{n+2}d'$. Hence
$$\nu(c-\sum_{i=0}^{n+1}\beta_i b^i)=\nu(b^{n+2}d')=(n+2)\nu(b)+\nu(d')\geq(n+2)\nu(b).$$
This completes the induction, then we are done because $\nu(c-\sum_{i\geq 0}\beta_i b^i)=\infty$ and so
$c=\sum_{i\geq 0}\beta_i b^i\in k[\![b]\!]$.
\end{proof}

\subsection{Bergman's centralizer theorem}\label{secberg}

Now since $k\langle X\rangle\subset K\llangle X\rrangle$, it follows from the above theorem that if $a\in k\langle
X\rangle$ is not constant, then $C(a;k\langle X\rangle)$ is commutative because $C(a;k\llangle X\rrangle)$ is
commutative. The next theorem is our main goal which shows that there is a similar result for $C(a;k\langle
X\rangle)$.

\begin{theorem}[Bergman's Centralizer Theorem, \cite{berg1969}]\label{berg}
If $a\in k\langle X\rangle$ is not constant, then the centralizer $C(a;k\langle X\rangle)\cong k[x]$, where $k[x]$ is
the polynomial algebra in one variable $x$.
\end{theorem}

We will not fully recover the original proof of Bergman's centralizer theorem since this is not our main idea.
However, we would take some necessary result in his original proof \cite{berg1969} which helps us to finish the
proof of that the centralizer is integrally closed. This will be shown in the Subsection \ref{bergori}.

\medskip

First of all, we need to emphasize that the proof of Cohn's centralizer theorem is included. Here is a sketch of the
proof.

\medskip

For simplicity, we denote by $C:=C(a;k\langle X\rangle)$ the centralizer of $a$ which from now on is not a
constant. Recall that the centralizer $C$ is also commutative. Moreover, $C$ is finitely generated, as module over
$k[a]$ or as algebra. Then since $k\langle X\rangle$ is a 2-fir (free ideal rings, cf. Lemma 1.5 in \cite{berg1969}),
and the center of a 2-fir is integrally closed, we obtain that the centralizer of $a$ is integrally closed in its field of
fractions after using the lifting to $k\langle X\rangle\otimes k(x)$ (where $x$ is a free variable). Then our aim is to
show that $C$ is a polynomial ring over $k$. In order to get this fact we shall study homomorphisms of $C$ into
polynomial rings. By using ``infinite" words, we obtained an embedding from $C$ into polynomial rings by
lexicographically ordered semigroup algebras, which completes this sketch of the proof. Indeed, any subalgebra not
equal to $k$ of a polynomial algebra $k[x]$ that is integrally closed in its own field of fractions is of form $k[y]$ (by
L\"{u}roth's theorem).

We conclude this section by pointing out that the method of ``infinite" words inspires us to find a possibility to prove
Bergman's centralizer theorem by deformation quantization. In the next section, we will establish this new approach
of quantization for generic matrices.

\section{Reduction to generic matrix}
In this section, we will establish an important theorem which gives a relation of commutative subalgebras in the free
associative algebra and the algebra of generic matrices. Let $k\langle X\rangle$ be the free associative algebra over
a field $k$ generated by a finite set $X=\{x_1,\dots,x_s\}$ of $s$ indeterminates, and let $k\langle
X_1,\dots,X_s\rangle$ be the algebra of $n\times n$ generic matrices generated by the matrices $X_{\nu}$. The
canonical homomorphism $\pi: k\langle x_1,\dots,x_s\rangle\to k\langle X_1,\dots,X_s\rangle$ shows in last section.

\medskip

We claim that if we have a commutative subalgebra of rank two in the free associative algebra $k\langle X\rangle$,
then we also have a commutative subalgebra of rank two if we consider a reduction to generic matrices of big
enough order $n$. We also call two elements of a free algebra \textit{algebraically independent} if the subalgebra
generated by these two elements is a free algebra of rank two. Otherwise we will call them \textit{algebraically
dependent}.
\medskip

In other words, if we have a commutative subalgebra $k[f,g]$ of rank two in the free associative algebra, then we
have to prove that its projection to generic matrices of some order also has rank two. i.e. $\pi(f),\pi(g)$ do not have
any relations.

\medskip

We need following theorem:

\begin{theorem}\label{generic}
Let $k\langle X\rangle$ be the free associative algebra over a field $k$ generated by a finite set $X$ of
indeterminates. If $k\langle X\rangle$ has a commutative subalgebra with two algebraically independent generators
$f,g\in k\langle X\rangle$, then the subalgebra of $n$ by $n$ generic matrices generated by reduction of $f$ and $g$
in $k\langle X_1,\dots,X_s\rangle$ also has rank two for big enough $n$.
\end{theorem}

\begin{proof}
Assume $k[f,g]$ be a commutative subalgebra generated by $f, g\in k\langle X\rangle\setminus k$ with rank two.
We denote $\bar{f},\bar{g}\in k\langle X_1,\dots,X_s\rangle$ to be the generic matrices of $f$ and $g$
respectively after reduction \ref{eqngeneric} of algebra of generic matrices with $n\times n$. The rank of
$k\langle\bar{f},\bar{g}\rangle$ must be $\leq2$ (i.e. it must be 1 or 2). Suppose the rank is 1, then for any two
elements $a,b\in k\langle f,g\rangle$, there exists a minimal polynomial $P(x,y)\in k[x,y]$ ($x, y$ are two free
variables) with degree $m$ such that $P(\bar{a},\bar{b})=0$ because the algebra of generic matrices is a domain
by Theorem \ref{amitthm}. On the other hand, by Amitsur-Levitzki Theorem \ref{amit-lev}, there exists no
polynomial with degree less than $2n$, such that $P(\bar{a},\bar{b})=0$. This leads to be a contradiction if we
choose $n>[m/2]$.
\end{proof}

Recall from the section \ref{secberg} that the centralizer $C:=C(a;k\langle X\rangle)$ of $a\in k\langle
X\rangle\setminus k$ is a commutative subalgebra of $k\langle X\rangle$, so from the above theorem, we conclude
that if the centralizer is a subalgebra in $k\langle X\rangle$ of rank two then the $\pi$-image subalgebra of $C$ has
also rank two.
\medskip

However, we prefer discussing this general case of subalgebras instead of just consider a centralizer subalgebra.
Furthermore, we want to prove that there is no commutative subalgebras of the free associative algebra $k\langle
X\rangle$ of rank greater than or equal to two.



\section{Quantization proof of rank one}

Up to our knowledge, there is no new proofs has been appeared after Bergman \cite{berg1969} for almost fifty
years. We are using a method of deformation quantization presented by M. Kontsevich to give an alternative proof
of Bergman's centralizer theorem. In this section \cite{zhang2017}, we got that the centralizer is a commutative
domain of transcendence degree one.

\medskip

Let $k\langle X\rangle$ be the free associative algebra over a field $k$ generated by $s$ free variables
$X=\{x_1,\dots,x_s\}$. Now, we concentrate on our proof that there is no commutative subalgebras of rank
greater than or equal to two. From the homomorphism $\pi: k\langle x_1,\dots,x_s\rangle\to k\langle
X_1,\dots,X_s\rangle$ and Theorem \ref{generic}, we are moving our goal from the elements of $k\langle
X\rangle$ to the algebra of generic matrices $k\langle X_1,\dots,X_s\rangle$, and we consider the quantization of
this algebra and its subalgebras.
\medskip

Let $A, B$ be two commuting generic matrices in $k\langle X_1,\dots,X_s\rangle$ which are algebraically
independent, i.e. $\rank k\langle A,B\rangle=2$. We have the following theorem.

\begin{theorem}\label{quan}
Let $A, B$ be two commuting generic matrices in $k\langle X_1,\dots,X_s\rangle$ with $\rank k\langle
A,B\rangle=2$, and let $\hat{A}$ and $\hat{B}$ be quantized images (by sending multiplications to star products
by means of Kontsevich's formal quantization) of $A$ and $B$ respectively by considering lifting $A$ and $B$ in
$k\langle X_1,\dots,X_s\rangle[\![\mathfrak{h}]\!]$. Then $\hat{A}$ and $\hat{B}$ do not commute. Moreover,
	\begin{equation}
	\frac{1}{\mathfrak{h}}[\hat{A},\hat{B}]_\star\equiv\begin{pmatrix} \frac{1}{\mathfrak{h}}\{\lambda_1,\mu_1\}& & 0 \\
	& \ddots & \\
	0 & &\frac{1}{\mathfrak{h}}\{\lambda_n,\mu_n\}\end{pmatrix}\mod{\mathfrak{h}}
	\end{equation}
	where $\lambda_i$ and $\mu_i$ are eigenvalues(weights) of $A$ and $B$ respectively.
\end{theorem}

To prove this theorem, we need some preparations. It is not easy to directly compute such two generic matrices with
order $n$. However, if we can diagonalize those matrices, then computation will be easier. So first of all, we should
show the possibilities. Without loss of generality, we may assume that one of the generic matrices $B$ is diagonal if
we have a proper choice of basis of the algebra of generic matrices. Now consider the other generic matrix $A$
which we mentioned above.


\begin{remark}
The generic matrix $A$ may not be diagonalizable over $k[x_{ij}^{(\nu)}]$, but it can be diagonalized over some
integral extension of the algebra $k[x_{ij}^{(\nu)}]$ with $i,j=1,\dots,n; \nu=1,\dots,s$.
\end{remark}


\begin{remark}
Any non-scalar element $A$ of the algebra of generic matrices must have distinct eigenvalues. In fact, by Amitsur's
Theorem \ref{amitthm}, namely, the algebra of generic matrices is an domain, if the minimal polynomial is not a
central polynomial, then the algebra can be embedded to a skew field. Hence, the minimal polynomial is irreducible,
and the eigenvalues are pairwise different.
\end{remark}

\begin{lemma}\label{lemquan}
Let $\hat{A}\equiv A_0+\mathfrak{h} A_1 (\mathrm{mod}\,\mathfrak{h}^2)$ be the quantized image of a generic
matrix $A\in k\langle X_1,\dots,X_s\rangle$, where $A_0$ is diagonal with distinct eigenvalues. Then, the quantized
images $\hat{A}$ can be diagonalized over some finite extension of $k[x_{ij}^{(\nu)}]$.
\end{lemma}

\begin{proof}
Without loss of generality, suppose $A_0$ is a diagonal generic matrix with distinct eigenvalues. We want to show
that there exists an invertible generic matrix $P$, such that $PAP^{-1}$ is diagonal. Now we consider their images
on $k\langle X_1,\dots,X_s\rangle[\![\mathfrak{h}]\!]$, we may assume $\hat{P}=I+\mathfrak{h}T$ and the
conjugation inverse $\hat{P}^{-1}=I-\mathfrak{h}T \mod \mathfrak{h}^2$ (where $I$ is the identity matrix ).
Then we have
$$(I+\mathfrak{h}T)(A_0+\mathfrak{h} A_1)(I-\mathfrak{h} T)=A_0+\mathfrak{h}([T,A_0]+A_1) \mod{\mathfrak{h}^2},$$
and we need to solve the equation $[T,A_0]=-A_1$.

This is clear since $A_0$ is diagonal. Let $A_0=\mathrm{diag}\{\lambda_1,\dots,\lambda_n\}$,
$T=(t_{ij})_{n\times n}$ and $A_1=(a_{ij})_{n\times n}$, then we have
$[T,A_0]=\left(\left(\lambda_i-\lambda_j\right) t_{ij}\right)_{n\times n}$. Hence, $$T=(t_{ij})_{n\times
n}=\left(-\frac{a_{ij}}{\lambda_i-\lambda_j}\right)_{n\times n}.$$

So far, we have determined the $\mathfrak{h}$ term of the matrix $\hat{A}$. Hence, we may assume
$\hat{A}\equiv A_0+\mathfrak{h}^2 A_2 \mod \mathfrak{h}^3$, then we continue to cancel the
$\mathfrak{h}^2$ term. Let $\hat{P}_2=I+\mathfrak{h}^2T_2$, and the conjugation inverse
$\hat{P}^{-1}_2=I-\mathfrak{h}^2 T_2$. Then, we have

$$(I+\mathfrak{h}^2T_2)(A_0+\mathfrak{h}^2 A_2)(I-\mathfrak{h}^2 T_2)=A_0+\mathfrak{h}^2([T_2,A_0]+A_2) \mod{\mathfrak{h}^3},$$

Hence, $T_2$ is determined by equation $[T_2,A_0]=-A_2=(a_{ij}^{(2)})_{n\times n}$. Similar computation
give all entries of $T_2$, namely $$T_2=\left(-\frac{a_{ij}^{(2)}}{\lambda_i-\lambda_j}\right)_{n\times n}.$$
Continue this process to cancel the term of $\mathfrak{h}^3$ etc., we obtain equations $[T_i,A_0]=-A_i$ for
$i=3,4,5\dots$. This leads to the result that $A$ could be diagonalized over the extension
$k[x_{ij}^{(\nu)}][\frac{1}{\lambda_i-\lambda_j}]$.
\end{proof}

Now $A, B$ are two algebraically independent but commuting generic matrices in $k\langle X_1,\dots,X_s\rangle$.
From previous discussion, we may assume $A$ and $B$ can be both diagonalized over an integral extension of
$k[x_{ij}^{(\nu)}]$. Consider result of diagonalization in $k\langle X_1,\dots,X_s \rangle[\![\mathfrak{h}]\!]$ and
then we compute the quantization commutator of two quantized generic matrices over $k\langle X_1,\dots,X_s
\rangle[\![\mathfrak{h}]\!]$. Now we can complete the proof of Theorem \ref{quan}.

\begin{proof}[Proof of Theorem \ref{quan}]
We have shown that $A,B$ can be both diagonalized over some finite extension of $k[x_{ij}^{(\nu)}]$, then
consider result of diagonalization with the quantization form in $k\langle X_1,\dots,X_s \rangle[\![\mathfrak{h}]\!]$,
i.e., we can write them into specific forms modulo $\mathfrak{h}^2$ as follows:
	\begin{equation*}
	\hat{A}\equiv \begin{pmatrix}
	\lambda_1 & & 0 \\
	& \ddots & \\
	0 & &\lambda_n\\
	\end{pmatrix}+\mathfrak{h} \begin{pmatrix}
	 \delta_1 & & * \\
	& \ddots & \\
	* & &\delta_n\\
	\end{pmatrix}\mod{\mathfrak{h}^2}
	\end{equation*}
	
	\begin{equation*}
	\hat{B}\equiv\begin{pmatrix}
	\mu_1& & 0 \\
	& \ddots & \\
	0 & &\mu_n\\
	\end{pmatrix}+\mathfrak{h} \begin{pmatrix}
	\nu_1 & & * \\
	& \ddots & \\
	* & & \nu_n\\
	\end{pmatrix}\mod{\mathfrak{h}^2}.
	\end{equation*}
	
	Then we can compute the quantization commutator,
	\begin{align*}
	[\hat{A},\hat{B}]_{\star}&:=\hat{A}\star \hat{B}-\hat{B}\star\hat{A}\equiv\begin{pmatrix}
	\{\lambda_1,\mu_1\}& & 0 \\
	& \ddots & \\
	0 & &\{\lambda_n,\mu_n\}\end{pmatrix}+\mathfrak{h} \vec{\lambda}\star\begin{pmatrix}
	0 & & * \\
	& \ddots & \\
	* & &0\\
	\end{pmatrix}\\
	&-\mathfrak{h}\begin{pmatrix}
	0 & & * \\
	& \ddots & \\
	* & &0\\
	\end{pmatrix}\star\vec{\lambda} +\mathfrak{h}\begin{pmatrix}
	0 & & * \\
	& \ddots & \\
	* & &0\\
	\end{pmatrix}\star\vec{\mu}-\mathfrak{h}\vec{\mu}\star\begin{pmatrix}
	0 & & * \\
	& \ddots & \\
	* & &0\\
	\end{pmatrix}\\
	&+ \mathfrak{h}^2 \left\{\begin{pmatrix}
	0 & & * \\
	& \ddots & \\
	* & &0\\
	\end{pmatrix},\begin{pmatrix}
	0 & & * \\
	& \ddots & \\
	* & &0\\
	\end{pmatrix} \right\} \mod{\mathfrak{h}^2}
	\end{align*}
	
Note that all terms have empty diagonals except the first term, and hence the quantization commutator
$[\hat{A},\hat{B}]_{\star}\neq 0 \mod \mathfrak{h^2}$, which completes the proof of the theorem by multiplying
$\frac{1}{\mathfrak{h}}$ on two sides of above equation.
\end{proof}

\begin{remark}

Suppose $\lambda_i$ and $\delta_i$, $i=1,\dots,n$ are algebraically dependent. Then there are polynomials $P_i$
in two variables such that $P_i(\lambda_i,\delta_i)=0$. Put $$P(x,y)=\prod_{i=1}^nP_i(x,y).$$ Then $P(A,B)$ is
diagonal matrix having zeros on the main diagonal, i.e. $P(A,B)=0$. It means that if $\rank k\langle A,B\rangle=2$
then $\lambda_i, \delta_i$ are algebraically independent for some $i$.
\end{remark}


Let us conclude this section by pointing out the whole process of this proof. Recall that we have the free associative
algebra $k\langle X\rangle$ over a field $k$, if we have a commutative subalgebra of rank two generated by $a,b\in
k\langle X\rangle$, then we may have a commutative subalgebra of the algebra of generic matrices $k\langle
X_1,\dots,X_s\rangle$ of rank two generated by $A,B$ (they are images of a homomorphism $\pi: k\langle
X\rangle\to k\langle X_1,\dots,X_s\rangle$). Consider the element $0=[a,b]$ of the free associative algebra
$k\langle X\rangle$, homomorphism $\pi$ and canonical quantization homomorphism $q$ sending multiplications to
star products, we yield that $$0=q\pi([a,b])=q[A,B]=[\hat{A},\hat{B}]_{\star}.$$ This leads a contradiction to
Theorem \ref{quan} which shows that $[\hat{A},\hat{B}]_{\star}\neq0$. So we obtain the following result.

\begin{theorem}
There is no commutative subalgebras of rank $\geq 2$ in the free associative algebra $k\langle X\rangle$. \qed
\end{theorem}

The centralizer ring is commutative from our discussion in section \ref{secberg}, and from the above theorem, it is of
rank 1. So it is a commutative subalgebra with form $k[x]$ for some $x\in k\langle X\rangle\setminus k$. We will
show it implies Bergman's centralizer theorem \ref{berg} in the next section.

\section{Centralizers are integrally closed}

We have shown that the centralizer $C$ is a commutative domain of transcendence degree one. For us, it was the
most interesting part of the proof of the Bergman's centralizer theorem. However, we have to prove the fact that
$C$ is integrally closed in order to complete the proof of Bergman's Centralizer Theorem. In our this work
\cite{zhang2018berg}, our proofs are based on the characteristic free instead of very rich and advanced P. Cohn
and G. Bergman's noncommutative divisibility theorem, we use generic matrices reduction, the invariant theory of
characteristic zero by C. Procesi \cite{Procesi} and the invariant theory of positive characteristic by A. N. Zubkov
\cite{Zubkov,Zubkov2} and S. Donkin \cite{Donkin,Donkin2}.

\medskip

\subsection{Invariant theory of generic matrices}

We will try to review some useful facts in the invariant theory of generic matrices.

Consider the algebra $\mathbb{A}_{n,s}$ of $s$-generated generic matrices of order $n$ over the ground field
$k$. Let $a_{\ell}=(a_{ij}^{\ell}), 1\leq i,j\leq n, 1\leq \ell \leq s$ be its generators. Let $R=k[a_{ij}^{\ell}]$ be
the ring of entries coefficients. Consider an action of matrices $M_n(k)$ on matrices in $R$ by conjugation, namely
$\varphi_B: B\mapsto M B M^{-1}$. It is well-known (refer to \cite{Procesi, Procesi2, Zubkov2}) that the
invariant function on this matrix can be expressed as a polynomial over traces $\tr(a_{i1},\dots,a_{is})$. Any
invariant on $\mathbb{A}_{n,s}$ is a polynomial of $\tr(a_{i1},\dots,a_{in})$. Note that the conjugation on $B$
induces an automorphism $\varphi_B$ of the ring $R$. Namely, $M(a_{ij})^{\ell}M^{-1}=(a'_{ij}{}^{\ell})$,
and $\varphi_B(M)$ of $R$ induces automorphism on $M_n(R)$. And for any $x\in \mathbb{A}_{n,s}$, we have
$$\varphi_B(x)=M x M^{-1}=\mathrm{Ad}_M(x).$$ Consider
$\varphi_B(x)=\mathrm{Ad}_M^{-1}\varphi_M(x)$. Then any element of the algebra of generic matrices is
invariant under $\varphi_M(x)$.

\medskip
When dealing with matrices in characteristic 0, it is useful to think that they form an algebra with a further unary
operation the $\textit{trace}$, $x\mapsto tr(x)$. One can formalize this as follows \cite{BPS1}:

\begin{definition}
An algebra with trace is an algebra equipped with an additional trace structure, that is a linear map $\tr: R\to R$
satisfying the following properties
$$\tr(ab) = \tr(ba), a\tr(b) = \tr(b)a, \tr(\tr(a)b) = \tr(a)\tr(b) \quad \text{for all}\; a,b\in R.$$
\end{definition}

\medskip

There is a well-known fact as follows.

\begin{theorem}
The algebra of generic matrices with trace is an algebra of concomitants, i.e. subalgebra of $M_n(R)$ is an invariant
under the action $\varphi_M(x)$.
\end{theorem}

This theorem was first proved by C. Procesi in \cite{Procesi} for the ground field $k$ of characteristic zero. If $k$
is a field of positive characteristic, we have to use not only traces, but also characteristic polynomials and their
linearization (refer to \cite{Donkin,Donkin2}). Relations between these invariants are discovered by C. Procesi
\cite{Procesi, Procesi2} for characteristic zero and A. N. Zubkov \cite{Zubkov, Zubkov2} for characteristic $p$.
C. de Concini and C. Procesi also generalized a characteristic free approach to invariant theory \cite{deConcProc}.

\medskip

Let us denote by $k_T\{X\}$ the algebra of generic matrices with traces. After above discussions, we have the
following proposition.

\begin{proposition}\label{ICwtr}
Let $n$ be a prime number, then the centralizer of $A\in k_T\{X\}$ is rationally closed in $k_T\{X\}$ and integrally
closed in $k_T\{X\}$.
\end{proposition}

\subsection{Centralizers are integrally closed}

Let $k\langle X\rangle$ be the free associative algebra as noted. Here we will prove the following theorem.

\begin{theorem}\label{integral}
The centralizer $C$ of non-trivial element $f$ in the free associative algebra is integrally closed.
\end{theorem}

Let $g,P,Q\in C:=C(f;F_z)$, and suppose $g Q^m=P^m$ for some positive integer $m$, i.e. in localization
$g=\frac{P^m}{Q^m}$. Then there exists $h\in C$, such that $h^m=g$. This means that the centralizer $C$ is
integral closed.

Consider the homomorphism $\pi$ from the free associative algebra $F_s$ to the algebra of generic matrices with
traces $k_T\{X\}$. Let us denote by $\bar{g}$ the image $\pi(g)$. Then we have following proposition.

\begin{proposition}
Consider the homomorphism $\pi: F_s\to k_T\{X\}$. Let the order of matrices be a prime number $p\gg 0$.
$\overline{g}=\pi(g), \overline{P}=\pi(P)$ and $\overline{Q}=\pi(Q)$. Then there exists $\overline{h}\in
k_T\{X\}$ such that
\begin{enumerate}
\item[1)] $\bar{h}^m=\overline{g}$,
\item[2)] $\bar{h}=\frac{\overline{P}}{\overline{Q}}$,
\item[3)] $\bar{h}\in \overline{C}$, where $\overline{C}=\pi(C)$.
\end{enumerate}
\end{proposition}

\begin{proof}
1) and 2) follows from Proposition \ref{ICwtr} that the algebra of generic matrices with traces of form is integral
closed. We need to prove 3). Note that all eigenvalues of $\bar{g}$ are pairwise different due to Proposition
\ref{propeigen}. So is $\bar{f}$. Hence $\bar{f},\bar{g}$ are diagonalizable and $\bar{h}$ can be diagonalized in
the same eigenvectors basis. Hence, by Proposition \ref{ProPcomm}, $\bar{h}$ commutes with $\bar{f}$, i.e.
$\bar{h}\in \overline{C}$.
\end{proof}

Now we have to prove that $\bar{h}$ in fact belongs to the algebra of generic matrices without trace. We use the
local isomorphism to get rid of traces.



\medskip


\begin{definition}[Local isomorphism]
Let $\mathbb{A}$ be an algebra with generators $a_1,\dots,a_s$ homogeneous respect this set of generators, and
let $\mathbb{A}'$ be an algebra with generators $a_1',\dots,a_s'$ homogeneous respect this set of generators. We
say that $\mathbb{A}$ and $\mathbb{A}'$ are locally $L$-isomorphic if there exist a linear map $\varphi:a_i\to
a_i'$ on the space of monomials of degree $\leq 2L$, and in this case for any two elements $b_1,b_2\in A$ with
highest term of degree $\leq L$, we have $$b_i=\sum_{j}M_{ij}(a_1,\dots,a_s),
b_i'=\sum_{j}M_{ij}(a_1',\dots,a_s'),$$ where $M_{ij}$ are monomials, and for $b=b_1\cdot b_2, b'=b_1'\cdot
b_2'$, we have $\varphi(b)=b'$.
\end{definition}

We need following lemmas, and propositions:

\begin{lemma}[Local isomorphism lemma]\label{lil}
For any $L$, if $s$ is  big enough prime, then the algebra of generic upper triangular matrices $\mathbb{U}_s$ is
locally $L$-isomorphic to the free associative algebra. Also reduction on the algebra of generic matrices of degree
$n$ provides an isomorphism up to degree $\leq 2s$.
\end{lemma}

Let us remind a well-known and useful fact.
\begin{proposition}
The trace of every element in $\mathbb{U}_s$ of any characteristic is zero.
\end{proposition}

In fact, we also proved

\begin{proposition}
If $n>n(L)$, then the algebra of generic matrices (without traces) is $L$-locally integrally closed.
\end{proposition}

\begin{lemma}
Consider the projection $\overline{\pi}$ of the algebra of generic matrices with trace to $\mathbb{U}_s$, sending
all traces to zero. Then we have $$\overline{\pi}(\overline{h})^m=\overline{\pi}(g).$$
\end{lemma}

\begin{proof}[Proof of Theorem \ref{integral}]
Let $p$ be a big enough prime number. For example, we can set $p\geq  2(\deg(f)+\deg(g)+\deg(P)+\deg(Q))$.
Because space of $k_T\{X\}$ of degree $\leq p$ is isomorphic to space of free associative algebra. We have
element $h$ corresponding to $\bar{h}$ up to this isomorphism. Due to local isomorphism, $h^m=g$, $h=P/Q$,
i.e. $hQ=P$.  Also we have $h$ commutes with $f$, i.e. $h\in C$.
\end{proof}


\subsection{Completion of the proof}\label{bergori}

From last two subsections, we have the following proposition:

\begin{proposition}\label{propratint}
Let $p$ be a big enough prime number, and $k\{X\}$  the algebra of generic matrices of order $p$. For any $A\in
k\{X\}$, the centralizer of $A$ is rationally closed and integrally closed in $k\{X\}$ over the center of $k\{X\}$.
\end{proposition}

In our previous paper \cite{zhang2017}, we establish that the centralizer in the algebra of generic matrices is a
commutative ring of transcendence degree one. According to Proposition \ref{propratint}, $C(A)$ is rationally
closed and integrally closed in $k\{X\}$. If $p$ is big enough, then $k\{X\}$ is $L$-locally integrally closed.

\medskip

Now we need one fact from the Bergman's paper \cite{berg1969}. Let $X$ be a totally ordered set, $W$ be the
free semigroup with identity $1$ on set $X$. We have the following lemma.
\begin{lemma}[Bergman]\label{lemword}
Let $u,v\in W\setminus\{1\}.$ If $u^{\infty}>v^{\infty}$, then we have
$u^{\infty}>(uv)^{\infty}>(vu)^{\infty}>v^{\infty}$.
\end{lemma}
\begin{proof}[Proof (Bergman)]
It suffices to show that the whole inequality is implied by $(uv)^{\infty}>(vu)^{\infty}$. Suppose
$(uv)^{\infty}>(vu)^{\infty}$, then we have following
$$(vu)^{\infty}=v(uv)^{\infty}>v(vu)^{\infty}=v^2(uv)^{\infty}>v^2(vu)^{\infty}=\cdots v^{\infty}.$$
Similarly, we obtain $(uv)^{\infty}<u^{\infty}.$
\end{proof}

Similarly, we also have inequalities with "$\geq$" replaced by "=" or "$\leq$".

\medskip

\begin{remark}
Similar constructions are used in \cite{BBL} for Burnside type problems or the height theorem of Shirshov.
\end{remark}

\medskip

Now let $R$ be the semigroup algebra on $W$ over field $k$, i.e. $R=F_s$ is the free associative algebra.
Consider $z\in \overline{W}$ be an infinite period word, and we denote $R_{(z)}$ be the $k$-subspace of $R$
generated by words $u$ such that $u=1$ or $u^{\infty}\leq z$. Let $I_{(z)}$ be the $k$-subspace spanned by
words $u$ such that $u\neq 1$ and $u^{\infty}<z$. Using Lemma \ref{lemword}, we can get that $R_{(z)}$ is a
subring of $R$ and $I_{{z}}$ is a two-sided ideal in $R_{(z)}$. It follows that $R_{(z)}/I_{(z)}$ will be
isomorphic to a polynomial ring $k[v]$.

\begin{proposition}[Bergman]\label{propberg}
If $C\neq k$ is a finitely generated subalgebra of $F_s$, then there is a homomorphism $f$ of $C$ in to polynomial
algebra over $k$ in one variable, such that $f(C)\neq k$.
\end{proposition}

\begin{proof}[Proof (Bergman)]
First let us totally order $X$. Let $G$ be a finite set of generators for $C$ and let $z$ be maximum over all
monomials $u\neq 1$ with nonzero coefficient in elements of $G$ of $u^{\infty}$. Then we have $G\subseteq
R_{(z)}$ and hence $C\subseteq R_{(z)}$, and the quotient map $f: R_{(z)}\to R_{(z)}/I_{(z)}\cong k[v]$ is
nontrivial on $C$.
\end{proof}

Now we can \textbf{finish the proof of Bergman's centralizer theorem}.

\begin{proof}
Consider homomorphism from the Proposition \ref{propberg}. Because $C$ is centralizer of $F_s$, it has
transcendence degree 1. Consider homomorphism $\rho$ send $C$ to the ring of polynomial. The homomorphism
has kernel zero, otherwise $\rho(C)$ will have smaller transcendence degree. Note that $C$ is integrally closed and
finitely generated, hence it can be embedded into polynomial ring of one indeterminate. Since $C$ is integrally
closed, it is isomorphic to polynomial ring of one indeterminate.

\medskip

Consider the set of system of $C_{\ell}$, $\ell$-generated subring of $C$ such that $C=\cup_{\ell}C_{\ell}$. Let
$\overline{C_{\ell}}$ be the integral closure of $C_{\ell}$. Consider set of embedding of $C_{\ell}$ to ring of
polynomial, then $\overline{C_{\ell}}$ are integral closure of those images, $\overline{C_{\ell}}=k[z_{\ell}]$,
where $z_{\ell}$ belongs to the integral closure of $C_{\ell}$. Consider sequence of $z_{\ell}$. Because
$k[z_{\ell}]\subseteq k[z_{\ell+1}]$, and degree of $z_{\ell+1}$ is strictly less than the degree of $z_{\ell}$.
Hence this sequence stabilizes for some element $x$. Then $k[z]$ is the needed centralizer.
\end{proof}

\subsection{On the rationality of subfields of generic matrices}   \label{Discourtions}

We will discuss some approaches to the following open problem.

\medskip

\begin{problem}
Consider the algebra of generic matrices $k\{X\}$ of order $s$. Consider $\Frac(k\{X\})$, and $K$ a subfield of
$\Frac(k\{X\})$ of transcendence degree one over the base field $k$. Is it true that $K$ is isomorphic to a rational
function over $k$, namely $K\cong k(t)$?
\end{problem}


\medskip

Let $k\{X\}$ be the algebra of generic matrices of a big enough prime order $s:=p$. Let $\Lambda$ be the
diagonal generic matrix $\Lambda=\diag{(\lambda_1,\dots,\lambda_s)}$ in $k\{X\}$, where transcendence
degrees satisfy in $\Trdeg k[\lambda_i]=1$. Let $N$ be another generic matrix, whose coefficients are algebraically
independent from $\lambda_1,\dots,\lambda_s$. It means that if $R$ is a ring of all coefficients of $N$, with
$\Trdeg(R)=s^2$, then $$\Trdeg R[\lambda_1,\dots,\lambda_s]=s^2+\Trdeg{k[\lambda_1,\dots,\lambda_s]}.$$

\begin{proposition}\label{ideal}
We consider the conjugation of generic matrices $\overline{f}$ and $\overline{g}$.
\begin{enumerate}
\item[a)] Let $k[f_{ij}]$ be a commutative ring , and $I=\langle f_{1i}\rangle\lhd k[f_{ij}]$ ($i>1$) be an ideal of
    $k[f_{ij}]$. Then we have $k[f_{11}]\cap I=0$.
\item[b)] Let $k[f_{ij},g_{ij}]$ be a commutative ring, and $J=\langle f_{1j},g_{1j}\rangle\lhd k[f_{ij},g_{ij}]$
    ($i,j>1$). For any algebraic function $P$ satisfies $P(f_{11},g_{11})=0$, which means $f,g$ are algebraic
    dependent on the $e_1$, then $k[f_{11},g_{11}]\cap J=0$.
\end{enumerate}
\end{proposition}

\begin{corollary}\label{coro32}
Let $\mathbb{A}$ be an algebra of generic matrices generated by $a_1,\dots,a_{s},a_{s+1}$. Let $f\in
k[a_1,\dots,a_s]$, $\varphi=a_{s+1}fa_{s+1}^{-1}$. Let $I=\langle \varphi_{1i}\rangle\lhd
k[a_{1},\dots,a_{s+1}]$. Then $k[\varphi_{11}]\cap I=0$.
\end{corollary}

\begin{proof}
Note that $f=\tau \Lambda \tau^{-1}$ for some $\tau$ and a diagonal matrix $\Lambda$ by proposition
\ref{ideal}. Then $\varphi=(a_{n+1}\tau)\Lambda(a_{n+1}\tau)^{-1}$ and we can treat $(a_{n+1}\tau)$ as a
generic matrix.
\end{proof}

\begin{theorem}\label{thmover}
Let $C:=C(f;F_n)$ be the centralizer ring of $f\in F_n\setminus k$. $\overline{C}$ is the reduction of generic
matrices, and $\overline{\overline{C}}$ is the reduction on first eigenvalue action. Then
$\overline{\overline{C}}\cong C$.
\end{theorem}

\begin{proof}
Let us recall that we already have $\overline{C}\cong C$ in \cite{zhang2017}. If we have $P(g_1,g_2)=0$, then
clearly $P(\lambda_1(g_1),\lambda_2(g_2))=0$ in the reduction on first eigenvalue action. Suppose
$\overline{\overline{P(g_1,g_2)}}=0$. Then $P(g_1,g_2)$ is an element of generic matrices with at least one zero
eigenvalue. Because minimal polynomial is irreducible, that implies that $P(g_1,g_2)=0$. It means any reduction
satisfying $\lambda_1$ satisfies completely. That what we want to prove.
\end{proof}

\medskip

Consider $\overline{C}$, for any $\overline{g}=(g_{ij})\in\overline{C}$. Investigate $g_{11}$. Suppose there is a
polynomial $P$ with coefficients in $k$, such that $P(f,g)=0$. We can make a proposition about intersection of the
ideals even sharper.

\medskip

Let $J=\langle f_{1j},g_{1j}\rangle$ ($j>1$) be an ideal of the commutative subalgebra $k[f_{ij},g_{ij}]$, then
$k[f_{11},g_{11}]\cap J=0$.

\medskip

From the discussion above and the theorem \ref{thmover}, we have the following proposition.

\begin{proposition}
$$k[f_{11},g_{11}] \mod J\cong k[f,g]$$
\end{proposition}

\begin{proof}
We have mod $J$ matrices from the following form:
$$\overline{\overline{f}}=\begin{pmatrix}
   \lambda_1 & 0 & \dots & 0 \\
   * & * & * & * \\
   * & * & * & * \\
   * & * & * & *
  \end{pmatrix},
  \overline{\overline{g}}=\begin{pmatrix}
   \lambda_2 & 0 & \dots & 0 \\
   * & * & * & * \\
   * & * & * & * \\
   * & * & * & *
  \end{pmatrix}$$
Then for any $H(\overline{\overline{f}},\overline{\overline{g}}) \mod J$, we have
$$\overline{\overline{f}}=\begin{pmatrix}
   H(\lambda_1,\lambda_2) & 0 & \dots & 0 \\
   * & * & * & * \\
   * & * & * & * \\
   * & * & * & *
  \end{pmatrix}$$
\end{proof}

Now we present an approach as follows. Consider $k(f,g)$. Let us extend the algebra of generic matrices by new
matrix $T$, independent from all others. Consider conjugation of $k(f,g)$ by $T$, $T k(f,g) T^{-1}$, and consider
$\tilde{f}=TfT^{-1}$ and $\tilde{g}=TgT^{-1}$. By Corollary \ref{coro32}, we have $$P(g_{11},f_{11})=0
\mod J.$$ On the other hand, we have $$k[f_{11},g_{11}]\cap J=0,$$ which means that $$P(f_{11},g_{11})=0
\mod J.$$ Put $f_{11}$ and $g_{11}$ be polynomial over commutative ring generated by all entries of $k[f,g]$
and $T$. Hence $\Frac (k(f,g))$ can be embedded into fractional field of rings of polynomials. According to
L\"{u}roth theorem, $\Frac(k(f,g))$ (hence $\Frac(C)$) is isomorphic to fields of rational functions in one variable.

\medskip

This will not guarantee rationality of our field, and there are counter examples in this situation. However, this
approach seems to be useful for highest term analysis.

\chapter{Automorphisms, augmentation topology, and approximation}
\label{Chapter2a} \lhead{Chapter 3. \emph{Automorphisms of $\Ind$-schemes, augmentation topology, and
approximation}}

\section{Introduction and main results}
This chapter is dedicated to the review of results of Kanel-Belov, Yu and Elishev on the geometry of the
$\Ind$-schemes
$$
\Aut \mathbb{K}[x_1,\ldots, x_n]
$$
and
$$
\Aut \mathbb{K}\langle x_1,\ldots, x_n\rangle
$$
of automorphisms of the polynomial algebra and the free associative algebra over an algebraically closed field, with
the number of generators $>2$. The inner character of $\Ind$ automorphisms of these $\Ind$-schemes, together
with the negative resolution of the automorphism group lifting problem, was established in \cite{KBYu}.

\subsection{Automorphisms of $K[x_1,\dots,x_n]$ and $K\langle x_1,\dots,x_n\rangle$}

Let $K$ be a field. The main objects of this study are the $K$-algebra automorphism groups $\Aut K[x_1,\dots,x_n]$ and \\
$\Aut K\langle x_1,\dots,x_n\rangle$ of the (commutative) polynomial algebra and the free associative algebra with
$n$ generators, respectively. The former is equivalent to the group of all polynomial one-to-one mappings of the
affine space $\mathbb{A}^n_{K}$. Both groups admit a representation as a colimit of algebraic sets of
automorphisms filtered by total degree (with morphisms in the direct system given by closed embeddings) which
turns them into topological spaces with Zariski topology compatible with the group structure. The two groups carry a
power series topology as well, since every automorphism $\varphi$ may be identified with the $n$-tuple
$(\varphi(x_1),\ldots,\varphi(x_n))$ of the images of generators. This topology plays an especially important role in
the applications, and it turns out -- as reflected in the main results of this study -- that approximation properties
arising from this topology agree well with properties of combinatorial nature.

$\Ind$-groups of polynomial automorphisms play a central part in the study of the Jacobian conjecture of O. Keller
as well as a number of problems of similar nature. One outstanding example is provided by a recent conjecture of
Kanel-Belov and Kontsevich (B-KKC), \cite{K-BK2,K-BK1}, which asks whether the group
$$
\Sympl({\mathbb C}^{2n})\subset \Aut(\mathbb{C}[x_1,\dots,x_{2n}])
$$
of complex polynomial automorphisms preserving the standard Poisson bracket $$\lbrace x_i,\;x_j\rbrace =
\delta_{i,n+j}-\delta_{i+n,j}$$ is isomorphic\footnote{In fact, the conjecture seeks to establish an isomorphism
 $\Sympl(K^{2n})\simeq \Aut(W_n(K))$ for any field $K$ of characteristic zero in a functorial manner.} to the group of automorphisms of the $n$-th Weyl algebra $W_n$
\begin{gather*}
W_n(\mathbb{C})=\mathbb{C}\langle x_1,\ldots,x_n,y_1,\ldots,y_n\rangle/I,\\
I= \left( x_ix_j-x_jx_i,\; y_iy_j-y_jy_i,\; y_ix_j-x_jy_i-\delta_{ij}\right).
\end{gather*}

The physical meaning of Kanel-Belov and Kontsevich conjecture is the invariance of the polynomial
symplectomorphism group of the phase space under the procedure of deformation quantization.

The B-KKC was conceived  during a successful search for a proof of stable equivalence of the Jacobian conjecture
and a well-known conjecture of Dixmier stating that $\Aut(W_n)=\End(W_n)$ over any field of characteristic zero.
In the papers \cite{K-BK2,K-BK1} a particular family of homomorphisms (in effect, monomorphisms)
$\Aut(W_n(\mathbb{C}))\rightarrow\Sympl({\mathbb C}^{2n})$ was constructed, and a natural question whether
those homomorphisms were in fact isomorphisms was raised. The aforementioned morphisms, independently studied
by Tsuchimoto to the same end, were in actuality defined as restrictions of morphisms of the saturated model of
Weyl algebra over an algebraically closed field of {\it positive} characteristic - an object which contains
$W_n(\mathbb{C})$ as a proper subalgebra. One of the defined morphisms turned out to have a particularly simple
form over the subgroup of the so-called tame automorphisms, and it was natural to assume that morphism was the
desired B-KK isomorphism (at least for the case of algebraically closed base field). Central to the construction is the
notion of infinitely large prime number (in the sense of hyperintegers), which arises as the sequence
$(p_m)_{m\in\mathbb{N}}$ of positive characteristics of finite fields comprising the saturated model. This leads to
the natural problem (\cite{K-BK1}):

\medskip
\noindent {\bf Problem.} Prove that the B-KK morphism is independent of the choice of the infinite prime
$(p_m)_{m\in\mathbb{N}}$.
\medskip

A general formulation of this question in the paper \cite{K-BK1} goes as follows:

For a commutative ring $R$ define
\begin{equation*}
R_{\infty}=\lim_{\rightarrow}\left( \prod_{p} R'\otimes \mathbb{Z}/p\mathbb{Z}\;/\;\bigoplus_{p} R'\otimes \mathbb{Z}/p\mathbb{Z}\right),
\end{equation*}
where the direct limit is taken over the filtered system of all finitely generated subrings $R'\subset R$ and the product
and the sum are taken over all primes $p$. This larger ring possesses a unique "nonstandard Frobenius"
endomorphism $\Fr:R_{\infty}\rightarrow R_{\infty}$ given by
\begin{equation*}
(a_p)_{\text{primes\;}p}\mapsto (a_p^p)_{\text{primes\;}p}.
\end{equation*}

The Kanel-Belov and Kontsevich construction returns a morphism
\begin{equation*}
\psi_R: \Aut (W_n(R))\rightarrow \Sympl R_{\infty}^{2n}
\end{equation*}

such that there exists a unique homomorphism
$$\phi_R: \Aut(W_n)(R)\rightarrow \Aut(P_n)(R_\infty)$$
obeying $\psi_R=\Fr_*\circ \phi_R$. Here $\Fr_*:\Aut(P_n)(R_\infty)\rightarrow \Aut(P_n)(R_\infty)$
 is the $\Ind$-group homomorphism induced by the Frobenius endomorphism of the coefficient ring, and $P_n$ is the commutative Poisson algebra, i.e. the polynomial algebra in $2n$ variables equipped with additional Poisson structure (so that $\Aut( P_n(R))$ is just $\Sympl( R^{2n})$ - the group of Poisson structure-preserving automorphisms).

\medskip
\noindent {\bf Question.}  {\it In the above formulation, does the image of $\phi_R$ belong to
$$\Aut(P_n)(i(R)\otimes {\mathbb Q})\,\,,$$
where $i:R\rightarrow R_\infty$ is the tautological inclusion? In other words, does there exist a unique
homomorphism
$$\phi_R^{can}:\Aut(P_n)(R)\rightarrow \Aut(P_n)(R\otimes {\mathbb Q})$$
 such that $\psi_R=\Fr_*\circ i_*\circ \phi_R^{can}$.}

\medskip

Comparing the two morphisms $\phi$ and $\varphi$ defined using two different free ultrafilters, we obtain a "loop"
element $\phi\varphi^{-1}$ of
 $\Aut_{\Ind}(\Aut(W_n))$,
(i.e. an automorphism which preserves the structure of infinite dimensional algebraic group). Describing this group
would provide a solution to this question.

\bigskip

In the spirit of the above we propose the following

\medskip
\noindent {\bf Conjecture.} {\it All automorphisms of the $\Ind$-group $\Sympl({\mathbb C}^{2n})$ are inner.}

A similar conjecture may be put forward for $\Aut(W_n(\mathbb{C}))$.

\medskip

Automorphism groups of Weyl algebras and their generalizations, as well as automorphisms of certain algebras of
vector fields, were studied in the works of Bavula \cite{Bav1, Bav2, Bav3}. Reduction to positive characteristic has
proven both fruitful and essential in the context of Weyl algebra. One of the precursors to the study of these algebras
in characteristic $p$ was the paper \cite{Bav8}.

\medskip

We are focused on the  investigation of the group $\Aut(\Aut(K[x_1,\dots,x_n]))$ and the corresponding
noncommutative (free associative algebra) case. This way of thinking has its roots in the realm of universal algebra
and universal algebraic geometry and was conceived in the pioneering work of Boris Plotkin. A more detailed
discussion can be found in \cite{KBLBerz}.



\

\noindent {\bf Wild automorphisms and the lifting problem.} In 2004, the celebrated Nagata conjecture over a field
$K$ of characteristic zero was proved by Shestakov and Umirbaev \cite{SU1, SU2} and a stronger version of the
conjecture was proved by Umirbaev and Yu \cite{UY}. Let $K$ be a field of characteristic zero. Every wild
$K[z]$-automorphism (wild $K[z]$-coordinate) of $K[z][x,y]$ is wild viewed as a $K$-automorphism
($K$-coordinate) of $K[x,y,z]$. In particular, the Nagata automorphism $(x-2y(y^2+xz)-(y^2+xz)^2z,
y+(y^2+xz)z, z)$ (Nagata coordinates $x-2y(y^2+xz)-(y^2+xz)^2z$ and $y+(y^2+xz)z$) are wild. In \cite{UY}, a
related question was raised:

\

\noindent {\bf The lifting problem.} \emph{Can an arbitrary wild automorphism (wild coordinate) of the polynomial
algebra $K[x,y,z]$ over a field $K$ be lifted to an automorphism (coordinate) of the free associative algebra
$K\langle x,y,z\rangle$?}

\medskip

In the paper \cite{BelovYuLifting}, based on the degree estimate \cite{MLY, YuYungChang}, it was proved that
any wild $z$-automorphism including the Nagata automorphism cannot be lifted as a $z$-automorphism (moreover,
in \cite{BKY} it is proved that every $z$-automorphism of $K\langle x, y,z\rangle$ is stably tame and becomes
tame after adding at most one variable). It means that  if every automorphism can be lifted, then it provides an
obstruction $z'$ to $z$-lifting and the question to estimate  such an obstruction is naturally raised.

In view of the above, we may ask the following:

\

 \noindent {\bf The automorphism group lifting problem.}
\emph{Is $\Aut(K[x_1,\dots,x_n])$ isomorphic to a subgroup of $\Aut(K\langle x_1,\dots,x_n\rangle)$ under the
natural abelianization?}

\medskip

The following examples show this problem is interesting and non-trivial.

\medskip
\noindent{\bf Example 1.} There is a surjective homomorphism (taking the absolute value) from $\mathbb C^*$
onto $\mathbb R^+$. But $\mathbb R^+$ is isomorphic to the subgroup $\mathbb R^+$ of $\mathbb C^*$ under
the homomorphism.

\medskip

\noindent{\bf Example 2.} There is a surjective homomorphism (taking the determinant) from $\text{GL}_n(\mathbb
R)$ onto $\mathbb R^*$. But obviously $\mathbb R^*$ is isomorphic to the subgroup $\mathbb R^*I_n$ of
$\text{GL}_n(\mathbb R)$.

\medskip

In this paper we prove that the automorphism group lifting problem has a negative answer.

The lifting problem  and the automorphism group lifting problem are closely related to the Kanel-Belov and
Kontsevich Conjecture (see Section \ref{KBConj}).

Consider a symplectomorphism $\varphi: x_i\to P_i,\; y_i\to Q_i$. It can be lifted to some automorphism
$\widehat{\varphi}$ of the quantized algebra $W_\hbar [[\hbar]]$: $$\widehat{\varphi}: x_i\to
P_i+P_i^1\hbar+\cdots+P_i^m\hbar^m;\ y_i\to Q_i+Q_i^1\hbar+\cdots+Q_i^m\hbar^m.$$ The point is to choose
a lift $\widehat{\varphi}$ in such a way that the degree of all $P_i^m, Q_i^m$ would be bounded. If that is true,
then the B-KKC follows.
\bigskip

\subsection{Main results}
The main results of this paper are as follows.

 \begin{theorem}        \label{ThAutAut}
Any $\Ind$-scheme automorphism $\varphi$ of $\NAut(K[x_1,\dots,x_n])$ for $n\ge 3$ is inner, i.e. is a
conjugation via some automorphism.
\end{theorem}

 \begin{theorem}        \label{ThAutAutFreeAssoc}
Any $\Ind$-scheme automorphism $\varphi$ of $\NAut(K\langle x_1,\dots,x_n\rangle)$ for $n\ge 3$ is semi-inner
(see definition \ref{DfSemi}).
\end{theorem}

$NAut$ denotes the group of {\em nice} automorphisms, i.e. automorphisms which can be approximated by tame
ones (definition \ref{DfNiceAprox}). In characteristic zero case every automorphism is nice.

For the group of  automorphisms of a semigroup a number of similar results on set-theoretical level was obtained
previously by Kanel-Belov, Lipyanski and Berzinsh \cite{BelovLiapiansk2, KBLBerz}. All these questions
(including $\Aut(\Aut)$ investigation) take root in the realm of Universal Algebraic Geometry and were proposed by
Boris Plotkin. Equivalence of two algebras having the same generalized identities and isomorphism of first order
means semi-inner properties of automorphisms (see \cite{BelovLiapiansk2, KBLBerz} for details).

\

\noindent {\bf Automorphisms of tame automorphism groups.} Regarding the tame automorphism group, something
can be done on the group- theoretic level. In the paper of H. Kraft and I. Stampfli \cite{KraftStampfli} the
automorphism group of the tame automorphism group of the polynomial algebra was thoroughly studied. In that
paper, conjugation of elementary automorphisms via translations played an important role. The results of our study
are different. We describe the group $\Aut(\TAut_0)$ of the group $\TAut_0$ of tame automorphisms preserving
the origin (i.e. taking the augmentation ideal onto an ideal which is a subset of the augmentation ideal). This is
technically more difficult, and will be universally and systematically done for both commutative (polynomial algebra)
case and noncommutative (free associative algebra) case. We observe a few problems in the shift conjugation
approach for the noncommutative (free associative algebra) case, as it was for commutative case in
\cite{KraftStampfli}. Any evaluation on a ground field element can return zero, for example in Lie polynomial
$[[x,y],z]$. Note that the calculations of $\Aut(\TAut_0)$ (resp. $\Aut_{\Ind}(\TAut_0)$, $\Aut_{\Ind}(\Aut_0)$)
imply also the same results for $\Aut(\TAut)$ (resp. $\Aut_{\Ind}(\TAut)$, $\Aut_{\Ind}(\Aut)$) according to the
approach of this article via stabilization by the torus action.

\begin{theorem}        \label{ThAutTAut}
Any  automorphism $\varphi$ of $\TAut_0(K[x_1,\dots,x_n])$ (in the group-theoretic sense) for $n\ge 3$ is inner,
i.e. is a conjugation via some automorphism.
\end{theorem}

 \begin{theorem}        \label{ThAutTAut0}
The group $\TAut_0(K[x_1,\dots,x_n])$ is generated by the automorphism $$x_1\to x_1+ x_2x_3,\; x_i\to x_i,
\;\;i\ne 1$$ and linear substitutions if $\Ch(K)\ne 2$ and $n>3$.
\end{theorem}

Let $G_N\subset \TAut(K[x_1,\dots,x_n])$, $E_N\subset \TAut(K\langle x_1,\dots,x_n\rangle)$ be tame
automorphism subgroups preserving the $N$-th power of the augmentation ideal.

 \begin{theorem}        \label{ThAutTAutOr}
Any  automorphism $\varphi$ of $G_N$ (in the group-theoretic sense) for $N\ge 3$ is inner, i.e. is given by a
conjugation via some automorphism.
\end{theorem}

\begin{definition}   \label{DfSemi}
An {\em anti-automorphism} $\Psi$ of a $K$-algebra $B$ is a vector space automorphism such that
$\Psi(ab)=\Psi(b)\Psi(a)$. For instance, transposition of matrices is an anti-automorphism. An anti-automorphism of
the free associative algebra $A$ is a {\em mirror anti-automorphism} if it sends $x_ix_j$ to $x_jx_i$ for some fixed
$i$ and $j$. If a mirror anti-automorphism $\theta$ acts identical on all generators $x_i$, then for any monomial
$x_{i_1}\cdots x_{i_k}$ we have
$$
\theta(x_{i_1}\cdots x_{i_k})=x_{i_k}\cdots x_{i_1}.
$$
Such an anti-automorphism will be generally referred to as \emph{the} mirror anti-automorphism.

An automorphism of $\Aut(A)$ is {\em semi-inner} if it can be expressed as a composition of an inner
automorphism and a conjugation by a mirror anti-automorphism.
\end{definition}

\begin{theorem}        \label{ThAutTAutFreass}
a) Any  automorphism $\varphi$ of $\TAut_0(K\langle
x_1,\dots,x_n\rangle)$ and also\\
 $\TAut(K\langle
x_1,\dots,x_n\rangle)$ (in the group-theoretic sense) for $n\ge 4$ is semi-inner, i.e. is a conjugation via some
automorphism and/or mirror anti-automorphism.

b) The same is true for $E_n$, $n\ge 4$.
\end{theorem}

The case of $\TAut(K\langle x,y,z\rangle)$ is substantially more difficult. We can treat it only on $\Ind$-scheme
level, but even then it is the most technical part of the paper (see section \ref{SbSc3VrbFreeAss}). For the
two-variable case a similar proposition is probably false.

\begin{theorem}            \label{ThTAss3Ind}
a) Let $\Ch(K)\ne 2$. Then $\Aut_{\Ind}(\TAut(K\langle
x,y,z\rangle))$ (resp. \\
$\Aut_{\Ind}(\TAut_0(K\langle x,y,z\rangle))$) is generated by conjugation by an automorphism or a mirror
anti-automorphism.

b) The same is true for $\Aut_{\Ind}(E_3)$.
\end{theorem}

By $\TAut$ we denote the tame automorphism group, $\Aut_{\Ind}$ is the group of $\Ind$-scheme automorphisms
(see section \ref{ScIndShme}).

Approximation allows us to formulate  the celebrated Jacobian conjecture for any characteristic.

\

{\bf Lifting of the automorphism groups.} In this article we prove that the automorphism group of
 polynomial algebra over an arbitrary field $K$ cannot be embedded
 into the automorphism group of free associative algebra induced by the  natural abelianization.

\begin{theorem}     \label{ThGroupLifting}
Let $K$ be an arbitrary field, $G=\Aut_0(K[x_1,\dots,x_n])$ and $n>2$. Then $G$ cannot be isomorphic to any
subgroup $H$ of $\Aut(K\langle x_1,\dots,x_n\rangle)$ induced by the natural abelianization. The same is true for
$\NAut(K[x_1,\dots,x_n])$.
\end{theorem}


\section{Varieties of automorphisms}

\subsection{Elementary and tame automorphisms}

Let $P$ be a polynomial that is independent of $x_i$ with $i$ fixed. An automorphism
$$x_i\to x_i+P,\; x_j\to x_j\;\; \mbox{for}\ i\ne j$$ is called {\em
elementary}. The group generated by linear automorphisms and elementary ones for all possible $P$ is called the
{\em tame automorphism group (or subgroup) $\TAut$} and elements of $\TAut$ are {\em tame automorphisms}.

\subsection{$\Ind$-schemes and $\Ind$-groups} \label{ScIndShme}

\begin{definition}
An {\it $\Ind$-variety} $M$ is the direct limit of algebraic varieties $M=\varinjlim \lbrace M_1\subseteq
M_2\cdots\rbrace$.
 An {\it $\Ind$-scheme} is an
$\Ind$-variety which is a group such that the group inversion is a morphism $M_i\rightarrow M_{j(i)}$ of algebraic
varieties, and the group multiplication induces a morphism from $M_i\times M_j$ to $M_{k(i,j)}$. A map $\varphi$
is a {\it morphism} of an $\Ind$-variety $M$ to an $\Ind$-variety $N$, if $\varphi(M_i)\subseteq N_{j(i)}$ and the
restriction $\varphi$ to $M_i$ is a morphism for all $i$. Monomorphisms, epimorphisms and isomorphisms are
defined similarly in a natural way.
\end{definition}

{\bf Example.} $M$ is the group of automorphisms of the affine space, and $M_j$ are the sets of all automorphisms
in $M$ with degree $\le j$.
\medskip

There is an interesting

\medskip
{\bf Problem.} \emph{Investigate growth functions of $\Ind$-varieties. For example, the dimension of varieties of
polynomial automorphisms of degree $\le n$.}

\medskip

Note that coincidence of growth functions of $\Aut(W_n(\mathbb{C}))$ and $\Sympl({\mathbb C}^{2n})$ would
imply the Kanel-Belov -- Kontsevich conjecture \cite{K-BK1}.

\begin{definition}             \label{DfAugm}
The ideal $I$ generated by variables $x_i$ is called the {\em augmentation ideal}.  For a fixed positive integer
$N>1$, the {\em augmentation subgroup $H_N$} is the group of all automorphisms $\varphi$ such that
$\varphi(x_i)\equiv x_i \mod I^N$. The larger group $\hat{H}_N\supset H_N$ is the group of automorphisms
whose linear part is scalar, and $\varphi(x_i)\equiv \lambda x_i \mod I^N$ ($\lambda$ does not depend on $i$).
We often say an arbitrary element of the group $\hat{H}_N$ is an automorphism that is homothety modulo (the
$N$-th power of) the augmentation ideal.
\end{definition}

\section{The Jacobian conjecture in any characteristic,
Kanel-Belov -- Kontsevich conjecture, and approximation}

\subsection{Approximation problems and Kanel-Belov -- Kontsevich Conjecture}    \label{KBConj}

Let us give formulation of the Kanel-Belov -- Kontsevich Conjecture:

\medskip
{\bf $B-KKC_n$}: $\Aut(W_n)\simeq\Sympl({\mathbb C}^{2n})$.
\medskip

A similar conjecture can be stated for endomorphisms

\medskip
{\bf $B-KKC_n$}: $\End(W_n)\simeq\Sympl\End({\mathbb C}^{2n})$.
\medskip

If the Jacobian conjecture $JC_{2n}$ is true, then the respective conjunctions over all $n$ of the two
conjectures are equivalent. 

It is natural to  approximate  automorphisms by tame ones. There exists such an approximation up to terms of any
order for polynomial automorphisms as well as Weyl algebra automorphisms, symplectomorphisms etc. However,
the naive approach fails.


It is known  that $\Aut(W_1)\equiv \Aut_1(K[x,y])$ where $\Aut_1$ stands for the subgroup of automorphisms of
Jacobian determinant one. However, considerations from \cite{Shafarevich} show that Lie algebra of the first group
is the algebra of derivations of $W_1$ and thus possesses no identities apart from the ones of the free Lie algebra,
another coincidence of the vector fields which diverge to zero, and has polynomial identities. These cannot be
isomorphic \cite{K-BK2,K-BK1}. In other words, this group has two coordinate system non-smooth with respect
to one another (but integral with respect to one another). One system is built from the coefficients of differential
operators in a fixed basis of generators, while its counterpart is provided by the coefficients of polynomials, which
are images of the basis
$\tilde{x}_i, \tilde{y}_i$. 

In the paper \cite{Shafarevich} functionals on ${\mathfrak m}/{\mathfrak m}^2$ were considered in order to define
the Lie algebra structure. In the spirit of that we have the following

\medskip
{\bf Conjecture.} The natural limit of ${\mathfrak m}/{\mathfrak m}^2$ is zero.
\medskip

It means that the definition of  the Lie algebra admits some sort of functoriality problem and it depends on the
presentation of (reducible) $\Ind$-scheme.

In his remarkable paper,  Yu. Bodnarchuk \cite{Bodnarchuk} established Theorem \ref{ThAutAut} by using
Shafarevich's results for the tame automorphism subgroup and for the case when the $\Ind$-scheme automorphism
is regular in the sense that it sends coordinate functions to coordinate functions. In this case the tame approximation
works (as well as for the symplectic case), and the corresponding method is similar to ours. We present it here in
order to make the text more self-contained, as well as for the purpose of tackling the noncommutative (that is, the
free associative algebra) case. Note that in general, for regular functions, if the Shafarevich-style approximation were
valid, then the Kanel-Belov -- Kontsevich conjecture would follow directly, which is absurd.


In the sequel, we do not assume regularity in the sense of \cite{Bodnarchuk} but only assume that the restriction of a
morphism on any subvariety is a morphism again. Note that morphisms of $\Ind$-schemes
$\Aut(W_n)\rightarrow\Sympl({\mathbb C}^{2n})$ have this property, but are not regular in the sense of
Bodnarchuk \cite{Bodnarchuk}.

We use the idea of singularity  which allows us to prove the augmentation subgroup structure preservation, so that
the approximation works in this case.

Consider the isomorphism $\Aut(W_1)\cong \Aut_1(K[x,y])$. It has a strange property. Let us add a small
parameter $t$. Then an element arbitrary close to zero with respect to $t^k$ does not go to zero  arbitrarily, so it is
impossible to make tame limit! There is a sequence of convergent product of elementary automorphisms, which is
not convergent under this isomorphism. Exactly the same situation happens  for $W_n$. These effects cause
problems in perturbative quantum field theory.

\medskip

\subsection{The Jacobian conjecture in any characteristic}
Recall that the Jacobian conjecture in characteristic zero states that any polynomial endomorphism
\begin{equation*}
\varphi:K^n\rightarrow K^n
\end{equation*}
with constant Jacobian is globally invertible.

A naive attempt to directly transfer this formulation to positive characteristic fails because of the counterexample
$x\mapsto x-x^p$ ($p=\Ch K$), whose Jacobian is everywhere $1$ but which is evidently not invertible.
Approximation provides a way to formulate a suitable generalization of the Jacobian conjecture to any characteristic
and put it in a framework of other questions.


\begin{definition}        \label{DfNiceAprox}
An endomorphism $\varphi\in\End(K[x_1,\dots,x_n])$ is {\em good} if\\
for any $m$ there exist $\psi_m\in\End(K[x_1,\dots,x_n])$ and\\
$\phi_m\in\Aut(K[x_1,\dots,x_n])$ such that
\begin{itemize}
    \item $\varphi=\psi_m\phi_m$
    \item $\psi_m(x_i)\equiv x_i\mod (x_1,\dots,x_n)^m$.
\end{itemize}

An automorphism $\varphi\in\Aut(K[x_1,\dots,x_n])$ is {\em nice} if for any $m$ there exist
$\psi_m\in\Aut(K[x_1,\dots,x_n])$ and $\phi_m\in\TAut(K[x_1,\dots,x_n])$ such that
\begin{itemize}
    \item $\varphi=\psi_m\phi_m$
    \item $\psi_m(x_i)\equiv x_i\mod (x_1,\dots,x_n)^m$, i.e. $\psi_m\in H_m$.
\end{itemize}
\end{definition}

Anick \cite{An} has shown that if $\Ch(K)=0$, any automorphism is nice. However, this is unclear in positive
characteristic.

\

{\bf Question.}\ {\it Is any automorphism over arbitrary field nice?}
\medskip

Ever good automorphism has Jacobian $1$, and all such automorphisms are good - and even nice - when
$\Ch(K)=0$. This observation allows for the following question to be considered a generalization of the Jacobian
conjecture to positive characteristic.

\medskip
{\bf The Jacobian conjecture in any characteristic:}\ {\it Is any good endomorphism over arbitrary field an
automorphism?}
\medskip

Similar notions can be formulated for the free associative algebra. That justifies the following

\medskip
{\bf Question.}\ {\it Is any automorphism of free associative algebra over arbitrary field nice?}
\medskip

\

{\bf Question (version of free associative positive characteristic case of JC).}\ {\it Is any good endomorphism of the
free associative algebra over arbitrary field an automorphism?}
\medskip

\subsection{Approximation for the  automorphism group of
 affine spaces}
\label{SbScAprox}

\medskip

Approximation is the most important tool utilized in this paper. In order to perform it, we have to prove that
$\varphi\in \Aut_{\Ind}(\Aut_0(K[x_1,\dots,x_n])$ preserves the structure of the augmentation subgroup.

The proof method utilized in theorems below works for commutative associative and free associative case. It is a
problem of considerable interest to develop similar statements for automorphisms of other associative algebras, such
as the commutative Poisson algebra (for which the $\Aut$ functor returns the group of polynomial
symplectomorphisms); however, the situation there is somewhat more difficult.

\bigskip

Suppose that $\varphi$ is an $\Ind$-automorphism (in either commutative of free associative case) such that it
stabilizes point-wise the set $T$ of automorphisms corresponding to the standard diagonal action of the maximal
torus (in the next section we will see that this implies that $\varphi$ also stabilizes every tame automorphism). The
following two continuity theorems, for the commutative and the free associative cases, respectively, constitute the
foundation of the approximation technique.

\begin{theorem}   \label{ThMainTechn}
Let $\varphi\in\Aut_{\Ind}(\Aut_0(K[x_1,\dots,x_n]))$ and let $H_N\subset \Aut_0(K[x_1,\dots,x_n])$ be the
subgroup of automorphisms which are identity modulo the ideal
 $(x_1,\dots,x_n)^N$ ($N>1$). Then $\varphi(H_N)\subseteq H_N$.

\end{theorem}

\begin{theorem}   \label{ThMainTechnFreeAss}
Let $\varphi\in\Aut_{\Ind}(\Aut_0(K\langle x_1,\dots,x_n\rangle ))$ and let $H_N$ be again the subgroup of
automorphisms which are identity modulo the ideal
 $(x_1,\dots,x_n)^N$. Then $\varphi(H_N)\subseteq H_N$.
\end{theorem}

\begin{corollary}
In both commutative and free associative cases under the assumptions above one has $\varphi= \Id$.
\end{corollary}

{\bf Proof.} Every automorphism can be approximated via the tame ones, i.e. for any $\psi$ and any $N$ there
exists a tame  automorphism $\psi'_N$ such that $\psi\psi'_N{}^{-1}\in H_N$.


The main point therefore is why $\varphi(H_N)\subseteq H_N$ whenever $\varphi$ is and $\Ind$-automorphism.

\medskip

{\bf Proof of Theorem \ref{ThMainTechn}.}

The method of proof is based upon the following useful fact from algebraic geometry:

\begin{lemma} \label{oneparamlem} Let

$$
\varphi: X\rightarrow Y
$$
be a morphism of affine varieties, and let $A(t)\subset X$ be a curve (or rather, a one-parameter family of points) in
$X$. Suppose that $A(t)$ does not tend to infinity as $t\rightarrow 0$. Then the image $\varphi A(t)$ under
$\varphi$ also does not tend to infinity as $t\rightarrow 0$.
\end{lemma}

The proof is straightforward and is left to the reader.

We now put the above fact to use. For $t>0$ let
$$
\hat{A}(t):\mathbb{A}^n_{K}\rightarrow \mathbb{A}^n_{K}
$$
be a one-parameter family of invertible linear transformations of the affine space preserving the origin. To that
corresponds a curve $A(t)\subset \Aut_0(K[x_1,\dots,x_n])$ of polynomial automorphisms whose points are linear
substitutions. Suppose that, as $t$ tends to zero, the $i$-th eigenvalue of $A(t)$ also tends to zero as $t^{k_i}$,
$k_i\in\mathbb{N}$. Such a family will always exist.

Suppose now that the degrees $\lbrace k_i,\;i=1,\ldots n\rbrace$ of singularity of eigenvalues at zero are such that
for every pair $(i,j)$, if $k_i\neq k_j$, then there exists a positive integer $m$ such that
$$
\text{either\;\;} k_im\leq k_j\;\;\text{or\;\;}k_jm\leq k_i.
$$

The largest such $m$ we will call the order of $A(t)$ at $t=0$. As $k_i$ are all set to be positive integer, the order
equals the integer part of $\frac{k_{\text{max}}}{k_{\text{min}}}$.

Let $M\in \Aut_0(K[x_1,\dots,x_n])$ be a polynomial automorphism.
\begin{lemma}    \label{Lm2} The curve $A(t)MA(t)^{-1}$ has no singularity at zero for any diagonalizable $A(t)$ of order $\leq N$ if and only if $M\in \hat{H}_N$, where $\hat{H}_N$ is the subgroup of automorphisms which are homothety modulo the
augmentation ideal.
\end{lemma}






{\bf Proof.} The  `If' part is elementary, for if $M\in \hat{H}_N$, the action of $A(t)MA(t)^{-1}$ upon any
generator $x_i$ (with $i$ fixed)\footnote{Without loss of generality we may assume that the coordinate functions
$x_i$ correspond to the principal axes of $\hat{A}(t)$.} is given by
\begin{eqnarray*}
A(t)MA(t)^{-1}(x_i)=\lambda x_i+t^{-k_i}\sum_{l_1+\cdots +l_n=N}a_{l_1\ldots l_n}t^{k_1l_1+\cdots+k_nl_n}x_1^{l_1}\cdots x_n^{l_n}+\\
+S_i(t,x_1,\ldots,x_n),
\end{eqnarray*}
where $\lambda$ is the homothety ratio of (the linear part of) $M$ and $S_i$ is polynomial in $x_1,\ldots,x_n$ of
total degree greater than $N$. Now, for any choice of $l_1,\ldots,l_n$ in the sum, the expression
$$
k_1l_1+\cdots+k_nl_n-k_i\geq k_{\text{min}}\sum l_j-k_i=k_{\text{min}}N-k_i\geq 0
$$
for every $i$, so whenever $t$ goes to zero, the coefficient will not blow up to infinity. Obviously the same argument
applies to higher-degree monomials within $S_i$.

\medskip

The other direction is slightly less elementary; assuming that $M\notin \hat{H}_N$, we need to show that there is a
curve $A(t)$ such that conjugation of $M$ by it produces a singularity at zero. We distinguish between two cases.

{\bf Case 1.} The linear part $\bar{M}$ of $M$ is not a scalar matrix. Then -- after a suitable basis change (see the
footnote) - it is not a diagonal matrix and has a non-zero entry in the position $(i,j)$. Consider a diagonal matrix
$A(t)=D(t)$ such that on all positions on the main diagonal except $j$-th it has $t^{k_i}$ and on $j$-th position it
has $t^{k_j}$. Then $D(t)\bar{M}D^{-1}(t)$ has $(i,j)$ entry with the coefficient $t^{k_i-k_j}$  and if
$k_j>k_i$ it has a singularity at $t=0$.

Let also $k_i<2k_j$. Then the non-linear part of $M$ does not produce singularities and cannot compensate the
singularity of the linear part.

{\bf Case 2.} The linear part $\bar{M}$ of $M$ is  a scalar matrix. Then conjugation cannot produce singularities in
the linear part and we as before are interested in the smallest non-linear term. Let $M\in H_N\backslash
H_{N+1}$. Performing a basis change if necessary, we may assume that
$$
\varphi(x_1)=\lambda x_1+\delta
x_2^N+S,
$$
where $S$ is a sum of monomials of degree $\ge N$ with coefficients in $K$.

Let $A(t)=D(t)$ be a diagonal matrix of the form $(t^{k_1}, t^{k_2},t^{k_1},\dots,t^{k_1})$ and let
$(N+1)\cdot k_2>k_1>N\cdot k_2$. Then in $A^{-1}MA$ the term $\delta x_2^N$ will be transformed into
$\delta x_2^N t^{Nk_2-k_1}$, and all other terms are multiplied by $t^{lk_2+sk_1-k_1}$ with $(l,s)\ne (1,0)$
and $ l, s>0$. In this case $lk_2+sk_1-k_1>0$ and we are done with the proof of  Lemma \ref{Lm2}.

\medskip

The next lemma is proved by direct computation. Recall that for $m>1$, the group $G_m$ is defined as the group
of all tame automorphisms preserving the $m$-th power of the augmentation ideal.

\begin{lemma}     \label{Lm3}

\

a) $[G_m,G_m]\subset H_m$, $m>2$. There exist elements \\
$\varphi\in H_{m+k-1}\backslash H_{m+k},\;\;\psi_1\in G_k,\;\;\psi_2\in G_m$, such that
$\varphi=[\psi_1,\psi_2]$.

b) $[H_m,H_k]\subset H_{m+k-1}$.

c) Let $\varphi\in G_m\backslash H_{m}$, $\psi\in H_k\backslash H_{k+1}$, $k>m$. Then $[\varphi,\psi]\in
H_k\backslash H_{k+1}$.

\end{lemma}

{\bf Proof.} a) Consider elementary automorphisms
\begin{gather*}
\psi_1: x_1\mapsto x_1+x_2^k,\;\;
 x_2\mapsto x_2,\;\;
 x_i\mapsto x_i,\;i> 2;
\end{gather*}
\begin{gather*}
\psi_2: x_1\mapsto x_1,\;\;
x_2\mapsto x_2+x_1^m, \;\;
x_i\mapsto x_i,\;i>2.
\end{gather*}

Set $\varphi=[\psi_1,\psi_2]=\psi_1^{-1}\psi_2^{-1}\psi_1\psi_2$. \\
Then
$$
\varphi: x_1\mapsto x_1-x_2^k+(x_2-(x_1-x_2^k)^m)^k,$$ $$x_2\mapsto
x_2-(x_1-x_2^k)^m+(x_1-x_2^k+(x_2-(x_1-x_2^k)^m)^k)^m, \;\; x_i\mapsto x_i,\;i>2.
$$
It is easy to see that if either $k$ or
$m$ is relatively prime with $\Ch(K)$, then not all terms of degree
$k+m-1$ vanish. Thus $\varphi\in H_{m+k-1}\backslash
H_{m+k}$.

Now suppose that $\Ch(K)\nmid m$, then obviously $m-1$ is relatively prime with $\Ch(K)$. Consider the
mappings
$$
\psi_1: x_1\mapsto x_1+x_2^k,\;\; x_2\mapsto x_2,\;\;x_i\mapsto x_i,\;i>2;
$$
$$
\psi_2: x_1\mapsto x_1,\;\; x_2\mapsto x_2+x_1^{m-1}x_3,\;\;x_i\mapsto x_i,\;i>2.
$$
Set again $\varphi'=[\psi_1,\psi_2]=\psi_1^{-1}\psi_2^{-1}\psi_1\psi_2$. Then $\varphi'$ acts as
\begin{gather*}
x_1\mapsto x_1-x_2^k+(x_2-(x_1-x_2^k)^{m-1}x_3)^k= \\
=x_1-k(x_1-x_2^k)^{m-1}x_2^{k-1}x_3+S, \\
x_2\mapsto x_2-(x_1-x_2^k)^{m-1}x_3+(x_1-x_2^k+(x_2-(x_1-x_2^k)^{m-1}x_3)^k)^{m-1}x_3, \\
x_i\mapsto x_i,\;\; i>2;
\end{gather*}
here $S$ stands for a sum of terms of degree $\ge m+k$. Again we see that $\varphi\in H_{m+k-1}\backslash
H_{m+k}$.

b) Let $$\psi_1: x_i\mapsto x_i+f_i;\ \psi_2: x_i\mapsto x_i+g_i,$$ for $i=1,\dots,n;$ here $f_i$ and $g_i$ do not
have monomials of degree less than or equal to $m$ and $k$, respectively. Then, modulo terms of degree $\ge
m+k$, we have $\psi_1\psi_2: x_i\mapsto x_i+f_i+g_i+\frac{\partial f_i}{\partial x_j} g_j$, so that modulo terms of
degree $\ge m+k-1$ we get $\psi_1\psi_2: x_i\mapsto x_i+f_i+g_i$ and $\psi_2\psi_1: x_i\to x_i+f_i+g_i$.
Therefore $[\psi_1,\psi_2]\in H_{m+k-1}$.

c) If $\varphi(I^m)\subseteq I^m$ and
$$
\psi: (x_1,\ldots,x_n)\mapsto (x_1+g_1,\ldots,x_n+g_n)
$$
is such that for some $i_0$ the polynomial $g_{i_0}$ contains a monomial of total degree $k$ (and all $g_i$ do not
contain monomials of total degree less than $k$), then, by evaluating the composition of automorphisms directly, one
sees that the commutator is given by
$$
[\varphi,\psi]: (x_1,\ldots,x_n)\mapsto (x_1+g_1+S_1,\ldots,x_n+g_n+S_n)
$$
with $S_i$ containing no monomials of total degree $<k+1$. Then the image of $x_{i_0}$ is $x_{i_0}$ modulo
polynomial of height $k$.

\begin{corollary}      \label{Cofinal}
Let $\Psi\in\Aut_{\Ind}(\NAut(K[x_1,\dots,x_n]))$. Then\\
$\Psi(G_n)=G_n$, $\Psi(H_n)=H_n$.
\end{corollary}

Corollary \ref{Cofinal} together with Proposition \ref{PrRips} of the next section imply Theorem
\ref{ThMainTechn}, for every nice automorphism, by definition, can be approximated by tame ones. Note that in
characteristic zero every automorphism is nice (Anick's theorem).

The proof of Corollary \ref{Cofinal} proceeds as follows (here for simplicity we put $\Char K = 0$, so that
$\NAut$ coincides with $\Aut$ thanks to Anick's theorem).

Let $\varphi$ be an $\Ind$-automorphism which stabilizes point-wise the standard action of the maximal torus.

1. We first note (and give a proof further along the text) that in this case $\varphi$ also stabilizes point-wise the set
of all tame automorphisms.

2. It follows from the singularity trick that:
$$
\varphi(\hat{H}_N)\subseteq \hat{H}_N
$$
(the reverse inclusion is also true due to the invertibility of $\varphi$). Namely, if $f=\varphi(g)$ is an automorphism
in $\varphi(\hat{H}_N)$ but not in $\hat{H}_N$ then there is a curve $A(t)$ of order $\leq N$ such that
$$
A(t)\circ f\circ A(t)^{-1}
$$
admits a singularity at $t=0$. But then
$$
\varphi^{-1}(A(t)\circ f\circ A(t)^{-1}) =A(t)\circ \varphi^{-1}(f)\circ A(t)^{-1}
$$
(this is thanks to the preservation of tame automorphisms) also admits a singularity at zero, which is a contradiction.

It is a fairly easy exercise to show that
$$
\varphi(\hat{H}_{N+1}\backslash \hat{H}_N) =\hat{H}_{N+1}\backslash \hat{H}_N
$$
for all $N$.

3. We now demonstrate that $\varphi(\hat{H}_N\backslash H_N) = \hat{H}_N\backslash H_N$ which together
with the preceding results will allow us to descend from homothety to identity modulo $N$.

A. Let $N>2$ first. Suppose $g\in\hat{H}_N\backslash H_N$. We take a tame automorphism $f$ which is given
by the sum of the identity map and a non-zero term of height two. Consider the automorphism
$$
g_f = f\circ g\circ f^{-1}.
$$
It is easy to see that $g_f\in \hat{H}_2$: as the linear part of $g$ is given by a scalar matrix not equal to the identity
matrix, the degree two component of $g\circ f^{-1}$ is proportional to the homothety ratio $\lambda\neq 1$,
therefore the composition with $f$ cannot compensate it.

On the other hand, if $\varphi(g)\in H_N$, i.e. the linear part of $\varphi(g)$ is the identity map, then the degree two
component of
$$
\varphi(g_f) = f\circ\varphi(g)\circ f^{-1}
$$
(this expression is again due to point-wise preservation of tame automorphisms) is equal to zero, which contradicts
$\varphi(\hat{H}_{N+1}\backslash \hat{H}_N) =\hat{H}_{N+1}\backslash \hat{H}_N$.

B. Let $N=2$. Suppose that $g\in\hat{H}_2\backslash H_2$ is a non-trivial homothety plus a term of height two.
The automorphism $g$ can be approximated by tame automorphisms, in particular there exists a tame automorphism
$\xi$ such that
$$
\xi\circ g
$$
is in $H_3$. From Case A it follows then that
$$
\varphi(\xi\circ g) = \xi\circ \varphi(g)
$$
is also in $H_3$. But then, since the linear part of $g$ is given by a non-trivial homothety, which means that $\xi$
scales it back to the identity matrix in order to approximate $g$ up to terms of height three, then the left action by
$\xi^{-1}$ reverses the scaling, so that the linear part of
$$
\xi^{-1}\circ\xi\circ \varphi(g) = \varphi(g)
$$
is given by a non-trivial homothety, which implies $\varphi(g) \in\hat{H}_2\backslash H_2$.

4. Finally, combining all of the results, we get $\varphi(H_N) = H_N$, $N>1$ as desired.

\subsection{Lifting of automorphism groups}

\subsubsection{Lifting of automorphisms from
$\Aut(K[x_1,\dots,x_n])$ to $\Aut(K\langle x_1,\dots,x_n\rangle)$}
\begin{definition} \label{defBialBir}
In the sequel, we call an action of the $n$-dimensional torus ${\mathbb T}^n$ on ${\mathbb K}\langle
x_1,\dots,x_n\rangle$ (the number of generators coincides with the dimension of the torus) \textbf{linearizable} if it
is conjugate to the standard diagonal action given by
$$
(\lambda_1,\ldots,\lambda_n)\;:\;(x_1,\ldots,x_n)\mapsto (\lambda_1x_1,\ldots,\lambda_nx_n).
$$
\end{definition}
The following result is a direct free associative analogue of a well-known theorem of Bia\l{}ynicki-Birula
\cite{BialBir2,BialBir1}. We will make frequent reference of the classical (commutative) case as well, which appears
as Theorem \ref{ThBialBir} in the text.

\begin{theorem}  \label{ThBialBirFree}
Any effective action of the $n$-torus on ${\mathbb K}\langle x_1,\dots,x_n\rangle$ is linearizable.
\end{theorem}

The proof is somewhat similar to that of Theorem \ref{ThBialBir}, with a few modifications. We provide it in
Chapter \ref{Chapter4}.

\medskip

As a corollary of the above theorem,  we get

\begin{proposition}
Let $T^n$ denote the standard torus action on $K[x_1,\ldots,x_n]$. Let $\widehat{T}^n$ denote its lifting to an
action on the free associative algebra $K\langle x_1,\ldots,x_n\rangle$. Then $\widehat{T}^n$ is also given by the
standard torus action.
\end{proposition}

{\bf Proof.} Consider the roots $\widehat{x_i}$ of this action. They are liftings of the coordinates $x_i$. We have
to prove that they generate the whole associative algebra.

Due to the reducibility of this action, all elements are product of eigenvalues of this action. Hence it is enough to
prove that eigenvalues of this action can be presented as a linear combination of this action. This can be done along
the lines of Bia\l{}ynicki-Birula \cite{BialBir1}. Note that all propositions of the previous section hold for the free
associative algebra. Proof of Theorem \ref{ThMainTechnFreeAss} is  similar. Hence we have the following

 \begin{theorem}        \label{ThAutAutFree}
Any $\Ind$-scheme automorphism $\varphi$ of $\Aut(K\langle x_1,\dots,x_n\rangle)$ for $n\ge 3$ is inner, i.e. is a
conjugation by some automorphism.
\end{theorem}

We therefore see that the group lifting (in the sense of isomorphism induced by the natural abelianization) implies the
analogue of Theorem \ref{ThMainTechn}.

This also implies that any automorphism group lifting, if exists, satisfies the approximation properties.

\begin{proposition}    \label{PrLiftStrct}
Suppose
$$
\Psi: \Aut(K[x_1,\dots,x_n])\rightarrow \Aut(K\langle
z_1,\dots,z_n\rangle)
$$
is a group homomorphism such that its composition with the natural map $\Aut(K\langle
z_1,\dots,z_n\rangle)\rightarrow \Aut(K[x_1,\dots,x_n])$ (induced by the projection $K\langle
z_1,\dots,z_n\rangle\rightarrow K[x_1,\dots,x_n]$) is the identity map. Then

\begin{enumerate}
    \item After a coordinate change $\Psi$ provides a correspondence between the standard torus actions
        $x_i\mapsto \lambda_ix_i$ and  $z_i\mapsto \lambda_iz_i$.
    \item Images of elementary automorphisms
    $$x_j\mapsto x_j,\; j\ne i,\;\; x_i\mapsto x_i+f(x_1,\dots,\widehat{x_i},\dots,x_n)$$
    are elementary automorphisms of the form
        $$z_j\mapsto z_j,\; j\ne i,\;\; z_i\mapsto z_i+f(z_1,\dots,\widehat{z_i},\dots,z_n).$$
    (Hence image of tame automorphism is tame automorphism).
    \item $\psi(H_n)=G_n$. Hence $\psi$ induces a map between the completion of the groups of
        $\Aut(K[x_1,\dots,x_n])$ and $\Aut(K\langle
z_1,\dots,z_n\rangle)$ with respect to the augmentation subgroup structure.
\end{enumerate}
\end{proposition}

\

\noindent {\bf Proof of Theorem \ref{ThGroupLifting}}

\noindent Any automorphism (including wild automorphisms such as the Nagata example) can be approximated by a
product of elementary automorphisms with respect to augmentation topology. In the case of the Nagata
automorphism corresponding to
$$\Aut(K\langle
x_1,\dots,x_n\rangle),$$ all such elementary automorphisms fix all
coordinates except  $x_1$ and $x_2$. Because of (2) and (3) of Proposition \ref{PrLiftStrct}, the lifted
automorphism would be an automorphism induced by an automorphism of $K\langle x_1,x_2,x_3\rangle$ fixing
$x_3$. However, it is impossible to lift the Nagata automorphism to such an automorphism due to the main result of
\cite{BelovYuLifting}. Therefore, Theorem \ref{ThGroupLifting} is proved.



\section{Automorphisms of the polynomial algebra and the approach of Bodnarchuk--Rips} \label{ScAutTameCoommN2}
 Let
$\Psi\in\Aut(\Aut(K[x_1,\dots,x_n]))$ (resp.
$\Aut(\TAut(K[x_1,\dots,x_n]))$,\\
$\Aut(\TAut_0(K[x_1,\dots,x_n]))$, $\Aut(\Aut_0(K[x_1,\dots,x_n]))$).

\subsection{Reduction to the case when $\Psi$ is identical on
$\SL_n$} We follow \cite{KraftStampfli} and \cite{Bodnarchuk} using the classical theorem of Bia\l{}ynicki-Birula
\cite{BialBir2,BialBir1}:

\begin{theorem}[Bia\l{}ynicki-Birula]  \label{ThBialBir}
Any effective action of torus ${\mathbb T}^n$ on ${\mathbb C}^n$ is linearizable (recall the definition
\ref{defBialBir}).
\end{theorem}

{\small

{\bf Remark.} An effective action of ${\mathbb T}^{n-1}$ on ${\mathbb C}^n$ is linearizable
\cite{BialBir1,BialBir2}. There is a conjecture whether any action of ${\mathbb T}^{n-2}$ on ${\mathbb C}^n$ is
linearizable, established for $n=3$. For codimension
 $>2$, there are positive-characteristic counterexamples \cite{Asanuma}.

\medskip

{\bf Remark.} Kraft and Stampfli \cite{KraftStampfli} proved (by considering periodic elements in $\mathbb T$)
that an effective action $T$ has the following property: if $\Psi\in\Aut(\Aut)$ is a group automorphism, then the
image of $T$ (as a subgroup of $\Aut$) under $\Psi$ is an algebraic group. In fact their proof is also applicable for
the free associative algebra case. We are going to use this result.

\medskip
}

Returning to the case of automorphisms $\varphi\in\Aut_{\Ind}\Aut$ preserving the $\Ind$-group structure, consider
now the standard action $x_i\mapsto\lambda_ix_i$ of the $n$-dimensional torus $\mathbb{T}\leftrightarrow
T^n\subset \Aut(\mathbb{C}[x_1,\ldots,x_n])$ on the affine space $\mathbb{C}^n$. Let $H$ be the image of
$T^n$ under $\varphi$. Then by Theorem \ref{ThBialBir} $H$ is conjugate to the standard torus $T^n$ via some
automorphism $\psi$. Composing $\varphi$ with this conjugation, we come to the case when $\varphi$ is the
identity on the maximal torus.  Then we have the following

\begin{corollary}
Without loss of generality, it is enough to prove Theorem \ref{ThAutAut} for the case when $\varphi|_{\mathbb
T}=\Id$.
\end{corollary}

Now we are in the situation when $\varphi$ preserves all linear mappings $x_i\mapsto \lambda_i x_i$.  We have to
prove that it is the identity.

\begin{proposition}[E. Rips, private communication] \label{PrRips}
Let $n>2$ and suppose $\varphi$ preserves the standard torus action on the commutative polynomial algebra. Then
$\varphi$ preserves all elementary transformations.
\end{proposition}

\begin{corollary}
Let $\varphi$ satisfy the conditions of Proposition \ref{PrRips}. Then $\varphi$ preserves all tame automorphisms.
\end{corollary}

{\bf Proof of Proposition \ref{PrRips}.} We state a few elementary lemmas.

\begin{lemma}  \label{LmTorComp}
Consider the diagonal action $T^1\subset T^n$ given by automorphisms: $\alpha: x_i\mapsto \alpha_i x_i$, $\beta:
x_i\mapsto \beta_i x_i$. Let $\psi: x_i\mapsto \sum_{i,J} a_{iJ}x^J,\; i=1,\dots,n$, where $J=(j_1,\dots,j_n)$ is the
multi-index, $x^J=x^{j_1}\cdots x^{j_n}$. Then

$$\alpha\circ\psi\circ\beta: x_i\mapsto \sum_{i,J} \alpha_i
a_{iJ}x^J\beta^J,$$

In particular,
$$\alpha\circ\psi\circ\alpha^{-1}: x_i\mapsto \sum_{i,J} \alpha_i
a_{iJ}x^J\alpha^{-J}.$$
\end{lemma}

Applying Lemma \ref{LmTorComp} and comparing the coefficients we get the following

\begin{lemma}  \label{LmLindiag}
Consider the diagonal $T^1$ action: $x_i\mapsto \lambda x_i$. Then the set of automorphisms commuting with this
action is exactly the set of linear automorphisms.
\end{lemma}

Similarly (using Lemma \ref{LmTorComp}) we obtain Lemmas \ref{LmMult1}, \ref{LmMult1n}, \ref{LmMult2}:

\begin{lemma}     \label{LmMult1}
a) Consider the following $ T^2$ action: $$x_1\mapsto \lambda\delta x_1,\; x_2\mapsto \lambda x_2,\; x_3\mapsto
\delta x_3,\; x_i\mapsto \lambda x_i,\; i>3.$$ Then the set $S$ of automorphisms commuting with this action is
generated by the following automorphisms: $$x_1\mapsto x_1+\beta x_2x_3,\; x_i\mapsto \varepsilon_i x_i,\; i>1,\;
(\beta,\varepsilon_i\in K).$$

b) Consider the following $T^{n-1}$ action:
$$x_1\mapsto
\lambda^I x_1,\; x_j\mapsto \lambda_j x_j,\; j>1\;(\lambda^I=\lambda_2^{i_2}\cdots\lambda_n^{i_n}).$$ Then the set $S$ of
automorphisms commuting with this action is generated by the following automorphisms: $$x_1\mapsto x_1+\beta
\prod_{j=2}^n x_j^{i_j},\ (\beta \in K).$$
\end{lemma}

{\bf Remark.} A similar statement for the free associative case is true, but one has to consider the set $\hat{S}$ of
automorphisms $x_1\mapsto x_1+h,\; x_i\mapsto \varepsilon_i x_i,\; i>1$, ($\varepsilon\in K$, and the polynomial
$h\in K\langle x_2,\dots,x_n\rangle$ has total degree $J$ - in the free associative case it is not just monomial
anymore).

\begin{corollary}      \label{CoActionToronElem}
Let $\varphi\in\Aut(\TAut(K[x_1,\dots,x_n]))$  stabilizing all elements from $\mathbb T$. Then $\varphi(S)=S$.
\end{corollary}

\begin{lemma}     \label{LmMult1n}
Consider the following $T^1$ action:
$$
x_1\mapsto \lambda^2 x_1,\;
x_i\mapsto \lambda x_i,\; i>1.
$$
Then the set $S$ of automorphisms commuting with this action is generated by the following automorphisms:
$$
x_1\mapsto x_1+\beta x_2^2,\; x_i\mapsto  \lambda_i x_i,\;
i>2,\; (\beta,\lambda_i\in K).
$$
\end{lemma}

\begin{lemma} \label{LmMult2}
Consider the set $S$ defined in the previous lemma. Then $[S,S]=\{uvu^{-1}v^{-1}\}$ consists of the following
automorphisms
$$x_1\mapsto x_1+\beta x_2x_3,\; x_2\mapsto x_2,\; x_3\mapsto  x_3,
\;(\beta\in K).$$
\end{lemma}

\begin{lemma}  \label{LmMult3}
Let $n\ge 3$. Consider the following set of automorphisms $$\psi_i: x_i\mapsto x_i+\beta_ix_{i+1}x_{i+2},\;
\beta_i\ne 0,\; x_k\mapsto x_k,\; k\ne i$$ for $i=1,\dots, n-1$. (Numeration is cyclic, so for example
$x_{n+1}=x_1$). Let $\beta_i\ne 0$ for all $i$. Then all of $\psi_i$ can be simultaneously conjugated by a torus
action to $$\psi_i': x_i\mapsto x_i+x_{i+1}x_{i+2},\; x_k\mapsto x_k,\; k\ne i$$ for $i=1,\dots, n$ in a unique
way.
\end{lemma}

{\bf Proof.} Let $\alpha: x_i\mapsto \alpha_i x_i$. Then by Lemma \ref{LmTorComp} we obtain
$$\alpha\circ\psi_i\circ\alpha^{-1}: x_i\mapsto x_i+\beta_i
x_{i+1}x_{i+2}\alpha_{i+1}^{-1}\alpha_{i+2}^{-1}\alpha_i$$ and
$$\alpha\circ\psi_i\circ\alpha^{-1}: x_k\mapsto x_k$$
for $k\ne i$.

Comparing the coefficients of the quadratic terms, we see that it is sufficient to solve the system:
$$\beta_i\alpha_{i+1}^{-1}\alpha_{i+2}^{-1}\alpha_i=1,\; i=1,\dots,n-1.$$
As $\beta_i\ne 0$ for all $i$, this system has a unique solution.

{\bf Remark.} In the free associative algebra case, instead of $\beta x_2x_3$ one has to consider $\beta
x_2x_3+\gamma x_3x_2$.

\subsection{The  lemma of Rips}

\begin{lemma}[E. Rips]  \label{LmRips}
Let $\Ch(K)\ne 2$, $|K|=\infty$. Linear transformations and $\psi_i'$ defined in Lemma \ref{LmMult3} generate
the whole tame automorphism group of $K[x_1,\ldots,x_n]$.
\end{lemma}

Proposition \ref{PrRips}  follows from Lemmas \ref{LmLindiag}, \ref{LmMult1}, \ref{LmMult1n}, \ref{LmMult2},
\ref{LmMult3}, \ref{LmRips}. Note that we have proved an analogue of Theorem \ref{ThAutAut} for tame
automorphisms.

\medskip
{\bf Proof of Lemma \ref{LmRips}.} Let $G$ be the group generated by elementary transformations as in Lemma
\ref{LmMult3}. We have to prove that is isomorphic to the tame automorphism subgroup fixing the augmentation
ideal. We are going to need some preliminaries.

\begin{lemma}  \label{LmR1}
Linear transformations of $K^3$ and $$\psi: x\mapsto x,\; y\mapsto y,\; z\mapsto z+xy$$ generate all mappings of
the form $$\phi_m^b(x,y,z): x\mapsto x,\; y\mapsto y,\; z\mapsto z+bx^m,\;\; b\in K.$$
\end{lemma}

{\bf Proof of Lemma \ref{LmR1}.} We proceed by induction. Suppose we have an automorphism
$$\phi^b_{m-1}(x,y,z): x\mapsto x,\; y\mapsto y,\; z\mapsto z+bx^{m-1}.$$  Conjugating by the linear
transformation ($z\mapsto y,\; y\mapsto z,\;x\mapsto x$), we obtain the automorphism
$$\phi_{m-1}^b(x,z,y): x\mapsto x,\; y\mapsto y+bx^{m-1},\; z\mapsto z.$$
Composing this on the right by $\psi$, we get the automorphism
$$\varphi(x,y,z): x\mapsto x,\; y\mapsto y+bx^{m-1},\; z\mapsto z+yx+x^m.$$ Note
that
$$\phi_{m-1}(x,y,z)^{-1}\circ\varphi(x,y,z): x\mapsto x,\; y\mapsto y,\;
z\mapsto z+xy+bx^m.$$
Now we see that
$$\psi^{-1}\phi_{m-1}(x,y,z)^{-1}\circ\varphi(x,y,z)=\phi^b_m$$ and
the lemma is proved.

\begin{corollary}  \label{CoLmR1}
Let $\Ch(K)\nmid n$ (in particular, $\Ch(K)\ne 0$) and $|K|=\infty$. Then $G$ contains all the transformations
$$z\mapsto z+bx^ky^l,\; y\mapsto y,\; x\mapsto x$$ such that $k+l=n$.
\end{corollary}

{\bf Proof.} For any invertible linear transformation
$$\varphi:
x\mapsto a_{11}x+a_{12}y,\; y\mapsto a_{21}x+a_{22}y,\; z\mapsto z; a_{ij}\in K$$
we have $$\varphi^{-1}\phi^b_{m}\varphi: x\mapsto x,\; y\mapsto y,\; z\mapsto z+b(a_{11}x+a_{12}y)^m.$$
Note that sums of such expressions contain all the terms of the form $bx^ky^l$. The corollary is proved.

\subsection{Generators of the tame automorphism group}

\begin{theorem}  \label{ThR2}
If $\Ch(K)\ne 2$ and $|K|=\infty$, then linear transformations and
$$\psi: x\mapsto x,\; y\mapsto y,\; z\mapsto z+xy$$ generate all mappings of the
form $$\alpha_m^b(x,y,z): x\mapsto x,\; y\mapsto y,\; z\mapsto z+byx^m,\;\; b\in K.$$
\end{theorem}

{\bf Proof of theorem \ref{ThR2}.} Observe that
$$\alpha=\beta\circ\phi^b_m(x,z,y): x\mapsto x+by^m,\; y\mapsto y+x+by^m,\;
z\mapsto z,$$ where $\beta: x\mapsto x,\; y\mapsto x+y,\;z\mapsto z$. Then
$$\gamma=\alpha^{-1}\psi\alpha: x\mapsto x,\; y\mapsto y,\; z\mapsto
z+xy+2bxy^m+by^{2m}.$$ Composing with $\psi^{-1}$ and
$\phi_{2m}^{2b}$ we get the desired $$\alpha_m^{2b}(x,y,z): x\mapsto x,\; y\mapsto y,\; z\mapsto z+2byx^m,\;\;
b\in K.$$

\begin{corollary}  \label{CoThR2}
Let $\Ch(K)\nmid n$ and $|K|=\infty$. Then $G$ contains all transformations of the form $$z\mapsto z+bx^ky^l,\;
y\mapsto y,\; x\mapsto x$$ such that $k=n+1$.
\end{corollary}

The {\bf proof} is similar to the proof of Corollary \ref{CoLmR1}. Note that either $n$ or $n+1$ is not a multiple
of $\Ch(K)$ so we have

\begin{lemma}  \label{LmR2n}
If $\Ch(K)\ne 2$ then linear transformations and $$\psi: x\mapsto x,\; y\mapsto y,\; z\mapsto z+xy$$ generate all
mappings of the form
$$\alpha_P: x\mapsto x,\; y\mapsto y,\; z\mapsto z+P(x,y),\;\; P(x,y)\in K[x,y].$$
\end{lemma}

We have proved Lemma \ref{LmRips} for the three variable case. In order to treat the case $n\ge 4$ we need one
more lemma.

\begin{lemma}  \label{LmR3Lin}
Let $M(\vec x)=a\prod x_i^{k_i},\;$ $a\in K,\; |K|=\infty$, $\Ch(K)\nmid k_i$ for at least one of $k_i$'s. Consider
the linear transformations denoted by $$f: x_i\mapsto y_i=\sum a_{ij}x_j,\; \det(a_{ij})\ne 0$$ and monomials
$M_f=M(\vec y)$. Then the linear span of $M_f$ for different $f$'s contains all homogenous polynomials of degree
$k=\sum k_i$ in $K[x_1,\dots,x_n]$.
\end{lemma}

{\bf Proof.} It is a direct consequence of the following fact. Let $S$ be a homogenous subspace of
$K[x_1,\dots,x_n]$ invariant  with respect to $GL_n$ of degree $m$. Then $S=S_{m/p^k}^{p^k}$,
$p=\Ch(K)$,\ $S_l$ is the space of all polynomials of degree $l$.

Lemma \ref{LmRips} follows from Lemma \ref{LmR3Lin} in a similar way as in the proofs of Corollaries
\ref{CoLmR1} and \ref{CoThR2}.

\

\subsection{$\Aut(\TAut)$ for general case}

 Now we consider the case when  $\Ch(K)$ is arbitrary, i.e. the remaining case
 $\Ch(K)=2$. Still $|K|=\infty$.
Although we are unable to prove the analogue of Proposition \ref{PrRips}, we can still play on the relations.

Let $$M=a\prod_{i=1}^{n-1} x_i^{k_i}$$ be a monomial, $a\in K$. For polynomial $P(x,y)\in K[x,y]$ we define
the elementary automorphism
$$\psi_P: x_i\mapsto x_i,\; i=1,\dots,n-1,\; x_n\mapsto
x_n+P(x_1,\dots,x_{n-1}).$$ We have $P=\sum M_j$ and $\psi_P$
naturally decomposes as a product of commuting $\psi_{M_j}$. Let $\Psi\in\Aut(\TAut(K[x,y,z]))$ stabilizing linear
mappings and $\phi$ (Automorphism $\phi$ defined in Lemma \ref{LmR1}). Then according to the corollary
\ref{CoActionToronElem} $\Psi(\psi_P)=\prod \Psi(\psi_{M_j})$. If $M=ax^n$ then due to Lemma \ref{LmR1}
$$\Psi(\psi_{M})=\psi_{M}.$$ We have to prove the same for other
type of monomials:

\begin{lemma}  \label{Le2nTestFreeAss}
Let $M$ be a monomial. Then
$$\Psi(\psi_{M})=\psi_{M}.$$
\end{lemma}

{\bf Proof.} Let $M=a\prod_{i=1}^{n-1} x_i^{k_i}$. Consider the automorphism $$\alpha: x_i\mapsto
x_i+x_1,\; i=2,\dots,n-1;\; x_1\mapsto x_1,\;x_n\mapsto x_n.$$ Then $$\alpha^{-1}\psi_M\alpha
=\psi_{x_1^{k_1}\prod_{i=2}^{n-1}(x_i+x_1)^{k_i}}=\psi_Q\psi_{ax_1^{\sum_{i=2}^{n-1} k_i}}.$$

Here the polynomial
$$Q=x_1^{k_1}\left(\prod_{i=2}^{n-1}(x_i+x_1)^{k_i}-ax_1^{\sum
k_i}\right).$$

It has the following form
$$Q=\sum_{i=2}^{n-1} N_i,$$
where $N_i$ are monomials such that none of them is proportional to a power of $x_1$.

According to Corollary \ref{CoActionToronElem}, $\Psi(\psi_M)=\psi_{bM}$ for some $b\in K$. We need only
to prove that $b=1$. Suppose the contrary, $b\ne 1$. Then
\begin{eqnarray*}
\Psi(\alpha^{-1}\psi_M\alpha)=\left(\prod_{[N_i,x_1]\ne 0}\Psi(\psi_{N_i})\right)\circ\Psi(\psi_{ax_1^{\sum_{i=2}^{n-1}
k_i}})=\\
\left(\prod_{[N_i,x_1]\ne
0}\psi_{b_iN_i}\right)\circ\psi_{ax_1^{\sum_{i=2}^{n-1} k_i}}
\end{eqnarray*}
for some $b_i\in K$.

On the other hand
$$\Psi(\alpha^{-1}\psi_M\alpha)=\alpha^{-1}\Psi(\psi_M)\alpha=\alpha^{-1}\psi_{bM}\alpha =
\left(\prod_{[N_i,x_1]\ne 0}
\psi_{bN_i}\right)\circ\psi_{ax_1^{\sum_{i=2}^{n-1} k_i}}$$

Comparing the factors $\psi_{ax_1^{\sum_{i=2}^{n-1} k_i}}$ and $\psi_{ax_1^{\sum_{i=2}^{n-1} k_i}}$ in
the last two products we get $b=1$. Lemma \ref{Le2nTestFreeAss} and hence Proposition \ref{PrRips} are
proved.

\section{The approach of Bodnarchuk--Rips to automorphisms of $\TAut(K\langle
x_1,\dots,x_n\rangle)$ ($n>2$)}

Now consider the free associative case. We treat the case $n>3$ on group-theoretic level and the case $n=3$ on
$\Ind$-scheme level. Note that if $n=2$ then $\Aut_0(K[x,y])=\TAut_0(K[x,y])\simeq \TAut_0(K\langle
x,y\rangle)=\Aut_0(K\langle x,y\rangle)$ and description of automorphism group of such objects is known due to J.
D\'{e}serti.

\subsection{The automorphisms of the tame automorphism \\
group of $K\langle x_1,\dots,x_n\rangle$, $n\ge 4$}

\begin{proposition}[E. Rips, private communication] \label{PrRipsFreeAss}
 Let $n>3$ and let $\varphi$ preserve the standard torus action
on the free associative  algebra $K\langle x_1,\dots,x_n\rangle$. Then $\varphi$ preserves all elementary
transformations.
\end{proposition}

\begin{corollary}
Let $\varphi$ satisfy the conditions of the proposition \ref{PrRipsFreeAss}. Then $\varphi$ preserves all tame
automorphisms.
\end{corollary}

For free associative algebras, we note that any automorphism preserving the torus action preserves also the
symmetric $$x_1\mapsto x_1+\beta(x_2x_3+x_3x_2),\; x_i\mapsto x_i,\; i>1$$ and the skew symmetric
$$x_1\mapsto x_1+\beta(x_2x_3-x_3x_2),\; x_i\mapsto x_i,\; i>1$$ elementary
automorphisms. The first property follows from Lemma \ref{LmMult1n}. The second one follows from the fact that
skew symmetric automorphisms commute with automorphisms of the following type $$x_2\mapsto x_2+x_3^2,\;
x_i\mapsto x_i,\; i\ne 2$$ and this property distinguishes them from elementary automorphisms of the form
$$x_1\mapsto x_1+\beta x_2x_3+\gamma x_3x_2,\; x_i\mapsto x_i,\; i>1.$$

Theorem \ref{ThAutAutFreeAssoc} follows from the fact that the forms $\beta x_2x_3+\gamma x_3x_2$
corresponding to general bilinear multiplication $$*_{\beta,\gamma}:(x_2,\;x_3)\to \beta x_2x_3+\gamma
x_3x_2$$ lead to associative multiplication if and only if $\beta=0$ or $\gamma=0$; the approximation also applies
(see section \ref{SbScAprox}).

\medskip
Suppose at first that $n=4$ and we are dealing with $K\langle x,y,z,t\rangle$.

\begin{proposition}       \label{PrGTame2}
The group $G$ containing all linear transformations and mappings
$$x\mapsto x,\; y\mapsto y,\; z\mapsto z+xy,\; t\mapsto t$$ contains also all
transformations of the form $$x\mapsto x,\; y\mapsto y,\; z\mapsto z+P(x,y),\; t\mapsto t.$$
\end{proposition}

{\bf Proof.} It is enough to prove that $G$ contains all transformations of the following form $$x\mapsto x,\;
y\mapsto y,\; z\mapsto z+aM,\; t\mapsto t,\;\; a\in K,$$ where $M$ is a monomial.

\medskip
{\bf Step 1.} Let $$M=a\prod_{i=1}^m x^{k_i}y^{l_i}\;\;\; \mbox{or}\;\;\; M=a\prod_{i=1}^m
y^{l_0}x^{k_i}y^{l_i}$$ or $$M=a\prod_{i=1}^m x^{k_i}y^{l_i}\;\;\; \mbox{or}\;\;\; M=a\prod_{i=1}^m
x^{k_i}y^{l_i}x^{k_{m+1}}.$$ Define the height of $M$, $H(M)$, to be the number of segments comprised of a
specific generator - such as $x^k$ - in the word $M$. (For instance, $H(a\prod_{i=1}^m
x^{k_i}y^{l_i}x^{k_{m+1}})=2m+1$.) Using induction on $H(M)$, one can reduce to the case when
$M=yx^k$. Let $M=M'x^k$ such that $H(M')<H(M)$. (Case when $M=M'y^l$ is obviously similar.) Let
$$\phi: x\to x,\; y\to y,\; z\to z+M',\; t\to t.$$
$$\alpha: x\to x,\; y\to y,\; z\to z,\; t\to t+zx^k.$$
Then $$\phi^{-1}\circ\alpha\circ\phi: x\to x,\; y\to y,\; z\to z,\; t\to t-M+ zx^k.$$

The automorphism $\phi^{-1}\circ\alpha\circ\phi$ is the composition of automorphisms $$\beta: x\to x,\; y\to y,\;
z\to z,\; t\to t-M$$ and
$$\gamma: x\to x,\; y\to y,\; z\to z, \;t\to t+zx^k.$$ Observe that $\beta$ is
conjugate to the automorphism $$\beta': x\to x,\; y\to y,\; z\to z-M,\; t\to t$$ by a linear automorphism $$x\to x,\;
y\to y,\; z\to t, \;t\to z.$$ Similarly, $\gamma$ is conjugate to the automorphism $$\gamma': x\to x,\; y\to y,\; z\to
z+yx^k,\; t\to t.$$ We have thus reduced to the case when $M=x^k$ or $M=yx^k$.

\medskip
{\bf Step 2.} Consider automorphisms $$\alpha: x\to x,\; y\to y+x^k,\; z\to z,\; t\to t$$ and $$\beta: x\to x,\; y\to y,\;
z\to z,\; t\to t+azy.$$ Then $$\alpha^{-1}\circ\beta\circ\alpha: x\to x,\; y\to y,\; z\to z,\; t\to t+azx^k+azy.$$ It is a
composition of the automorphism $$\gamma: x\to x,\; y\to y,\; z\to z,\; t\to t+azx^k$$ which is conjugate to the
needed automorphism $$\gamma': x\to x,\; y\to y,\; z\to z+yx^k,\; t\to t$$ and an automorphism $$\delta: x\to x,\;
y\to y,\; z\to z,\; t\to t+azy,$$ which is conjugate to the automorphism $$\delta': x\to x,\; y\to y,\; z\to z+axy,\; t\to
t$$ and then to the automorphism $$\delta'': x\to x,\; y\to y,\; z\to z+xy,\; t\to t$$ (using similarities). We have
reduced the problem to proving the statement $$G\ni\psi_M,\;\; M=x^k$$ for all $k$.

\medskip
{\bf Step 3.} Obtain the automorphism $$x\to x,\; y\to y+x^n,\; z\to z,\; t\to t.$$ This problem is similar to the
commutative case of $K[x_1,\dots,x_n]$ (cf. Section \ref{ScAutTameCoommN2}).

Proposition \ref{PrGTame2} is proved.

Returning to the general case $n\geq 4$, let us formulate the remark made after Lemma \ref{LmMult1} as follows:

\begin{lemma}     \label{LmMult1Ass}
 Consider the following $T^{n-1}$ action: $$x_1\to \lambda^I
x_1,\; x_j\to \lambda_j x_j,\;\; j>1; \;\;\lambda^{I}=\lambda_2^{i_2}\cdots\lambda_n^{i_n}.$$ Then the set $S$
of automorphisms commuting with this action is generated by the following automorphisms: $$x_1\to x_1+H,\; x_i\to
x_i;\;\; i>1,$$ where $H$ is any homogenous polynomial of total degree $i_2+\cdots+i_n$.
\end{lemma}

Proposition \ref{PrGTame2} and Lemma \ref{LmMult1Ass} imply

\begin{corollary}                  \label{CoLmMult1Ass}
Let $\Psi\in\Aut(\TAut_0(K\langle x_1,\dots,x_n\rangle))$ stabilize all elements of torus and linear automorphisms,
$$\phi_P: x_n\to x_n+P(x_1,\dots,x_{n-1}),\; x_i\to x_i,\;
i=1,\dots,n-1.$$ Let $P=\sum_IP_I$, where $P_I$ is the homogenous component
of $P$ of multi-degree $I$. Then

a) $\Psi(\phi_P): x_n\to x_n+P^\Psi(x_1,\dots,x_{n-1}),\; x_i\to x_i,\; i=1,\dots,n-1$.

b) $P^\Psi=\sum_I P_I^\Psi$; here $P_I^\Psi$ is homogenous  of multi-degree $I$.

c) If $I$ has positive degree with respect to one or two variables, then $P_I^\Psi=P_I$.
\end{corollary}

Let $\Psi\in\Aut(\TAut_0(K\langle x_1,\dots,x_n\rangle))$ stabilize all elements of torus and linear automorphisms,
$$\phi: x_n\to x_n+P(x_1,\dots,x_{n-1}),\; x_i\to x_i,\;
i=1,\dots,n-1.$$

Let $\varphi_Q: x_1\to x_1,\;x_2\to x_2,\; x_i\to x_i+Q_i(x_1,x_2),\; i=3,\dots,n-1,\; x_n\to x_n;\;$
$Q=(Q_3,\dots,Q_{n-1})$. Then $\Psi(\varphi_Q)=\varphi_Q$ by Proposition \ref{PrGTame2}.

\begin{lemma}        \label{LmQPsi12}
a) $\varphi_Q^{-1}\circ\phi_P\circ\varphi_Q =\phi_{P_Q}$, where
$$P_Q(x_1,\dots,x_{n-1})=P(x_1,x_2,x_3+Q_3(x_1,x_2),\dots,x_{n-1}+Q_{n-1}(x_1,x_2)).$$

b) Let $P_Q=P_Q^{(1)}+P_Q^{(2)}$, $P_Q^{(1)}$ consist of all terms containing one of the variables
$x_3,\dots,x_{n-1}$, and let $P_Q^{(1)}$ consist of all terms containing just $x_1$ and $x_{2}$. Then
$$P^\Psi_Q=P_Q^\Psi=P_Q^{(1)\Psi}+P_Q^{(2)\Psi}=P_Q^{(1)\Psi}+P_Q^{(2)}$$.
\end{lemma}

\begin{lemma}        \label{LmQPSub}
If $P_Q^{(2)}=R_Q^{(2)}$ for all $Q$ then $P=R$.
\end{lemma}

{\bf Proof.} It is enough to prove that if $P\ne 0$ then $P^{(2)}_Q\ne 0$ for appropriate
$Q=(Q_3,\dots,Q_{n-1})$. Let $m=\deg(P), \;Q_i=x_1^{2^{i+1}m}x_2^{2^{i+1}m}$. Let $\hat{P}$ be the
highest-degree component of $P$, then $\hat{P}(x_1,x_2,Q_3,\dots,Q_{n-1})$ is the highest-degree component
of $P^{(2)}_Q$. It is enough to prove that
$$\hat{P}(x_1,x_2,Q_3,\dots,Q_{n-1})\ne 0.$$ Let $x_1\prec x_2\prec
x_2\prec\cdots\prec x_{n-1}$ be the standard lexicographic order. Consider the lexicographically minimal term
$M$ of $\hat{P}$. It is easy to see that the term
$$M|_{Q_i\to x_i},\; \;i=3,\;n-1$$ cannot cancel with any other term
$$N|_{Q_i\to x_i},\;\; i=3,\;n-1$$ of
$\hat{P}(x_1,x_2,Q_3,\dots,Q_{n-1})$. Therefore $\hat{P}(x_1,x_2,Q_3,\dots,Q_{n-1})\ne 0$.

Lemmas \ref{LmQPsi12} and \ref{LmQPSub} imply

\begin{corollary}                  \label{CoSemifinal}
Let $\Psi\in\Aut(\TAut_0(K\langle x_1,\dots,x_n\rangle))$ stabilize all elements of torus and linear automorphisms.
Then $P^\Psi=P$, and $\Psi$ stabilizes all elementary automorphisms and therefore the entire group
$\TAut_0(K\langle x_1,\dots,x_n\rangle)$.
\end{corollary}

We obtain the following

\begin{proposition}            \label{PrFinalNGEQ4}
Let $n\ge 4$ and let $\Psi\in\Aut(\TAut_0(K\langle x_1,\dots,x_n\rangle))$ stabilize all elements of torus and
linear automorphisms. 
Then either $\Psi=\Id$ or $\Psi$ acts as conjugation by the mirror anti-automorphism.
\end{proposition}

Let $n\ge 4$. Let $\Psi\in\Aut(\TAut_0(K\langle x_1,\dots,x_n\rangle))$ stabilize all elements of torus and linear
automorphisms. Denote by $EL$ an elementary automorphism
$$
EL:x_1\to x_1,\;\ldots,\;x_{n-1}\to x_{n-1},\;x_n\to x_n+x_1x_2
$$
(all other elementary automorphisms of this form, i.e. $x_k\to x_k+x_ix_j,\;x_l\to x_l$ for $l\neq k$ and $k\neq
i,\;k\neq j,\;i\neq j$, are conjugate to one another by permutations of generators).

We have to prove that $\Psi(EL)=EL$ or $\Psi(EL): x_i\to x_i;\; i=1,\dots,x_{n-1},\; x_n\to x_n+x_2x_1$. The
latter corresponds to $\Psi$ being the conjugation with the mirror anti-automorphism of $K\langle
x_1,\dots,x_n\rangle$.

Define for some $a,b\in K$ $$x*_{a,b}y=axy+byx.$$

Then, in any of the above two cases, $$\Psi(EL): x_i\to x_i;\; i=1,\dots,x_{n-1},\; x_n\to x_n+x_1*_{a,b}x_2$$
for some $a,b$.

The following lemma is elementary:

\begin{lemma}                     \label{LmAssoc}
The operation $*=*_{a,b}$ is associative if and only if  $ab=0$. 
\end{lemma}

The associator of $x,\;y,\;\text{and}\;z$ is given by
\begin{eqnarray*}
\lbrace x,y,z\rbrace_{*}\equiv (x*y)*z-x*(y*z)=\\
ab(zx-xz)y+aby(xz-zx)=ab[y,[x,z]].
\end{eqnarray*}


Now we are ready to prove Proposition \ref{PrFinalNGEQ4}. For simplicity we treat only the case $n=4$ -- the
general case is dealt with analogously. Consider the automorphisms $$\alpha: x\to x,\; y\to y,\; z\to z+xy,\; t\to t,$$
$$\beta: x\to x,\; y\to y,\; z\to z,\; t\to t+xz,$$ $$h:x\to x,\; y\to y,\;
z\to z,\; t\to t-xz.$$ (Manifestly $h=\beta^{-1}$.) Then $$\gamma=h\alpha^{-1}\beta\alpha=[\beta,\alpha]: x\to x,\;
y\to y,\; z\to z,\; t\to t-x^2y.$$ Note that $\alpha$ is conjugate to $\beta$ via a generator permutation
$$
\kappa:x\to x,\;y\to z,\;z\to t,\;t\to y,\;\;\kappa\circ\alpha\circ\kappa^{-1}=\beta
$$
and $$\Psi(\gamma): x\to x,\; y\to y,\; z\to z,\; t\to t-x*(x*y).$$

Let $$\delta: x\to x,\; y\to y,\; z\to z+x^2,\; t\to t,$$
$$\epsilon: x\to
x,\; y\to y,\; z\to z,\; t\to t+zy.$$ Let
$\gamma'=\epsilon^{-1}\delta^{-1}\epsilon\delta$. Then $$\gamma': x\to x,\; y\to y,\; z\to z,\; t\to t-x^2y.$$ On the
other hand we have
$$\varepsilon=\Psi(\epsilon^{-1}\delta^{-1}\epsilon\delta): x\to x,\;
y\to y, \;z\to z,\; t\to t-(x^2)*y.$$ We also have
$\gamma=\gamma'$. Equality $\Psi(\gamma)=\Psi(\gamma')$ is equivalent to the equality $x*(x*y)=x^2*y$. This
implies $x*y=xy$ and we are done.

\subsection{The group $\Aut_{\Ind}(\TAut(K\langle x,y,z\rangle))$}  \label{SbSc3VrbFreeAss}
This is the most technically loaded part of the present study. At the moment we are unable to accomplish the
objective of describing the entire group $\Aut\TAut(K\langle x,y,z\rangle)$. In this section we will determine only its
subgroup
 $\Aut_{\Ind}\TAut_0(K\langle
x,y,z\rangle)$, i.e. the group of $\Ind$-scheme automorphisms, and prove Theorem \ref{ThTAss3Ind}. We use the
approximation results of Section \ref{SbScAprox}. In what follows we suppose  that $\Ch(K)\ne 2$. As in the
preceding chapter, $\{x,y,z\}_*$ denotes the associator of $x,y,z$ with respect to a fixed binary linear operation
$*$, i.e.
$$\{x,y,z\}_*=(x*y)*z-x*(y*z).$$ 

\begin{proposition}                      \label{LeTAss3Ind}
Let $\Psi\in\Aut_{\Ind}(\TAut_0(K\langle x,y,z\rangle))$ stabilize all linear automorphisms. Let $$\phi: x\to x,\; y\to
y,\; z\to z+xy.$$ Then either $$\Psi(\phi): x\to x,\; y\to y,\; z\to z+axy$$ or
$$\Psi(\phi): x\to x,\; y\to y,\; z\to z+byx$$
for some $a,b\in K$.
\end{proposition}

{\bf Proof.}\ Consider the automorphism $$\phi: x\to x,\; y\to y,\; z\to z+xy.$$ Then $$\Psi(\phi): x\to x,\; y\to y,\;
z\to z+x*y,$$ where $x*y=axy+byx$. Let $a\neq 0$. We can make the star product $*=*_{a,b}$ into
$x*y=xy+\lambda yx$ by conjugation with the mirror anti-automorphism and appropriate linear substitution. We
therefore need to prove that $\lambda=0$, which implies $\Psi(\phi)=\phi$.

The following two lemmas are proved by straightforward computation.

\begin{lemma}     \label{LmAssociator}
Let $A=K\langle x,y,z\rangle$. Let $f*g=fg+\lambda fg$. Then $\{f,g,h\}_*=\lambda[g,[f,h]]$.

In particular $\{f,g,f\}_*=0$,
$f*(f*g)-(f*f)*g=-\{f,f,g\}_*=\lambda [f,[f,g]]$,\\
$(g*f)*f-g*(f*f)=\{g,f,f\}_*=\lambda [f,[f,g]]$.
\end{lemma}


\begin{lemma}     \label{LmxyyzCommutant}
Let $\varphi_1: x\to x+yz,\; y\to y,\; z\to z$; $\varphi_2: x\to x,\; y\to y,\; z\to z+yx$;
$\varphi=\varphi_2^{-1}\varphi_1^{-1}\varphi_2\varphi_1$. Then modulo terms of order $\ge 4$ we have:

$$\varphi: x\to x+y^2x,\; y\to y,\; z\to z-y^2z$$ and $$\Psi(\varphi):
x\to x+y*(y*x),\; y\to y,\; z\to z-y*(y*z).$$
\end{lemma}


\begin{lemma}     \label{LmElemStepxxy}
a) Let  $\phi_l: x\to x,\; y\to y,\; z\to z+y^2x$.  Then $$\Psi(\phi_l): x\to x,\; y\to y,\; z\to z+y*(y*x).$$

b) Let  $\phi_r: x\to x,\; y\to y,\; z\to z+xy^2$.  Then $$\Psi(\phi_r): x\to x,\; y\to y,\; z\to z+(x*y)*y.$$

\end{lemma}

{\bf Proof.} According to the results of the previous section we have
$$\Psi(\phi_l): x\to x,\; y\to y,\; z\to z+P(y,x)$$ where $P(y,x)$ is
homogenous of degree 2 with respect to $y$ and degree 1 with respect to $x$. We have to prove that
$H(y,x)=P(y,x)-y*(y*x)=0$.

Let $\tau: x\to z,\; y\to y,\; z\to x;\; \tau=\tau^{-1},\;\; \phi'=\tau\phi_l\tau^{-1}: x\to x+y^2z,\; y\to y,\; z\to z$. Then
$\Psi(\phi_l'): x\to x+P(y,z),\; y\to y,\; z\to z$.

Let $\phi_l''=\phi_l\phi_l': x\to x+P(y,z),\; y\to y,\; z\to z+P(y,x)$ modulo terms of degree $\ge 4$.

Let $\tau: x\to x-z,\; y\to y,\; z\to z$ and let $\varphi_2$, $\varphi$~ be the automorphisms described in Lemma
\ref{LmxyyzCommutant}.

Then $$T=\tau^{-1}\phi_l^{-1}\tau\phi_l'':x\to x,\;y\to y,\;z\to z$$ modulo terms of order $\ge 4$.

On the other hand $$\Psi(T): x\to x+H(y,z)-H(y,x),\; y\to y,\; z\to z+P$$ modulo terms of order $\ge 4$. Because
$\deg_y(H(y,x)=2,\;\deg_x(H(y,x))=1$ we get $H=0$.

 Proof of  b) is  similar.

\begin{lemma}     \label{LmElemSquare}
a)  Let  $$\psi_1: x\to x+y^2,\; y\to y,\; z\to z;\;\; \psi_2: x\to x,\; y\to y,\; z\to z+x^2.$$  Then
$$[\psi_1,\psi_2]=\psi_2^{-1}\psi_1^{-1}\psi_2\psi_1: x\to x,\; y\to
y,\; z\to z+y^2x+xy^2,$$ $$\Psi([\psi_1,\psi_2]): x\to x,\;y\to y,\;
z\to z+(y*y)*x+x*(y*y).$$

b) $$\phi_l^{-1}\phi_r^{-1}[\psi_1,\psi_2]: x\to x,\;y\to y,\;z\to z$$ modulo terms of order $\ge 4$ but
\begin{gather*}
\Psi\left(\phi_l^{-1}\phi_r^{-1}[\psi_1,\psi_2]\right): x\to x,\;y\to y,\\
z\to
z+(y*y)*x+x*(y*y)-(x*y)*y-y*(y*x)=\\
=z+4\lambda [x[x,y]]
\end{gather*}
modulo terms of order $\ge 4$.
\end{lemma}

{\bf Proof.} a) can be obtained by direct computation. b) follows from a) and the lemma \ref{LmAssociator}.


Proposition \ref{LeTAss3Ind} follows from Lemma \ref{LmElemSquare}.

We need a few auxiliary lemmas. The first one is an analogue of the hiking procedure from
\cite{BelovUzyRowenSerdicaBachtur,BelovIAN}.

\begin{lemma}                  \label{LeHiking}
Let $K$ be  algebraically closed, and let \ $n_1,\dots,n_m$ be positive integers. Then there exist  $k_1,\dots,k_s\in
{\mathbb Z}$ and $\lambda_1,\dots,\lambda_s\in K$ such that

\begin{itemize}
    \item $\sum k_i=1$ modulo $\Ch(K)$ (if $\Ch(K)=0$ then $\sum k_i=1$).
    \item $\sum_i k_i^{n_j}\lambda_i=0$ for all $j=1,\dots,m$.
\end{itemize}
\end{lemma}

For $\lambda\in K$ we define an automorphism $\psi_\lambda: x\to x,\; y\to y,\; z\to\lambda z$.

The next lemma provides for some translation between the language of polynomials and the group action language. It
is similar to the hiking process \cite{BelovIAN,BelovUzyRowenSerdicaBachtur}.

\begin{lemma}         \label{LeForHiking}
Let $\varphi\in K\langle x,y,z\rangle$. Let $\varphi(x)=x,\; \varphi(y)=y+\sum_i R_i+R',\; \varphi(z)=z+Q$. Let
$\deg(R_i)=N$, let also the degree of all monomials in $R'$ be greater than $N$, and let the degree of all
monomials in $Q$ be greater than or equal to $N$. Finally, assume $\deg_z(R_i)=i$ and the $z$-degree of all
monomials of $R_1$ greater than $0$.

Then

a) $\psi_\lambda^{-1}\varphi\psi_\lambda: x\to x,\; y\to y+\sum_i \lambda^iR_i+R'',\; z\to z+Q'$. Also the total
degree of all monomials comprising $R'$ is greater than $N$, and the degree of all monomials of $Q$ is greater than
or equal to $N$.

b) Let $\phi=\prod \left(\psi_{\lambda^{-1}_i}\varphi\psi_{\lambda_i}\right)^{k_i}$. Then

$$\phi: x\to x,\; y\to y+\sum_i R_i\lambda_i^{k_i}+S,\; z\to z+T$$
where the degree of all monomials of $S$ is greater than $N$ and the degree of all monomials of $T$ is greater than
or equal to $N$.
\end{lemma}

{\bf Proof.} a) By direct computation. b) is a consequence of a).

{\bf Remark.} In the case of characteristic zero, the condition of $K$ being algebraically closed can be dropped.
After hiking for several steps, we need to prove just

\begin{lemma}                  \label{LeHikingCh0}
Let $\Ch(K)=0$, let $n$ be a positive integer. Then there exist $k_1,\dots,k_s\in {\mathbb Z}$ and
$\lambda_1,\dots,\lambda_s\in K$ such that

\begin{itemize}
    \item $\sum k_i=1$.
    \item $\sum_i k_i^{n}\lambda_i=0$.
\end{itemize}
\end{lemma}

Using this lemma we can cancel out all the terms in the product in the Lemma \ref{LeForHiking} except for the
constant one. The proof of Lemma \ref{LeHikingCh0} for any field of zero characteristic can be obtained through
the following observation:

\begin{lemma}     \label{LeInclExcl}
$$\left(\sum_{i=1}^n \lambda_i\right)^n-\sum_j
\left(\lambda_1+\cdots+\widehat{\lambda_j}+\cdots+\lambda_n\right)^n+\cdots+
$$
$$+(-1)^{n-k}\sum_{i_1<\cdots<i_k}\left(x_{i_1}+\cdots+x_{i_k}\right)^n+
\cdots+(-1)^{n-1}\left(x_1^n+\cdots+x_n^n\right)=n!\prod_{i=1}^nx_i
$$
and if $m<n$ then
$$\left(\sum_{i=1}^n \lambda_i\right)^m-\sum_j
\left(\lambda_1+\cdots+\widehat{\lambda_j}+\cdots+\lambda_n\right)^m+\cdots+
$$
$$+(-1)^{n-k}\sum_{i_1<\cdots<i_k}\left(x_{i_1}+\cdots+x_{i_k}\right)^m+
\cdots+(-1)^{n-1}\left(x_1^m+\cdots+x_n^m\right)=0.
$$
\end{lemma}

The lemma \ref{LeInclExcl} allows us to replace the $n$-th powers by product of constants, after that the statement
of Lemma \ref{LeHikingCh0} becomes transparent.

\begin{lemma}               \label{Le2xyAprox}
Let $\varphi: x\to x+R_1,\; y\to y+R_2,\; z\to z'$, such that the total degree of all monomials in $R_1,\; R_2$ is
greater than or equal to $N$. Then for $\Psi(\varphi): x\to x+R_1',\; y\to y+R_2',\; z\to z''$ with the total degree of
all monomials in $R_1', R_2'$ also greater than or equal to $N$.
\end{lemma}

{\bf Proof.} Similar to the proof of Theorem \ref{ThMainTechn}.

Lemmas \ref{Le2xyAprox}, \ref{LeForHiking}, \ref{LeHiking} imply the  following statement.

\begin{lemma}                  \label{LmTechInvar}
Let $\varphi_j\in \Aut_0(K\langle x,y,z\rangle),\; j=1,2$, such that
$$\varphi_j(x)=x,\ \varphi_j(y)=y+\sum_i R_i^j+R'_j,\;
\varphi_j(z)=z+Q_j.$$ Let $\deg(R_i^j)=N$, and suppose that the degree of all monomials
in $R'_j$ is greater than $N$, while the degree of all monomials in $Q$ is greater than or equal to $N$;
$\deg_z(R_i)=i$, and the $z$-degree of all monomials in $R_1$ is positive. Let $R_0^1=0,\; R_0^2\ne 0$.

Then $\Psi(\varphi_1)\ne\varphi_2$.
\end{lemma}

Consider the automorphism

$$\phi: x\to x,\; y\to y,\; z\to z+P(x,y).$$
Let $\Psi\in\Aut_{\Ind}\TAut_0(k\langle x,y,z\rangle)$ stabilize the standard torus action pointwise. Then
$$\Psi(\phi): x\to x,\; y\to y,\; z\to z+Q(x,y).$$ We denote
$$\bar{\Psi}(P)=Q.$$ Our goal is to prove that $\bar{\Psi}(P)=P$
for all $P$ if $\Psi$ stabilizes all linear automorphisms and $\bar{\Psi}(xy)=xy$. We proceed by strong induction on
total degree. The base case corresponds to $k=1$ and $l=1$ and is assumed. We then heave

\begin{lemma}
$$\bar{\Psi}(x^ky^l)=x^ky^l$$
provided that $\bar{\Psi}(P)=P$ for all monomials $P(x,y)$ of total degree $< k+l$.
\end{lemma}

{\bf Proof.}

Let $$\phi: x\to x,\; y\to y,\; z\to z+x^ky^l,$$
$$\varphi_1: x\to x+y^l,\; y\to y,\; z\to z,$$ $$\varphi_2: x\to x,\; y\to
y+x^k,\; z\to z,$$ $$\varphi_3: x\to x,\; y\to y,\; z\to z+xy,$$ $$h: x\to x,\; y\to y,\; z\to z-x^{k+1}.$$ Then, for
$k>1$ and $l>1$

\begin{gather*}
g=h\varphi_3^{-1}\varphi_1^{-1}\varphi_2^{-1}\varphi_3\varphi_1\varphi_2:\\
x\mapsto x-y^l+(y-(x-y^l)^k)^l,\\
y\mapsto y-(x-y^l)^k+(x-y^l+(y-(x-y^l)^k)^l)^k,\\
z\mapsto z-xy-x^{k+1}+(x-y^l)(y-(x-y^l)^k).
\end{gather*}
Observe that the height of $g(x)-x$, $g(y)-y$ and $g(z)-z$ is at least $k+l-1$, when $k>1$ or $l>1$. We then use
Theorem \ref{ThMainTechn} and the induction step. Applying $\Psi$ yields the result because
$\Psi(\varphi_i)=\varphi_i,\; i=1,2,3$ and $\varphi(H_N)\subseteq H_N$ for all $N$. The lemma is proved.

Let $$M_{k_1,\dots,k_s}=x^{k_1}y^{k_2}\cdots y^{k_s}$$ for even $s$ and
$$M_{k_1,\dots,k_s}=x^{k_1}y^{k_2}\cdots x^{k_s}$$ for odd
$s$, $k=\sum_{i=1}^n k_i$. Then
$$M_{k_1,\dots,k_s}=M_{k_1,\dots,k_{s-1}}y^{k_s}$$ for even $s$
and
$$M_{k_1,\dots,k_s}=M_{k_1,\dots,k_{s-1}}x^{k_s}$$ for odd
$s$.

We have to prove that $\bar{\Psi}(M_{k_1,\dots,k_{s}})=M_{k_1,\dots,k_{s}}$. By induction we may assume
that $\bar{\Psi}(M_{k_1,\dots,k_{s-1}})=M_{k_1,\dots,k_{s-1}}$.

For any monomial $M=M(x,y)$ we define an automorphism
$$\varphi_M: x\to x,\; y\to y,\; z\to z+M.$$

We also define the automorphisms $$\phi_k^e: x\to x,\; y\to y+zx^k,\; z\to z$$ and $$\phi_k^o: x\to x+zy^k,\; y\to
y,\; z\to z.$$ We will present the case of even $s$ - the odd $s$ case is similar.

Let $D_{zx^k}^e$ be a derivation of $K\langle x,y,z\rangle$ such that $D_{zx^k}^e(x)=0$,
 $D_{zx^k}^e(y)=zx^k$, $D_{zx^k}^e(z)=0$.
Similarly, let $D_{zy^k}^o$ be a derivation of $K\langle x,y,z\rangle$ such that $D_{zy^k}^o(y)=0$,
$D_{zx^k}^o(x)=zy^k$, $D_{zy^k}(z)^o=0$.

The following lemma is proved by direct computation:

\begin{lemma}
Let
$$u=\phi_{k_{s}}^e{}^{-1}\varphi(M_{k_1,\dots,k_{s-1}})^{-1}\phi_{k_{s}}^e\varphi(M_{k_1,\dots,k_{s-1}})$$
for even $s$ and
$$u=\phi_{k_{s}}^o{}^{-1}\varphi(M_{k_1,\dots,k_{s-1}})^{-1}\phi_{k_{s}}^o\varphi(M_{k_1,\dots,k_{s-1}})$$
for odd $s$. Then

$$u:x\to x,\; y\to y+M_{k_1,\dots,k_s}+N',\; z\to
z+D^e_{zx^k}(M_{k_1,\dots,k_{s-1}})+N$$ for even $s$ and
$$u:x\to x+M_{k_1,\dots,k_s}+N',\; y\to y,\; z\to
z+D^o_{zx^k}(M_{k_1,\dots,k_{s-1}})+N$$ for odd $s$,
where $N$, $N'$ are sums of terms of degree $>k=\sum_{i=1}^sk_i$.
\end{lemma}

Let $\psi(M_{k_1,\dots,k_{s}}): x\to x,\; y\to y,\; z\to z+M_{k_1,\dots,k_{s}}$, $$\alpha_e: x\to x,\; y\to
y-z,\;z\to z, \alpha_o: x\to x-z,\; y\to y,\;z\to z,$$ Let $P_M=\Psi(M)-M$. Our goal is to prove that $P_M=0$.

Let
$$v=\psi(M_{k_1,\dots,k_{s}})^{-1}\alpha_{e}\psi(M_{k_1,\dots,k_{s}})u\alpha_{e}^{-1}$$
for even $s$ and
$$v=\psi(M_{k_1,\dots,k_{s}})^{-1}\alpha_{o}\psi(M_{k_1,\dots,k_{s}})u\alpha_{o}^{-1}$$
for odd $s$.

The next  lemma  is also proved by direct computation:

\begin{lemma}      \label{LmFinalType}
a) $$v: x\to x,\; y\to y+H,\; z\to z+H_1+H_2$$ for even $s$ and
$$v: x\to x+H,\; y\to y,\; z\to z+H_1+H_2$$
for odd $s$

b)
$$\Psi(v): x\to x, y\to y+P_{M_{k_1,\dots,k_{s}}}+\widetilde{H}, z\to z+\widetilde{H_1}+\widetilde{H_2}$$
for even $s$ and
$$\Psi(v): x\to x+P_{M_{k_1,\dots,k_{s}}}+\widetilde{H}, y\to y, z\to z+\widetilde{H_1}+\widetilde{H_2}$$
for odd $s$, where $H_2$, $\widetilde{H_2}$ are sums of terms of degree greater than $k=\sum_{i=1}^s k_i$,
$H$, $\widetilde{H}$ are sums of terms of degree $\ge k$ and positive $z$-degree, $H_1$, $\widetilde{H_1}$
are sums of terms of degree $k$ and positive $z$-degree.
\end{lemma}

%
%
%
%

{\bf Proof of Theorem \ref{ThTAss3Ind}.} Part b) follows from part a). In order to prove a) we are going to show
that $\bar{\Psi}(M)=M$ for any monomial $M(x,y)$ and for any $\Psi\in\Aut_{\Ind}(\TAut(\langle x,y,z\rangle))$
stabilizing the standard torus action $T^3$ and $\phi$. The automorphism $\Psi(\Phi_M)$ has the form described in
Lemma \ref{LmFinalType}. But in this case Lemma \ref{LmTechInvar} implies $\bar{\Psi}(M)-M=0$.

%
%
%
%
%

\section{Some open questions concerning the tame automorphism group}

As the conclusion of the paper, we would like to raise the following questions.

\begin{enumerate}
    \item Is it true that any  automorphism $\varphi$ of $\Aut(K\langle x_1,\dots,x_n\rangle)$ (in the group-theoretic
        sense - that is, not necessarily an automorphism preserving the $\Ind$-scheme structure) for $n=3$ is
        semi-inner, i.e. is a conjugation by some automorphism or mirror anti-automorphism?
    \item Is it true that $\Aut(K\langle x_1,\dots,x_n\rangle)$ is generated by affine automorphisms and
        automorphism $x_n\to x_n+x_1x_2,\; x_i\to x_i,\; i\ne n$? For $n\ge 5$ it seems to be easier and the answer
        is probably positive, however for $n=3$ the answer is known to be negative, cf. Umirbaev
        \cite{umirbaev1995ext} and Drensky and Yu \cite{DYuStrongAnik}. For $n\ge 4$ we believe the answer is
        positive.
    \item  Is it true that $\Aut(K[x_1,\dots,x_n])$ is generated by linear automorphisms and automorphism $x_n\to
        x_n+x_1x_2,\; x_i\to x_i,\; i\ne n$? For $n=3$ the answer is negative: see the proof of the Nagata conjecture
        \cite{SU1,SU2,UY}. For $n\ge 4$ it is plausible that the answer is positive.
    \item Is any  automorphism $\varphi$ of $\Aut(K\langle x,y,z\rangle)$ (in the group-theoretic sense) semi-inner?
    \item Is it true that the conjugation in Theorems \ref{ThAutTAut} and \ref{ThAutTAutFreass} can be done by
        some tame automorphism? Suppose $\psi^{-1}\varphi\psi$ is tame for any tame $\varphi$. Does it follow
        that $\psi$ is tame?
    \item Prove Theorem \ref{ThTAss3Ind} for $\Ch(K)=2$. Does it hold on the set-theoretic level, i.e.
        $\Aut(\TAut(K\langle
x,y,z\rangle))$ are generated by conjugations by an automorphism or the mirror anti-automorphism?
\end{enumerate}

Similar questions can be formulated for nice automorphisms.

\chapter{Approximation by tame automorphisms and the Kontsevich conjecture} 
\label{Chapter3} \lhead{Chapter 4. \emph{Approximation by tame automorphisms}}

The first four sections of this chapter is based on the paper \cite{KGE}. In light of this approximation method, A.
Elishev, A. Kanel-Belov, and  J.-T. Yu \cite{K-BE4} proposed a proof of the Kontsevich Conjecture.

\section{Endomorphisms of \texorpdfstring{$\mathbb{K}[X]$}{K[X]}, \texorpdfstring{$W_n(\mathbb{K})$}{Wn(K)} and \texorpdfstring{$P_n(\mathbb{K})$}{Pn(K)}}

\subsection{Definitions and notation}

The $n$-th Weyl algebra $W_{n}(\mathbb{K})$ over $\mathbb{K}$ is by definition the quotient of the free
associative algebra
\begin{equation*}
\mathbb{K}\langle a_1,\ldots,a_n,b_1,\ldots,b_n\rangle
\end{equation*}
by the two-sided ideal generated by elements
\begin{equation*}
b_ia_j-a_jb_i-\delta_{ij},\;\;a_ia_j-a_ja_i,\;\;b_ib_j-b_jb_i,
\end{equation*}
with $1\leq i,j\leq n$. One can think of $W_{n}(\mathbb{K})$ as the algebra
\begin{equation*}
\mathbb{K}[x_1,\ldots,x_n,y_1,\ldots,y_n]
\end{equation*}
with two sets of $n$ mutually commuting generators (images of the free generators under the canonical projection)
which interact according to $[y_i,x_j]=y_ix_j-x_jy_i=\delta_{ij}$, although, unless the context necessitates
clarification, we would like to denote the Weyl algebra henceforth by $W_{n}(\mathbb{K})$ in order to avoid
confusion with $\mathbb{K}[X]$ -- notation reserved for the ring of polynomials in commuting variables.

\medskip

The polynomial algebra $\mathbb{K}[x_1,\ldots,x_N]$ itself is the quotient of the free associative algebra by the
congruence that makes all its generators commutative. When $N=2n$ is even, the  algebra $A_{2n}$ carries an
additional structure of the Poisson algebra -- namely, a bilinear map $$\lbrace
\;,\;\rbrace:\mathbb{K}[x_1,\ldots,x_N]\otimes \mathbb{K}[x_1,\ldots,x_N]\rightarrow
\mathbb{K}[x_1,\ldots,x_N]$$ that turns $\mathbb{K}[x_1,\ldots,x_N]$ into a Lie algebra and acts as a
derivation with respect to polynomial multiplication. Under a fixed choice of generators, this map is given by the
canonical Poisson bracket
\begin{equation*}
\lbrace x_i,x_j\rbrace = \delta_{i,n+j}-\delta_{i+n,j}.
\end{equation*}
We denote the pair $(\mathbb{K}[x_1,\ldots,x_{2n}],\lbrace \;,\;\rbrace)$ by $P_n(\mathbb{K})$. In our
discussion the coefficient ring $\mathbb{K}$ is a field of characteristic zero, and for later purposes (Proposition 4.3)
we require $\mathbb{K}$ to be algebraically closed. Thus one may safely assume $\mathbb{K}=\mathbb{C}$ in
the sequel.

\medskip

Throughout we assume all homomorphisms to be unital and preserving all defining structures carried by the objects in
question. Thus, by a Weyl algebra endomorphism we always mean a $\mathbb{K}$-linear ring homomorphism
$W_{n}(\mathbb{K})$ into itself that maps $1$ to $1$. Similarly, the set $\End \mathbb{K}[x_1,\ldots,x_n]$
consists of all $\mathbb{K}$-endomorphisms of the polynomial algebra, while $\End P_{n}$ is the set of
polynomial endomorphisms preserving the Poisson structure. We will call elements of the group $\Aut P_{n}$
\textbf{polynomial symplectomorphisms}, due to the fact that they can be identified with polynomial one-to-one
mappings $\mathbb{A}^{2n}_{\mathbb{K}}\rightarrow\mathbb{A}^{2n}_{\mathbb{K}}$ of the affine space
$\mathbb{A}^{2n}_{\mathbb{K}}$ which preserve the symplectic form
\begin{equation*}
\omega=\sum_{i}dp_i\wedge dx_i.
\end{equation*}

\medskip

Any endomorphism $\varphi$ of $\mathbb{K}[x_1,\ldots,x_N]$, $P_n(\mathbb{K})$ or $W_n(\mathbb{K})$
can be identified with the ordered set
\begin{equation*}
(\varphi(x_1),\;\varphi(x_2),\;\ldots)
\end{equation*}
of images of generators of the corresponding algebra. For $\mathbb{K}[x_1,\ldots,x_N]$ and
$P_n(\mathbb{K})$, the polynomials $\varphi(x_i)$ can be decomposed into sums of homogeneous components;
this means that the endomorphism $\varphi$ may be written as a formal sum
\begin{equation*}
\varphi = \varphi_0+\varphi_1+\cdots,
\end{equation*}
where $\varphi_k$ is a string (of length $N$ and $2n$, respectively) whose entries are homogeneous polynomials of
total degree $k$.\footnote{We set $\Deg x_i=1$.} Accordingly, the height $\Ht(\varphi)$ of the endomorphism is
defined as
\begin{equation*}
\Ht(\varphi)=\inf\lbrace k\;|\;\varphi_k\neq 0\rbrace,\;\;\Ht(0)=\infty.
\end{equation*}
This is not to be confused with the degree of endomorphism, which is defined as $\Deg(\varphi)=\sup\lbrace
k\;|\;\varphi_k\neq 0\rbrace$.\footnote{For $W_n$ the degree is well defined, but the height depends on the
ordering of the generators.} The height $\Ht(f)$ of a polynomial $f$ is defined quite similarly to be the minimal
number $k$ such that the homogeneous component $f_k$ is not zero. Evidently, for an endomorphism
$\varphi=(\varphi(x_1),\;\ldots,\;\varphi(x_N))$ one has
\begin{equation*}
\Ht(\varphi)=\inf\lbrace \Ht(\varphi(x_i))\;|\;1\leq i\leq N\rbrace.
\end{equation*}

\medskip

The function
\begin{equation*}
d(\varphi,\psi)=\exp(-\Ht(\varphi-\psi))
\end{equation*}
is a metric on $\End \mathbb{K}[x_1,\ldots,x_N]$. We will refer to the corresponding topology on $\End$ (and on
subspaces such as $\Aut$ and $\TAut$) as the formal power series topology.

\medskip

\subsection{Tame automorphisms}
We call an automorphism $\varphi\in\Aut\mathbb{K}[x_1,\ldots,x_N]$ \textbf{elementary} if it is of the form
\begin{equation*}
\varphi = (x_1,\ldots,\;x_{k-1},\;ax_k+f(x_1,\ldots,x_{k-1},\;x_{k+1},\;\ldots,\;x_N),\;x_{k+1},\;\ldots,\;x_N)
\end{equation*}
with $a\in\mathbb{K}^{\times}$. Observe that linear invertible changes of variables -- that is, transformations of the
form
\begin{equation*}
(x_1,\;\ldots,\;x_N)\mapsto (x_1,\;\ldots,\;x_N)A,\;\;A\in\GL(N,\mathbb{K})
\end{equation*}
are realized as compositions of elementary automorphisms.

The subgroup of $\Aut\mathbb{K}[x_1,\ldots,x_N]$ generated by all elementary automorphisms is the group
$\TAut \mathbb{K}[x_1,\ldots,x_N]$ of so-called \textbf{tame automorphisms}.

\medskip

Let $P_{n}(\mathbb{K})=\mathbb{K}[x_1,\ldots,x_{n},p_1,\ldots,p_n]$ be the polynomial algebra in $2n$
variables with Poisson structure. It is clear that for an elementary
$\varphi\in\Aut\mathbb{K}[x_1,\ldots,x_{n},p_1,\ldots,p_n]$ to be a symplectomorphism, it must either be a linear
symplectic change of variables -- that is, a transformation of the form
\begin{equation*}
(x_1,\;\ldots,\;x_n,\;p_1,\;\ldots,\;p_n)\mapsto (x_1,\;\ldots,\;x_n,\;p_1,\;\ldots,\;p_n)A
\end{equation*}
with $A\in\Sp(2n,\mathbb{K})$ a symplectic matrix, or be an elementary transformation of one of two following
types:
\begin{equation*}
(x_1,\;\ldots,\;x_{k-1},\;x_k+f(p_1,\;\ldots,\;p_n),\;x_{k+1},\;\ldots,\;x_n,\;p_1,\;\ldots,\;p_n)
\end{equation*}
or
\begin{equation*}
(x_1,\;\ldots,\;x_{n},\;p_1,\;\ldots,\;p_{k-1},\;p_k+g(x_1,\;\ldots,\;x_n),\;p_{k+1},\;\ldots,\;p_n).
\end{equation*}
Note that in both cases we do not include translations of the affine space into our consideration, so we may safely
assume the polynomials $f$ and $g$ to be at least of height one.

The subgroup of $\Aut P_{n}(\mathbb{K})$ generated by all such automorphisms is the group $\TAut
P_{n}(\mathbb{K})$ of \textbf{tame symplectomorphisms}. One similarly defines the notion of tameness for the
Weyl algebra $W_n(\mathbb{K})$, with tame elementary automorphisms having the exact same form as for
$P_n(\mathbb{K})$.

\medskip

The automorphisms which are not tame are called \textbf{wild}. It is unknown at the time of writing whether the
algebras $W_n$ and $P_n$ have any wild automorphisms in characteristic zero for $n>1$, however for $n=1$ all
automorphisms are known to be tame \cite{Jung, VdK, ML1, ML2}. On the other hand, the celebrated example of
Nagata
\begin{equation*}
(x+(x^2-yz)x,\;y+2(x^2-yz)x+(x^2-yz)^2z,\;z)
\end{equation*}
provides a wild automorphism of the polynomial algebra $\mathbb{K}[x,y,z]$.

\medskip

It is known due to Kanel-Belov and Kontsevich \cite{K-BK1,K-BK2} that for $\mathbb{K}=\mathbb{C}$ the
groups
\begin{equation*}
\TAut W_n(\mathbb{C})\;\;\text{and}\;\;\TAut P_n(\mathbb{C})
\end{equation*}
are isomorphic. The homomorphism between the tame subgroups is obtained by means of non-standard analysis and
involves certain non-constructible entities, such as free ultrafilters and infinite prime numbers. Recent effort
\cite{K-BE, K-BE2} has been directed to proving the homomorphism's independence of such auxiliary objects,
with limited success.

\section{Approximation by tame automorphisms}

Let $\varphi\in\Aut \mathbb{K}[x_1,\ldots,x_N]$ be a polynomial automorphism. We say that $\varphi$ is
approximated by tame automorphisms if there is a sequence
\begin{equation*}
\psi_1,\;\psi_2,\ldots,\;\psi_k,\ldots
\end{equation*}
of tame automorphisms such that
\begin{equation*}
\Ht((\psi_k^{-1}\circ\varphi)(x_i)-x_i)\geq k
\end{equation*}
for $1\leq i\leq N$ and all $k$ sufficiently large. Observe that any tame automorphism $\psi$ is approximated by
itself -- that is, by a stationary sequence $\psi_k=\psi$.

\medskip

The following two theorems are the main results of this chapter.

\begin{thm}\label{theorem31}
Let $\varphi=(\varphi(x_1),\;\ldots,\;\varphi(x_N))$ be an automorphism of the polynomial algebra
$\mathbb{K}[x_1,\ldots,x_N]$ over a field $\mathbb{K}$ of characteristic zero, such that its Jacobian
\begin{equation*}
\J(\varphi)=\Det \left[\frac{\partial \varphi(x_i)}{\partial x_j}\right]
\end{equation*}
is equal to $1$. Then there exists a sequence $\lbrace \psi_k\rbrace\subset \TAut \mathbb{K}[x_1,\ldots,x_N]$ of
tame automorphisms approximating $\varphi$.
\end{thm}

\begin{thm}\label{theorem32}
Let $\sigma=(\sigma(x_1),\;\ldots,\;\sigma(x_n),\;\sigma(p_1),\;\ldots,\;\sigma(p_n))$ be a symplectomorphism of
$\mathbb{K}[x_1,\ldots,x_n,p_1,\ldots,p_n]$ with unit Jacobian.
Then there exists a sequence $\lbrace \tau_k\rbrace\subset \TAut P_n(\mathbb{K})$ of tame symplectomorphisms
approximating $\sigma$.
\end{thm}

Theorem \ref{theorem31} is a special case of a classical result of Anick \cite{An} (Anick proved approximation for
all \'etale maps, not just automorphisms). We give here a slightly simplified proof suitable for our context. The
second theorem first appeared in \cite{KGE} and is essential in our approach to the lifting problem in deformation
quantization.

\medskip

The proof of Theorem \ref{theorem31} consists of several steps each of which amounts to composing a given
automorphism $\varphi$ with a tame transformation of a specific type -- an operation which allows one to dispose in
$\varphi(x_i)$ ($1\leq i\leq N$) of monomial terms of a given total degree, assuming that the lower degree terms
have already been dealt with. Thus the approximating sequence of tame automorphisms is constructed by induction.
As it was mentioned before, we disregard translation automorphisms completely: all automorphisms discussed here
are origin-preserving, so that the polynomials $\varphi(x_i)$ have zero free part. This of course leads to no loss of
generality.

\medskip

The process starts with the following straightforward observation.

\smallskip

\begin{lem}
There is a linear transformation $A\in\GL(N,\mathbb{K})$
\begin{equation*}
(x_1,\;\ldots,\;x_N)\mapsto (x_1,\;\ldots,\;x_N)A
\end{equation*}
such that its composition $\varphi_A$ with $\varphi$ fulfills
\begin{equation*}
\Ht(\varphi_A(x_i)-x_i)\geq 2
\end{equation*}
for all $i\in\lbrace 1,\ldots,N\rbrace$.
\end{lem}

\begin{proof}
Consider $$A_1=\left[\frac{\partial \varphi(x_i)}{\partial x_j}\right]\left( 0,\ldots,0\right)$$ -- the linear part of
$\varphi$. Its determinant is equal to the value of $\J(\varphi)$ at zero, and $\J(\varphi)$ is a non-zero constant.
Composing $\varphi$ with the linear change of variables induced by $A_1^{-1}$ (on the left) results in an
automorphism $\varphi_A$ that is identity modulo $\mathit{O}(x^2)$.
\end{proof}

Using the above lemma, we may replace $\varphi$ with $\varphi_A$ (and suppress the $A$ subscript for
convenience), thus considering automorphisms which are close to the identity in the formal power series topology.

\smallskip

The next lemma justifies the inductive step: suppose we have managed, by tame left action, to eliminate the terms of
degree $2,\ldots, k-1$, then there is a sequence of elementary automorphisms such that their left action eliminates
the term of degree $k$. This statement translates into the following lemma.

\begin{lem}
Let $\varphi$ be a polynomial automorphism such that
$$
\varphi(x_1)=x_1+f_1(x_1,\ldots,x_n)+r_1,\;\;\ldots,\;\;\varphi(x_n)=x_n+f_n(x_1,\ldots,x_n)+r_n
$$
and $f_i$ are homogeneous of degree $k$ and $r_i$ are the remaining terms (thus $\Ht(r_i)>k$). Then one can find
a sequence $\sigma_1,\ldots, \sigma_m$ of tame automorphisms whose composition with $\varphi$ is given by
\begin{equation*}
\sigma_m\circ\ldots\circ\sigma_1\circ\varphi:\;x_1\mapsto x_1+F_i(x_1,\ldots,x_n)+R_1,\;\;\ldots,\;\;x_n\mapsto x_n+F_n(x_1,\ldots,x_n)+R_n
\end{equation*}
with $F_i$ homogeneous of degree $k+1$ and $\Ht(R_i)>k+1$.
\end{lem}
\begin{proof}
We will first show how to get rid of degree $k$ monomials in the images of all but one generator and then argue that
the remaining image is rectified by an elementary automorphism. Let $N\leq n$ be the number of images
$\varphi(x_i)$ such that $f_i\neq 0$, and let $x_1$ and $x_2$ be two generators\footnote{Evidently, no loss of
generality results from such explicit labelling.} corresponding to non-zero term of degree $k$. The image of $x_1$
admits the following presentation as an element of the polynomial ring $\mathbb{K}[x_3,\ldots,x_n][x_1,x_2]$:
\begin{equation*}
\varphi(x_1)= x_1+\sum_{d}\sum_{p+q=d}\lambda_{p,q}x_1^p x_2^q + r_i
\end{equation*}
where the coefficients $\lambda_{p,q}$ are polynomials of the remaining variables (thus the double sum above is
just a way to express $f_1$ as a polynomial in $x_1$ and $x_2$ with coefficients given by polynomials in the rest of
the variables).

Consider the transformation $\Phi_{\lambda\mu}$ of the following form
\begin{equation*}
x_1\mapsto x_1-\lambda (x_1+\mu x_2)^d,\;\;x_2\mapsto x_2-\lambda \mu^{-1}(x_1+\mu x_2)^d,\;\; x_3\mapsto x_3,\;\;\ldots,\;\; x_n\mapsto x_n,
\end{equation*}
with $\lambda \in \mathbb{K}[x_3,\ldots,x_n]$   and $\mu \in\mathbb{K}$. This mapping is equal to the
composition $\psi_{\mu}\circ\phi_{\lambda\mu}\circ\psi_{\mu}^{-1}$ with
\begin{equation*}
\psi_{\mu}:x_1\mapsto x_1+\mu x_2,\;\;x_2\mapsto x_2
\end{equation*}
and
\begin{equation*}
\phi_{\lambda\mu}: x_1\mapsto x_1,\;\;x_2\mapsto x_2+\lambda\mu^{-1}x_1^d
\end{equation*}
and so is a tame automorphism. As the ground field $\mathbb{K}$ has characteristic zero, it is infinite, so that we
can find numbers $\mu_1,\;\ldots,\; \mu_{l(d)}$ such that the polynomials $$(x+\mu_1 y)^d,\;\ldots,\;(x+\mu_{l(d)}
y)^d$$ form a basis of the $\mathbb{K}$-module of homogeneous polynomials in $x$ and $y$ of degree $d$ (this
is an easy exercise in linear algebra). Therefore, by selecting $\Phi_{\lambda\mu}$ with appropriate polynomials
$\lambda_{p,q}$ and $\mu_i$ corresponding to the basis, we eliminate, by acting with $\Phi_{\lambda\mu}$ on the
left, the degree $d$ terms in the double sum. Iterating for all $d$, we dispose of $f_1$ entirely.

The above procedure yields a new automorphism $\tilde{\varphi}$ which is a composition of the initial
automorphism $\varphi$ with a tame automorphism. The number $\tilde{N}$ of images of $x_i$ under
$\tilde{\varphi}$ with non-zero term of degree $k$ equals $N-1$; therefore, the procedure can be repeated a finite
number of times to give an automorphism $\varphi_1$, such that the image under $\varphi_1$ of only one generator
contains a non-zero term of degree $k$. Let
\begin{equation*}
\varphi_1(x_n)=x_n+g_n(x_1,\ldots,x_n)+\tilde{r}_n
\end{equation*}
be the image of that generator (again, no loss of generality results from us having labelled it $x_n$). We claim now
that the polynomial $g_n$ does not depend on $x_n$.

Indeed, otherwise the Jacobian of $\varphi_1$ (which must be a constant and is in fact equal to $1$ in our setting)
would have a degree $k-1$ component given by
$$
\partial_{x_n} g_n(x_1,\ldots,x_n)\neq 0
$$
(remember that by construction $g_1=\ldots=g_{n-1}=0$), which yields a contradiction. Note that another way of
looking at this condition is that if a polynomial mapping of the form
\begin{equation*}
x_1\mapsto x_1 + H_1(x_1,\ldots,x_n),\;\;\ldots,\;\;x_n\mapsto x_n + H_n(x_1,\ldots,x_n),\;\;\Ht(H_i)>1
\end{equation*}
is an automorphism, the higher-degree part $(H_1,\ldots, H_n)$ must have traceless Jacobian:
$$
\Tr \left(\frac{\partial H_i}{\partial x_j}\right)=0.
$$

Finally, since $g_n$ does not contain $x_n$, an elementary automorphism
\begin{equation*}
x_1\mapsto x_1,\;\;\ldots,\;\;x_{n-1}\mapsto x_{n-1},\;\;x_n\mapsto x_n-g_n(x_1,\ldots,x_n)
\end{equation*}
eliminates this term. The lemma is proved.
\end{proof}

The last lemma concludes the proof of Theorem \ref{theorem31} by induction. The proof of the inductive step is
essentially a statement that a certain vector space invariant under a linear group action is, in a manner of speaking,
big enough to allow for elimination by elements of the group. More precisely, let $T_{n,k}(\mathbb{K})$ be the
vector space of all \textbf{traceless} $n$ by $n$ matrices whose entries are homogeneous of degree $k$
polynomials from $\mathbb{K}[x_1,\ldots,x_n]$, and let the group $\GL(n,\mathbb{K})$ act on $T_{n,k}$ as
follows: for $A\in\GL(n,\mathbb{K})$ and $v\in T_{n,k}$, the image $A(v)$ is obtained by taking the product
matrix $v A^{-1}$ and then performing (entry-wise in $vA^{-1}$) the linear change of variables induced by $A$.
Then one has the following

\begin{prop}
If $V\subset T_{n,k}(\mathbb{K})$ is a $\mathbb{K}$-submodule invariant under the defined above action of
$\GL(n,\mathbb{K})$, then either $V=0$ or $V= T_{n,k}(\mathbb{K})$.
\end{prop}

Properties of similar nature played an important role in \cite{K-BMR1, K-BMR2}. The invariance under linear
group action will become somewhat more pronounced in the symplectomorphism case.

\section{Approximation by tame symplectomorphisms and lifting to Weyl algebra}
We turn to the proof of the more relevant Theorem \ref{theorem32} on the symplectic tame approximation. The
strategy is analogous to the proof of approximation for polynomial automorphisms with unit Jacobian, with a few
more elaborate details which we now consider.

\medskip

The first step of the proof copies the polynomial automorphism case and takes the following form.

\begin{lem}
There is a linear transformation $A\in\Sp(2n,\mathbb{K})$
\begin{equation*}
(x_1,\;\ldots,\;x_n,\;p_1,\;\ldots,\;p_n)\mapsto (x_1,\;\ldots,\;x_n,\;p_1,\;\ldots,\;p_n)A
\end{equation*}
such that its composition $\sigma_A$ with $\sigma$ fulfills
\begin{equation*}
\Ht(\sigma_A(x_i)-x_i)\geq 2,\;\;\Ht(\sigma_A(p_i)-p_i)\geq 2
\end{equation*}
for all $i\in\lbrace 1,\ldots,n\rbrace$.
\end{lem}

\medskip

We now proceed to formulate the inductive step in the proof as the following main lemma.

\begin{lem}
Let $\sigma$ be a polynomial symplectomorphism such that
\begin{equation*}
\sigma(x_i)=x_i+U_i,\;\;\sigma(p_i)=p_i+V_i
\end{equation*}
and $U_i$ and $V_i$ are of height at least $k$. Then there exists a tame symplectomorphism $\sigma_k$ such that
the polynomials $\tilde{U}_i=(\sigma_k^{-1}\circ\sigma)(x_i)-x_i$ and
$\tilde{V}_i=(\sigma_k^{-1}\circ\sigma)(p_i)-p_i$ are of height at least $k+1$.
\end{lem}
\begin{proof}
In order to establish the inductive step, we are going to need the following classical result (which is a particular case
of Corollary 17.21 in Fulton and Harris \cite{FulHar}).

\begin{lem}
Suppose $\mathbb{K}$ is an infinite field, $A=\mathbb{K}[x_1,\ldots,x_n,p_1,\ldots,p_n]$ is the polynomial
algebra with standard $\mathbb{Z}$-grading according to the total degree
\begin{equation*}
A=\bigoplus_{d\geq 0} A_d, \;\; A_d=\lbrace \text{homogeneous polynomials of total degree d}\rbrace.
\end{equation*}
Let $V$ be a $\mathbb{K}$-submodule of $A$ invariant under the action of $\Sp(2n, \mathbb{K})$ (given by
linear symplectic changes of variables). Suppose $V$ is contained in a given homogeneous component $A_d$. If
$V \neq 0$ then $V = A_d$.
\end{lem}

We now turn to the proof of the inductive step.  Suppose that
\begin{equation*}
\sigma: x_i\mapsto x_i + f_i + P_i,\;\;p_j\mapsto p_j + g_j + Q_j
\end{equation*}
is a polynomial symplectomorphism, where $f_i$ and $g_j$ are degree $k$ components and the height of $P_i$
and $Q_j$ is greater than $k$. The preservation of the symplectic structure by $\sigma$ means that the $k$-th
component obeys the following identities:
\begin{equation*}
\lbrace x_i, f_j\rbrace - \lbrace x_j, f_i\rbrace=0
\end{equation*}
and
\begin{equation*}
\lbrace p_i, f_j\rbrace - \lbrace p_j, f_i\rbrace=0
\end{equation*}
where $\lbrace\;,\;\rbrace$ is the Poisson bracket corresponding to the symplectic form. In the case of standard
symplectic structure these identities translate into
\begin{equation*}
\frac{\partial f_i}{\partial p_j}-\frac{\partial f_j}{\partial p_i}=0,\;\; \frac{\partial g_i}{\partial x_j}-\frac{\partial g_j}{\partial x_i}=0,
\end{equation*}
in which one recognizes the condition for an appropriate differential form to be closed. The triviality of affine space
cohomology then implies that there exists a polynomial $F(x_1,\ldots,x_n,p_1,\ldots,p_n)$, homogeneous of degree
$k+1$, such that
\begin{equation*}
\frac{\partial F}{\partial p_i}=f_i,\;\; \frac{\partial F}{\partial x_i}=g_i;
\end{equation*}
in this way the $k$-component of a symplectomorphism is generated by a homogeneous polynomial. The tame
symplectomorphism group acts on the space of all such generating polynomials (the image of a polynomial is the
polynomial corresponding to the $k$-component of the composition with the tame symplectomorphism), and the
orbit of this tame action carries the structure of a $\mathbb{K}$-module (one may easily come up with a
symplectomorphism corresponding to the sum of two generating polynomials). Therefore this space fulfills the
conditions of the previous lemma, which in this case implies that one can, by a composition with a tame
symplectomorphism, eliminate the $k$-component. The main lemma, and therefore the Theorem 3.2, is proved.
\end{proof}

\medskip

Once the approximation for the case of symplectomorphisms has been established, we can investigate the problem
of lifting symplectomorphisms to Weyl algebra automorphisms. More precisely, one has the following
\begin{prop}
Let $\mathbb{K}=\mathbb{C}$ and let $\sigma:P_n(\mathbb{C})\rightarrow P_n(\mathbb{C})$ be a
symplectomorphism over complex numbers. Then there exists a sequence
\begin{equation*}
\psi_1,\;\psi_2,\;\ldots,\;\psi_k,\;\ldots
\end{equation*}
of tame automorphisms of the $n$-th Weyl algebra $W_n(\mathbb{C})$, such that their images $\sigma_k$ in
$\Aut P_n(\mathbb{C})$ approximate $\sigma$.
\end{prop}
\begin{proof}
This is an immediate corollary of Theorem \ref{theorem32} and the existence of tame subgroup isomorphism
\cite{K-BK1}.

\end{proof}
A few comments are in order. First, the quantization of elementary symplectomorphisms is a very simple procedure:
one needs only replace the $x_i$ and $p_i$ by their counterparts $\hat{x}_i$ and $\hat{p}_i$ in the Weyl algebra
$W_n$. Because the transvection polynomials $f$ and $g$ (in the expressions for elementary symplectomorphisms)
depend, as it has been noted, on one type of generators (resp. $p$ and $x$), the quantization is well defined.

Second, as the tame automorphism groups $\TAut W_n(\mathbb{C})$ and $\TAut P_n(\mathbb{C})$ are
isomorphic, the correspondence between sequence of tame symplectomorphisms converging to
symplectomorphisms and sequences of tame Weyl algebra automorphisms is one to one. The main question is how
one may interpret these sequences as endomorphisms of $W_n(\mathbb{C})$.

\smallskip

Our construction shows that these sequences of tame automorphisms may be thought of as (vectors of) power series
-- that is, elements of
\begin{equation*}
\mathbb{C}[[\hat{x}_1,\ldots,\hat{x}_n,\hat{p}_1,\ldots,\hat{p}_n]]^{2n}.
\end{equation*}

The main problem therefore consists in verifying that these vectors have entries polynomial in generators -- that is,
that the limits of lifted tame sequences are Weyl algebra endomorphisms.

\medskip

One could take a more straightforward (albeit an equivalent) approach to the lifting of symplectomorphisms by
following the prescription of deformation quantization: starting with a symplectic automorphism of the polynomial
algebra $A=\mathbb{K}[x_1,\ldots,x_n,p_1,\ldots,p_n]$, one constructs a map of $A[[\hbar]]$, the algebra of
formal power series (in Planck's constant $\hbar$), which preserves the star product satisfying Weyl algebra
identities. The approximation theory as developed in this text is then a property of the $\hbar$-adic topology. The
(algebraically closed version of) Kontsevich Conjecture would then follow if one were to establish a cutoff theorem.

Unfortunately, this naive approach is deficient in the sense that the resulting lifting is not generally canonical with
respect to the choice of the converging sequence. One needs a more elaborate strategy to construct the lifting map.
One such strategy will be discussed in the next chapter.

\section{Conclusion}

We have developed tame approximation theory for symplectomorphisms in formal power series topology. By virtue
of the known correspondence between tame automorphisms of the even-dimensional affine space and tame
automorphisms of the Weyl algebra, which is the object corresponding to the affine space in terms of deformation
quantization, we have arrived at the lifting property of symplectomorphisms. This line of research may yield new
insights into endomorphisms of the Weyl algebra, the Dixmier conjecture, and the Jacobian conjecture.

\medskip

Inspired by this approximation idea, A. Elishev, A. Kanel-Belov and J.-T. Yu \cite{K-BE4} provide the augmented
automorphisms to prove the Belov-Kontsevich Conjecture. We will show it in the next section.


\section{Augmented Weyl algebra structure}

This section, which is credited to A. Elishev, A. Kanel-Belov and J.-T. Yu \cite{K-BE4}, is the main idea of
augmented Weyl algebra structure and the solution of the Belov-Kontsevich Conjecture.

We first state the following theorem.
\begin{thm} \label{thmgabber}
The mappings
$$
\Phi_N:\Aut^{\leq N} W_{n,\mathbb{C}}\rightarrow \Aut^{\leq N} P_{n,\mathbb{C}}
$$
induced by $\Phi$ are morphisms of algebraic varieties.
\end{thm}

The proof can be found in \cite{K-BE2}. This theorem has an exact (and crucial to our approach) analogue in the
setting of the quantized algebras $W_n$ and $P_n$.

Thus, this section is devoted to the study of some conjectures arising in connection with Jacobian conjecture
(namely, Kontsevich conjecture and related questions), as well as the study of geometric and topological properties
of $\Ind$-schemes of automorphisms of polynomial algebras playing a certain role in approaches to solving the
above conjecture. The results of this study, in addition to their importance in the context of the recovery problem,
Kontsevich conjecture and related issues, are of independent interest.

\medskip

In order to resolve the symplectomorphism lifting problem and construct the inverse to the homomorphism $\Phi$,
we introduce the augmented and skew augmented Weyl and Poisson algebras.

The \emph{augmented}, or $h$-\emph{augmented} Weyl algebra $W^h_{n, \mathbb{C}}$ is defined as the
quotient of the free algebra on $(2n+1)$ indeterminates $\mathbb{C}\langle a_1,\ldots, a_n, b_1,\ldots, b_n,
c\rangle$ by the two-sided ideal generated by elements
$$
a_ia_j-a_ja_i,\;\;b_ib_j-b_jb_i\;\;,
b_ia_j-a_jb_i-\delta_{ij}c,\;\;
a_ic-ca_i,\;\;
b_ic-cb_i.
$$
The algebra $W^h_{n, \mathbb{C}}$, in other words, differs from $W_{n, \mathbb{C}}$ in the form of the
commutation relations -- in the case of $W^h_{n, \mathbb{C}}$, the coordinate-momenta pairs of generators
commute into $h$ (which is added as a central variable to the algebra; the variable $h$ thus somewhat resembles the
Planck constant) -- and one can return to the non-augmented algebra $W_{n, \mathbb{C}}$ by
\emph{specializing} the augmentation parameter to $h=1$. The augmented Poisson algebra, denoted by $P^h_{n,
\mathbb{C}}$, is defined similarly: one adds the variable $h$ to the commutative polynomial algebra of $2n$
generators and endows it with the Poisson bracket defined as:
$$
\lbrace p_i, x_j\rbrace = h\delta_{ij}.
$$

\smallskip

It can be verified that these new algebras behave in a way almost identical to the one we described in the prequel; in
particular, the notions of tame automorphism, tame (modified) symplectomorphism and homomorphism
\begin{equation*}
\Phi:\Aut^{\leq N}(W^h_{n,\mathbb{C}})\rightarrow \Aut^{\leq N}(P^h_{n,\mathbb{C}}).
\end{equation*}
(defined for a fixed infinite prime) which is identical on the tame points, are present. We also note that the action of
any $h$-augmented automorphism (or, correspondingly, symplectomorphism) on $h$ is necessarily a dilation
$$
h\mapsto \lambda h
$$
where $\lambda$ is a constant. Indeed, the image of $h$ cannot contain monomials proportional to $x_i$ or $p_j$
(otherwise the commutation relations will not hold), and it cannot be a polynomial in $h$ of degree greater than one,
as in that case Also, the proof of the counterpart of the Theorem \ref{thmgabber} is established in a similar fashion.
\begin{thm} \label{thmgabberh}
The mappings
$$
\Phi^h_N:\Aut^{\leq N} W^h_{n,\mathbb{C}}\rightarrow \Aut^{\leq N} P^h_{n,\mathbb{C}}
$$
induced by $\Phi^h$ are morphisms of normalized algebraic varieties.
\end{thm}

\smallskip

As we shall see, one can prove the counterpart to the Conjecture \ref{mainconj} for these augmented algebras, and
then demonstrate that the specialization to $h=1$ yields the isomorphism between automorphism groups of the
non-augmented algebras. The construction of the augmented version of the isomorphism, however, requires to
further modify the algebras by making the commutators between $d_i$ and $x_j$ nonzero for $i\neq j$.

This pair of auxiliary, \emph{skew augmented} algebras, denoted by $W_{n,\mathbb{C}}^h[k_{ij}]$ and
$P_{n,\mathbb{C}}^h[k_{ij}]$ (which correspond to augmented Weyl and Poisson algebras, respectively), are
defined as follows. Let the augmented Poisson generators be denoted by $\xi_i$ with $1\leq i\leq 2n$, which we will
call the main generators, (the passage from $x_i$ and $p_j$ to $\xi_i$ is made for the sake of uniformity of notation
-- and in fact, from the viewpoint of the singularity trick which we use below in order to establish canonicity of the
lifting, all of these generators are on equal footing, unlike the standard form $x_i$ and $p_j$, for which only
symplectic transformations are permitted), and let $[k_{ij}]$ be a skew-symmetric array (a skew matrix) of central
variables. The algebra $P_{n,\mathbb{C}}^h[k_{ij}]$ is generated by $2n$ commuting variables $\xi_i$, the
augmentation variable $h$ and the variables $[k_{ij}]$ (thus being the polynomial algebra in these variables); the
Poisson bracket is defined on the generators $\xi_i$:
$$
\lbrace \xi_i,\xi_j\rbrace = hk_{ij}.
$$
The bracket of any element with $h$ or with any of the $k_{ij}$ is zero.

The skew version of the algebra $W_{n,\mathbb{C}}$ is defined analogously.

\smallskip

It is easily seen that the new algebras essentially share the positive-characteristic properties with $W_n$ and $P_n$,
from which it follows that a mapping
$$
\Phi^{hk}:\Aut W_{n,\mathbb{C}}^h[k_{ij}]\rightarrow \Aut P_{n,\mathbb{C}}^h[k_{ij}]
$$
analogous to $\Phi$ and $\Phi^h$ can be defined for every infinite prime $[p]$. In a manner identical to the previous
section it can be established that this mapping consists of a system of morphisms of the normalized varieties
$\Aut^{\leq N}$, thus yielding the skew augmented analogue of Theorems \ref{thmgabber} and \ref{thmgabberh}.

\begin{thm}\label{thmgabberskew}
The mappings $\Phi^{hk}_N$ are morphisms of normalized varieties.
\end{thm}

\smallskip


The plan of the proof of the main theorem goes as follows. Given that in all three considered cases -- the
non-augmented, the $h$-augmented and the skew augmented case -- the $\Ind$-morphism between (the
normalizations of) $\Ind$-varieties of automorphisms is well defined, we will examine its properties. In particular, we
are going to establish the continuity of the morphism $\Phi^{hk}$ -- or, to be more precise, its restriction to a
certain subspace -- in the power series topology (defined by the choice of grading below). That result, together with
the tame approximation and $\Phi^{hk}$ being the identity map on the tame automorphisms -- a property which
also holds in all three considered cases -- will allow us to prove that the lifted limits of tame sequences are
independent of the choice of the converging sequence (canonicity of lifting) and then demonstrate that the lifted limits
are given by polynomials and not power series, i.e. that the lifted limits are (skew augmented Weyl algebra)
automorphisms. Effectively we will establish the skew augmented version of the Kontsevich conjecture, or rather the
more relevant isomorphism between subgroups $\Aut_{k} $ of automorphisms which act linearly on the auxiliary
variables $k_{ij}$. Most of the conceptually non-trivial topological machinery -- namely, the singularity trick
mentioned in the introduction, are employed at this first stage. In fact, the good behavior of the skew augmented
algebras with respect to the singularity trick is the sole reason for introducing these algebras in our proof.

Next, we will connect the skew augmented algebras with the $h$-augmented algebras by means of a localization
argument. Once this is done, the establishing of canonicity of lifting and polynomial nature of the lifted limits in the
$h$-augmented case becomes a fairly straightforward affair.

Lastly, in order to demonstrate that the results for the $h$-augmented algebras imply the Kontsevich isomorphism,
we will need to specialize to $h=1$. The procedure requires some effort, and in fact extension of the domain for the
constructed inverse morphism will be needed. The procedure will finalize the proof.

\subsection{Continuity of \texorpdfstring{$\Phi^{hk}$}{Phi} and the singularity trick}

We will now study the power series topology induced on subgroups of skew augmented algebra automorphisms
which are linear on $k_{ij}$. As usual, the topology is metric, and it is induced by the grading specified according to
the following assignment of degrees to the generators:

$$\Deg h = 0,$$ $$\Deg k_{ij} = 2,$$ $$\Deg \xi_i = 1.$$ Note that since $h$ and $k_{ij}$ appear as products in the commutation relations, one could assign degree two to the augmentation parameter $h$ and degree zero to the skew-form variables $k_{ij}$, in analogy with the case of augmented algebra $P^h_n$, while essentially preserving the $\Ind$-scheme structure of $\Aut$.

\smallskip

The metric which induces the power series topology is defined as
$$
\rho(\varphi, \psi) = \exp(-\Ht(\varphi-\psi))
$$
where
$$
\varphi-\psi = (\varphi(\xi_1)-\psi(\xi_1),\ldots,\varphi(\xi_{2n})-\psi(\xi_{2n}),\ldots)
$$
is the algebra endomorphism defined by its images (on $\xi_i$, $h$ and $k_{ij}$), and the \textbf{height}
$\Ht(\varphi)$ of an endomorphism is defined as the minimal total degree $m$ such that in one of the generator
images under $\varphi$ a non-zero homogeneous component of degree $m$ exists.

Symbolically, we say that the power series topology is defined via the powers of the augmentation ideal $I$
$$
I=(\xi_1,\ldots, \xi_{2n}, h, \lbrace k_{ij}\rbrace)
$$
just as it is so in the commutative case, when every variable carries degree one.

The system of neighborhoods $\lbrace H_N\rbrace$ of the identity automorphism in $\Aut
P_{n,\mathbb{C}}^h[k_{ij}]$ is defined by setting
$$
H_N = \lbrace g\in \Aut P_{n,\mathbb{C}}^h[k_{ij}]\;:\; g(\eta) \equiv \eta \;(\text{mod}\;I^N)\rbrace
$$
(here $\eta$ denotes any generator in the set $\lbrace \xi_1,\ldots, \xi_{2n},h, \lbrace k_{ij}\rbrace\rbrace$, so that
elements of $H_N$ are precisely those automorphisms which are identity modulo terms which lie in $I^N$; again,
the phrase "mod $I^N$" is short-hand for the distance as defined above).

Similar notions of grading, topology, and system of standard neighborhoods of a point, are valid for the algebra
$W_{n,\mathbb{C}}^h[k_{ij}]$ (once the proper ordering of the generators in the chosen set is fixed). The
neighborhoods of the identity point for this algebra will be denoted by $G_N$.

The point of introducing the skew algebras $W_{n,\mathbb{C}}^h[k_{ij}]$ and $P_{n,\mathbb{C}}^h[k_{ij}]$
is that a certain singularity analysis procedure (the singularity tricks mentioned in the introduction) can be
implemented for these algebras in full analogy with the case of the commutative polynomial algebra processed in our
preceding study \cite{KBYu}, while on the other hand there seems to be no straightforward way to execute the
singularity trick for the algebras $W^h_n$ and $P^h_n$. Furthermore, after adjunction of $k_{ij}^{-1}$ (together
with the entries of the inverse matrix) and extension of scalars (the localization procedure) one can embed the
$h$-augmented $\mathbb{C}$-algebras $W^h_n$ and $P^h_n$ in the skew augmented algebras over the larger
coefficient ring, thus connecting the $h$-augmented automorphisms with the skew augmented ones which are linear
on $k_{ij}$, as we shall see below.

\smallskip

We will establish the continuity of the direct morphism $\Phi^{hk}$ and perform the singularity trick in the following
stable form.

Consider the algebra $P_{n+1,\mathbb{C}}^h[k_{ij}]$ with $(2n+2)$ main generators $\lbrace \xi_1,\ldots,
\xi_{2n}, u,v\rbrace$. Let
$$
\Aut_{u,v,k} P_{n+1,\mathbb{C}}^h[k_{ij}]
$$
denote the set of all automorphisms $\varphi$ of $P_{n+1,\mathbb{C}}^h[k_{ij}]$ such that:

1. $\varphi(\xi_i) = \xi_i + S_i$, where $S_i$ is a polynomial (in $\xi_i$, $u$, $v$, $h$ and $k_{ij}$) such that its
height with respect to $\lbrace \xi_1,\ldots, \xi_{2n}, u,v\rbrace$ is at least two.

2. $\varphi(u) = u$, $\varphi(v) = v$.

3. $\varphi(k_{ij})$ is a $\mathbb{C}$-linear combination of $k_{ij}$, i.e.
$\varphi\in\Aut_{k}P_{n+1,\mathbb{C}}^h[k_{ij}]$.

Define the grading as before: $\xi_i$, $u$, $v$ carry degree one, $h$ carries degree zero, and $k_{ij}$ carry
degree two.

Denote by $H_N^{u,v,k}$ the subgroups of $\Aut_{u,v,k} P_{n+1,\mathbb{C}}^h[k_{ij}]$ consisting of
elements which are the identity map modulo terms of height $N$ with respect to the grading defined above. Also,
the definition is repeated for the skew augmented Weyl algebra $W_{n+1,\mathbb{C}}^h[k_{ij}]$; the resulting
subgroups are denoted by $G_N^{u,v,k}$.

\smallskip

The purpose of the singularity trick set up below is the proof of the following result, which establishes continuity of
the direct morphism.

\begin{prop}\label{skewthetaprop}
If $\Phi^{hk}$ is the restriction of the direct morphism to $\Aut_k$, then
$$
\Phi^{hk}(G_N^{u,v,k}) \subseteq H_N^{u,v,k}
$$
for every $N$.
\end{prop}

The singularity trick is essentially a criterion for an automorphism $\varphi$ to be an element of $H_N^{u,v,k}$,
expressed in terms of \emph{asymptotic behavior of certain parametric families} associated to it. The parametric
families of automorphisms are constructed from $\varphi$ by conjugating it with $\mathbb{C}$-linear changes of the
main generators (the latter are given by the set $\lbrace \xi_1,\ldots, \xi_{2n}\rbrace$). Such parameterized variable
changes are given by $(2n+2)$ by $(2n+2)$ matrices $\Lambda(t)$ with
$$
(\xi_1,\ldots, \xi_{2n},u,v)\mapsto (\xi_1,\ldots, \xi_{2n},u,v)\Lambda(t)
$$
representing the action (such transformations of the main generators induce appropriate mappings of $[k_{ij}]$).
Note that if $\varphi$ is in $H_N^{u,v,k}$, then the conjugation by $\Lambda(t)$ is also in $H_N^{u,v,k}$, as the
action upon $u$ and $v$ is that of $\Lambda(t)\circ \Lambda(t)^{-1}$.

\smallskip

We are going to examine the behavior of such one-parameter families near singularities of $\Lambda(t)$.

Suppose that, as $t$ tends to zero, the $i$-th eigenvalue of $\Lambda(t)$ also tends to zero as $t^{m_i}$, $m_i\in\mathbb{N}$. 

Let $\lbrace m_i,\;i=1,\ldots 2n+2\rbrace$ be the set of degrees of singularity of eigenvalues of $\Lambda(t)$ at
zero. Suppose that for every pair $(i,j)$ the following holds: if $m_i\neq m_j$, then there exists a positive integer
$M$ such that
$$
\text{either\;\;} m_iM\leq m_j\;\;\text{or\;\;}m_jM\leq m_i.
$$
We will call the largest such $M$ the \textbf{order} of $\Lambda(t)$ at $t=0$. As $m_i$ are all set to be positive
integer, the order equals the integer part of $\frac{m_{\text{max}}}{m_{\text{min}}}$.

We now formulate the criterion
\begin{prop}[Singularity trick]\label{singtrick}
An element $\varphi\in\Aut_{u,v,k} P_{n+1,\mathbb{C}}^h[k_{ij}]$ belongs to $H_N^{u,v,k}$ if and only if for
every linear matrix curve $\Lambda(t)$ of order $\leq N$ the curve
$$
\Lambda(t)\circ\varphi\circ\Lambda(t)^{-1}
$$
does not have a singularity (a pole) at $t=0$.
\end{prop}
\begin{proof}
Suppose $\varphi\in H_N^{u,v,k}$ and fix a one-parametric family $\Lambda(t)$. Without loss of generality, we
may assume that the first $2n$ main generators $\lbrace \xi_1,\ldots, \xi_{2n}\rbrace$ correspond to eigenvectors
of $\Lambda(t)$. If $\xi_i$ denotes any of these main generators, then the action of $\Lambda(t)\circ\varphi\circ
\Lambda(t)^{-1}$ upon it reads
$$
\Lambda(t)\circ\varphi\circ \Lambda(t)^{-1}(\xi_i) = \xi_i + t^{-m_i}\sum_{l_1+\cdots+l_{2n} = N}a_{l_1\ldots l_{2n}}t^{m_1l_1+\cdots+m_{2n}l_{2n}}P_{i}(\xi_1,\ldots, \xi_{2n},h,k_{ij}) + S_i
$$
where $P_i$ is homogeneous of total degree $N$ (in the previously defined grading) and the height of $S_i$ is
greater than $N$. One sees that for any choice of $l_1,\ldots,l_{2n}$ in the sum, the expression
$$
m_1l_1+\cdots+m_{2n}l_{2n} - m_i\geq m_{\text{min}}\sum l_j-m_i=m_{\text{min}}N-m_i\geq 0,
$$
so whenever $t$ goes to zero, the coefficient will not go to infinity. The same argument applies to higher-degree
monomials within $S_i$.

The other direction is established by contraposition. Assuming $\varphi\notin H_N^{u,v,k}$, we need to prove the
existence of linear curves with suitable eigenvalue behavior near $t=0$ which create singularities via conjugation with
the given automorphism.

Suppose first that the image of $\xi_1$ under $\varphi$ possesses a monomial which is not divisible by $\xi_1$ or
any $k_{1j}$ ($j\neq 1$). Then one can take $m_1$ and $m_2<m_1$ such that
$$
(N+1)m_2\geq m_1\geq Nm_2
$$
and set the curve $\Lambda(t)$ to be given by a diagonal matrix with entries
$t^{m_1},\;t^{m_2},\;t^{m_2},\ldots$. It is easily checked that conjugation of $\varphi$ by this curve creates a
pole at the coefficient of the chosen monomial.

The general case can be reduced to this special case by means of transformations of the form ($\lambda$ and
$\delta$ are suitable constants)
\begin{gather*}
\xi_1\mapsto \xi_1 + \lambda u + \delta v,\\
k_{ij}\mapsto k_{ij},\;\;1<i,j\leq 2n,\\
k_{1j}\mapsto k_{1j} + \lambda k_{2n+1,j} + \delta k_{2n+2,j},\\
k_{1,2n+1}\mapsto k_{1,2n+1} + \delta k_{2n+2,2n+1},\\
k_{1,2n+2}\mapsto k_{1,2n+2} + \lambda k_{2n+1,2n+2}.
\end{gather*}
Conjugation with these transformations create in the image of $\xi_1$ under the resulting automorphism a monomial
from the previous case. In order to obtain the curve $\Lambda(t)$ from the diagonal curve acting on the conjugated
automorphism, one needs only conjugate it with the inverse of the above transform. The singularity trick is proved.

\end{proof}
The skew augmented Weyl algebra counterpart of the singularity trick is valid.
\begin{cor}\label{singtrickweyl}
An element $\varphi\in\Aut_{u,v,k} W_{n+1,\mathbb{C}}^h[k_{ij}]$ belongs to $G_N^{u,v,k}$ if and only if for
every linear matrix curve $\Lambda(t)$ of order $\leq N$ the curve
$$
\Lambda(t)\circ\varphi\circ\Lambda(t)^{-1}
$$
does not have a singularity at $t=0$.
\end{cor}
The proof of this statement is essentially the same as that of Proposition \ref{singtrick}. Note that thanks to the
choice of grading -- the one in which the degree of all $k_{ij}$ is two -- the reordering of the non-commuting
variables in a word cannot produce monomials of smaller total degree.

\smallskip

The implementation of the singularity trick in the proof of Proposition \ref{skewthetaprop} requires also the
following general fact.
\begin{lem} \label{lem1} Let
$$
\Phi: X\rightarrow Y
$$
be a morphism of affine algebraic sets, and let $\varphi(t)$ be a curve (more simply, a one-parameter family of
points) in $X$. Suppose that $\varphi(t)$ does not tend to infinity as $t\rightarrow 0$. Then the image $\Phi
\varphi(t)$ under $\Phi$ also does not tend to infinity as $t\rightarrow 0$.
\end{lem}
The proof of the Lemma is an easy exercise and is left to the reader.

Proposition \ref{skewthetaprop} is now an elementary consequence of the above Lemma together with the
singularity trick (Proposition \ref{singtrick} and Corollary \ref{singtrickweyl}). Indeed, let us assume the contrary --
i.e. that for some $N$
$$
\Phi^{hk}(G_N^{u,v,k}) \nsubseteq H_N^{u,v,k}.
$$
Then there exists an element $\varphi\in G_N^{u,v,k}$ such that its image $\Phi^{hk}(\varphi)\notin
H_N^{u,v,k}$. By Proposition \ref{singtrick}, there is a linear automorphism (matrix) curve $\Lambda(t)$ of order
$\leq N$ such that the curve
$$
\Lambda(t)\circ \Phi^{hk}(\varphi)\circ \Lambda(t)^{-1}
$$
has a pole at $t=0$. Since $\Phi^{hk}$ is point-wise stable on linear variable changes, the latter curve is the image
under $\Phi^{hk}$ of the curve
$$
\Lambda(t)\circ \varphi\circ \Lambda(t)^{-1}.
$$
By our assumption, $\varphi\in G_N^{u,v,k}$; therefore, by Corollary \ref{singtrickweyl}, the curve above has no
singularity at $t=0$. But then the statement that the curve $$\Lambda(t)\circ \Phi^{hk}(\varphi)\circ
\Lambda(t)^{-1}$$ -- which is the image of the former curve under the morphism $\Phi^{hk}$ -- has a singularity
at $t=0$ yields a contradiction with Lemma \ref{lem1}. Proposition \ref{skewthetaprop} is proved.

The immediate consequence of Proposition \ref{skewthetaprop} is the following result.
\begin{thm}\label{skewcontthm}
The mapping
$$
\Phi^{hk}: \Aut_{u,v,k} W_{n,\mathbb{C}}^h[k_{ij}]\rightarrow \Aut_{u,v,k} P_{n,\mathbb{C}}^h[k_{ij}]
$$
is continuous in the power series topology defined at the start of the section.
\end{thm}
This was the main objective of the singularity trick, and this result will provide the means to establish the canonicity of
the symplectomorphism lifting procedure we define in the next subsection.

\subsection{Lifting in the \texorpdfstring{$h$}{h}-augmented and skew augmented cases}

We now proceed with the resolution of the symplectomorphism lifting problem for both augmented and skew
augmented algebras. More specifically, we will show how the results of the singularity analysis procedure conducted
in the previous subsection provide for a way to construct the inverse to the homomorphisms $\Phi^h$ and
$\Phi^{hk}$.


\smallskip

Suppose given an automorphism $\varphi\in \Aut P^h_{n,\mathbb{C}}$ of the $h$-augmented Poisson algebra.
Without loss of generality, we may assume that the linear part of $\varphi$ is the identity matrix: indeed, one can
compose $\varphi$ with tame automorphisms (tame approximation of automorphisms of $P^h_{n,\mathbb{C}}$ is
valid according to an argument similar to that of \cite{KGE}), so that the linear part of the resulting automorphism is
the identity map; also the morphism $\Phi^h$ is point-wise stable on tame automorphisms.

We add two more $h$-Poisson variables (and lift $\varphi$ to an automorphism of the new algebra by demanding it
be stable on the new generators) and, correspondingly, consider the skew Poisson version -- the algebra
$P_{n+1,\mathbb{C}}^h[k_{ij}]$ with the last two variables denoted by $u$ and $v$. Our objective is to realise
the algebra $P^h_{n,\mathbb{C}}$ as a subalgebra in an appropriate localization of
$P_{n+1,\mathbb{C}}^h[k_{ij}]$. To that end, we consider the algebra $P_{n+1,\mathbb{C}}^h[k_{ij}]$  and
transform the main generators
$$
\lbrace \xi_1,\ldots, \xi_{2n},u,v\rbrace
$$
to
$$
\lbrace x_1,\ldots, x_{2n},u,v\rbrace
$$
with $\lbrace x_i, u\rbrace = 0$ and $\lbrace x_i, v\rbrace = 0$. The change of the generating set is required to
properly define the action of $\varphi$, so that it will be an automorphism and will be in agreement with the
conditions of Proposition \ref{skewthetaprop}. The variable change is done according to
$$
x_i = \xi_i - \alpha_i u - \beta_i v
$$
with $\alpha_i = k_{i,2n+2}k_{2n+1,2n+2}^{-1}$ and $\beta_i = - k_{i,2n+1}k_{2n+1,2n+2}^{-1}$ for
$i=1,\ldots, 2n$. We extend the coefficient ring by adding the necessary variables. The new generators $\lbrace
x_1,\ldots, x_{2n}\rbrace$ commute according to
\begin{gather*}
\lbrace x_i, x_j\rbrace = h(k_{ij} - \alpha_jk_{i,2n+1}+\alpha_ik_{j,2n+1}-\beta_jk_{i,2n+2}+\\
\beta_ik_{j,2n+2}+ (\alpha_i\beta_j - \alpha_j\beta_i)k_{2n+1,2n+2}) = h\tilde{k}_{ij}.
\end{gather*}
Note that the new commutation relation matrix, which we denote by $[\tilde{k}_{ij}]$\footnote{We exclude $u$
and $v$, so that $i$ and $j$ run from $1$ to $2n$.}, is again skew-symmetric, and that its entries are
$\mathbb{C}$-polynomial in the entries of the initial matrix and their inverses.

We now reduce the matrix $[\tilde{k}_{ij}]$ to the standard form (corresponding to the algebra $P^h_n$) by
transforming $\lbrace x_1,\ldots, x_{2n}\rbrace$ to $\lbrace q_1,\ldots, q_n,p_1,\ldots, p_n\rbrace$ with
$$
\lbrace p_i, q_j \rbrace = h\delta_{ij}.
$$
The new variables $p_i$ and $q_j$ are expressed as linear combinations of $x_1,\ldots, x_{2n}$ with coefficients
in the appropriate polynomial ring.

The algebra $P^h_{n,\mathbb{C}}$ is therefore a subalgebra of the algebra generated by $$\lbrace q_1,\ldots,
q_n,p_1,\ldots, p_n, u,v\rbrace$$ (together with $h$ as the augmentation variable), as the Poisson bracket takes its
proper form after the standard form reduction, while $\mathbb{C}$ is a subring of the coefficient ring.

We extend our automorphism $\varphi$ to act on this algebra: on $p_i$ and $q_j$ its action is given by definition,
and we impose $\varphi(u) = u$, $\varphi(v) = v$ and $\varphi(k_{ij}) = k_{ij}$. Thus, starting from $\varphi$ we
have arrived at an automorphism $\bar{\varphi}$ of the localized skew Poisson algebra.

With respect to the skew augmented algebra generator set  $\lbrace \xi_1,\ldots, \xi_{2n},u,v\rbrace$ this
automorphism is generally not polynomial in $k_{ij}$, although it always will be polynomial in $h$. In order to
construct from it an automorphism of the skew algebra, we need to get rid of the denominators first. This is
accomplished by the following lemma.
\begin{lem}\label{conjugationlemma}
For every $\bar{\varphi}$ constructed as above, there is a polynomial $P$ in $k_{ij}$, such that conjugation of
$\bar{\varphi}$ with the transformation
$$
(\xi_1,\ldots, \xi_{2n},u,v)\mapsto (P\xi_1,\ldots, P\xi_{2n},Pu,Pv),\;\;h\mapsto P^2h
$$
is polynomial in $k_{ij}$. The polynomial $P$ depends only on the two systems of algebra generators.
\end{lem}
\begin{proof}
Indeed, the denominators in the expression for $\bar{\varphi}$ are polynomial in $k_{ij}$ coming from the
separation of the $(u,v)$-plane and the standard form reduction (at which point the determinant of $[\tilde{k}_{ij}]$
makes its contribution). One can therefore find appropriate $P(k_{ij})$ to cancel these denominators. Furthermore,
the polynomial $P$ depends only on the two generator systems -- more specifically, on the transformation matrix
between those systems.
\end{proof}
We denote the result of the conjugation of Lemma \ref{conjugationlemma} by $\varphi^P$. The images of the main
generators (both in the cases of the initial -- skew -- generators as well as those which correspond to the standard
form) under $\varphi^P$ are, by Lemma \ref{conjugationlemma}, polynomial in $k_{ij}$, and are also by
construction polynomial in $h$.

The automorphism $\varphi^P$, when acting upon the standard form generators
$$
\lbrace q_1,\ldots, q_n,p_1,\ldots, p_n\rbrace
$$
can be viewed as an automorphism of the $h$-augmented Poisson algebra $P^h_{n,\mathbb{C}[\lbrace
hk_{ij}\rbrace]}$ over the polynomial ring $\mathbb{C}[\lbrace hk_{ij}\rbrace]$. The $\mathbb{Z}$-grading of
this algebra is specified by assigning degree $1$ to the main generators and degree $0$ to $h$ and all $k_{ij}$. As
an automorphism of this Poisson algebra, $\varphi^P$ admits, by an argument virtually identical to the main result of
\cite{KGE}, a tame automorphism (symplectomorphism) sequence converging to it in the power series topology
induced by the above grading. Let us fix such a sequence and denote it by $\lbrace \psi_m\rbrace$.

Every element $\psi_k$ of the tame sequence is such that the images under $\psi_m$ of the generators $\lbrace
q_1,\ldots, q_n,p_1,\ldots, p_n\rbrace$ are polynomial in $h$ and $k_{ij}$. Importantly, the tame sequence
$\lbrace \psi_m\rbrace$ converging to the $\mathbb{C}[\lbrace hk_{ij}\rbrace]$-automorphism $\varphi^P$ can
be connected to the original $h$-augmented symplectomorphism $\varphi$ by inversion of the procedure which
leads to the definition of $\varphi^P$. Precisely, we have the following statement.
\begin{prop}\label{goodbehaviorprop}
For every $\varphi$ and every sequence $\lbrace \psi_m\rbrace$ converging to $\varphi^P$ as above, there is a
sequence $\lbrace \sigma_m\rbrace$ of tame $\mathbb{C}$-symplectomorphisms of $P_{n,\mathbb{C}}$
converging to $\varphi$ with respect to the topologies with $\Deg h = 0$ and $\Deg h = 2$.
\end{prop}
\begin{proof}
The sequence $\lbrace \sigma_m\rbrace$ is constructed from $\lbrace\psi_m\rbrace$ by reversing the conjugation
by $P$ and disposing of the stable variables $u$ and $v$. We note that the conjugation is a group homomorphism,
which means that it suffices to prove that the reverse conjugation disposes of $k_{ij}$ in every elementary tame
automorphism (as $\psi_m$ are composition of elementary automorphisms). The latter property, however, is
obvious.

The convergence of $\lbrace \sigma_m\rbrace$ follows immediately from Lemma \ref{conjugationlemma} and
Lemma \ref{convergencelemma} below.
\end{proof}

When acting upon the localized skew Poisson algebra generators $\lbrace \xi_1,\ldots, \xi_{2n}\rbrace$, the
augmented symplectomorphisms $\psi_m$ need not be polynomial in $k_{ij}$, and therefore $\psi_m$ are not in
general images of automorphisms of the skew Poisson algebra under localization. This is remedied by application of
Lemma \ref{conjugationlemma}: one can find a polynomial $P_1$ in the variables $k_{ij}$, such that the
conjugation of every element $\psi_m$ of the tame sequence with the mapping
$$
(\xi_1,\ldots, \xi_{2n},u,v)\mapsto (P_1\xi_1,\ldots, P_1\xi_{2n},P_1u,P_1v),\;\;h\mapsto P_1^2h
$$
yields an automorphism of the localized skew Poisson algebra polynomial in $k_{ij}$. Again, as in Proposition
\ref{goodbehaviorprop}, one can return to a sequence of tame symplectomorphisms of $P_{n,\mathbb{C}}$ by
reversing the conjugation.

We then have the following statement.
\begin{lem}\label{convergencelemma}
The sequence $\lbrace \psi^{P_1}_m\rbrace$ converges to the conjugated automorphism $(\varphi^P)^{P_1}$ in
the power series topology with $\Deg h = \Deg k_{ij} = 0$ as well as in the power series topology with $\Deg h =
0$, $\Deg k_{ij} = 2$.
\end{lem}
\begin{proof}
The first half of the statement follows from the construction of the tame sequence and from the observation that, due
to the fact that the two coordinate systems are connected by a transformation that has zero free term, the height of
the polynomials $P$ and $P_1$ is at least one.

One then obtains convergence in the power series topology relevant to the singularity trick (Proposition
\ref{skewthetaprop}) from that in the approximation power series topology by noting that giving a non-zero degree
to $k_{ij}$ may only make the consecutive approximations closer to the limit in the corresponding metric.

\end{proof}

The sequence $\lbrace \psi^{P_1}_m\rbrace$ can, due to its polynomial character with respect to $k_{ij}$, be
thought of as a sequence of tame automorphisms of the skew Poisson algebra $P_{n+1,\mathbb{C}}^h[k_{ij}]$
over $\mathbb{C}$ converging to $(\varphi^P)^{P_1}$. We now take the pre-image of the sequence $\lbrace
\psi^{P_1}_m\rbrace$ under the morphism $\Phi^{hk}$ to obtain the sequence
$$
\lbrace \sigma'_m\rbrace=\lbrace (\Phi^{hk})^{-1}(\psi^{P_1}_m)\rbrace
$$
of automorphisms of the skew augmented Weyl algebra $W_{n+1,\mathbb{C}}^h[k_{ij}]$. We now may take the
formal limit of this sequence, which is ostensibly dependent on the choice of the convergent tame sequence $\lbrace
\psi_m\rbrace$, apart from the point $\varphi$ itself. This limit, which we denote by
$$
\Theta^h_P(\varphi, \lbrace \psi_m\rbrace)
$$
to reflect the dependence on the sequence, is given by formal power series in the skew augmented Weyl generators.
Applying the inverse to the conjugations performed earlier and disposing of the stable variables (whose presence is
justified by the form of the singularity trick and is therefore needed in the proof of independence of the choice of the
convergent sequence, as we shall see below), we arrive at a vector of formal power series (the entries of which
correspond to images of the generators) in the generators of the $h$-augmented Weyl algebra
$W^h_{n,\mathbb{C}}$.

\smallskip

The most important consequence of Theorem \ref{skewcontthm} is the independence of the lifted sequence's formal
limit of the choice of the approximating tame sequence $\lbrace \psi_m\rbrace$. We have the following proposition.

\begin{prop} \label{propcanon}
Let $\varphi$ be an automorphism of $P^h_{n,\mathbb{C}}$ and let
$$
\psi_1,\;\ldots,\; \psi_m,\;\ldots
$$
and
$$
\psi'_1,\;\ldots,\; \psi'_m,\;\ldots
$$
be two sequences of tame automorphisms which converge to $\varphi^P$ as in the construction above. Then the
lifted sequences
$$
\lbrace (\Phi^{hk})^{-1}(\psi^{P_1}_m)\rbrace\;\;\text{and}\;\;\lbrace (\Phi^{hk})^{-1}(\psi'^{P_1}_m)\rbrace
$$
converge to the same automorphism of the power series completion of the skew augmented Weyl algebra. This
means that one must have
$$
\left((\Phi^{hk})^{-1}(\psi^{P_1}_m)\right)^{-1}\circ (\Phi^{hk})^{-1}(\psi'^{P_1}_m)\equiv \Id\;(\text{mod}\;I^{N(k)})
$$
with $N(k)\rightarrow \infty$ as $k\rightarrow \infty$.
\end{prop}
\begin{proof}
This result follows immediately from the continuity of $\Phi^{hk}$ established in the previous subsection.
\end{proof}

The meaning of this Proposition to the symplectomorphism lifting problem is clear: in our construction, the formal
limit
$$
\Theta^h_P(\varphi, \lbrace \psi_m\rbrace)
$$
is independent of $\lbrace \psi_m\rbrace$ and is therefore a well-defined function of the point $\varphi$.
Furthermore, as can be inferred directly from the tame approximation, this function is homomorphic -- it preserves
the group structure given by composition of automorphisms. As the conjugations are also homomorphisms, we
conclude that $h$-augmented symplectomorphisms $\varphi \in \Aut P^h_{n,\mathbb{C}}$ are lifted
homomorphically to endomorphisms of the power series completion $\hat{W}^h_{n,\mathbb{C}}$ of the
$h$-augmented Weyl algebra. For an augmented symplectomorphism $\varphi$, we denote its image with respect
to the lifting map by
$$
\Theta^h(\varphi).
$$

Our next objective is to demonstrate that for every symplectomorphism $\varphi$, the image $\Theta^h(\varphi)$
under the lifting map is in fact an automorphism of the $h$-augmented Weyl algebra. In fact, it remains to show only
that the generator images with respect to $\Theta^h(\varphi)$ cannot be given by infinite series: indeed, that would
imply that $\Theta^h(\varphi)$ is an $h$-augmented Weyl endomorphism; the invertibility of the lifted mapping
follows from the canonicity of lifting: indeed, $\Theta^h$ not only preserves compositions but also maps inverses to
inverses, therefore for any symplectomorphism $\varphi$ the mapping $\Theta^h(\varphi^{-1})$ will be the inverse
of $\Theta^h(\varphi)$.

Alternatively, one can arrive at the invertibility after one shows the polynomial character of the lifted endomorphism:
it is known \cite{Tsu2, K-BK2} that the direct morphism $\Phi$ -- and hence, by a straightforward extension of the
argument in the aforementioned work, its augmented analogue $\Phi^h$ -- distinguishes automorphisms, i.e. the
image of a non-automorphism cannot be an automorphism.

\smallskip

The main properties of the mapping $\Theta^h$ can now be summarized in the following way.
\begin{prop}\label{thetaproperties}

1. There exists a well-defined mapping $\Theta^h$ whose domain is $\Aut P^h_{n,\mathbb{C}}$ and whose
codomain lies in the set of automorphisms of the power series completion of $W^h_{n,\mathbb{C}}$.

2. $\Theta^h$ is a group homomorphism.

3. For a fixed $\varphi$, the coordinates of $\Theta^h(\varphi)$ -- i.e. the coefficients with respect to the fixed
generator basis decomposition -- are given by polynomials in the coordinates of $\varphi$.

4. The skew augmented analogue $\Theta^{hk}$ of $\Theta^{h}$ is continuous in the power series topology.

\end{prop}

\begin{proof}
The first two statements follow immediately from the construction. The third statement follows from the fact that the
lifting is independent of the approximating sequence: indeed, that implies that the coefficients of the lifted limit are
read off any valid approximating sequence, or more precisely its finite subset (consisting of first several elements, as
the limit symplectomorphism is polynomial). But then the coefficients are polynomial in the coordinates of the lifted
tame elements, and since only a finite number of them suffices, they are also polynomial in the coordinates of the
initial symplectomorphism.

The continuity of $\Theta^{hk}$ is established in a manner identical to that of $\Phi^{hk}$ -- namely with the help
of Proposition \ref{singtrick} a proof similar to that of Theorem \ref{skewcontthm} can be executed. Note that it
follows from the third statement of this Proposition that $\Theta^{hk}$ fulfills the conditions of Lemma \ref{lem1},
therefore by combining the previously proved statements with the analogous steps for the lifting map, one can show
that
$$
\Phi^{hk}(G_N^{u,v,k}) = H_N^{u,v,k}.
$$
\end{proof}

\subsection{The lifted limit is polynomial}
We proceed with establishing the polynomial character of the image $\Theta^h(\varphi)$.
\begin{thm}\label{thmpolynomial}
Let
$$
\Theta^h: \Aut P^h_{n,\mathbb{C}}\rightarrow \Aut \hat{W}^h_{n,\mathbb{C}}
$$
be the lifting homomorphism constructed in the previous section and let, as before, $x_1,\ldots, x_n, d_1,\ldots,
d_n,h$ denote the generators of $W^h_{n,\mathbb{C}}$ together with its deformation parameter. Then, for every
augmented symplectomorphism $\varphi\in \Aut P^h_{n,\mathbb{C}}$, the images
$$
\Theta^h(\varphi)(x_1),\ldots, \Theta^h(\varphi)(d_n)
$$
are polynomials in $x_i$, $d_i$ and $h$.
\end{thm}
\begin{proof}
Suppose that, contrary to the statement of the theorem, for a fixed $\varphi$ there is an index $i$ such that,
say,\footnote{The case for $d_i$ is processed analogously.} $\Theta^h(\varphi)(x_i)$ is a true infinite series of Weyl
monomials.

Let $\lambda$ be a parameter and let
$$
\tau_\lambda: (x_1,\ldots, d_n, h)\mapsto (\lambda x_1,\ldots,\lambda  d_n,\lambda^2 h)
$$
denote the family of dilation transformations parameterized by $\lambda$. For fixed $\varphi$, define
$$
\varphi_\lambda = \tau_\lambda^{-1}\circ\varphi\circ \tau_\lambda
$$
to be the parametric family of $h$-augmented symplectomorphisms constructed by conjugating $\varphi$ with the
dilations.

We introduce a pair of auxiliary variables $u$ and $v$, $\lbrace v,u\rbrace = h$ (their Weyl counterparts will be
also denoted by $u$ and $v$ -- obviously it does not create any ambiguities) and, for a fixed large enough positive
integer $k$, define the following parametric family of linear transformations
$$
\psi_\lambda: u\mapsto u + \lambda^kx_i,\;\;p_i\mapsto p_i - \lambda^k v
$$
(all other generators are unchanged). As always, we extend the action of $\varphi$ to the auxiliary variables by
setting $\varphi(u) = u$ and $\varphi(v) = v$ (while the dilation extends to $(u,v)\mapsto (\lambda u, \lambda v)$).
Consider the following parametric family of $h$-augmented symplectomorphisms:
$$
\varphi_{t,\lambda}= \varphi_\lambda\circ \psi_\lambda\circ\varphi_{\lambda}^{-1}.
$$
The conjugation of $\varphi$ by the inverse to the dilation $\tau_\lambda$ amounts to multiplying each homogeneous
component of degree $m$ by $\lambda^{1-m}$. Therefore (as $\varphi$ is polynomial) for large enough $k$, the
curve $\varphi_{t,\lambda}$ can be continuously extended by its limit at $\lambda = 0$ -- namely, by the identity
symplectomorphism. Continuity is understood in the sense of continuous dependence of coordinates of the
symplectomorphism on the parameter.

Now, the image of $u$ under the lifted curve $\Theta^h(\varphi_{t,\lambda})$ is
$$
u + \lambda^k\Theta^h(\varphi_\lambda)(x_i)
$$
and, as by assumption $\Theta^h(\varphi)(x_i)$ is an infinite series, is itself an infinite series. But since $\Theta^h$ is
identical on linear transformations, we have
$$
\Theta^h(\varphi_\lambda) = \tau_\lambda^{-1}\circ\Theta^h(\varphi) \circ \tau_\lambda
$$
from which it follows that in the image $\Theta^h(\varphi_{t,\lambda})(u)$ there will be monomials with coefficients
proportional to $\lambda^{-m}$ for $m$ greater than any fixed arbitrary natural number.

However, it follows from the third statement in Proposition \ref{thetaproperties} that the coordinates (i.e.
coefficients of Weyl monomials) of the lifted symplectomorphism are continuous functions of the coordinates of the
symplectomorphism, therefore since the curve $\varphi_{t,\lambda}$ is regular at $\lambda = 0$, then so must also
be its image under $\Theta^h$ -- a contradiction.
\end{proof}

We can now combine this theorem with the results of the previous subsection.
\begin{thm} \label{inversethm}
The lifting homomorphism $\Theta^h$ is the inverse to the direct homomorphism $\Phi^h$.
\end{thm}
\begin{proof}
Indeed, Theorem \ref{thmpolynomial} shows that the compositions $\Phi^h\circ\Theta^h$ and
$\Theta^h\circ\Phi^h$ are well defined. In order to prove that, say, $\Phi^h\circ\Theta^h=\Id$ one changes the basis
of generators (and handles the extension of the base ring as in the prequel) to that of the skew augmented algebra
and uses the fact that $\Phi^{hk}\circ\Theta^{hk}$ coincides with the identity map on the dense subset of tame
symplectomorphisms and hence must be the identity map everywhere (the spaces in question are metric spaces, in
particular they are Hausdorff).
\end{proof}
This theorem has one important corollary.
\begin{cor}\label{correduction}
The lifting map $\Theta^h$ and the direct map $\Phi^h$ are consistent with modulo infinite prime reductions. That is,
for any symplectomorphism $\varphi$, almost all its modulo $p_m$ ($m$ in the index set in the ultraproduct
decomposition) reductions coincide with the (twisted by inverse Frobenius) restrictions to the center of the modulo
$p_m$ reductions of its lifting $\Theta^h(\varphi)$.
\end{cor}
\begin{proof}
This statement is essentially a reformulation of Theorem \ref{inversethm}, if one takes into account the construction
of $\Phi^h$ -- the image of the Weyl algebra automorphism being reconstructed from the ultraproduct of positive
characteristic automorphisms restricted to the center.

\end{proof}

\begin{remark}
Alternatively, one could come up with a line of reasoning more conforming to the combinatorial side of the
constructions employed thus far. If there were a symplectomorphism $\varphi$ which, after lifting and subsequent
direct homomorphism action (ultraproduct decomposition followed by restriction to the center) produces a different
symplectomorphism $\varphi'$, then one may take a sequence of elementary symplectomorphism whose total action
on $\varphi$ maps it to an element which is the identity map modulo terms of degree $N$. If $\varphi'\neq \varphi$,
then the action of the same elementary sequence on $\varphi'$ will produce an element which admits a term of
degree $1<M<N$. The two objects are then mapped to automorphisms of the skew augmented algebra and lifted
(we also note that the mapping induced by the change of basis is continuous, as it is in essence a dilation of
generators by a polynomial of height at least one). By the singularity trick (Proposition \ref{singtrick}) such an object
can be conjugated by an appropriate linear variable change in order to produce a singular curve. Now, by
construction, the skew version of $\varphi$ lifts to a skew Weyl automorphism, and again by the singularity trick
(essentially by the continuity of $\Theta^{hk}$) the lifting of the partially approximated automorphism (i.e. after the
action of the elementary automorphisms) cannot have a singularity of order $\leq N$. However, as the skew version
of $\varphi$ and $\varphi'$ (as well as its partial approximation) corresponds in the ultraproduct to restrictions to the
center, the restriction of an object which is not singular of order $\leq N$ must also be non-singular of order $\leq
N$, in contradiction with the existence of $\varphi'$.
\end{remark}

The $h$-augmented counterpart of the Kontsevich conjecture follows at once from Theorem \ref{inversethm}.
\begin{thm}\label{mainthmh}
The homomorphism
$$
\Phi^h :\Aut W^h_{n,\mathbb{C}}\rightarrow \Aut P^h_{n,\mathbb{C}}
$$
is an isomorphism.
\end{thm}

As the map $\Phi^h$ is closely related to the morphism $\Phi^{hk}$ of the skew augmented case, and,
correspondingly, as $\Theta^h$ is related to the lifting map for the skew augmented symplectomorphisms (essentially
given by $\Theta^h_P$), we obtain another important consequence of Theorem \ref{thmpolynomial}.
\begin{thm}\label{mainthmskew}
Let $\Aut_{k}P_{n,\mathbb{C}}^h[k_{ij}]$ and $\Aut_{k}W_{n,\mathbb{C}}^h[k_{ij}]$ denote the
automorphism subgroups of the skew Poisson and Weyl algebras consisting of those automorphisms that map
$k_{ij}$ to $\mathbb{C}$-linear combinations of $k_{ij}$. Then the mapping
$$
\Phi^{hk}:\Aut_{k}W_{n,\mathbb{C}}^h[k_{ij}]\rightarrow \Aut_{k}P_{n,\mathbb{C}}^h[k_{ij}]
$$
is an isomorphism.
\end{thm}
Theorem \ref{mainthmskew} is significant to the proof of the independence of the morphism $\Phi$ of the choice of
infinite prime given in \cite{K-BE2}.

\subsection{Specialization}

Now that we have established the isomorphism between the automorphism groups of the $h$-augmented Weyl and
Poisson algebras, the proof of the Main Theorem reduces to specializing to $h=1$. That, however, is by no means a
trivial affair, as one needs to take care not only of the automorphisms polynomial in $h$ (for which the existence of
lifting has been established), but also of those which are polynomial in $h^{-1}$.

The necessity of extension of the domain of the lifting map can be seen from the following argument. Suppose
$\varphi^h$ is an automorphism of the $h$-augmented algebra $P^h_{n,\mathbb{C}}$ which acts as the identity
map on $h$. Since it is stable on $h$, it corresponds to an automorphism of the $\mathbb{C}[h]$-algebra (where
$h$ is a parameter and not a generator, which the can be effectively adjoined to the ground field) generated by
$x_i,\;p_j$ with the Poisson bracket containing $h$. This object, after appropriate localization, maps to an
automorphism of the (augmented) Poisson algebra $P^h_n$ with the ground field $\mathbb{C}(h)$. On the other
hand, any automorphism $\varphi$ of $P_{n,\mathbb{C}}$ can be made into a $\mathbb{C}(h)$-automorphism
$\varphi^h$ by introducing a scalar $h$ and conjugating $\varphi$ with a mapping
$$
x_i' = hx_i,\;\;p_j' = p_j.
$$
The resulting transformation will be an automorphism of the Poisson $\mathbb{C}(h)$-algebra with the bracket as in
the augmented algebra $P^h_n$,\emph{ however in general the images of the generators under this automorphism
will contain negative powers of $h$.  } Its specialization to $h=1$ returns it to $\varphi$. Therefore, every
polynomial symplectomorphism has a pre-image under specialization of the $\mathbb{C}(h)$-algebra
automorphisms. The conclusion is that Theorem \ref{mainthmh} does not immediately imply that $\Phi$ is an
isomorphism; rather, the domain of the lifting map $\Theta^h$ needs to be extended to the points with rational
dependency on the augmentation parameter, at which point the claim that $\Phi$ has an inverse given by the
specialization of the extended lifting map $\Theta^h$ becomes valid.

\smallskip

The extension of the domain is accomplished in the following way. For a symplectomorphism $\varphi$ which is
rational in $h$, we will construct images $\Theta^h(\varphi)(x_i)$ and $\Theta^h(\varphi)(d_i)$ one by one by
introducing auxiliary variables and twisting the symplectomorphism in order to create an object polynomial in $h$ --
using the fact that the action of $\Theta^h$ is well defined -- from which the form of the corresponding lifted
generator image may be extracted. As the procedure yields not all of the images simultaneously, we will need to
check its canonical nature as well as verify the commutation relations.

We fix $i$, $1\leq i\leq n$, which corresponds to the image of $x_i$, and introduce a pair $u$, $v$ of auxiliary
variables which are extra $x$ and $p$ with respect to the augmented Poisson bracket. We also add the
corresponding augmented Weyl variables which we denote by $\hat{u}$ and $\hat{v}$. Let
$$
\lambda = h^k
$$
be $k$-th power of the augmentation parameter, for large enough $k$. Define the automorphism
$$
\psi_{\lambda}: u\mapsto u + \lambda x_i,\;\;p_i\mapsto p_i - \lambda v.
$$
We extend $\varphi$ to the new algebra by its identical action on the auxiliary variables and denote the extended
map by $\varphi_a$. Consider the following twisted automorphism:
$$
\varphi_{t,\lambda} = \varphi_a\circ\psi_{\lambda}\circ\varphi_a^{-1}.
$$
As $k$ can be taken arbitrarily large, the mapping $\varphi_{t,\lambda}$ will be polynomial in $h$ for all $k>k_0$
(where $k_0$ depends on $\varphi$ but is finite for the fixed automorphism). We can now read off the expression
for the image of $x_i$ under $\varphi$ from the action of $\varphi_{t,\lambda}$ on the auxiliary variable $u$:
$$
\varphi_{t,\lambda}(u) = u + h^k\varphi(x_i);
$$
the expression is polynomial in $h$, and $\varphi_{t,\lambda}$ thus admits lifting to an automorphism of the
$h$-augmented Weyl algebra $W^h_{n,\mathbb{C}}$. As we will show in a moment, the action of the lifted
automorphism on $\hat{u}$ will be given by the expression
$$
\hat{u} + h^k P_i(x_1,\ldots, d_n,h)
$$
($P_i$ is polynomial in $x_1,\ldots, d_n$ and rational -- or, more precisely, Laurent-polynomial -- in $h$) so that
one can set
$$
\hat{\varphi}(x_i) = P_i(x_1,\ldots, d_n,h)
$$
and thus, for all $x_i$, obtain the action of the lifted symplectomorphism. Switching the roles of $x_i$ and $d_i$
allows for reconstruction of the images of $d_i$. As a result, we get a mapping
$$
\hat{\varphi}: (x_1,\ldots, d_n)\mapsto (P_1(x_1,\ldots, d_n,h),\ldots, Q_n(x_1,\ldots, d_n,h)).
$$
Note, however, that as the lifting is not defined for the components in the composition (with the exception of
$\psi_{\lambda}$), one cannot immediately conclude that the image of $\hat{u}$ under the lifted map will be of the
form as above, or that the parts which depend on $x_1,\ldots, d_n$ will combine to a well defined automorphism --
these properties need to be verified.

\smallskip

The first step is to ensure the constructed mapping $\hat{\varphi}$ is well defined (is canonical with respect to
$\varphi$). This property in fact follows from the consistency with modulo infinite prime reductions given by
Corollary \ref{correduction}, and is manifested in the form of the next two lemmas.

\begin{lem}\label{twistedliftcanon}
Suppose $\theta$ is an $h$-augmented polynomial symplectomorphism over $\mathbb{C}$. Denote by $\lbrace
\theta_p\rbrace$ the sequence of characteristic $p$ symplectomorphisms representing its modulo $[p]$ reduction.
For a generic element $p$ in a sequence representing $[p]$, denote the Weyl generators by $x_1,\ldots, x_n,
d_1,\ldots, d_n$ and the corresponding $p$-th powers generating the center of the Weyl algebra over
$\mathbb{F}_p$ by  $\xi_1,\ldots, \xi_n, \eta_1,\ldots, \eta_n$. Then, for almost all $p$ in $[p]$ (in the sense of the
ultrafilter), the image under $\theta_p$ of every central generator admits a unique pre-image Weyl polynomial
$\hat{H}$ with respect to taking the $p$-th power and pulling back the coefficients by the inverse Frobenius
automorphism.
\end{lem}
\begin{proof}
We prove the statement for $H = \theta_p(\xi_i)$ -- the case $\eta_j$ is identical.

Suppose first that
$$
\theta_p(\xi_i) = \xi_i = x_i^p.
$$
Then the Newton polyhedron of the image $\theta_p(\xi_i)$ has only one vertex, therefore -- as taking the $p$-th
power only dilates the Newton polyhedron -- the polynomial $\hat{\theta}_p(x_i)$ must be equal to $x_i$.

The general case uses Corollary \ref{correduction}, which states that modulo $[p]$ reductions of $\theta$ and its
lifting $\hat{\theta}$ are consistent -- that is, for almost all $p$ in $[p]$, the restriction of $\hat{\theta}_p$ to the
center (twisted by the inverse Frobenius acting on the coefficients) coincides with $\theta_p$. The application is as
follows. Suppose
$$
H = \theta_p(\xi_i)
$$
is the image of $\xi_i$. From Corollary \ref{correduction} we know that
$$
H = \Fr_*^{-1}\hat{\theta}_p(x_i^p)
$$
where $\Fr_*^{-1}$ is the action of the inverse Frobenius automorphism on the coefficients of the polynomial. The
last equation is equivalent to
$$
\hat{\theta}_p^{-1}(\Fr_*(H)) = x_i^p.
$$
By the special case above, there exists a unique Weyl polynomial $\hat{G}$ such that
$$
\hat{G}^p = \hat{\theta}_p^{-1}(\Fr_*(H)).
$$
But then
$$
H = \Fr_*^{-1}(\hat{\theta}_p(\hat{G}^p))
$$
which is exactly what we wanted.
\end{proof}

It will be convenient to denote the one-to-one correspondence between modulo $p$ reductions of central
polynomials coming from characteristic zero symplectomorphisms with their Weyl liftings by $\Phi^h_p$ (for this
correspondence is, as evidenced by Lemma \ref{twistedliftcanon}, shares essential nature with the characteristic
zero direct homomorphism $\Phi^h$).

We now apply the above lemma in order to establish the form of the pre-image Weyl polynomial in the case of
auxiliary variables $u,v$ and the central polynomial of a special type.
\begin{lem}\label{formlemma}
Let $u,v$ denote the extra Poisson variables, and let
$$
H = u + h^k\varphi(x_i)
$$
be the image of $u$ under the twisted automorphism coming from $\varphi$ as above ($\varphi(x_i)$ is rational in
$h$ but $h^k\varphi(x_i)$ is polynomial in $h$). Then the unique pre-image $\hat{H}$ of $H$ with respect to the
correspondence $\Phi^h_p$ of the previous lemma has the form
$$
\hat{H} = \hat{u} + h^k P_i(x_1,\ldots, d_n,h)
$$
where $P_i$ is rational in $h$.
\end{lem}
\begin{proof}
We establish the statement in several elementary steps. Firstly, as $H$ does not contain the auxiliary variable $v$,
$\hat{H}$ does not contain its Weyl counterpart $\hat{v}$: indeed, otherwise the Newton polyhedron of $H$
would contain (in the case of $v$ carrying great enough weight to make the corresponding monomial the
highest-order term) a vertex corresponding to the monomial containing $\hat{v}$. \footnote{Note that $\Phi^h_p$
behaves toward the Newton polyhedra as the homomorphism taking the $p$-th power does.}

Now let
$$
\hat{H} = Q(\hat{u}) + R
$$
where every monomial in $R$ is proportional to generators other than $u$. Then
$$
H = \Phi^h_p(\hat{H}) = \Phi^h_p(Q) + \Phi^h_p(R),
$$
as the two differential operators $Q$ and $R$ commute with each other and therefore taking the $p$-th power is
executed as in the commutative case. By Lemma \ref{twistedliftcanon}, we must have
$$
Q(\hat{u}) = \hat{u}.
$$

Finally, we show that if $\hat{H}$ contains monomials which are products of $\hat{u}$ with other generators, then
$\Phi^h_p(\hat{H})\neq H$.  Indeed, if such a monomial had a non-zero coefficient in  $\hat{H}$, then there would
exist a grading under which this monomial would be the highest-order term (corresponding to a vertex in the Newton
polyhedron). Then the image $\Phi^h_p(\hat{H})$ would also have a monomial corresponding to this highest-order
term with non-zero coefficient, as taking the $p$-th power dilates the polyhedron and therefore maps the extremal
points to extremal points.

The conclusion is that the polynomial $\hat{H}$ has the form
$$
\hat{u} + \tilde{P}_i(x_1,\ldots, d_n,h).
$$
Taking out $h^k$ from $\tilde{P}_i$ leaves us with the form we needed.
\end{proof}

Lemma \ref{formlemma} provides a canonical way to relate the images $\varphi(x_i)$ and $\varphi(d_j)$ of the
initial symplectomorphism with the Weyl pre-images. Therefore, as an array of differential operators, the lifting
$\hat{\varphi}$ is well defined. We denote the polynomials in Weyl generators in the correspondence by
$$
(\hat{\varphi}(x_1),\ldots, \hat{\varphi}(d_n)).
$$

We now need to verify the commutation relations in order to establish its homomorphic character. Again we have
two lemmas.
\begin{lem}\label{commutationlemma1}
\begin{gather*}
[\hat{\varphi}(x_i), \hat{\varphi}(x_j)] = [\hat{\varphi}(d_i), \hat{\varphi}(d_j)]=0,\\
[\hat{\varphi}(x_i), \hat{\varphi}(d_j)]=0,\;\;i\neq j.
\end{gather*}
\end{lem}
\begin{proof}
It suffices to prove
$$
[\hat{\varphi}(x_1), \hat{\varphi}(x_2)]=0
$$
thanks to the variable re-labelling and the existence of the "Fourier transform" -- the automorphism
$$
x_i\mapsto d_i, \;\;d_i\mapsto -x_i.
$$
We introduce two pairs of auxiliary Poisson variables,  $u_1, u_2, v_1, v_2$, and for $\lambda = h^k$ and $k$
large enough consider the automorphism $\psi$:
\begin{gather*}
u_1\mapsto u_1 + \lambda x_1,\;\;u_2\mapsto u_2+\lambda x_2\\
p_1\mapsto p_1-\lambda v_1,\;\;p_2\mapsto p_2-\lambda v_2
\end{gather*}
($\psi$ acts an the identity map on the rest of the generators).

We take the twisted automorphism
$$
\varphi_{t,\lambda} = \varphi_a\circ\psi\circ\varphi_a^{-1}
$$
with $\psi$ now being the chosen linear transformation and take $k$ to be large enough so that the twisted
automorphism is a polynomial symplectomorphism. We then lift it with $\Theta^h$ to the $h$-augmented Weyl
algebra as before.

By Lemma \ref{formlemma}, the images of the Weyl counterparts $\hat{u}_i$ ($i=1,2$) of $u_1, u_2$ under the
lifted twisted automorphism will have the form
$$
\hat{u}_i + \lambda T_i,
$$
the polynomials $T_i$ do not contain the auxiliary variables and $\Phi^h_p(T_i) = \varphi(x_i),\;\;i=1,2$.

Now, as $\varphi_{t,\lambda}$ and its lifting are automorphisms, we must have
$$
[\hat{u}_1 + \lambda T_1,\hat{u}_2 + \lambda T_2]=0
$$
so that
$$
[\hat{u}_1, \hat{u}_2] + \lambda([\hat{u}_1, T_2] + [T_1, \hat{u}_2]) + \lambda^2 [T_1,T_2] = 0
$$
from which it follows immediately that
$$
[T_1,T_2] = 0
$$
as desired.
\end{proof}

\begin{lem}\label{commutationlemma2}
$$[\hat{\varphi}(d_i), \hat{\varphi}(x_i)] = h.$$
\end{lem}
\begin{proof}
We proceed in an manner analogous to the previous lemma: we construct the appropriate twisting from whose lifting
the relevant images may be read off and then evaluate the commutator.

Let $u,v$ be auxiliary Poisson variables and let
$$
\psi_1: u\mapsto u+\lambda x_i,\;\;p_i\mapsto p_i - \lambda v,
$$
$$
\psi_2: v\mapsto v+\mu p_i,\;\; x_i\mapsto x_i - \mu u
$$
(in both cases the other generators are mapped to themselves). Consider the composition
$$
\theta = \psi_1\circ\psi_2.
$$
Then
$$
\theta(u) = u+\lambda x_i,\;\; \theta(v) = v+\mu p_i - \lambda\mu v
$$
and
$$
\theta(x_i) = x_i - \mu u - \lambda\mu x_i,\;\;\theta(p_i) = p_i - \lambda v.
$$
Take
$$
\varphi_{t,\lambda\mu} = \varphi_a\circ\theta\circ\varphi_a^{-1}
$$
where as before $\varphi_a$ extends from $\varphi$ by the identical action on $u,v$. The images of $u,v$ under
$\varphi_{t,\lambda\mu}$ read:
$$
\varphi_{t,\lambda\mu}(u) = u + \lambda \varphi(x_i),\;\;\varphi_{t,\lambda\mu}(v) = (1-\lambda\mu)v +\mu\varphi(p_i).
$$
By properly selecting $\lambda$ and $\mu$ as polynomials in $h$, we can make $\varphi_{t,\lambda\mu}$ into a
polynomial $h$-augmented symplectomorphism and therefore lift it with $\Theta^h$. Again, by Lemma
\ref{formlemma}, the action of the lifted automorphism on $\hat{u}$ and $\hat{v}$ will have the needed form (with
the part dependent on $x_i$, $d_j$ given by the images under $\hat{\varphi}$). Now, the commutator of the
images of $\hat{u}$ and $\hat{v}$ must be equal to $h$. We therefore have the following:
$$
h = [(1-\lambda\mu)\hat{v} + \mu\hat{\varphi}(d_i), \hat{u}+\lambda\hat{\varphi}(x_i)]=h(1-\lambda\mu) + \lambda\mu[\hat{\varphi}(d_i),\hat{\varphi}(x_i)]
$$
from which the statement follows directly.
\end{proof}

The conclusion is that augmented symplectomorphisms rational in $h$ are lifted to endomorphisms of the augmented
Weyl algebra (also rational in $h$) by a homomorphism whose restriction to points polynomial in $h$ coincides with
$\Theta^h$. What remains to show is that the lifted mappings are automorphisms, however, this is accomplished by
an argument similar to that for points polynomial in $h$ (cf. discussion immediately preceding Proposition
\ref{thetaproperties}). Also, thanks to Lemma \ref{formlemma}, we known that the lifting of points rational in $h$ is
also the inverse mapping to the extension to these points of the direct homomorphism $\Phi^h$. The specialization to
$h=1$ may now be safely executed, and the Main Theorem follows.


\chapter{Torus actions on free associative algebras,
lifting and Bialynicki-Birula type theorems} \label{Chapter4} \lhead{Chapter 5. \emph{Torus actions on free
associative algebras}}

We first prove that every maximal torus action on the free algebra is conjugate to a linear action. This statement is
the free algebra analogue of a classical theorem of A. Bia\l{}ynicki-Birula. This chapter is based on two papers
\cite{TA1, TA2}.

\section{Actions of algebraic tori}

In this section we recall basic definitions of the theory of torus actions, as formulated by Bia\l{}ynicki-Birula
\cite{BialBir1, BB3} and others.


\smallskip

Let $\mathbb{K}$ be the ground field. Let $I$ be a finite or a countable index set and let $$Z=\{z_i: i\in I\}$$ be
the set of variables, which is sometimes referred to as the alphabet.

\begin{def} \label{unnecessarydef1}
The free associative algebr $F_I(\mathbb{K})=\mathbb{K}\left\langle  Z \right\rangle $ is the algebra generated by
words in the alphabet $Z$ (as usually, word concatenation gives the multiplication of monomials and extends linearly
to define the multiplication in the algebra).
\end{def}

Any element of $\mathbb{K}\left\langle Z \right\rangle $ can be written uniquely in the form
$$\sum\limits_{k=0}^{\infty} \sum\limits_{i_1,\ldots,i_k \in I}^{}a_{i_1,i_2,\ldots,i_k}z_{i_1}z_{i_2}\ldots z_{i_k},$$
where the coefficients $a_{i_1,i_2,\ldots,i_k}$ are elements of the field $\mathbb{K}$ and all but finitely many of
these elements are zero.

\smallskip

In our context, the alphabet $Z$ is the same as the set of algebra generators, therefore the terms "monomial" and
"word" will be used interchangeably.

\smallskip

In the sequel, we employ the following short-hand notation for a free algebra monomial. For an element $z$, its
powers are defined intuitively. Any monomial $z_{i_1}z_{i_2}\ldots z_{i_k}$ can then be written in a reduced form
with subwords $zz\ldots z$ replaced by powers.

We then write
$$
z^I = z_{j_1}^{i_1}z_{j_2}^{i_2}\ldots z_{j_k}^{i_k}
$$
where by $I$ we mean an assignment of $i_k$ to $j_k$ in the word $z^I$. Sometimes we refer to $I$ as a
multi-index, although the term is not entirely accurate. If $I$ is such a multi-index, its abosulte value $|I|$ is defined
as the sum $i_1+\cdots+ i_k$.

\smallskip

For a field $\mathbb{K}$, let $\mathbb{K}^{\times}=\mathbb{K}\backslash \{0\}$ denote the multiplicative
group of its non-zero elements viewed.

\begin{Def} \label{defgroup}
An $n$-dimensional algebraic $\mathbb{K}$-torus is a group
$$
\mathbb{T}_n\simeq (\mathbb{K}^{\times})^n
$$
(with obvious multiplication).
\end{Def}
Denote by $\mathbb{A}^n$ the affine space of dimension $n$ over $\mathbb{K}$.
\begin{Def} \label{defaction}
A (left, geometric) torus action is a morphism
$$
\sigma: \mathbb{T}_n\times \mathbb{A}^n\rightarrow \mathbb{A}^n.
$$
that fulfills the usual axioms (identity and compatibility):
$$
\sigma(1,x)=x,\;\;\sigma(t_1,\sigma(t_2,x))=\sigma(t_1t_2,x).
$$

The action $\sigma$ is \textbf{effective} if for every $t\neq 1$ there is an element $x\in \mathbb{A}^n$ such that
$\sigma(t,x)\neq x$.
\end{Def}

In \cite{BialBir1}, Bia\l{}ynicki-Birula proved the following two theorems, for $\mathbb{K}$ algebraically closed.

\begin{thm} \label{BBthm1}
Any regular action of $\mathbb{T}_n$ on $\mathbb{A}^n$ has a fixed point.
\end{thm}

\begin{thm} \label{BBthm2}
Any effective and regular action of $\mathbb{T}_n$ on $\mathbb{A}^n$ is a representation in some coordinate
system.
\end{thm}

The term "regular" is to be understood here as in the algebro-geometric context of regular function
(Bia\l{}ynicki-Birula also considered birational actions).

\smallskip

In the following section (dedicated to the proof of the free algebra version of Theorems \ref{BBthm1} and
\ref{BBthm2}), the ground field is algebraically closed.

\medskip

As was mentioned in the introduction, an algebraic group action on $\mathbb{A}^n$ is the same as the
corresponding action by automorphisms on the algebra
$$
\mathbb{K}[x_1,\ldots,x_n]
$$
of coordinate functions. In other words, it is a group homomorphism
$$
\sigma: \mathbb{T}_n\rightarrow \Aut \mathbb{K}[x_1,\ldots,x_n].
$$
An action is effective if and only if $\Ker\sigma = \lbrace 1\rbrace$.

The polynomial algebra is a quotient of the free associative algebra
$$
F_n = \mathbb{K}\langle z_1,\ldots,z_n\rangle
$$
by the commutator ideal $I$ (it is the two-sided ideal generated by all elements of the form $fg-gf$). The definition
of torus action on the free algebra is thus purely algebraic.

In this chapter we establish the free algebra version of the Bia\l{}ynicki-Birula theorem. The latter is formulated as
follows.

\begin{thm} \label{BBfree}
Suppose given an action $\sigma$ of the algebraic $n$-torus $\mathbb{T}_n$ on the free algebra $F_n$. If
$\sigma$ is effective, then it is linearizable.
\end{thm}

\smallskip

The linearity (or linearization) problem, as it has become known since Kambayashi, asks whether all (effective,
regular) actions of a given type of algebraic groups on the affine space of given dimension are conjugate to
representations. According to Theorem \ref{BBfree}, the linearization problem extends to the noncommutative
category. Several known results concerning the (commutative) linearization problem are summarized below.

\begin{enumerate}
	\item Any effective regular torus action on $\mathbb{A}^2$ is linearizable (Gutwirth \cite{Gutwi}).
	\item Any effective regular torus action on $\mathbb{A}^n$ has a fixed point (Bia\l{}ynicki-Birula \cite{BialBir1}).
	\item Any effective regular action of $\mathbb{T}_{n-1}$ on $\mathbb{A}^n$ is linearizable (Bia\l{}ynicki-Birula \cite{BialBir2}).
    \item Any (effective, regular) one-dimensional torus action (i.e., action of $\mathbb{K}^{\times})$ on
        $\mathbb{A}^3$ is linearizable (Koras and Russell \cite{KoRu2}).
	\item If the ground field is not algebraically closed, then a torus action on $\mathbb{A}^n$ need not be linearizable. In \cite{Asanuma}, Asanuma proved that over any field $\mathbb{K}$, if there exists a non-rectifiable closed embedding from $\mathbb{A}^{m}$ into $\mathbb A^{n}$, then there exist non-linearizable effective actions of $(\mathbb{K}^{\times})^r$ on $\mathbb A^{1+n+m}$ for $1\le r\le 1+m$.
	\item When $\mathbb {K}$ is infinite and has positive characteristic, there are examples of non-linearizable torus actions on $\mathbb{A}^{n}$ (Asanuma \cite{Asanuma}).
\end{enumerate}

\begin{remark} \label{nonessential1}
A closed embedding $\iota:\mathbb A^m\to\mathbb A^n$ is said to be rectifiable if it is conjugate to a linear
embedding by an automorphism of $\mathbb A^n$.
\end{remark}

As can be inferred from the review above, the context of the linearization problem is rather broad, even in the case
of torus actions. The regulating parameters are the dimensions of the torus and the affine space. This situation is due
to the fact that the general form of the linearization conjecture (i.e., the conjecture that states that any effective
regular torus action on any affine space is linearizable) has a negative answer.

\smallskip

Transition to the noncommutative geometry presents the inquirer with an even broader context: one now may vary
the dimensions as well as impose restrictions on the action in the form of preservation of the PI-identities. Caution is
well advised. Some of the results are generalized in a straightforward manner -- the proof in the next section being
the typical example, others require more subtlety and effort. Of some note to us, given our ongoing work in
deformation quantization (see, for instance, \cite{KGE}) is the following instance of the linearization problem, which
we formulate as a conjecture.

\begin{conj} \label{BBsympln}
For $n\geq 1$, let $P_n$ denote the commutative Poisson algebra, i.e. the polynomial algebra
$$
\mathbb{K}[z_1,\ldots,z_{2n}]
$$
equipped with the Poisson bracket defined by
$$
\lbrace z_i, z_j\rbrace = \delta_{i,n+j}-\delta_{i+n,j}.
$$
Then any effective regular action of $\mathbb{T}_n$ by automorphisms of $P_n$ is linearizable.
\end{conj}

A version of Theorem \ref{commautautthm} for the commutative Poisson algebra is a conjecture of significant
interest. It turns out that the algebra $P_n$ admits a certain augmentation by central variables which distort the
Poisson structure, such that the automorphism group of the resulting algebra admits the property of Theorem
\ref{commautautthm}. The case is studied in the paper \cite{K-BE4}.

\section{Maximal torus action on the free algebra}

In this section, we provide proof to the free algebra version (Theorem \ref{BBfree}) of the Bia\l{}ynicki-Birula
theorem \cite{BialBir1}.

The proof proceeds along the lines of the original commutative case proof of Bia\l{}ynicki-Birula.

If $\sigma$ is the effective action of Theorem \ref{BBfree}, then for each $t\in \mathbb{T}_n$ the automorphism
$$
\sigma(t): F_n\rightarrow F_n
$$
is given by the $n$-tuple of images of the generators $z_1,\ldots,z_n$ of the free algebra:
$$
(f_1(t,z_1,\ldots,z_n),\ldots,f_n(t,z_1,\ldots,z_n)).
$$
Each of the $f_1,\ldots, f_n$ is a polynomial in the free variables.

\begin{lem} \label{fixedorigin}
There is a translation of the free generators
$$
(z_1,\ldots,z_n)\rightarrow (z_1-c_1,\ldots,z_n-c_n),\;\;(c_i\in\mathbb{K})
$$
such that (for all $t\in\mathbb{T}_n$) the polynomials $f_i(t,z_1-c_1,\ldots,z_n-c_n)$ have zero free part.
\end{lem}
\begin{proof}
This is a direct corollary of Theorem \ref{BBthm1}. Indeed, any action $\sigma$ on the free algebra induces, by
taking the canonical projection with respect to the commutator ideal $I$, an action $\bar{\sigma}$ on the
commutative algebra $\mathbb{K}[x_1,\ldots,x_n]$. If $\sigma$ is regular, then so is $\bar{\sigma}$. By Theorem
\ref{BBthm1}, $\bar{\sigma}$ (or rather, its geometric counterpart) has a fixed point, therefore the images of
commutative generators $x_i$ under $\bar{\sigma}(t)$ (for every $t$) will be polynomials with trivial degree-zero
part. Consequently, the same will hold for $\sigma$.
\end{proof}

We may then suppose, without loss of generality, that the polynomials $f_i$ have the form
$$
f_i(t,z_1,\ldots,z_n)=\sum_{j=1}^{n}a_{ij}(t)z_j + \sum_{j,l=1}^{n}a_{ijl}(t)z_jz_l + \sum_{k=3}^{N}\sum_{J,|J|=k}a_{i,J}(t)z^J
$$
where by $z^J$ we denote, as in the introduction, a particular monomial
$$
z_{i_1}^{k_1}z_{i_2}^{k_2}\ldots
$$
(a word in the alphabet $\lbrace z_1,\ldots, z_n\rbrace$ in the reduced notation; $J$ is the multi-index in the sense
described above); also, $N$ is the degree of the automorphism (which is finite) and $a_{ij}, a_{ijl},\ldots$ are
polynomials in $t_1,\ldots, t_n$.

As $\sigma_t$ is an automorphism, the matrix $[a_{ij}]$ that determines the linear part is non-singular. Therefore,
without loss of generality we may assume it to be diagonal (just as in the commutative case \cite{BialBir1}) of the
form
$$
\diag(t_1^{m_{11}}\ldots t_n^{m_{1n}},\ldots, t_1^{m_{n1}}\ldots t_n^{m_{nn}}).
$$

Now, just as in \cite{BialBir1}, we have the following
\begin{lem} \label{lem1-Bia}
The power matrix $[m_{ij}]$ is non-singular.
\end{lem}
\begin{proof}
Consider a linear action $\tau$ defined by
$$
\tau(t):(z_1,\ldots, z_n)\mapsto (t_1^{m_{11}}\ldots t_n^{m_{1n}}z_1,\ldots, t_1^{m_{n1}}\ldots t_n^{m_{nn}}z_n),\;\; (t_1,\ldots,t_n)\in\mathbb{T}_n.
$$
If $T_1\subset T_n$ is any one-dimensional torus, the restriction of $\tau$ to $\mathbb{T}_1$ is non-trivial.
Indeed, were it to happen that for some $\mathbb{T}_1$,
$$
\tau(t)z=z,\;\;t\in \mathbb{T}_1,\;\;(z=(z_1,\ldots,z_n))
$$
then our initial action $\sigma$, whose linear part is represented by $\tau$, would be identity modulo terms of degree
$>1$:
$$
\sigma(t)(z_i) = z_i + \sum_{j,l}a_{ijl}(t)z_jz_l+\cdots.
$$
Now, equality $\sigma(t^2)(z)=\sigma(t)(\sigma(t)(z))$ implies
\begin{align*}
\sigma(t)(\sigma(t)(z_i))&=\sigma(t)\left(z_i+\sum_{jl}a_{ijl}(t)z_jz_l+\cdots\right) \\&=
z_i+ \sum_{jl}a_{ijl}(t)z_jz_l+\sum_{jl}a_{ijl}(t)(z_j+\sum_{km}a_{jkm}(t)z_kz_m+\cdots)\\&(z_l+\sum_{k'm'}a_{lk'm'}(t)z_{k'}z_{m'}+\cdots)+\cdots\\&=
z_i+\sum_{jl}a_{ijl}(t^2)z_jz_l+\cdots
\end{align*}
which means that
$$
2a_{ijl}(t)=a_{ijl}(t^2)
$$
and therefore $a_{ijl}(t)=0$. The coefficients of the higher-degree terms are processed by induction (on the total
degree of the monomial). Thus
$$
\sigma(t)(z) = z,\;\;t\in\mathbb{T}_1
$$
which is a contradiction since $\sigma$ is effective. Finally, if $[m_{ij}]$ were singular, then one would easily find a
one-dimensional torus such that the restriction of $\tau$ were trivial.
\end{proof}

Consider the action
$$
\varphi(t) = \tau(t^{-1})\circ\sigma(t).
$$
The images under $\varphi(t)$ are
$$
(g_1(z,t),\ldots, g_n(z,t)),\;\;(t = (t_1,\ldots,t_n))
$$
with
$$
g_i(z,t) = \sum g_{i,m_1\ldots m_n}(z)t_1^{m_1}\ldots t_n^{m_n},\;\;m_1,\ldots, m_n\in\mathbb{Z}.
$$
Define $G_i(z) = g_{i,0\ldots 0}(z)$ and consider the map $\beta:F_n\rightarrow F_n$,
$$
\beta:(z_1,\ldots,z_n)\mapsto (G_1(z),\ldots, G_n(z)).
$$
\begin{lem} \label{lem2}
$\beta\in \Aut F_n$ and
$$
\beta = \tau(t^{-1})\circ\beta\circ\sigma(t).
$$
\end{lem}
\begin{proof}
This lemma mirrors the final part in the proof in \cite{BialBir1}. The conjugation is straightforward, since for every
$s,t\in\mathbb{T}_n$ one has
$$
\varphi(st) = \tau(t^{-1}s^{-1})\circ\sigma(st) = \tau(t^{-1})\circ\tau(s^{-1})\circ\sigma(s)\circ\sigma(t)=\tau(t^{-1})\circ\varphi(s)\circ\sigma(t).
$$

Denote by $\hat{F}_n$ the power series completion of the free algebra $F_n$, and let $\hat{\sigma}$,
$\hat{\tau}$ and $\hat{\beta}$ denote the endomorphisms of the power series algebra induced by corresponding
morphisms of $F_n$. The endomorphisms $\hat{\sigma}$, $\hat{\tau}$, $\hat{\beta}$ come from (polynomial)
automorphisms and therefore are invertible.

Let
$$
\hat{\beta}^{-1}(z_i) \equiv B_i(z) = \sum_{J}b_{i,J}z^J
$$
(just as before, $z^J$ is the monomial with multi-index $J$). Then
$$
\hat{\beta}\circ\hat{\tau}(t)\circ\hat{\beta}^{-1}(z_i) = B_i(t_1^{m_{11}}\ldots t_n^{m_{1n}}G_1(z),\ldots,t_1^{m_{n1}}\ldots t_n^{m_{nn}}G_n(z)).
$$
Now, from the conjugation property we must have
$$
\hat{\beta}=\hat{\sigma}(t^{-1})\circ\hat{\beta}\circ\hat{\tau}(t),
$$
therefore $\hat{\sigma}(t) = \hat{\beta}\circ\hat{\tau}(t)\circ\hat{\beta}^{-1}$ and
$$
\hat{\sigma}(t)(z_i) = \sum_{J}b_{i,J}(t_1^{m_{11}}\ldots t_n^{m_{1n}})^{j_1}\ldots (t_1^{m_{n1}}\ldots t_n^{m_{nn}})^{j_n}G(z)^J;
$$
here the notation $G(z)^J$ stands for a word in $G_i(z)$ with multi-index $J$, while the exponents $j_1,\ldots, j_n$
count how many times a given index appears in $J$ (or, equivalently, how many times a given generator $z_i$
appears in the word $z^J$).

Therefore, the coefficient of $\hat{\sigma}(t)(z_i)$ at $z^J$ has the form
$$
b_{i,J}(t_1^{m_{11}}\ldots t_n^{m_{1n}})^{j_1}\ldots (t_1^{m_{n1}}\ldots t_n^{m_{nn}})^{j_n}+S
$$
with $S$ a finite sum of monomials of the form
$$
c_L (t_1^{m_{11}}\ldots t_n^{m_{1n}})^{l_1}\ldots (t_1^{m_{n1}}\ldots t_n^{m_{nn}})^{l_n}
$$
with $(j_1,\ldots,j_n)\neq (l_1,\ldots,l_n)$. Since the power matrix $[m_{ij}]$ is non-singular, if $b_i,J\neq 0$, we
can find a $t\in\mathbb{T}_n$ such that the coefficient is not zero. Since $\sigma$ is an algebraic action, the degree
$$
\sup_{t}\Deg (\hat{\sigma})
$$
is a finite integer $N$. With the previous statement, this implies that
$$
b_{i,J} = 0,\;\;\text{whenever}\;\;|J|>N.
$$
Therefore, $B_i(z)$ are polynomials in the free variables. What remains is to notice that
$$
z_i = B_i(G_1(z),\ldots,G_n(z)).
$$
Thus $\beta$ is an automorphism.
\end{proof}
From Lemma \ref{lem2} it follows that
$$
\tau(t) = \beta^{-1}\circ \sigma(t)\circ\beta
$$
which is the linearization of $\sigma$. Theorem \ref{BBfree} is proved.

\section{Discussion}
The noncommutative toric action linearity property has several useful applications. In the work \cite{KBYu}, it is
used to investigate the properties of the group $\Aut F_n$ of automorphisms of the free algebra. As a corollary of
Theorem \ref{BBfree}, one gets
\begin{cor} \label{cor1}
Let $\theta$ denote the standard action of $\mathbb{T}_n$ on $K[x_1,\ldots,x_n]$ -- i.e., the action
$$
\theta_t: (x_1,\ldots,x_n)\mapsto (t_1x_1,\ldots,t_nx_n).
$$
 Let $\tilde{\theta}$ denote its lifting to an action on the free associative algebra $F_n$. Then
$\tilde{\theta}$ is also given by the standard torus action.
\end{cor}

This statement plays a part, along with a number of results concerning the induced formal power series topology on
$\Aut F_n$, in the establishment of the free associative analogue of Theorem \ref{commautautthm}.

\smallskip

The proofs in this paper, for the most part, were based upon the techniques from the commutative category. It is,
however, a problem of legitimate interest to try and obtain proofs for various linearity statements using tools specific
to the category of associative algebras, bypassing the known commutative results. As one outstanding example of
this problem, we expect the free associative analogue of the second Bia\l{}ynicki-Birula theorem to hold and
formulate it here as a conjecture.

\begin{conj} \label{conjlin}
Any effective action of $\mathbb{T}_{n-1}$ on $F_n$ is linearizable.
\end{conj}

Also of independent interest is the following instance of the linearity problem.

\begin{conj} \label{BBsymp2}
For $n\geq 1$, let $P_n$ denote the commutative Poisson algebra, i.e. the polynomial algebra
$$
\mathbb{K}[z_1,\ldots,z_{2n}]
$$
equipped with the Poisson bracket defined by
$$
\lbrace z_i, z_j\rbrace = \delta_{i,n+j}-\delta_{i+n,j}.
$$
Then any effective regular action of $\mathbb{T}_n$ by automorphisms of $P_n$ is linearizable.
\end{conj}
This problem is loosely analogous to the Bia\l{}ynicki-Birula theorem, in the sense of maximality of torus with
respect to the dimension of the configurations space (spanned by $x_i$). There seems to be no straightforward way
of finding the linearizing canonical coordinates on the phase space, however. For the $\Ind$-variety $\Aut P_n$, a
version of Theorem \ref{commautautthm} may be stated. The geometry of $\Aut P_n$ is relevant to problems of
deformation quantization.

\chapter{Jacobian conjecture, Specht and Burnside type problems}
\label{Chapter6} \lhead{Chapter 6. \emph{Jacobian conjecture, Specht and Burnside type problems}}

This chapter explores an approach to polynomial mappings and the Jacobian Conjecture and related questions,
initiated by A.V.~Yagzhev, whereby these questions are translated to identities of algebras, leading to a solution  in
\cite{Yag1} of the version  of the Jacobian Conjecture for free associative algebras. (The first version, for two
generators, was obtained by Dicks and J.~Levin \cite{Di, DiLev}, and the full version by Schofield~\cite{Schof}.)
We start by laying out the basic framework in this introduction. Next, we set up Yagzhev's correspondence to
algebras in \S\ref{operad}, leading to the basic notions of weak nilpotence and Engel type. In
\S\ref{SbScJCArbVar0} we discuss the Jacobian Conjecture in the context of various varieties, including the free
associative algebra.

Given any  polynomial endomorphism $\phi$ of the $n$-dimensional affine space $A^n_\k=\spec
\k[x_1,\dots,x_n]$ over a field $\k$, we define its {\it Jacobian matrix} to be the matrix
$$ \left(
\partial \phi^*(x_i)/ \partial x_j\right)_{1\le i,j\le n}.$$
The determinant of the Jacobian matrix is called  the {\it Jacobian} of $\phi$. The celebrated {\bf Jacobian
Conjecture} $\JC_n$ in dimension $n\ge 1$ asserts that {\it for any field~$\k$ of characteristic zero, any polynomial
endomorphism $\phi$ of $\A^n_\k $ having Jacobian $1$ is an
 automorphism.} Equivalently, one can say that $\phi$ preserves the
 standard top-degree differential form $dx_1\wedge\dots\wedge
dx_n\in\Omega^n(\A^n_\k)$. References to this well known problem and related questions can be found
in~\cite{BCW}, \cite{Kulik01}, and \cite{vdE}.
 By the Lefschetz principle it is sufficient to consider the case  $\k=\C$;
obviously, $\JC_n$ implies $\JC_m$ if $n>m$.
  The conjecture $\JC_n$ is obviously true  in the case
$n=1$, and it is open for $n\ge 2$.

The {\bf Jacobian Conjecture}, denoted as $\JC $, is the conjunction of the conjectures $\JC_n$ for
 all finite $n$.  The Jacobian
Conjecture has many reformulations (such as the Kernel Conjecture and the Image Conjecture,
cf.~\cite{vdE,VDEssenImage,VDEssenImageA,VDEssenWenhua,VDEssenWenhua1} for details) and is closely
related to questions concerning quantization. It is stably equivalent to the following conjecture of Dixmier, concerning
automorphisms of the Weyl algebra $W_n$, otherwise known as the {\it quantum affine algebra}.

\medskip
{\bf Dixmier Conjecture $DC_n$:}\ {\it Does $\End(W_n)=\Aut(W_n)$}?
\medskip

The implication $DC_n\to JC_n$ is well known,  and the inverse implication $JC_{2n}\to DC_n$ was recently
obtained independently by Tsuchimoto \cite{Tsu2} (using $p$-curvature) and Belov and Kontsevich \cite{K-BK2},
\cite{K-BK1} (using Poisson brackets on the center of the Weyl algebra). Bavula~\cite{Bavula1new} has obtained
a shorter proof, and also obtained a positive solution of an analog of the Dixmier Conjecture for integro differential
operators, cf.~\cite{Bavula1}. He also proved that every monomorphism of the  Lie algebra of triangular polynomial
derivations is an automorphism \cite{Bavula2} (an analog of Dixmier's conjecture).

The Jacobian Conjecture is closely related to many questions of affine algebraic geometry concerning affine space,
such as the Cancellation Conjecture  (see Section \ref{SbScRelQuestJC}). If we replace the variety of commutative
associative algebras (and the accompanying affine spaces) by an arbitrary algebraic variety
\label{PgAlgVar}\footnote{Algebraic geometers use  word {\it variety}, roughly speaking, for objects whose local
structure is obtained from the solution of system of algebraic equations. In the framework of universal algebra, this
notion is used for subcategories of algebras defined by a  given set of identities. A deep analog of these notions is
given in \cite{BelovIAN}.}, one easily gets a counterexample to the  $\JC$. So, strategically these questions deal
with some specific properties of affine space which we do not yet understand, and for which we do not have the
appropriate formulation apart from these very difficult questions.

 It seems that these properties do indicate some sort of quantization. From that
perspective, noncommutative analogs of these problems (in particular, the Jacobian Conjecture and the analog of the
Cancellation Conjecture)  become interesting for free associative algebras, and more generally, for arbitrary
varieties of algebras.

We work in the language of universal algebra, in which an algebra is defined in terms of a set of operators, called its
{\it signature}. This approach enhances the investigation of the Yagzhev correspondence between endomorphisms
and algebras. We work with deformations and so-called {\it packing properties} to be introduced in
Section~\ref{SbScJCArbVar0} and Section~\ref{SbSbScPacking}, which denote specific noncommutative
phenomena which enable one to solve   the  $\JC$ for the free associative algebra.

 From the viewpoint of universal
algebra,  the Jacobian conjecture becomes a problem of ``Burnside type,'' by which we mean the question of
whether a given finitely generated algebraic structure satisfying given periodicity conditions is necessarily finite,
cf.~Zelmanov~\cite{Zelmanov}.
 Burnside originally posed the question
of the finiteness of a finitely generated  group satisfying the identity $x^n=1$. (For odd $n\ge 661,$
counterexamples were found by Novikov and Adian, and quite recently  Adian reduced the estimate from $661$ to
$101$).  Another class of counterexamples was discovered by Ol'shanskij~\cite{Ol}. Kurosh posed the question of
local finiteness of algebras whose elements are algebraic over the base field. For algebraicity of bounded degree, the
question has a positive solution, but otherwise there are the Golod-Shafarevich counterexamples.

 Burnside type
problems play an important role in algebra. Their solution in the associative case is closely tied to Specht's problem
of whether any set of polynomial identities can be deduced from a finite subset. The $\JC$ can be  formulated in the
context of whether one system of identities implies another, which also relates to Specht's problem.

In the Lie algebra case there is a similar notion. An element $x\in L$ is called {\it Engel of degree $n$} if
$[\dots[[y,x],x]\dots, x]=0$ for any $y$ in the Lie algebra~$L$. Zelmanov's result that any finitely generated Lie
algebra of bounded Engel degree is nilpotent yielded his solution of the Restricted Burnside Problem for groups.
Yagzhev introduced the notion of {\it Engelian} and {\it weakly nilpotent} algebras of arbitrary signature (see
Definitions \ref{Engtyp}, \ref{DfWeaklyNilp}), and proved that the $\JC$ is equivalent to the question of weak
nilpotence of algebras of Engel type satisfying a system of Capelli identities, thereby showing the relation of the
$\JC$ with problems of Burnside type.

\medskip

 {\bf A negative approach.} Let us mention
a way of constructing counterexamples. This approach, developed by Gizatullin, Kulikov,  Shafarevich, Vitushkin,
and others,
is related to
decomposing polynomial mappings into the composition of $\sigma$-processes
\cite{Gizatullin,Kulik01,Shafarevich,Vit1,Vitushkin1,Vit2}. It allows one to solve some polynomial automorphism
problems, including tameness problems, the most famous of which is  {\em Nagata's Problem} concerning the
wildness of  Nagata's automorphism
 $$  (x,y,z)\mapsto(x-2(xz+y^2)y-(xz+y^2)^2z,\, y+(xz+y^2)z,\, z),$$
cf.~\cite{Nag}. Its solution by Shestakov and Umirbaev~\cite{SU2} is the major advance in this area in the last
decade. The Nagata automorphism can be constructed as a product of automorphisms of $K(z)[x,y]$, some of
them having non-polynomial coefficients (in $K(z)$).
 The following theorem of
Abhyankar-Moh-Suzuki \cite{AM,Su} and \cite{M} can be viewed in this context:

\medskip
\textbf{AMS Theorem.}  If $f$ and $g$ are polynomials in $K[z]$ of degrees $n$ and $m$ for which $K[f, g] =
K[z]$, then $n$ divides $m$ or $m$ divides $n$.\medskip

Degree estimate theorems are polynomial analogs to Liouville's approximation theorem in algebraic number theory
(\cite{BonnetVenerau,Kuroda,YuYungChang,MLY}). T.~Kishimoto has proposed  using a program of Sarkisov, in
particular for Nagata's Problem. Although difficulties remain in applying $``\sigma$-processes'' (decomposition of
birational mappings into standard blow-up operations) to the affine case, these may provide new insight. If we
consider affine transformations of the plane, we have relatively simple singularities at infinity, although for bigger
dimensions they can be more complicated. Blow-ups provide some understanding of birational mappings with
singularities. Relevant information may be provided in the affine case.  The paper~\cite{BEW} contains some deep
 considerations about singularities.

\section{The Jacobian Conjecture and Burnside type problems, via algebras}\label{operad}

In this section we translate the  Jacobian Conjecture to the language of algebras and their identities. This can be done
at two levels:  At the level of the algebra obtained from a polynomial mapping, leading to the notion of {\it weak
nilpotence} and {\it Yagzhev algebras} and at the level of the differential and the algebra arising from the Jacobian,
leading to the  notion of {\it Engel type}. The Jacobian Conjecture   is  the link between these two notions.

\subsection{The Yagzhev correspondence}\label{Yagcor}

\subsubsection{Polynomial mappings in universal algebra}

Yagzhev's approach is to pass from algebraic geometry to universal algebra. Accordingly, we work in the
framework of a universal algebra $A$ having signature $\Omega$. $A ^{(m)}$ denotes $A \times \dots \times A $,
taken $m$ times.

We fix a commutative, associative base ring~$C$, and consider $C$-modules equipped with extra {\em operators}
$A ^{(m)}\to A,$ which we call $m$-{\em ary}. Often one of these operators will be (binary) multiplication. These
operators will be multilinear, i.e., linear with respect to each argument. Thus, we can define the {\em degree} of an
operator to be its number of arguments. We say an operator $\Psi(x_1, \dots, x_m)$ is {\em symmetric} if
$\Psi(x_1, \dots, x_m) = \Psi(x_{\pi(1)}, \dots, x_{\pi(m)})$ for all permutations~ $\pi$.

\begin{definition}\label{DefPolMap}
A {\em string} of operators is defined inductively. Any  operator $\Psi(x_1, \dots, x_m)$ is a string of degree $m$,
and if $s_j$ are strings of degree $d_j,$ then $\Psi(s_1, \dots, s_m)$ is a string of degree $\sum_{j=1}^m d_j$.  A
mapping
 $$\alpha: A ^{(m)} \to A$$
 is called {\em polynomial} if it can be
expressed as a sum of strings of operators of the algebra $A$. The {\em degree} of the mapping is the maximal
length of these strings. \end{definition}

{\bf Example.} Suppose an algebra $A$ has two extra operators: a binary operator $\alpha(x,y)$ and a tertiary
operator $\beta(x,y,z)$. The mapping $F: A\to A$ given by $\ x\to x+\alpha(x,x)+\beta(\alpha(x,x),x,x)$ is a
polynomial mapping of~$A$, having degree $4$. Note that if $A$ is finite dimensional as a vector space, not every
polynomial mapping of $A$ as an affine space is a polynomial mapping of $A$ as an algebra.
\medskip

\subsubsection{Yagzhev's
correspondence  between polynomial mappings and algebras}

Here we associate an algebraic structure to each polynomial map. Let $V$ be an $n$-dimensional vector space
over the field $\k$, and $F:V\to V$  be a polynomial mapping of degree $m$.
 Replacing $F$ by the composite $TF$,  where $T$ is a translation such that
$TF(0)=0$, we may assume that $F(0) =0$. Given a base $\{\vec{e}_i\}_{i=1}^n$ of $V$, and for an element
$v$ of $V$ written uniquely as a sum $\sum x_i\vec{e}_i$, for $x_i \in \k$, the coefficients of $\vec{e}_i$ in
$F(v)$ are (commutative) polynomials in the~ $x_i$. Then $F$ can be written in the following form:
$$
x_i\mapsto F_{0i}(\vec{x})+F_{1i}(\vec{x})+\dots+F_{mi}(\vec{x})
$$
where each $F_{\alpha i}(\vec{x})$ is a homogeneous form of degree $\alpha$, i.e.,
$$
F_{\alpha i}(\vec{x})=\sum_{j_1 + \dots + j_n =\alpha}\kappa_J
x_1^{j_1}\cdots x_n^{j_n},
$$
with $F_{0i}=0$ for all $i$,  and $F_{1i}(\vec x)=\sum_{k=1}^n \mu_{ki}x_k$.

 We are
interested in invertible mappings that have a nonsingular Jacobian matrix $(\mu_{ij})$. In particular, this matrix is
nondegenerate at the origin. In this case $\det(\mu_{ij})\ne 0$, and by composing $F$ with an affine transformation
we arrive at the situation for which $\mu_{ki}=\delta_{ki}$. Thus,  the mapping $F$ may be taken to have  the
following form:

\begin{equation}   \label{Eq1Fk}
x_i\to x_i-\sum_{k=2}^m F_{ki}.
\end{equation}

Suppose we have a mapping as in \eqref{Eq1Fk}. Then the Jacobi matrix can be written as $E - G_1 - \dots -
G_{m-1}$ where $G_i$ is an $n\times n$ matrix with entries which are homogeneous polynomials of degree $i$. If
the Jacobian is $1$, then it is invertible with inverse a polynomial matrix (of  homogeneous degree at most
$(n-1)(m-1)$, obtained via the adjoint matrix).

If we write the inverse as a formal power series, we   compare the homogeneous components and get:
 \begin{equation}  \label{EqMJ}
\sum_{j_im_{j_i}=s}M_J=0,
\end{equation}
where $M_J$ is the sum of products $a_{\alpha_1} a_{\alpha_q}$ in which the factor $a_j$ occurs $m_j$~times,
and $J$ denotes the multi-index $(j_1,\dots,j_q)$.

Yagzhev considered the cubic homogeneous mapping $\vec{x}\to \vec{x}+(\vec{x},\vec{x},\vec{x}),$ whereby
the Jacobian matrix becomes $E-G_3.$ We return to  this case in~Remark~\ref{Jac3}.
  The   slightly more general approach
 given here presents the
Yagzhev correspondence more clearly and also provides tools for investigating deformations and packing properties
(see Section \ref{SbSbScPacking}). Thus, we consider not only the cubic case (i.e. when the mapping has the form
$$x_i\to x_i+P_i(x_1,\dots,x_n);\ i=1,\dots, n,$$ with $P_i$
 cubic homogenous polynomials), but the more general situation of arbitrary degree.

For any $\ell$, the set of (vector valued) forms $\{F_{\ell,i}\}_{i=1}^n$ can be interpreted as a homogeneous
mapping $\Phi_\ell: V\to V$ of degree $\ell$.  When  $\ch(\k)$ does not divide $\ell$, we take instead the
polarization of this mapping, i.e.~the multilinear symmetric mapping
$$
\Psi_\ell:V^{\otimes\ell}\to V
$$
 such that
$$
(F_{\ell ,i}(x_1),\dots,F_{\ell,
i}(x_n))=\Psi_\ell(\vec{x},\dots,\vec{x}) \cdot \ell!
$$
Then Equation (\ref{Eq1Fk}) can be rewritten as

\begin{equation}   \label{Eq2Fk}
\vec{x}\to \vec{x}-\sum_{\ell=2}^m
\Psi_{\ell}(\vec{x},\dots,\vec{x}).
\end{equation}

We define the  algebra  $(A, \{ \Psi_{ \ell}\}),$ where $A$ is the vector space $V$ and the~$ \Psi_{ \ell}$ are
viewed as operators $A^{\ell}\to A.$

\begin{definition}
 The \textbf{Yagzhev correspondence} is the correspondence from
 the polynomial mapping $(V,F)$ to the   algebra  $(A, \{ \Psi_{ \ell}\}).$
\end{definition}

\subsection{Translation of the invertibility condition to the language of
identities}

The next step is to bring in algebraic varities, defined in terms of identities.

\begin{definition}
 A {\em polynomial identity} (PI) of $A$ is a
polynomial mapping of~$A,$ all of whose values are identically zero.

The {\it algebraic variety} generated by an algebra $A$, denoted as $\Var(A)$, is the class of all algebras satisfying
the same PIs as $A$.
\end{definition}

Now we come to a crucial idea of Yagzhev:

\medskip
{\it The invertibility of $F$ and the invertibility of the Jacobian of $F$ can be expressed via (2)  in the language of
polynomial identities.}
\medskip

%

Namely, let $y=F(x)=x-\sum_{\ell =2}^m\Psi_\ell (x)$. Then

\begin{equation}   \label{EqTerm}
F^{-1}(x)=\sum_t t(x),
\end{equation}

 \noindent where each $t$ is a {\it
term}, a formal expression in the mappings $\{\Psi_\ell \}_{\ell =2}^m$ and the symbol $x$. Note that the
expressions $\Psi_2(x,\Psi_3(x,x,x))$ and $\Psi_2(\Psi_3(x,x,x),x)$ are different although they represent same
element of the algebra. Denote by $|t|$  the number of occurrences of variables, including multiplicity, which are
included in $t$.

The invertibility of $F$ means that, for all $q\ge q_0$,
\begin{equation}    \label{EqTerm2}
\sum_{|t|=q}t(a)=0, \quad \forall a \in A.
\end{equation}

Thus we have translated  invertibility of the mapping $F$ to the language of identities. (Yagzhev had an analogous
formula, where the terms only involved $\Psi_3$.)
%

\begin{definition}     \label{DfWeaklyNilp}
  An element $a\in A$ is called {\em nilpotent} of   index $\le n$
if $$M(a, a, \dots, a)= 0$$  for each monomial $M(x_1, x_2, \dots )$ of degree $\ge n$.   $A$~is {\em weakly
nilpotent} if each element of $A$ is nilpotent.  $A$~is {\em weakly nilpotent of class~$k$} if each element of $A$
is nilpotent of index $k$. (Some authors use the terminology {\it index} instead of {\em class}.) Equation
\eqref{EqTerm2} means $A$ is weakly nilpotent.

To stress this fundamental notion of Yagzhev, we define a
 {\em Yagzhev algebra} of   {\em order} $q_0$ to be a weakly nilpotent algebra, i.e.,
 satisfying the
identities  \eqref{EqTerm2}, also called the {\it system of Yagzhev identities} arising from $F$.
\end{definition}

Summarizing, we get the following fundamental translation from conditions on the endomorphism $F$ to identities of
algebras.

\begin{theorem}
The endomorphism $F$ is invertible if and only if the corresponding algebra  is a Yagzhev algebra of high enough
order.
\end{theorem}

\subsubsection{Algebras of Engel type}

The analogous procedure can be carried out for the differential mapping. We recall that $\Psi_\ell$ is a symmetric
multilinear mapping of degree $\ell$.  We denote the mapping $y\to\Psi_\ell (y,x,\dots,x)$ as $\Ad_{\ell-1}(x)$.

\begin{definition}\label{Engtyp}
An algebra  $A$ is of {\em Engel type}  $s$ if it satisfies a system of identities
\begin{equation}  \label{EqMJ2}
\sum_{\ell m_\ell =s}  \quad \sum_ {\alpha _{1}+ \cdots + \alpha
_{q}=m_\ell }\Ad_{\alpha_1}(x)\cdots \Ad_{\alpha_q}(x)= 0.
\end{equation}
$A$  is of {\em Engel type}    if  $A$   has Engel type $s$  for some $s$.
\end{definition}

.

\begin{theorem}\label{corresconj}
The endomorphism $F$ has Jacobian $1$ if and only if the corresponding algebra has Engel type $s$  for some $s$.
\end{theorem}
\begin{proof}
Let $x'=x+dx$. Then
\begin{equation}\begin{aligned}
\Psi_\ell (x')= & \Psi_\ell (x)+\ell \Psi_\ell (dx,x,\dots,x)
\\ \qquad \qquad + & \quad \mbox{forms\ containing \ more\ than\ one occurence\ of }
dx.\end{aligned}
\end{equation}

Hence the differential of the mapping
$$
F:\vec{x}\mapsto \vec{x}-\sum_{\ell =2}^m \Psi_\ell (\vec{x}, \dots, \vec{x})
$$
is
$$
\left(E-\sum_{\ell =2}^m \ell \Ad_{\ell-1} (x)\right)\cdot dx
$$
The identities (\ref{EqMJ}) are equivalent to the   system of identities  (\ref{EqMJ2}) in the signature
$\Omega=(\Psi_2,\dots,\Psi_m)$,
 taking $a_{\alpha_j} = \Ad_{\alpha_j}$ and $m_j = \deg \Psi_\ell -1$.


\end{proof}

Thus, we have reformulated the condition of invertibility of the Jacobian in the language of identities.

As explained in \cite{vdE}, it is well known from \cite{BCW} and \cite{Yag4}   that the Jacobian Conjecture can
be reduced to the cubic homogeneous case; i.e., it is enough to consider mappings of type
$$
x\to x+\Psi_3(x,x,x).
$$
In this case the Jacobian assumption is equivalent to the {\it Engel condition} ~-- nilpotence of the mapping
$\Ad_3(x)[y]$ (i.e. the mapping $y\to (y,x,x)$). Invertibility, considered in \cite{BCW}, is equivalent to
  weak nilpotence, i.e., to the identity $\sum_{|t|=k}t=0$
   holding for all sufficiently large $k$.

\medskip
\begin{rem}\label{Jac3}  In the cubic homogeneous case, $j=1$, $\alpha _j=2$ and $m_j = s$,
and we define the linear map
$$\Ad_{xx}: y\to (x,x,y)$$ and the index set $T_j\subset
\{1,\dots,q\}$ such that $i\in T_j$ if and only if $\alpha_i=j$.

Then the equality (\ref{EqMJ2}) has  the following form:
$$
\Ad_{xx}^{s/2}=0.
$$
 Thus, for a ternary symmetric algebra, Engel type
means that the operators $\Ad_{xx}$ for all $x$ are nilpotent. In other words,  the mapping $$Ad_3(x): y\to
(x,x,y)$$ is nilpotent. Yagzhev called this the  {\em Engel condition}. (For Lie algebras the nilpotence of the
operator $\Ad_x: y\to (x,y)$ is the usual Engel condition. Here we have a generalization for arbitrary signature.)

Here  are Yagzhev's original  definitions, for edification. A binary algebra $A$ is
 {\em Engelian} if for any element $a\in A$ the subalgebra
$<\!R_a,L_a\!>$ of vector space endomorphisms of $A$  generated by
 the left multiplication operator $L_a$ and the right
 multiplication operator
~$R_a$ is nilpotent, and  {\em weakly Engelian} if for any element
$a\in A$ the operator $R_a+L_a$ is nilpotent. 
\end{rem}

This leads us to the {\it Generalized Jacobian Conjecture}:

\medskip
{\bf Conjecture.} {\it Let $A$ be an algebra with symmetric $\k$-linear operators $\Psi_\ell$, for $ \ell
=1,\dots,m$. In any variety of
 Engel type, $A$ is a Yagzhev algebra. }
\medskip

By Theorem~\ref{corresconj}, this conjecture would yield  the Jacobian Conjecture.

\subsubsection{The case of binary algebras }
When $A$ is a binary algebra, {\it Engel type} means that the left and right multiplication mappings are both
nilpotent.

 A well-known result of
S.~Wang~\cite{BCW} shows that the Jacobian Conjecture holds for quadratic mappings
$$
\vec{x}\to \vec{x}+\Psi_2(\vec{x},\vec{x}).
$$

If  two different points $(x_1,\dots,x_n)$ and $(y_1,\dots,y_n)$ of an affine space  are mapped to the same point
by $(f_1,\dots,f_n)$, then the fact that the vertex of a parabola is in the middle of the interval whose endpoints are at
the roots shows that all $f_i(\vec{x})$ have gradients at this midpoint $P=(\vec{x}+\vec{y})/2$ perpendicular to
the line segment $[\vec{x},\vec{y}]$. Hence   the  Jacobian is zero at the midpoint $P$. This fact holds in any
characteristic $\ne 2$.

In Section~\ref{LiftY}  we prove the following theorem of Yagzhev, cf.~Definition~\ref{Cap1} below:

\begin{theorem}[Yagzhev]\label{Yag2}
Every symmetric binary Engel type algebra of order $k$  satisfying the
  system of Capelli identities  of order $n$  is
weakly nilpotent, of weak nilpotence index  bounded by some function $F(k,n)$.
\end{theorem}

\begin{rem} Yagzhev 
formulates his theorem in the following way:\medskip

 {\it Every binary weakly Engel
algebra of order $k$ satisfying   the system of Capelli identities of order $n$   is weakly nilpotent, of   index
bounded by some function $F(k,n)$}.

\medskip
We obtain this  reformulation, by replacing the algebra $A$ by the algebra~$A^+$ with multiplication given by
$(a,b)=ab+ba$.
\end{rem}

The following problems may help us understand the situation:

\medskip
{\bf Problem.}\ {\it Obtain a straightforward proof of this theorem and deduce from it the Jacobian
 Conjecture for quadratic
mappings. }
\medskip

{\bf Problem. (Generalized Jacobian Conjecture for quadratic mappings)}\ {\it Is every symmetric binary algebra  of
Engel type
 $k$, a Yagzhev algebra?}
\medskip

\subsubsection{The case of ternary algebras}

As we have observed, Yagzhev reduced the Jacobian Conjecture over a field of characteristic zero to the question:

\medskip
{\it Is every  finite dimensional ternary Engel algebra  a Yagzhev algebra? }
\medskip

\noindent  Dru\'zkowski \cite{Druz2,Dru2} reduced this to the case when all cubic forms $\Psi_{3i}$ are cubes of
linear forms. Van den Essen and his school reduced the  $\JC$ to the symmetric case; see
\cite{VDEssenBondt,VDEssenBondt1} for details.
  Bass,  Connell, and  Wright \cite{BCW} use other methods including inversions. Yagzhev's approach matches that of
\cite{BCW}, but using identities instead.   

\subsubsection{An example  in nonzero characteristic of an Engel algebra that is not  a Yagzhev algebra}
\label{SbScEngNoWNilpCharpos}

Now we give an example,   over an arbitrary field $\k$ of characteristic $p>3$, of a finite  dimensional Engel
algebra that is not  a Yagzhev algebra, i.e., not weakly nilpotent. This means that the situation for binary algebras
   differs intrinsically from that for ternary algebras,
and it would be worthwhile to understand why.

\begin{theorem}      \label{ThEngWeakNoNilp}
If\ $\ch(\k)=p>3$, then there exists a finite dimensional $\k$-algebra that is Engel  but not weakly nilpotent.
\end{theorem}
\Proof
 Consider the noninvertible mapping $F: \k[x]\to \k[x]$ with
Jacobian~1:
$$
F: x\to x+x^p.
$$
We introduce  new commuting indeterminates $\{y_i\}_{i=1}^n$ and extend this mapping to $k[x,y_1,\dots,y_n]$
by sending $y_i\mapsto y_i$. If $n$ is big enough, then it is possible to find  tame automorphisms $G_1$ and
$G_2$ such that $G_1\circ F\circ G_2$ is a cubic mapping $\vec{x}\to \vec{x}+\Psi_3(\vec{x})$, as follows:

Suppose we have a mapping
$$
F:x_i\to P(x)+M
$$
where $M=t_1t_2t_3t_4$ is a monomial of degree at least $4$. Introduce two new commuting indeterminates $z,
y$ and take $F(z)=z$, $F(y)=y$.

Define the mapping $G_1$ via $G_1(z)=z+t_1t_2$, $G_1(y)=y+t_3t_4$ with $G_1$ fixing all other
indeterminates; define $G_2$ via $G_2(x)=x-yz$ with $G_2$ fixing all other indeterminates.

The composite mapping $G_1\circ F\circ G_2$ sends $x$ to $P(x)-yz-yt_1t_2-zt_3t_4$, $y$ to $y+t_3t_4$, $z$
to $z+t_1t_2$, and agrees with $F$ on all other indeterminates.

Note that we have removed the  monomial $M=t_1t_2t_3t_4$ from the image of $F$, but instead have obtained
various monomials of smaller degree ($t_1t_2$ , $t_3t_4$, $zy$, $zt_3t_4$, $yt_1t_2$). It is easy to see that this
process terminates.

Our new mapping $H(x)=x+\Psi_2(x)+\Psi_3(x)$ is noninvertible and has Jacobian~1. Consider its blowup
$$R: x\mapsto x+T^2y+T\Psi_2(x),\  y\mapsto y-\Psi_3(x),\  T\mapsto T.$$


This mapping $R$ is invertible if and only if the initial mapping is invertible, and has   Jacobian~1 if and only if the
initial mapping has Jacobian~1, by~\cite[Lemma 2]{Yag4}. This mapping is also cubic homogeneous. The
corresponding ternary algebra is Engel, but not weakly nilpotent. \Endproof

This example shows that a direct combinatorial approach to the Jacobian Conjecture encounters difficulties, and in
working with related Burnside type problems (in the sense of Zelmanov \cite{Zelmanov}, dealing with  nilpotence
properties of Engel algebras, as indicated in the introduction), one should take into account specific properties arising
in characteristic zero.

\begin{defn}
An algebra $ A$ is  {\em nilpotent} of   {\em class} $\le n$ if $M(a_1, a_2, \dots )= 0$  for each monomial
$M(x_1, x_2, \dots )$ of degree $\ge n$. An  ideal $I$ of $ A$ is  {\em strongly nilpotent} of {\em class} $\le n$ if
$M(a_1, a_2, \dots )= 0$  for each monomial $M(x_1, x_2, \dots )$  in which indeterminates of total degree $\ge
n$ have been substituted to elements of $I$.
\end{defn}

 Although the notions of nilpotent and  strongly nilpotent coincide in the associative case,
they differ for ideals of nonassociative algebras. For example, consider the following algebra suggested by
Shestakov: $A$ is the algebra generated by $a,b,z$ satisfying the relations $a^2 =b,$ $ bz = a$ and all other
products 0. Then $I = Fa+Fb$ is nilpotent as a subalgebra,
 satisfying $I^3 = 0$ but not strongly nilpotent (as an ideal), since
$$b = ((a(bz))z)a \ne 0,$$
and one can continue indefinitely in this vein.
 Also, \cite{KuS} contains an example of a finite dimensional non-associative algebra without any ideal which is maximal witih respect to being nilpotent as a subalgebra.

 In connection with the Generalized Jacobian Conjecture
in characteristic~0, it follows  from results of Yagzhev ~\cite{Yagzhev9}, also cf.~\cite{GorniZampieri},  that there
exists a $20$-dimensional Engel algebra over $\mathbb Q$, not weakly nilpotent, satisfying the identities
\begin{align*}
x^2y=-yx^2, \quad (((yx^2)x^2)x^2)x^2=0,\cr (xy+yx)y=2y^2x, \quad
x^2y^2=0.
\end{align*}

  However, this algebra can be seen to be
Yagzhev (see Definition \ref{DfWeaklyNilp}).

 For associative
algebras, one uses the term ``nil'' instead of ``weakly nilpotent.'' Any nil subalgebra of a
  finite  dimensional associative algebra is nilpotent, by
Wedderburn's Theorem \cite{Wedderburn}). Jacobson generalized this result to other settings, cf.~\cite[Theorem
15.23]{Row}, and Shestakov~\cite{ShestakovWedderburn} generalized it to a wide class of
  Jordan algebras (not necessarily commutative).

\medskip
Yagzhev's investigation of weak nilpotence has applications to the Koethe Conjecture,  for algebras over
uncountable fields.  He reproved:

\medskip
{\bf *} {\it In every associative algebra over an uncountable field, the sum of every two nil right  ideals is a nil right
ideal} \cite{yagzhevKoethe}.
\medskip

(This was proved first by Amitsur~\cite{Am}. Amitsur's result is for affine algebras, but one can easily reduce to the
affine
case.) 
%
%

\subsubsection{Algebras satisfying systems of Capelli identities}\label{SbSbScCap}

\begin{definition}\label{Cap1} The
 {\it Capelli polynomial} $C_k$ of order $k$ is
 $$C_k := \sum _{\sigma\in S_k } (-1)^{\sigma}x_{\sigma(1)}y_1 \cdots
 x_{\sigma(k)}y_k.$$
 \end{definition}

  It is obvious that an associative algebra
satisfies the  Capelli identity  $c_k$ iff, for any monomial $M(x_{1},\ldots,x_{k},y_{1},\ldots,y_{r})$ multilinear in
the $x_i$, the following equation holds  identically in $A$:
\begin{equation}        \label{EquVarStar}
\sum_{\sigma\in S_k }(-1)^{\sigma} M(v_{\sigma(1)},\ldots,
  v_{\sigma(k)},
  y_{1},\ldots,y_{r}) = 0.
\end{equation}


However, this does not apply to nonassociative algebras, so we need to generalize this condition.

\begin{definition} The algebra $A$ satisfies a
 {\bf system of Capelli identities} of order $k$, if \eqref{EquVarStar}
holds identically in $A$ for any monomial $M(x_{1},\ldots,x_{k},y_{1},\ldots,y_{r})$ multilinear in the $x_i$.
\end{definition}

Any algebra of dimension $<k$ over a field satisfies a
 system of Capelli identities  of order $k$.
%
%
Algebras satisfying systems of Capelli identities behave much like finite dimensional algebras. They were introduced
and systematically studied by  Rasmyslov \cite{RazmyslovIANCapely}, \cite{RazmyslovBook}.

 Using Rasmyslov's method,
Zubrilin~\cite{Zubrilin4}, also see \cite{RazmZubrilin,Zubrilin1}, proved that if $A$ is an arbitrary algebra satisfying
the system of Capelli identities of order~$n$, then the chain of ideals defining the {\it solvable radical} stabilizes at
the $n$-th step. More precisely, we utilize a Baer-type radical, along the lines of Amitsur~\cite{am2}.

Given an algebra $A$, we define $\Nilp_1: = \Nilp_1(A)  = \sum \{$Strongly nilpotent  ideals of $A\},$ and
inductively, given $\Nilp_k $, define $\Nilp_{k+1} $ by $\Nilp_{k+1}/\Nilp_k = \Nilp_1(A/\Nilp_k).$ For a limit
ordinal $\alpha,$ define
$$\Nilp_\alpha  = \cup _{\beta < \alpha} \Nilp_\beta.$$ This must
stabilize at some ordinal $\alpha$, for which we define $\Nilp(A) = \Nilp_\alpha.$

  \medskip

Clearly $\Nilp(A/\Nilp(A)) = 0;$ i.e., $A/\Nilp(A)$ has no nonzero strongly nilpotent ideals.   Actually,
Amitsur~\cite{am2} defines $\zeta(A)$ as built up from ideals having trivial multiplication, and proves
\cite[Theorem~1.1]{am2} that $\zeta(A)$ is the intersection of the prime ideals of $A$.

We shall use the notion of {\it sandwich}, introduced by Kaplansky and Kostrikin, which is a powerful tool for
Burnside type problems \cite{Zelmanov}. An ideal $I$ is called a {\it sandwich ideal} if, for any $k$,
$$M(z_1,z_2,x_1,\dots,x_k)=0$$
for any $z_1, z_2\in I$, any set of elements $x_1,\dots,x_k$, and any multilinear monomial $M$ of degree $k+2$.
(Similarly, if the operations of an algebra have degree $\le \ell$, then it  is natural to use $\ell$-sandwiches, which by
definition satisfy the property that $$M(z_1,\dots,z_\ell,x_1,\dots,x_k)=0$$ for any $z_1, \dots,z_\ell \in I$, any set
of elements $x_1,\dots,x_k$, and any multilinear monomial $M$ of degree $k+\ell$.)

The next useful lemma follows from a result from \cite{Zubrilin4}:

\begin{lemma}     \label{LeSandvich}
If  an ideal $I$ is strongly  nilpotent of class  $\ell$, then there exists a decreasing sequence of ideals
$I=I_1\supseteq\dots\supseteq I_{l+1}=0$ such that $I_s/I_{s+1}$ is a sandwich ideal in $A/I_{s+1}$ for all
$s\le l$.
\end{lemma}

\begin{definition}     \label{DfRepresAlg}
An algebra $A$ is {\em representable}  if it can be embedded into an algebra finite dimensional over some extension
of the ground field.
\end{definition}

\begin{rem} \ Zubrilin \cite{Zubrilin4}, properly clarified, proved the more precise
statement, that if an algebra $A$ of arbitrary signature satisfies a system of Capelli identities $C_{n+1}$, then there
exists a sequence $B_0\subseteq B_1\subseteq\dots\subseteq B_n$  of strongly nilpotent ideals  such that:
\begin{itemize}
  \item The natural projection of $B_i$ in $A/B_{i-1}$ is a strongly nilpotent ideal.

  \item $A/B_n$ is representable.
  \item If $I_1\subseteq I_2\subseteq\dots\subseteq I_n$ is any sequence of ideals of $A$ such that
      $I_{j+1}/I_{j}$ is a sandwich ideal in $A/I_j$, then $B_n\supseteq I_n$.
\end{itemize}

 Such a sequence of ideals will be called a {\em Baer-Amitsur
  sequence}.
 In affine space the Zariski closure of the radical is radical, and hence the factor algebra is representable.
 (Although the radical  coincides with the linear closure if the base field is infinite (see \cite{BelovUzyRowen2}),
 this assertion holds for arbitrary  signatures and base fields.) Hence in representable algebras,
 the Baer-Amitsur sequence stabilizes after finitely many steps.
  Lemma \ref{LeSandvich} follows from these
considerations.
\end{rem}

Our next main goal is to prove Theorem \ref{ThSandvich} below, but first we need another notion.
\subsubsection{The tree associated to a monomial}         \label{SbSbScTreeMon}

 Effects of nilpotence have been used   by different authors in
another language. We  associate a {\it rooted labelled tree} to  any monomial: Any branching vertex indicates the
symbol of an operator, whose outgoing edges are the terms in the corresponding symbol. Here is the precise
definition.

\begin{definition}     \label{DfTreeMon}
Let $M(x_1,\dots,x_n)$ be a monomial in an algebra $A$ of arbitrary signature. One can associate the tree $T_M$
by an inductive procedure:

\begin{itemize}
  \item If $M$ is a single variable, then $T_M$ is just the vertex $\bullet$.
  \item Let $M=g(M_1,\dots,M_k)$, where $g$ is a $k$-ary operator. We assume  inductively    that the trees
      $T_i,$ $ i=1,\dots, k,$ are already defined.
Then the tree $T_M$ is the disjoint union of the $T_i$, together with the root $\bullet$ and arrows starting with
$\bullet$ and ending with the roots of the trees $T_i$.
\end{itemize}

\end{definition}

{\bf Remark.} Sometimes one  labels $T_M$ according to the operator $g$ and the positions inside $g$.

\medskip

If the outgoing degree of each vertex is 0 or 2, the tree is called {\it binary}. If the outgoing degree of each vertex is
either 0 or 3, the tree is called {\it ternary}. If each operator is binary, $T_M$ will be binary; if each operator is
ternary, $T_M$ will be ternary.

\subsection{Lifting Yagzhev algebras}\label{LiftY}

 Recall Definitions \ref{DfWeaklyNilp}  and   \ref{Engtyp}.

\begin{theorem}         \label{ThSandvich}
Suppose $A$ is an algebra of   Engel type, and let $I$ be a sandwich ideal of $A$. If $A/I$ is Yagzhev, then $A$ is
Yagzhev.
\end{theorem}

{\bf Proof.} The proof follows easily from the following two propositions.

\medskip

  Let $k$ be the class of weak  nilpotence of $A/I$. We call a
branch of the tree {\it fat} if it has more than $k$ entries.

\begin{statement}
a) The sum of all monomials of any degree $s>k$ belongs to $I$.

b) Let $x_1,\dots, x_n$ be  fixed indeterminates, and $M$ be an arbitrary monomial, with $s_1,\dots,s_\ell>k$.
Then

\begin{equation}      \label{EqSnvch}
\sum_{|t_1|=s_1,\dots, |t_\ell|=s_\ell} M(x_1,\dots,x_n,t_1,\dots, t_\ell)\equiv 0.
\end{equation}

c)  The sum of all monomials of degree $s$, containing at least $\ell$ non-intersecting fat branches, is zero.
\end{statement}

{\bf Proof.} a) is just a reformulation of the weak  nilpotence
 of~$A/I$; b) follows from a) and the sandwich property
of an ideal $I$. To get c) from~b),  it is enough to consider the highest non-intersecting fat branches.

\begin{statement}[Yagzhev]
The linearization of the sum of all terms with a fixed fat branch of length $n$ is the complete linearization of the
function
$$
\sum_{\sigma\in S_n}\prod(\Ad_{k_{\sigma(i)}})(z)(t).
$$
\end{statement}


Theorem 1.2, Lemma \ref{LeSandvich}, and  Zubrilin's  result give us the following major result:

\begin{theorem}\label{eqcond}
In characteristic zero, the Jacobian conjecture is equivalent to the following statement:

Any algebra of Engel type satisfying some system of Capelli identities is a Yagzhev algebra.
\end{theorem}

This theorem generalizes the following result of Yagzhev:

\begin{theorem}
The Jacobian conjecture is equivalent to the following statement:

Any ternary Engel algebra in characteristic 0 satisfying a system of Capelli identities is a Yagzhev algebra.
\end{theorem}

The Yagzhev correspondence and the results of this section
  (in particular, Theorem~\ref{eqcond}) yield the proof of
  Theorem~\ref{Yag2}.

\subsubsection{Sparse identities}

Generalizing  Capelli identities, we say that an algebra satisfies a system of {\it sparse identities} when there exist
$k$ and coefficients $\alpha_{\sigma}$ such that for any monomial $M(x_{1},\ldots,x_{k},y_{1},\ldots,y_{r})$
multilinear in $x_i$ the following equation holds:
\begin{equation}        \label{EquVarStarGM}
\sum_{\sigma }\alpha_{\sigma} M(c_{1}v_{\sigma(1)}d_{1},\ldots,
  c_{k}v_{\sigma(k)}d_{k},
  y_{1},\ldots,y_{r}) = 0.
\end{equation}

Note that one need only check  \eqref{EquVarStarGM} for monomials. The system of Capelli identities is a special
case of a system of sparse identities (when $\alpha_{\sigma}=(-1)^\sigma$). This concept ties in with the following
``few long branches'' lemma \cite{Zubrilin2}, concerning the structure of trees of monomials for algebras with sparse
identities:

\begin{lemma}[Few long branches]\label{spars}
Suppose an algebra $A$ satisfies a system of sparse identities of order $m$. Then any monomial is linearly
representable by monomials such that the corresponding tree has not more than  $m-1$ disjoint branches of length
$\ge m$.

\end{lemma}

 Lemma~\ref{spars} may   be useful in studying nilpotence of Engel algebras.
%
%

\subsection{Inversion formulas and problems of  Burnside type}
We have seen that the $\JC$  relates to problems of ``Specht type'' (concerning whether one set of polynomial
identities implies another), as well as problems of Burnside type.

Burnside type problems become more complicated in  nonzero characteristic; cf.~Zelmanov's review article
\cite{Zelmanov}.

  Bass, Connell, and Wright~\cite{BCW} attacked the $\JC$ by means
  of inversion formulas.
 D.~Wright \cite{Wright_B}
wrote an inversion formula for the symmetric case and related it to a combinatorial structure called the {\it
Grossman--Larson Algebra}. Namely, write $F=X-H$, and define $J(H)$ to be the Jacobian matrix of $H$.
Wright proved the $\JC$ for the case where $H$ is homogeneous and $J(H)^3=0$, and also for the case where
$H$ is cubic and $J(H)^4=0$; these correspond in Yagzhev's terminology to the cases of Engel type $3$ and $4$,
respectively. Also, the so-called {\it chain vanishing theorem} in~\cite{Wright_B} follows from  Engel type. Similar
results were obtained earlier by Singer~\cite{Singer} using tree formulas for formal inverses. The inversion formula,
introduced in \cite{BCW}, was investigated by D.~Wright and his school.  Many authors  use the language of
so-called  {\it tree expansion} (see \cite{Wright_B,Singer} for details). In view of Theorem
\ref{ThEngWeakNoNilp}, the tree expansion technique  should be highly nontrivial.

 The Jacobian Conjecture can be formulated
as a question of quantum field theory (see \cite{Abdelmalek}), in which tree expansions are seen to correspond to
Feynmann diagrams.

In the papers \cite{Singer} and \cite{Wright_B} (see also \cite{Wright1}),  trees with one label correspond  to
elements of the algebra $A$  built by Yagzhev, and $2$-labelled trees correspond to the elements of the operator
algebra $D(A)$ (the algebra generated by operators $x\to M(x,\vec{y})$, where $M$ is some monomial). These
authors deduce weak nilpotence from the Engel conditions of degree 3 and 4.   The inversion formula for
automorphisms of tensor product of Weyl algebras and the ring of polynomials was studied intensively in the papers
\cite{Bavula,Bavula1new}. Using  techniques from \cite{K-BK2}, this yields a slightly different proof of the
equivalence between the  $\JC$ and DC, by an argument similar to one given in \cite{YagzhevLast}. Yagzhev's
approach makes the situation much clearer, and the known approaches to the Jacobian Conjecture using inversion
formulas can be explained from this viewpoint.

\begin{rem}\label{Jac4} The most
recent inversion formula (and probably the most algebraically explicit one) was obtained by V.~Bavula
\cite{Bavula3}. The coefficient $q_0$ can be made explicit in \eqref{EqTerm2}, by means of  the Gabber
Inequality, which says that if
$$f: K^n\to K_n; \quad x_i\to f_i(\vec x)$$ is a polynomial
automorphism, with
 $\deg(f)=\max_i\deg(f_i)$, then $\deg(f^{-1})\le \deg(f)^{n-1}$)
\end{rem}

In fact, we are working with {\it operads}, cf.~the classical book~\cite{MSS}. A review of operad theory and its
relation with physics and $\PI$-theory in particular Burnside type problems, will appear in D.~Piontkovsky
\cite{Piontkovski}; see also \cite{Piontkovski1,PiontkovskiKhoroshkin}.  Operad theory provides a supply of
natural identities and varieties, but they also correspond  to geometric facts. For example, the Jacobi identity
corresponds to the fact that the altitudes of a triangle are concurrent. M.~Dehn's observations that the Desargue
property of a projective plane corresponds to associativity of its coordinate ring, and Pappus' property to its
commutativity, can be considered as a first step in operad theory. Operads are important in mathematical physics,
and formulas for the famous Kontsevich quantization theorem resemble formulas for the inverse mapping. The
operators considered here are operads.


\section{The Jacobian Conjecture for varieties, and deformations}\label{SbScJCArbVar0}

In this section we consider   analogs of the  $\JC$ for other varieties of algebras, partially with the aim on throwing
light on the classical  $\JC$ (for the commutative associative polynomial algebra).

\subsection{Generalization of the Jacobian Conjecture to arbitrary
varieties}                  \label{SbScJCArbVar}

J.~Birman~\cite{Bir} already proved the $\JC$ for free groups in 1973.
 The $\JC$ for free associative algebras (in two
generators) was established   in~1982 by~W.~Dicks and J.~Levin \cite{DiLev, Di}, utilizing Fox derivatives, which
we describe later on. Their result was reproved by Yagzhev~\cite{Yag2}, whose  ideas are sketched in this section.
Also see  Schofield~\cite{Schof}, who proved  the full version. Yagzhev then applied these ideas to other varieties
of algebras \cite{Yag1,Yagzhev9} including nonassociative commutative algebras and anti-commutative algebras;
U.U.~Umirbaev \cite{Umirbaev} generalized these to
``Schreier varieties,'' defined by the property that every subalgebra of a free algebra is free. 
The $\JC$ for free Lie algebras was proved by Reutenauer~\cite{Reu}, Shpilrain~\cite{Shp}, and Umirbaev
\cite{Umi}.

The Jacobian Conjecture for varieties generated by finite dimensional algebras, is closely related to the Jacobian
Conjecture in the usual commutative associative  case, which is the most important.

Let $\goth M$ be a variety of algebras of some signature~$\Omega$ over a given field~$\k$ of characteristic zero,
and $\k_{\goth M}\!<\vec{x}\!>$~the relatively free algebra in $\goth M$ with generators $\vec{x}=\{x_i :i\in I\}$.
We assume that $|\Omega|, |I|<\infty$, $I=1,\dots, n$.

Take a set $\vec{y}=\{y_i\}_{i=1}^n$ of new indeterminates. For any $f(\vec{x})\in k_{\goth M}\!<\!\vec{x}\!>$
one can define an element $\hat{f}(\vec{x},\vec{y})\in k_{\goth M}\!<\!\vec{x},\vec{y}\!>$ via the equation

\begin{equation}         \label{EqDerFn}
f(x_1+y_1,\dots,x_n+y_n)=f(\vec{y})+\hat{f}(\vec{x},\vec{y})+
R(\vec{x},\vec{y})
\end{equation}
\noindent where $\hat{f}(\vec{x},\vec{y})$ has degree $1$ with respect to $\vec x$, and $R(\vec{x},\vec{y})$ is
the sum of monomials of degree  $\ge 2$  with respect to $\vec x$; $\hat{f}$ is a generalization of the differential.

Let $\alpha\in\End(k_{\goth M}\!<\!\vec{x}\!>)$, i.e.,
\begin{equation}      \label{EqEnd}
\alpha: x_i\mapsto f_i(\vec{x});\ i=1,\dots,n.
\end{equation}

\begin{definition}\label{Jac1}
Define the Jacobi endomorphism $\hat{\alpha}\in\End(\k_{\goth M}\!<\!\vec{x},\vec{y}\!>)$ via the equality
\begin{equation}       \label{EqDif}
\hat{\alpha}:
\left\{
\begin{array}{cc}
x_i\to\hat{f}_i(\vec{x}),\cr y_i\to y_i.
\end{array}
\right.
\end{equation}
\end{definition}

The Jacobi mapping $f \mapsto \hat{ f}$ satisfies the chain rule, in the sense that it  preserves composition.

\begin{rem}\label{chain}
 It is not
difficult to check (and is well known) that if $\alpha\in\Aut(\k_{\goth M}\!<\!\vec{x}\!>)$ then
$\hat{\alpha}\in\Aut(\k_{\goth M}\!<\!\vec{x},\vec{y}\!>)$.
\end{rem}

The inverse implication is called the {\it Jacobian Conjecture for the variety~$\goth M$}. Here is an important
special case.

\begin{definition}
Let $A\in{\goth M}$ be a finite dimensional algebra, with base $\{\vec{e}_i\}_{i=1}^N$. Consider a set of
commutative indeterminates $\vec{\nu}=\{\nu_{si}|s=1,\dots,n; i=1,\dots,N\}$. The elements
$$z_j=\sum_{i=1}^N\nu_{ji} \vec{e}_i; \quad  j=1,\dots,n$$
are called  {\em generic elements of $A$}.

\end{definition}

Usually in the matrix algebra ${\mathbb M}_m(\k),$ the set of matrix units $\{e_{ij}\}_{i,j=1}^m$ is taken as the
base. In this case $e_{ij}e_{kl}=\delta_{jk}e_{il}$ and $z_l=\sum_{ij}\lambda_{ij}^le_{ij}$, $l=1,\dots,n$.

\begin{definition}
A   {\em generic matrix} is a matrix whose entries are distinct
   commutative indeterminates, and the so-called   {\em algebra of
generic matrices of order $m$} is generated by associative generic $m\times m$ matrices.
\end{definition}

The algebra of generic matrices is prime, and every prime,
  relatively free, finitely generated associative
$PI$-algebra is isomorphic to an algebra of generic matrices. If we include taking traces as an operator in the
signature, then we get the {\it algebra of generic matrices with trace}. That algebra is a Noetherian module over its
center.

Define the $\k$-linear mappings
$$
\Omega_i: \k_{\goth M}\!<\!\vec{x}\!>\to \k[\nu];\quad i=1,\dots,
n
$$
via the relation
$$
f(\sum_{i=1}^N\nu_{1i}e_i,\dots,\sum_{i=1}^N\nu_{ni}e_i)=
\sum_{i=1}^N (f\Omega_i)e_i.
$$
It is easy to see that the polynomials $f\Omega_i$ are uniquely determined by $f$.

One can define the mapping
$$
\varphi_A: \End(k_{\goth M}\!<\!\vec{x}\!>)\to\End(k[\vec{\nu}])
$$
as follows: If
$$
\alpha\in\End(k_{\goth M}\!<\!\vec{x}\!>): x_s\to f_s(\vec{x})\quad
s=1,\dots,n
$$
then $\varphi_A(\alpha)\in\End(k[\vec{\nu}])$ can be defined via the relation
$$
\varphi_A(\alpha): \nu_{si}\to P_{si}(\vec{\nu}); \quad s=1,\dots,n; \quad
i=1,\dots,n,
$$
where $P_{si}(\vec{\nu})=f_s\Omega_i$.

The following proposition is well known.

\begin{statement}[\cite{Yagzhev9}]
Let $A\in{\goth M}$ be a finite dimensional algebra, and $\vec{x}=\{x_1,\dots,x_n\}$ be a finite set of
commutative indeterminates. Then the mapping $\varphi_A$ is a semigroup homomorphism, sending $1$ to $1$,
and automorphisms to automorphisms. Also the mapping $\varphi_A$ commutes with the operation $\widehat{} $
of taking the Jacobi endomorphism, in the sense that $\widehat{\varphi_A(\alpha)}=\varphi_A(\hat{\alpha})$. If
$\varphi$ is invertible, then $\widehat{\varphi}$ is also invertible.
%
\end{statement}

This proposition is important, since as noted after Remark~\ref{chain}, the opposite direction is  the  $\JC$.

\subsection{Deformations and the Jacobian Conjecture for free associative
algebras}

\begin{definition} A {\it $T$-ideal} is a   completely characteristic ideal, i.e., stable under
any endomorphism.\end{definition}

\begin{statement}
Suppose $A$ is a relatively free algebra in the variety $\goth M$, $I$~ is a $T$-ideal in $A$, and ${\goth
M}'=\Var(A/I)$. Any polynomial mapping $F: A\to A$
 induces  a natural mapping $F':A/I\to A/I$, as well as a
mapping $\widehat{F'}$ in~ ${\goth M}'$. If $F$ is invertible, then $F'$ is invertible; if $\hat{F}$ is invertible, then
$\widehat{F'}$ is also invertible.
\end{statement}

For example, let $F$ be a polynomial endomorphism of the free associative algebra $k\!<\vec{x}\!>$, and $I_n$
be the $T$-ideal of the algebra of generic matrices of order $n$. Then $F(I_n)\subseteq I_n$ for all $n$. Hence
$F$ induces an endomorphism $F_{I_n}$ of $k\!<\vec{x}\!>/I_n$. In particular, this is a semigroup
homomorphism. Thus, if $F$ is invertible, then $F_{I_n}$ is invertible, but not vice versa.

The Jacobian mapping $\widehat{F_{I_n}}$ of the reduced endomorphism $F_{I_n}$ is   the reduction of the
Jacobian mapping of $F$.

\subsubsection{The Jacobian Conjecture and the packing property}   \label{SbSbScPacking}

\indent   This subsection is based on the {\it packing property} and deformations. Let us illustrate the main idea. It is
well known that the composite of ALL quadratic extensions of $\mathbb Q$ is infinite dimensional over $\mathbb
Q$.  Hence all such extensions cannot be embedded (``packed'')  into a single commutative finite dimensional
$\mathbb Q$-algebra. However, all of them can be packed into ${\mathbb M}_2(\mathbb Q)$. We formalize the
notion of packing in \S\ref{pack1subsec}.  Moreover, for ANY elements NOT in $\mathbb Q$ there is a
parametric family of embeddings (because it embeds non-centrally and thus can be deformed via conjugation by a
parametric set of matrices). Uniqueness thus means belonging to the center. Similarly, adjoining noncommutative
coefficients allows one to decompose polynomials, as to be elaborated below.

This idea allows us to solve equations via a finite dimensional extension, and to find  a parametric sets of solutions if
some solution does not belong  to the original algebra. That situation contradicts   local invertibility.

Let $F$ be an endomorphism of the free associative algebra having invertible Jacobian.  We suppose that $F(0)=0$
and
$$F(x_i)=x_i+\sum\mbox{terms\ of}\ \mbox{order}\ge 2 .$$
We intend to show how the invertibility of the Jacobian implies invertibility of the mapping $F$.

Let $ {Y}_1,\dots, {Y}_k$ be generic $m \times m$ matrices. Consider the system of equations
$$
\left\{F_i( {X}_1,\dots, {X}_n)=Y_i; \quad i=1,\dots,k
\right\}.
$$
This system has a solution over some finite extension of order $m$ of the field generated by the center of the algebra
of generic matrices {\it with trace}.

Consider the set of block diagonal $mn\times mn$ matrices:

\begin{equation}        \label{EqMatrType}
A=\left(
\begin{array}{lcccl}
A_1&0& &\dots&0\cr 0&A_2&0&\dots&0\cr \vdots& &\ddots&
&\vdots \cr 0& &\dots& &A_n
\end{array}
\right),
\end{equation}
where the  $A_j$ are   $m\times m$ matrices.

Next, we  consider the system of equations
\begin{equation}\label{EqMatnk}
    \left\{F_i(X_1,\dots,X_n)=Y_i; \quad i=1,\dots,k \right\},
\end{equation}
where the $mn\times mn$ matrices $Y_i$ have the form (\ref{EqMatrType}) with the $A_j$ generic matrices.

Any $m$-dimensional  extension of the base field $\k$ is embedded into~$\mathbb{M}_m(\k)$. But
$\mathbb{M}_{mn}(\k)\simeq \mathbb{M}_{m}(\k)\otimes\mathbb{M}_{n}(\k)$. It follows that for appropriate
$m$, the system~(\ref{EqMatnk}) has a unique solution in the matrix ring with traces. (Each    is given by a matrix
power series where the summands are matrices whose entries  are homogeneous forms, seen by  rewriting   $Y_i  =
X_i + \sum \text{terms of order2}$ as       $X_i  =  Y_i + \sum \text{terms of order  2}$, and iterating.)  The
solution is unique  since their entries are distinct commuting indeterminates.

If $F$ is invertible, then this solution must have   block diagonal form. However, if $F$ is not invertible, this solution
need not have  block diagonal form.  Now we translate invertibility of the Jacobian to the language of {\bf parametric
families} or {\bf deformations}.

Consider the matrices
$$
E_\lambda^\ell=
\left(
\begin{array}{lcccr}
E&0&& \dots&0\cr 0&\ddots&&\dots&0\cr 0&\dots&\lambda\cdot E&
&0\cr \vdots& &\dots&\ddots&\vdots\cr 0& &\dots&  &E
\end{array}
\right) $$ where $E$ denotes the identity matrix. (The index
$\ell$ designates the position of the  block $\lambda\cdot E$.) Taking  $X_j$ not to be  a block diagonal matrix,
then for some $\ell$ we obtain a non-constant parametric family $E_\lambda^\ell X_j(E_\lambda^\ell )^{-1}$
dependent on $\lambda$.

On the other hand, if $Y_i$ has form (\ref{EqMatrType}) then $E_\lambda^\ell Y_i(E_\lambda^\ell)^{-1}=Y_i$
for all $\lambda\ne 0$; $ \ell=1,\dots,k$.

Hence, if $F_{I_n}$ is not an automorphism, then we have a {\bf continuous parametric set of solutions}. But if the
Jacobian mapping is invertible, it is locally 1:1, a contradiction. This argument yields the following result:

\begin{theorem}
For $F\in\End(k\!<\vec{x}\!>)$, if the Jacobian of $F$ is invertible, then  the reduction $F_{I_n}$ of $F$, modulo
the $T$-ideal of the algebra of generic matrices, is invertible.
\end{theorem}

See \cite{Yag1} for further details of the proof. Because any  relatively free affine algebra of characteristic 0
  satisfies the set of identities of some
matrix algebra, it is the quotient of the algebra of generic matrices by some $T$-ideal $J$. But $J$ maps into itself
after any endomorphism of the algebra. We conclude:

\begin{corollary}
If $F\in\End(k\!<\!\vec{x}\!>)$ and the Jacobian of $F$ is invertible, then the reduction $F_{J}$ of $F$ modulo
any proper $T$-ideal $J$  is invertible.
\end{corollary}

In order to get invertibility of   $\vec F$ itself, Yagzhev used the additional ideas:
\begin{itemize}
\item The block diagonal technique works equally well on  skew fields. \item The above  algebraic constructions
    can be carried out on Ore extensions, in particular for the {\it Weyl algebras}
    $W_n=\k[x_1,\dots,x_n;\partial_1,\dots,\partial_n]$.

\item By a result of L.~Makar-Limanov, the free associative algebra   can be embedded into the ring of fractions of
    the Weyl algebra. This provides a nice presentation for mapping the free algebra.
\end{itemize}


\begin{definition}
Let $A$ be an algebra, $B\subset A$ a subalgebra, and $\alpha: A\to A$  a polynomial mapping of $A$ (and hence
$\alpha(B)\subset B$, see Definition \ref{DefPolMap}). $B$ is a {\it test algebra for $\alpha$}, if
$\alpha(A\backslash B)\ne A\backslash B$.
\end{definition}

The next theorem shows the universality of the notion of a test algebra. An endomorphism is called {\it rationally
invertible} if it is invertible over {Cohn's}
 skew field of fractions~\cite{Cohn1} of
$\k\!<\!\vec{x}\!>$.

\begin{theorem}[Yagzhev]
For any $\alpha\in \End(k\!<\vec{x}\!>)$, one of the two statements holds:
\begin{itemize}
\item $\alpha$ is rationally invertible, and its reduction to any finite dimensional factor  also is rationally invertible.
    \item There exists a test algebra for some finite dimensional reduction of $\alpha$.
\end{itemize}
\end{theorem}

This theorem implies the Jacobian conjecture for free associative algebras. We do not go into details,
 referring the reader to the
papers \cite{Yag1} and \cite{Yagzhev9}.

\medskip

{\bf Remark.} The same idea is used in quantum physics. The polynomial $x^2+y^2+z^2$ cannot be decomposed
for any commutative ring of coefficients. However, it can decomposed using noncommutative ring of coefficients
(Pauli matrices). The Laplace operator in $3$-dimensional space can be decomposed in such a manner.

\subsubsection{Reduction to nonzero characteristic}

 One can work with deformations equally well in nonzero
characteristic. However, the naive Jacobian condition does not give us parametric families, because of
 consequences  of inseparability. Hence it is interesting using
deformations to get a reasonable version of the $\JC$
 for characteristic $p> 0$, especially because of recent
progress in the    $\JC$  related to the reduction of holonomic modules to the case of characteristic $p$ and
investigation of the $p$-curvature or Poisson brackets on the center \cite{K-BK2}, \cite{K-BK1}, \cite{Tsu1}.

In his very last paper \cite{YagzhevLast} A.V.~Yagzhev approached the    $\JC$ using positive characteristics. He
noticed that the existence of a counterexample
 is
equivalent to the existence of an Engel, but not Yagzhev, finite dimensional ternary algebra in each positive
characteristic $p\gg 0$. (This fact is also used in the papers \cite{K-BK2,K-BK1,Tsu1}.)

If a counterexample to the    $\JC$ exists, then  such an algebra $A$ exists  even over a finite field, and hence  can
be finite. It generates a locally finite variety of algebras that are of Engel type, but not Yagzhev. This situation can be
reduced to the case of a locally semiprime variety. Any relatively free algebra of this variety is semiprime, and the
centroid of its localization is a finite direct sum of fields. The situation can be reduced to one field, and he tried to
construct an embedding which is not an automorphism. This would contradict the finiteness property.

Since a reduction of an endomorphism as a mapping on points of finite height may be an automorphism,  the issue of
injectivity also arises. However, this approach looks promising, and may involve new ideas, such as  in the papers
\cite{K-BK2,K-BK1,Tsu1}. Perhaps different infinitesimal conditions (like the Jacobian condition in characteristic
zero) can be found.

\subsection{The Jacobian Conjecture for other classes of
algebras} Although the Jacobian Conjecture remains open for commutative associative algebras, it has been
established for other classes of algebras, including free associative algebras, free Lie algebras, and free metabelian
algebras. See \S\ref{SbScJCArbVar} for further details.

 An
algebra is {\it metabelian} if it satisfies the identity $[x,y][z,t]=0$.

The case of free metabelian algebras, established by Umirbaev \cite{umirbaev1995ext}, involves some interesting
new ideas that we describe now. His method of proof is by means of co-multiplication, taken from the theory of
Hopf algebras and quantization. Let $A^{op}$ denote the opposite algebra of  the free associative algebra $A$,
with generators $t_i$. For $f\in A$ we denote the corresponding element of $A^{op}$ as  $f^*$. Put $\lambda:
A^{op}\otimes A\to A$ be the mapping such that $\lambda(\sum f_i^* \otimes g_i)=\sum f_ig_i$.
$I_A:=\ker(\lambda)$ is a free $A$ bimodule with generators $t_i^*\otimes 1-1\otimes t_i$. The mapping $d_A:
A\to I_A$ such that $d_A(a)=a^*\otimes 1-1\otimes a$ is called the {\it universal derivation} of $A$. The {\it Fox
derivatives} $\partial a/\partial t_i\in A^{op}\otimes A$ \cite{F} are defined via $d_A(a)=\sum_i (t_i^*\otimes
1-1\otimes t_i)\partial a/\partial t_i$, cf.~\cite{DiLev} and~\cite{umirbaev1995ext}.

Let $C=A/\Id([A,A])$, the free commutative associative algebra, and let $B=A/\Id([A,A])^2$, the free metabelian
algebra. Let
$$\partial(a)=(\partial a/\partial t_1,\dots,\partial a/\partial
t_n).$$ One can define the natural derivations

$$\bar{\partial}:A\to (A'\otimes A)^n\to(C'\otimes C)^n,$$

\begin{equation}\label{nat} \tilde{\partial}:A\to (C'\otimes C)^n\to C^n.\end{equation}
where the mapping $(C'\otimes C)^n)\to C^n$ is induced by $\lambda$.  Then
$\ker(\bar{\partial})=\Id([A,A])^2+F$ and $\bar{\partial}$ induces a derivation $B\to (C'\otimes C)^n$, whereas
$\tilde{\partial}$ induces the usual derivation $C\to C^n$. Let $\Delta: C\to C'\otimes C$ be the mapping induced
by $d_A$, i.e., $\Delta(f)=f^*\otimes 1-1\otimes f$, and let $z_i=\Delta(x_i)$.
The {\it Jacobi matrix} is defined in the  natural way, and provides the formulation of the
 $\JC$ for free metabelian algebras. One of the crucial steps in
proving the~$\JC$ for free metabelian algebras is the following homological lemma from~\cite{umirbaev1995ext}:

\begin{lemma}
Let $\vec{u}=(u_1,\dots,u_n)\in (C^{op}\otimes C)^n$. Then $\vec{u}=\bar{\partial}(\bar{w})$ for some
$w\in\Id([A,A])$ iff
$$\sum z_iu_i=0.$$
\end{lemma}

The proof also requires the following theorem:

\begin{theorem}
Let $\varphi\in\End(C)$. Then $\varphi\in\Aut(C)$ iff $\Id(\Delta( \varphi(x_i) ))_{i=1}^n=\Id(z_i)_{i=1}^n$.
\end{theorem}

The paper \cite{umirbaev1995ext} also includes the following result:

\begin{theorem}  \label{ThBClifting}
Any automorphism of $C$ can be extended to an automorphism of $B$, using the $\JC$ for the free metabelian
algebra $B$.
\end{theorem}

This is a nontrivial result, unlike the extension of an automorphism of $B$ to an automorphism of $A/\Id([A,A])^n$
for any $n>1$.
%
%
%
%
%
%
%
%
%

\subsection{Questions related to the Jacobian Conjecture}                            \label{SbScRelQuestJC}
Let us turn to other interesting questions which can be linked to the Jacobian Conjecture. The quantization
procedure is a bridge between the commutative and noncommutative cases and is deeply connected to the $\JC$
and related questions. Some of these questions also are discussed in the paper \cite{DrYuLift}.

Relations between the free associative algebra and the classical commutative situation are very deep. In particular,
Bergman's theorem that any commutative subalgebra of the free associative algebra is isomorphic to a polynomial
ring in one indeterminate is the noncommutative analog of Zak's theorem \cite{Zaks} that any integrally closed
subring  of a polynomial ring of Krull dimension $1$ is isomorphic to a polynomial ring in one indeterminate.

For example, Bergman's theorem is used to describe the automorphism group $\Aut(\End(\k\langle
x_1,\dots,x_n\rangle))$ \cite{KBLBerz}; Zak's theorem is used in the same way to describe the group
$\Aut(\End(\k[x_1,\dots,x_n]))$ \cite{BelovLiapiansk2}.

\medskip
{\bf Question.}\ {\it Can one prove Bergman's theorem via quantization?}
\medskip

Quantization could be a key idea for understanding Jacobian type problems in  other varieties of algebras.

\paragraph{1. Cancellation problems.}$ $

\medskip

We recall three classical problems.

\medskip

  {\bf 1.} {\it Let $K_1$ and $K_2$ be affine domains for
which $K_1[t]\simeq K_2[t]$. Is it true that $K_1\simeq K_2$? }

  {\bf 2.} {\it Let $K_1$ and $K_2$ be an affine fields
for which $K_1(t)\simeq K_2(t)$. Is it true that $K_1\simeq K_2$? In particular, if $K(t)$  is a field of rational
functions  over the field~$\k$, is it true that $K$ is also a field of rational functions over $\k$? }

  {\bf 3.} {\it If $K[t]\simeq \k [x_1,\dots,x_n]$, is it
true that $K\simeq \k [x_1,\dots,x_{n-1}]$? }
\medskip

The answers to Problems 1 and 2 are `No' in general (even if $\k ={\mathbb C}$); see the fundamental
paper~\cite{BCSS}, as well as \cite{BYML} and the references therein. However,  Problem 2 has a positive
solution in low dimensions. Problem 3 is currently called the {\it Cancellation Conjecture}, although Zariski's original
cancellation conjecture was for fields (Problem~2). See (\cite{MiySugie}, \cite{KZ}, \cite{Dan}, \cite{ShYu1}) for
Zariski's conjecture and related problems. For $n\ge 3$, the Cancellation Conjecture (Problem 3) remains open, to
the best of our knowledge, and it is reasonable to pose the Cancellation Conjecture for free associative rings and
ask the following:

\medskip
\noindent {\bf Question.} {\it If $K*\k [t]\simeq \k \!\!<\!x_1,\dots,x_n\!>$, then is $K\simeq \k
\!<\!x_1,\dots,x_{n-1}\!>? $}
\medskip

This question was solved  for $n=2$  by V.~Drensky and J.T.~Yu \cite{DrYu}.

\paragraph{2. The Tame Automorphism Problem.} Yagzhev utilized his
approach to study the tame automorphism problem. Unfortunately, these  papers are not preserved.

It is easy to see that every endomorphism $\phi$ of a commutative algebra can be lifted to some endomorphism of
the free associative algebra, and hence to some endomorphism of the algebra of generic matrices. However, it is not
clear that any automorphism $\phi$ can be lifted to an automorphism.

We recall that an automorphism  of $\k [x_1,\dots,x_n]$ is  {\em elementary} if it has the form
$$x_1 \mapsto x_1 + f(x_2, \dots, x_n), \quad x_i \mapsto x_i, \quad \forall i \ge 2.$$
 A  {\em tame automorphism} is a product of elementary automorphisms, and a non-tame
automorphism is called  {\em wild}. The ``tame automorphism problem'' asks whether any automorphism is tame.
Jung~\cite{Jung} and van der Kulk \cite{VdK} proved this for $n=2$, (also see \cite{Nie1,Nie2} for free groups,
\cite{Cohn1} for free Lie algebras, and \cite{ML2,Czer2} for free associative algebras), so one takes $n>2.$

 Elementary automorphisms can be
lifted to automorphisms of the free associative algebra; hence every tame automorphism can be so lifted. If an
automorphism $\varphi$ cannot be lifted to an automorphism of the algebra of generic matrices, it cannot be tame.
This give us approach to the tame automorphism problem.

We can lift an automorphism of $\k [x_1,\dots,x_n]$ to an endomorphism of $\k \!\!<x_1,\dots,x_n\!>$ in many
ways. Then replacing $x_1,\dots,x_n$ by $N\times N$ generic matrices induces a polynomial mapping $F_{(N)}:
\k ^{nN^2}\to \k ^{nN^2}$.

For each automorphism $\varphi$, the invertibility of this mapping can be transformed into  compatibility   of some
system of equations. For example, Theorem~10.5 of \cite{Peretz} says that the Nagata automorphism is wild,
provided that a certain system of five equations in $27$ unknowns has no solutions. Whether Peretz' method  can
effectively attack tameness questions remains to be seen. The wildness of the Nagata automorphism was established
by Shestakov and Umirbaev~\cite{SU2}. One important ingredient in the proof is
 {\it degree estimates} of an expression $p(f,g)$ of algebraically independent polynomials $f$ and $g$ in terms of the degrees of $f$ and $g$,
 provided neither leading term is proportional to a power of the other, initiated by Shestakov and Umirbaev~\cite{Shes1}.  An exposition based on their
method is given in Kuroda \cite{Kuroda}.

One of the most important tools is the degree estimation technique, which in the multidimensional case becomes the
analysis of leading terms, and is more complicated. We refer to the deep papers
\cite{BonnetVenerau,Kuroda,Marek}. Several papers of Kishimoto
 contain gaps, but also provide deep insights.

One can also ask the weaker question of ``coordinate tameness:''
 Is the image of $(x,y,z)$ under the Nagata automorphism the image  under some
 (other) tame automorphism? This also fails, by
 \cite{UY}.

An automorphism $\varphi$ is called {\em stably tame} if,  when several new indeterminates~$\{t_i\}$ are adjoined,
the extension of $\varphi$ given by  $\varphi'(t_i)=t_i$ is tame; otherwise it is called {\em stably wild}. Stable
tameness of automorphisms of $\k [x,y,z]$ fixing~$z$ is proved in \cite{BEW}; similar results for $\k \langle
x,y,z\rangle$ are given in \cite{BKY}.

Yagzhev tried to construct wild automorphisms via polynomial automorphisms of the Cayley-Dickson algebra with
base $\{\vec{e}_i\}_{i=1}^8$, and the set $\{\nu_i, \xi_i, \varsigma_i\}_{i=1}^8$ of commuting indeterminates.
Let
 $$x=\sum \nu_i\vec{e}_i,\ y=\sum \xi_i\vec{e}_i,\ z=\sum
 \varsigma_i\vec{e}_i.$$
 Let $(x,y,z)$ denote the associator $(xy)z - x (yz)$ of the elements $x,y,z$,
 and write
 $$(x,y,z)^2=\sum f_i(\vec{\nu},\vec{\xi},\vec{\varsigma})\vec{e}_i.$$

Then the endomorphism $G$ of  the polynomial algebra given by
 $$G:\
 \nu_i\to\nu_i+f_i(\vec{\nu},\vec{\xi},\vec{\varsigma}),\quad
 \xi_i\to\xi_i,\ \varsigma_i\to\varsigma_i,$$
is an automorphism, which likely is stably wild.

In the free associative case, perhaps it is possible to construct an example of an automorphism, the wildness of
which could be proved by considering its Jacobi endomorphism
(Definition~\ref{Jac1}). 
Yagzhev tried to construct
 examples of algebras $R=A\otimes A^{op}$ over which there are invertible
matrices that cannot decompose as products of elementary ones. Yagzhev conjectured that the automorphism
 $$x_1\to x_1+y_1(x_1y_2-y_1x_2),\ x_2\to x_2+(x_1y_2-y_1x_2)y_2,\
y_1\to y_1,\ y_2\to y_2$$
 of the free associative algebra is wild.

 Umirbaev~\cite{UmirbaevAnic}  proved  in characteristic 0 that
the {\it Anick automorphism} $x\to x + y(xy-yz)$, $y\to y$, $z\to z+(zy-yz)y$ is wild, by using metabelian algebras.
The proof uses description of the defining relations of 3-variable automorphism groups~\cite{U,U5,U6}. Drensky
and  Yu~\cite{DYuStrongAnik,DY4} proved in characteristic 0 that the image of $x$ under the Anick
Automorphism is not the image of any tame automorphism.

\medskip
{\bf Stable Tameness Conjecture}. {\it Every automorphism of the   polynomial
 algebra $\k [x_1, \dots, x_n]$, resp.~of the free  associative algebra $\k \!\!<\!x_1 ,\dots, x_n\!>$, is stably tame}.
\medskip

Lifting in the free associative case is related to quantization. It provides some light on the similarities and differences
between the commutative and noncommutative cases. Every tame automorphism of the polynomial ring can be lifted
to an automorphism of the free associative algebra. There was a conjecture that {\it any wild $z$-automorphism of
$\k [x,y,z]$ (i.e., fixing $z$) over an arbitrary field $\k$ cannot be lifted to a $z$-automorphism of
$\k\!\!<\!x,y,z\!>$.} In particular, the Nagata automorphism cannot be so lifted \cite{DrYuLift}. This conjecture
was solved by Belov and J.-T.Yu \cite{BelovYuLifting} over an arbitrary field. However, the general lifting
conjecture is still open. In particular, it is not known whether the Nagata automorphism can be lifted to an
automorphism of the free algebra. (Such a lifting could not fix $z$).

The paper \cite{BelovYuLifting}  describes all the $z$-automorphisms of $\k \!\!<\!x,y,z\!>$ over an arbitrary field
$\k $. Based on that work, Belov and J.-T.Yu \cite{BKY} proved that every $z$-automorphism of $\k
\!\!<\!x,y,z\!>$  is stably tame, for all fields~$\k $. A~similar result in the commutative case is proved by Berson,
van den Essen, and Wright~\cite{BEW}. These are    important first steps towards solving the stable tameness
conjecture in the noncommutative and commutative cases.

The free associative situation is much more rigid than the polynomial case. Degree estimates for the free associative
case are the same for prime
  characteristic   \cite{YuYungChang} as in characteristic~0~\cite{MLY}. The methodology is  different
 from the commutative case,
 for which degree estimates (as well as
examples of wild automorphisms) are not known in prime characteristic.

J.-T.Yu found some evidence of a connection between the lifting conjecture and the Embedding Conjecture of
Abhyankar and Sathaye. Lifting seems to be ``easier''.

\subsection{Reduction to simple algebras}\label{pack1subsec}
This subsection is devoted to finding test algebras.

Any prime algebra $B$   satisfying a system of Capelli identities of order $n+1$ ($n$ minimal such) is said to have
{\it rank} $n$. In this case,  its operator algebra is PI.
The localization of $B$ is a simple algebra of dimension $n$ over its centroid, which is a field. This is the famous {\it
rank theorem} \cite{RazmyslovBook}.

\subsubsection{Packing properties}\label{pack1subsubsec}

\begin{definition}
Let ${\cal M}=\{{\goth M}_i: i \in I\}$ be an arbitrary set of varieties of algebras. We say that $\cal M$ satisfies
the {\em packing property}, if for any $n\in{\mathbb N}$ there exists a prime algebra $A$ of rank $n$ in some
${\goth M}_j$ such that any prime algebra in any~${\goth M}_i$ of rank $n$ can be embedded into some central
extension  $K\otimes A$ of $A$.

$\cal M$ satisfies the {\em  finite packing property}  if, for any finite set of prime algebras $A_j\in{\goth M}_i$,
there exists a prime algebra $A$ in some ${\goth M}_k$ such that each $A_j$ can be embedded into $A$.
\end{definition}

The set of proper subvarieties of associative algebras satisfying a system of Capelli identities of some order $k$
satisfies the packing property (because any simple associative algebra is a matrix algebra over field).

 However,   the  varieties of
alternative algebras satisfying a system of Capelli identities of order $>8$,   or of Jordan algebras satisfying a system
of Capelli identities of order $>27$,    do not even satisfy the finite packing property. Indeed, the matrix algebra of
order $2$ and the Cayley-Dickson algebra cannot be embedded into a common prime alternative algebra. Similarly,
${\mathbb H}_3$ and the Jordan algebra of symmetric matrices cannot be embedded into a common Jordan prime
algebra. (Both of these assertions follow easily by considering their PIs.)

It is not known whether or not  the packing property holds for Engel algebras satisfying a system of Capelli identities;
knowing the answer would enable us to resolve the  $\JC$, as will be seen below.

\begin{theorem}
If the set of varieties of Engel algebras (of arbitrary fixed order) satisfying a system of Capelli identities of some
order satisfies the packing property, then the Jacobian Conjecture has a positive solution.
\end{theorem}

\begin{theorem}\label{ThWeakStsfy}
The set of varieties from the previous theorem satisfies the finite packing property.
\end{theorem}

Most of the remainder of this section is devoted to the proof of these two theorems.

\medskip
{\bf Problem.} {\it Using  the packing property and deformations, give a reasonable analog of the $\JC$ in nonzero
characteristic.} (The naive approach using only  the determinant of the Jacobian does not work.)
\medskip

 \subsubsection{Construction of simple Yagzhev algebras}   

Using the Yagzhev correspondence and composition of elementary automorphisms it is possible to construct a new
algebra of Engel type.

\begin{theorem}
Let $A$ be an algebra of Engel type. Then $A$ can be embedded into a prime algebra of Engel type.
\end{theorem}

\Proof Consider the  mapping $F:V\to V$  (cf.~\eqref{Eq1Fk}) given by
$$F:\
x_i\mapsto x_i+\sum_j\Psi_{ij};\quad i=1,\dots,n$$ (where the
$\Psi_{ij}$ are forms of homogenous degree $j$). Adjoining new indeterminates $\{t_i\}_{i=0}^n$, we put
$F(t_i)=t_i$ for $ i=0,\dots,n$.

Now we  take the transformation
$$
G:\ t_0\mapsto t_0,\quad x_i\mapsto x_i,\quad  \ t_i\mapsto
t_i+t_0x_i^2,\quad\mbox{for}\quad  i=1,\dots,n.
$$
The composite $F\circ G$ has invertible Jacobian (and hence the corresponding algebra has Engel type) and can be
expressed as follows:
$$
F\circ G:\ x_i\mapsto x_i+\sum_j\Psi_{ij},\quad \ t_0\mapsto
t_0,\quad t_i\mapsto t_i+t_0x_i^2\quad \mbox{for}\quad i=1,\dots,n.
$$
It is easy to see that the corresponding algebra $\widehat{A}$ also satisfies the following properties:

\begin{itemize}
\item $\widehat{A}$ contains $A$ as a subalgebra (for $t_0=0$).
\item If $A$ corresponds to a cubic homogenous mapping (and thus is Engel) then $\widehat{A}$ also
    corresponds to a cubic homogenous mapping (and thus is Engel). \item If some of the forms $\Psi_{ij}$ are not
    zero, then $A$ does not have nonzero ideals with product 0, and hence is prime (but its localization need not be
    simple!).
\end{itemize}

Any algebra $A$ with  operators  can be embedded, using the previous construction, to a prime algebra with
nonzero multiplication.   The theorem is proved. \Endproof

Embedding via the previous theorem preserves the cubic homogeneous case, but does not yet  give us an embedding
into a simple algebra of Engel type.

\begin{theorem}\label{3.8}
Any algebra $A$   of Engel type  can be embedded into a simple algebra  of Engel type.
\end{theorem}

\Proof We  start from the following observation:

\begin{lemma}
Suppose $A$ is a finite dimensional algebra, equipped with a base $\vec{e}_1,\dots,\vec{e}_n,\vec{e}_{n+1}$. If
for any $1\le i, j\le n+1$ there exist operators $\omega_{ij}$ in the signature $\Omega(A)$ such that
$\omega_{ij}(\vec{e}_i,\dots,\vec{e}_i,\vec{e}_{n+1})=\vec{e}_j$, with all other values on the base vectors
being zero, then $A$ is simple.
\end{lemma}

This lemma implies:

\begin{lemma}
Let $F$ be a polynomial endomorphism of ${\mathbb C}[x_1,\dots,x_n;t_1,t_2]$, where $$F(x_i)=\sum_j
\Psi_{ij}.$$ For notational convenience we put $x_{n+1}=t_1$ and $x_{n+2}=t_2$. Let $\{k_{ij}\}_{i=1,j}^s$
be a set of natural numbers such that
\begin{itemize}
\item For any $x_i $ there exists $k_{ij}$ such that among all $\Psi_{ij}$ there is {\em exactly one} term of
    degree $k_{ij}$, and it has the form $\Psi_{i,k_{ij}}=t_1x_j^{k_{ij}-1}$.

\item For $t_2$ and any $x_i$ there exists $k_{iq}$ such that among all $\Psi_{ij}$ there is {\em exactly one}
    term of degree $k_{iq}$, and it has the form $\Psi_{n+2,k_{iq}}=t_1x_j^{k_{iq}-1}$.

\item For $t_1$ and any $x_i$ there exists $k_{iq}$ such that among all $\Psi_{ij}$ there is {\em exactly one}
    term of degree $k_{iq}$, and it has the form $\Psi_{n+1,k_{iq}}=t_2x_j^{k_{iq}-1}$.
\end{itemize}
Then the corresponding algebra is simple.

\end{lemma}

\Proof Adjoin the term $t_\ell x_i^{k-1}$   to the $x_i$, for $\ell = 1,2$. Let $e_i$
  be the base vector  corresponding to $x_i$. Take the
  corresponding
$k_{ij}$-ary operator $$\omega: \omega(\vec{e}_i,\dots,\vec{e}_i,\vec{e}_{n+\ell}))=\vec{e}_j,$$ with all other
products   zero. Now we apply the previous lemma.
 \Endproof

\medskip
{\bf Remark.} In order to be flexible with constructions via the Yagzhev correspondence, we are working in the
general, not necessary cubic, case.
\medskip

Now we can conclude the proof of Theorem~\ref{3.8}. Let $F$ be the mapping corresponding to the algebra $A$:

$$F:\
x_i\mapsto x_i+\sum_j\Psi_{ij},\quad i=1,\dots,n,$$ where
$\Psi_{ij}$ are forms of homogeneous degree $j$. Let us adjoin new indeterminates $\{t_1, t_2\}$ and put
$F(t_i)=t_i,$ for $ i=1, 2$.

We choose all $ k_{\alpha,\beta}>\max(\deg(\Psi_{ij}))$ and assume that these numbers are sufficiently large. Then
we consider the mappings
$$
G_{k_{ij}}: x_i\mapsto x_i+x_j^{k_{ij}-1}t_1,\ i\le n;\quad
t_1\mapsto t_1;\quad t_2\mapsto t_2;\quad x_s\mapsto x_s\
\mbox{for}\ s\ne i.
$$
$$
G_{k_{i(n+2)}}: t_2\mapsto x_i^{k_{ij}-1}t_1;\quad t_1\mapsto
t_1;\quad x_s\mapsto x_s\ \mbox{for}\  1\le s\le n.
$$
$$
G_{k_{i(n+1)}}: t_1\mapsto x_i^{k_{ij}-1}t_2;\quad t_2\mapsto
t_2;\quad x_s\mapsto x_s\ \mbox{for}\  1\le s\le n.
$$
These mappings are elementary automorphisms.

Consider the mapping $H=\circ_{k_{ij}}G_{k_{ij}}\circ F$, where the composite is taken in order of ascending
$k_{\alpha\beta}$, and then with $F$.
 If the $k_{\alpha\beta}$ grow quickly enough, then the terms
 obtained in the previous step do not affect the
lowest term obtained at the next step, and this term will be as described in the lemma.  The theorem is proved.
\Endproof

\medskip
{\bf Proof of Theorem \ref{ThWeakStsfy}.} The direct sum of Engel type algebras  is also of Engel type, and by
Theorem~\ref{3.8} can be embedded
 into a simple algebra of Engel type.
\Endproof

\paragraph{The Yagzhev correspondence and algebraic extensions.}$ $

For notational simplicity, we  consider a cubic homogeneous mapping
$$
F:\ x_i\mapsto x_i+\Psi_{3i}(\vec{x}).
$$
We shall construct the Yagzhev correspondence of an algebraic extension.

Consider the equation
$$
t^s=\sum_{p=1}^s\lambda_pt^{s-p},
$$
where the $\lambda_p$ are formal parameters. If $m\ge s$, then for some $\lambda_{pm}$, which  can be
expressed as polynomials in $\{\lambda_p\}_{p=1}^{s-1}$,  we have
$$
t^m= \sum_{p=1}^s\lambda_{pm}t^{s-p}.
$$

Let $A$ be the algebra corresponding to the mapping $F$. Consider
$$A\otimes \k[\lambda_1,\dots,\lambda_s]$$ and its finite algebraic
extension $\hat{A}=A\otimes \k[\lambda_1,\dots,\lambda_s,t]$. Now we take the mapping corresponding (via the
Yagzhev correspondence) to the ground ring $R=\k [\lambda_1,\dots,\lambda_s]$ and algebra $\hat{A}$.

For $m=1,\dots,s-1$, we define new formal indeterminates, denoted as~$T^mx_i$.  Namely, we put $T^0x_i=x_i$
and for $m\ge s$, we identify $ T^mx_i$ with $ \sum_{p=1}^s\lambda_{pm}T^{s-p}x_i $, where
$\{\lambda_p\}_{p=1}^{s-1}$ are formal parameters in the
 centroid of some extension $R\otimes A$.
Now we extend the mapping~$F$, by putting
$$
F(T^mx_i)=T^mx_i+T^{3m}\Psi_{3i}(\vec{x}),\quad m=1,\dots,s-1.
$$
We get a natural mapping corresponding to the algebraic extension.

Now we can take more symbols $T_j$, $ j=1,\dots, s$, and equations
$$
T_j^s=\sum_{p=1}^s\lambda_{pj}T_j^{s-p}
$$
and a new set of indeterminates $x_{ijk}=T_j^kx_i$ for $ j=1,\dots,s$ and $ i=1,\dots, n$. Then we put
$$
x_{ijm}=T_j^mx_i= \sum_{p=1}^s\lambda_{jpm}T_j^{s-p}x_i
$$
and
$$
F(x_{ijm})=x_{ijm}+T_j^{3m}\Psi_{3i}(\vec{x}),\quad m=1,\dots,s-1.
$$
This yields an ``algebraic extension'' of $A$.

\paragraph{Deformations of algebraic extensions.} Let $m=2$. Let us
introduce new indeterminates $y_1,y_2$, put $F(y_i)=y_1, \ i=1,2$, and compose $F$ with the automorphism
$$G: T_1^1x_i\mapsto T_1^1x_i+y_1x_i, \quad T_1^1x_i\mapsto T_2^1x_i+y_1x_i, \quad x_i \mapsto x_i,
\quad i = 1,2,$$ $$ y_1\mapsto y_1+y_2^2y_1,\quad  y_2\mapsto y_2. $$ (Note that the $T_1^1x_i$ and
$T_2^1x_i$ are {\it new} indeterminates and not proportional to $x_i$!) Then compose $G$ with the
automorphism $H: y_2\mapsto y_2+y_1^2$, where $H$  fixes the other indeterminates. Let us call the
corresponding new algebra $\hat{A}$. It is easy to see that $\Var(A)\ne\Var(\hat{A})$.


 Define an  {\it identity of the pair} $(A,B)$, for
$A\subseteq B$ to be a polynomial in two sets of indeterminates $x_i, z_j$ that vanishes whenever the $x_i$ are
evaluated in $A$ and $z_j$ in $B$.) The {\it variety of the pair} $(A,B)$ is the class of  pairs of algebras   satisfying
the identities of $(A,B)$.

 Recall that by  the rank theorem, any prime algebra $A$ of rank $n$
can be embedded into an $n$-dimensional simple algebra $\hat{A}$. We consider the variety of the  pair
$(A,\hat{A})$.

Considerations of deformations yield the following:

\begin{statement}
Suppose for all simple $n$-dimensional pairs there exists a universal pair in which all of them can be embedded.
Then the Jacobian Conjecture has a positive solution.
\end{statement}

We see the relation with

\medskip
{\bf The Razmyslov--Kushkulei theorem}~\cite{RazmyslovBook}: \ {\it Over an algebraically closed field, any two
finite dimensional simple algebras satisfying the same identities are isomorphic.}
\medskip

The difficulty in applying this theorem is that the identities may depend on parameters. Also, the natural generalization
of the Rasmyslov--Kushkulei theorem for a variety and subvariety does not hold: Even if $\Var(B)\subset\Var(A)$,
where $B$ and $A$ are simple finite dimensional algebras over some algebraically closed field,
  $B$ need not be embeddable to $A$.

\chapter*{Conclusion}

The quantization program constitutes a substantial and well designed approach to the Jacobian conjecture, as well as
to various related topics in algebra and algebraic geometry. The recent developments presented in this review have
been instrumental in our investigation of Kontsevich conjecture as well as the establishment of results of independent
interest.

Furthermore, as can be seen from the discussion of the work of A.V. Yagzhev, there are substantial areas of the
theory which require further development and which might, conceivably, hold the insights necessary for the resolution
of the Jacobian conjecture.

While at present the quantization approach does not seem to be adequately developed for a successful attack on the
Jacobian problem to happen (as evidenced by our discussion of Kontsevich conjecture), and the rather substantial
critique of the general quantization and lifting philosophy (due to Orevkov and others) exists, further research and
development of the theory is well advised.


\begin{thebibliography} {99}


\bibitem{Abdelmalek} A.~Abdesselam, { The Jacobian Conjecture as a problem of Perturbative Quantum Field
    theory}, Ann. Henri Poincare {4} (2003), no.~2, 199--215.


\bibitem{AM} S.~Abhyankar and T.~Moh, { Embedding of the line in the plane}, J.~Reine Angew.~Math.~{
    276} (1975), 148--166.

\bibitem{Am} S.A.~Amitsur, {Algebras over infinite fields}, Proc. Amer. Math. Soc. { 7} (1956), 35--48.

\bibitem{am2} S.A.~Amitsur, { A general theory of radicals: III Applications}, Amer.~J.~Math.~{75}
    (1954), 126--136.

\bibitem{ALB} J.~Alev and L.~Le~Bruyn.
\newblock Automorphisms of generic 2 by 2 matrices.
\newblock In { Perspectives in Ring Theory}, pages 69--83. Springer, 1988.


\bibitem{AmLev1} A. S. Amitsur and J. Levitzki.  Minimal identities for algebras. Proc. Amer. Math.
    Soc.  1 (1950), 449-463.

\bibitem{AmLev2} A. S. Amitsur and J. Levitzki.  Remarks on minimal identities for algebras. Proc.
    Amer. Math. Soc. 2 (1951), 320-327.

\bibitem{An} D.~J. Anick.
\newblock Limits of tame automorphisms of $k [x_1,\ldots, x_n]$.
\newblock  Journal of Algebra, 82(2):459--468, 1983.

\bibitem{art1978} V.~A. Artamonov.
\newblock Projective metabelian groups and {Lie} algebras.
\newblock  Izvestiya: Mathematics, 12(2):213--223, 1978.

\bibitem{art1985} V.~A. Artamonov.
\newblock Projective modules over universal enveloping algebras.
\newblock Mathematics of the USSR-Izvestiya, 25(3):429, 1985.

\bibitem{art1991} V.~A. Artamonov.
\newblock Nilpotence, projectivity, decomposability.
\newblock Siberian Mathematical Journal, 32(6):901--909, 1991.

\bibitem{art1998} V.~A. Artamonov.
\newblock The quantum {Serre} problem.
\newblock Russian Math. Surveys, 53(4):3--77, 1998.

\bibitem{art2008} V.~A. Artamonov.
\newblock Quantum polynomials advances in algebra and combinatorics, vol.
  19-34, 2008.

\bibitem{Artin1999} M.~Artin.
\newblock Noncommutative {Rings}, 1999.

\bibitem{AKZ} I.~Arzhantsev, K.~Kuyumzhiyan, and M.~Zaidenberg.
\newblock Infinite transitivity, finite generation, and {Demazure} roots.
\newblock Advances in Mathematics, 351:1--32, 2019.

\bibitem{Asanuma} T.~Asanuma.
\newblock Non-linearizable algebraic $k^*$-actions on affine spaces.
\newblock  Preprint, 1996.

\bibitem{Backelin} E.~Backelin.
\newblock Endomorphisms of quantized {Weyl} algebras.
\newblock  Letters in Mathematical Physics, 97(3):317--338, 2011.

\bibitem{Bass84} H.~Bass.
\newblock A non-triangular action of $G_a$ on $A^3$.
\newblock Journal of Pure and Applied Algebra, 33(1):1--5, 1984.

\bibitem{BCW} H.~Bass, E.~H. Connell, and D.~Wright.
\newblock The {Jacobian} conjecture: reduction of degree and formal expansion
  of the inverse.
\newblock  Bulletin of the American Mathematical Society, 7(2):287--330,
  1982.


\bibitem{Bav8} V.V. Bavula. A question of Rentschler and the Dixmier problem. (English) Zbl.
    0995.16019. Ann. Math. (2) 154, No.3, 683-702 (2001).

\bibitem{Bav1} V.V. Bavula. Generalized Weyl algebras and diskew polynomial rings. arXiv:
    1612.08941.

\bibitem{Bav3} V.V. Bavula. The group of automorphisms of the Lie algebra of derivations of a
    polynomial algebra. J. Alg. Appl. 16 (2017), arXiv: 1304.3836.

\bibitem{Bav2} V.V. Bavula. The groups of automorphisms of the Lie algebras of formally analytic
    vector fields with constant divergence. Comptes Rendus Mathematique, 352 (2). pp. 85-88, arXiv:
    1311.2284.


\bibitem{Bavula3} V.V.~Bavula,  {The inversion formulae for automorphisms of Weyl algebras and  polynomial
    algebras,} J. Pure Appl.~Algebra {210} (2007), 147--159.

\bibitem{Bavula} V.V.~Bavula,  {The inversion formulae for automorphisms of polynomial algebras and rings of
    differential operators in prime characteristic,} J. Pure Appl. Algebra {212} (2008), no.~10, 2320--2337.



\bibitem{Bavula1} V.V.~Bavula,  { An analogue of the conjecture of Dixmier is true for the algebra of polynomial
    integro-differential operators,} J. of Algebra { 372} (2012),  237--250.

\bibitem{Bavula2} V.V.~Bavula,  {Every monomorphism of the Lie algebra of unitriangular polynomial
    derivations is an authomorphism,} C.R. Acad.Sci. Paris, Ser. 1, { 350}, (2012), no 11--12, 553--556.



\bibitem{Bavula1new} V.V.~Bavula,  {\it The ${\rm Jacobian \ Conjecture}_{2n}$ implies the  {\rm Dixmier} {\rm
    Problem}$_n$,} arXiv:math/0512250.

\bibitem{BCSS} A.~Beauville, J.-L.~Colliot-Thelene, J-J Sansuc, and P. Swinnerton-Dyer, {\it Varietes stablement
    rationnelles non rationnelles}, Ann. Math. { 121}, (1985) 283--318.

\bibitem{bayen1978} F.~Bayen, M.~Flato, C.~Fronsdal, A.~Lichnerowicz, and D.~Sternheimer.
\newblock Deformation theory and quantization. {I.} Deformations of symplectic
  structures.
\newblock  Annals of Physics, 111(1):61--110, 1978.

\bibitem{KBTrees} A.~Belov.
\newblock Linear recurrence equations on a tree.
\newblock  Mathematical Notes, 78(5):603--609, 2005.

\bibitem{BelovIAN} A. Belov. Local finite basis property and local representability of varieties of
    associative rings. Izv. Rus. Acad.  Sci. 74 (2010) 3-134. English transl.: Izvestiya:
    Mathematics, 74 (2010) 1-126.



\bibitem{BBRY} A.~Belov, L.~Bokut, L.~Rowen, and J.-T. Yu.
\newblock The {Jacobian} conjecture, together with {Specht} and {Burnside}-type
  problems.
\newblock In Automorphisms in Birational and Affine Geometry, pages
  249--285. Springer, 2014.


\bibitem{BYML} A.~Belov, L.~Makar-Limanov, and J.T.~Yu, {\it On the Generalised Cancellation Conjecture},
    J. of Algebra {\bf 281} (2004), 161--166.


\bibitem{BelovUzyRowen2}
   A.~Belov, L.H.~Rowen, and U.~Vishne,
   {\it Structure of Zariski-closed algebras,}
Trans. Amer. Math. Soc. {\bf 362} (2012), 4695--4734.

\bibitem{BelovYuLifting} A. Belov-Kanel and J.-T. Yu. On the lifting of the Nagata automorphism.
    Selecta Math. (N.S.) 17 (2011), 935-945.

\bibitem{KBLBerz} A. Kanel-Belov, A. Berzins, R. Lipyanski. Automorphisms of the semigroup of
    endomorphisms of free associative algebras. Int. Journ. of Algebra and Comp., Vol. 17,:5/6
    (2007), 923--939, arXiv:math/0512273.

\bibitem{K-BE} A.~Belov-Kanel and A.~Elishev.
\newblock On planar algebraic curves and holonomic {$D$}-modules in positive
  characteristic.
\newblock Journal of Algebra and Its Applications, 15(08):1650155, 2016.

\bibitem{K-BK1} A.~Belov-Kanel and M.~Kontsevich.
\newblock Automorphisms of the {Weyl} algebra.
\newblock Letters in Mathematical Physics, 74(2):181--199, 2005.

\bibitem{K-BK2} A.~Belov-Kanel and M.~Kontsevich.
\newblock The {Jacobian} conjecture is stably equivalent to the {Dixmier}
  conjecture.
\newblock  Moscow Mathematical Journal, 7(2):209--218, 2007.

\bibitem{BelovLiapiansk2} A.~Belov-Kanel and R.~Lipyanski.
\newblock Automorphisms of the endomorphism semigroup of a polynomial algebra.
\newblock  Journal of Algebra, 333(1):40--54, 2011.

\bibitem{BKY} A. Belov-Kanel and J.-T. Yu. Stable tameness of automorphisms of $F\langle
    x,y,z\rangle$ fixing $z$. Selecta Math. (N.S) 18 (2012), 799-802.

\bibitem{berg1969} G.~M. Bergman.
\newblock Centralizers in free associative algebras.
\newblock Transactions of the American Mathematical Society,
  137:327--344, 1969.

\bibitem{berg1978diamond} G.~M. Bergman.
\newblock The diamond lemma for ring theory.
\newblock Advances in Mathematics, 29(2):178--218, 1978.

\bibitem{BEW} J.~Berson, A.~van den Essen, and D.~Wright, {\em Stable tameness of two-dimensional
    polynomial automorphisms over a regular ring}, 2007 (rev.~2010), Advances in Mathematics {\bf 230} (2012),
    2176--2197.

\bibitem{Bir} J.~Birman, {\em An inverse function theorem for free groups}, Proc. Amer. Math. Soc. {\bf 41}
    (1973), 634--638.

\bibitem{BonnetVenerau} P.~Bonnet and S.~V\'en\'ereau. {\it Relations between the leading terms of a polynomial
    automorphism.}, J. of Algebra {\bf 322} (2009),
 no. 2, 579---599. 13 Aug 2008.

\bibitem{Berz} A.~Berzins.
\newblock The group of automorphisms of semigroup of endomorphisms of free
  commutative and free associative algebras.
\newblock arXiv preprint math/0504015, 2005.

\bibitem{BialBir1} A.~Bia\l{}ynicki-Birula.
\newblock Remarks on the action of an algebraic torus on $k^n$, {I}.
\newblock Bull. Acad. Polon. Sci. S{\'e}r. Sci. Math. Astronom. Phys,
  14:177--181, 1966.

\bibitem{BialBir2} A.~Bia\l{}ynicki-Birula.
\newblock Remarks on the action of an algebraic torus on $k^n$, {II}.
\newblock Bull. Acad. Polon. Sci. S{\'e}r. Sci. Math. Astronom. Phys,
  15:123--125, 1967.

\bibitem{BB3} A.~Bia\l{}ynicki-Birula.
\newblock Some theorems on actions of algebraic groups.
\newblock {\em Annals of mathematics}, pages 480--497, 1973.

\bibitem{Bit} T.~Bitoun.
\newblock The $p$-support of a holonomic {$D$}-module is lagrangian, for $p$
  large enough.
\newblock arXiv preprint arXiv:1012.4081, 2010.

\bibitem{Bodnarchuk} Yu. Bodnarchuk. Every regular automorphism of the affine Cremona group is
    inner. J. Pure Appl. Algebra 157 (2001), 115-119.

\bibitem{ShirshovMemory} L.~Bokut, V. Latyshev, I. Shestakov and E.~Zelmanov.
\newblock  Selected works of {A. I. Shirshov}.
\newblock Springer, 2009.

\bibitem{bokut62lie} L.~A. Bokut.
\newblock Embedding {Lie} algebras into algebraically closed {Lie} algebras.
\newblock Algebra i Logika, 1:47--53, 1962.

\bibitem{bokut62alg} L.~A. Bokut.
\newblock Embedding of algebras into algebraically closed algebras.
\newblock In Doklady Akademii Nauk, volume 145(5), pages 963--964.
  Russian Academy of Sciences, 1962.

\bibitem{bokut66} L.~A. Bokut.
\newblock Theorems of embedding in the theory of algebras.
\newblock In Colloq. math, volume~14, pages 349--353, 1966.

\bibitem{BPS1} M.~Bresar, C.~Procesi, and S.~Spenko.
\newblock Functional identities on matrices and the {Cayley-Hamilton}
  polynomial.
\newblock arXiv preprint arXiv:1212.4597, 2013.

\bibitem{Campbell} L.~A. Campbell.
\newblock A condition for a polynomial map to be invertible.
\newblock Mathematische Annalen, 205(3):243--248, 1973.

\bibitem{cohn1964} P.~M. Cohn.
\newblock Subalgebras of free associative algebras.
\newblock  Proceedings of the London Mathematical Society, 3(4):618--632,
  1964.

\bibitem{cohn74progress} P.~M. Cohn.
\newblock Progress in free associative algebras.
\newblock Israel Journal of Mathematics, 19(1-2):109--151, 1974.

\bibitem{cohn92skew} P.~M. Cohn.
\newblock A brief history of infinite-dimensional skew fields.
\newblock Math. Scient, 17:1--14, 1992.

\bibitem{Cohn1} P.M. Cohn. Free Rings and Their Relations. 2nd edition, Academic Press (1985).

\bibitem{Czer1} A.~J. Czerniakiewicz.
\newblock Automorphisms of a free associative algebra of rank 2. I.
\newblock Transactions of the American Mathematical Society,
  160:393--401, 1971.

\bibitem{Czer2} A.~J. Czerniakiewicz.
\newblock Automorphisms of a free associative algebra of rank 2. II.
\newblock Transactions of the American Mathematical Society,
  171:309--315, 1972.

 \bibitem{Dan}
W.~Danielewski, {\it On the cancellation problem and automorphism groups of affine algebraic varieties}, preprint,
Warsaw 1989.

\bibitem{dBvdE2} M.~De~Bondt and A.~Van~den Essen.
\newblock The {Jacobian} conjecture for symmetric Druzkowski mappings.
\newblock University of Nijmegen, Department of Mathematics, 2004.

\bibitem{dBvdE1} M.~De~Bondt and A.~Van~den Essen.
\newblock A reduction of the {Jacobian} conjecture to the symmetric case.
\newblock Proceedings of the American Mathematical Society,
  133(8):2201--2205, 2005.

\bibitem{deConcProc} C.~De~Concini and C.~Procesi.
\newblock A characteristic free approach to invariant theory.
\newblock In Young Tableaux in Combinatorics, Invariant Theory, and
  Algebra, pages 169--193. Elsevier, 1982.


\bibitem{Deserti} J. D\'{e}serti. Sur le groupe des automorphismes polynomiaux du plan affine. J.
    Algebra 297 (2006) 584-599.

\bibitem{Di} W. Dicks, Automorphisms of the free algebra of rank two. Group actions on rings
    (Brunswick, Maine, 1984), Contemp. Math., 43 (1985) 63-68.

\bibitem{DiLev} W. Dicks and J. Lewin. Jacobian conjecture for free associative algebras. Commun.
    Alg. 10, No. 12 (1982), 1285–1306.

\bibitem{Dix} J.~Dixmier.
\newblock {Sur les algebres de Weyl}.
\newblock  Bulletin de la Soci{\'e}t{\'e} math{\'e}matique de France,
  96:209--242, 1968.

\bibitem{Dodd} C.~Dodd.
\newblock The $p$-cycle of holonomic {$D$}-modules and auto-equivalences of the
  {Weyl} algebra.
\newblock arXiv preprint arXiv:1510.05734, 2015.

\bibitem{Donkin} S.~Donkin.
\newblock Invariants of several matrices.
\newblock Inventiones mathematicae, 110(1):389--401, 1992.

\bibitem{Donkin2} S.~Donkin.
\newblock Invariant functions on matrices.
\newblock Mathematical Proceedings of the Cambridge Philosophical
  Society, 113(1):23--43, 1993.

\bibitem{DrYu} V.~Drensky and J.-T. Yu.
\newblock A cancellation conjecture for free associative algebras.
\newblock Proceedings of the American Mathematical Society,
  136(10):3391--3394, 2008.

  \bibitem{DY4}  V.~Drensky and J.-T.~Yu,{\it The strong Anick
conjecture}, Proc. Natl. Acad. Sci. USA {\bf 103} (2006),  4836--4840.


\bibitem{DrYuLift}
 V.~Drensky and J.-T.~Yu, {\it Coordinates and automorphisms
of polynomial and free associative algebras of rank three,} Front. Math. China~{\bf 2 (1)} (2007),  13--46.

  \bibitem {DYuStrongAnik} V. Drensky and J.-T. Yu. The strong Anick conjecture is true. J. Eur.
    Math. Soc. (JEMS) 9 (2007), 659-679.

\bibitem{Druz1} L.~Dru{\.z}kowski.
\newblock An effective approach to {Keller's Jacobian} conjecture.
\newblock Mathematische Annalen, 264(3):303--313, 1983.

\bibitem{Druz2} L.~Dru{\.z}kowski.
\newblock The {Jacobian} conjecture: symmetric reduction and solution in the
  symmetric cubic linear case.
\newblock Annales Polonici Mathematici, 1(87):83--92, 2005.

\bibitem{Dru2} L.M.~Dru\'zkowski,
 {\it   New reduction in the Jacobian conjecture. Effective methods in
 algebraic and analytic geometry, 2000 (Krakow)}, Univ. Iagel. Acta Math. No. 39 (2001), 203--206.

\bibitem{Eli-phd} A.~Elishev.
\newblock Automorphisms of polynomial algebras, quantization and {Kontsevich}
  conjecture.
\newblock Moscow Institute of Physics and Technology, PhD Thesis, 2019.

\bibitem{TA1} A.~Elishev, A.~Kanel-Belov, F.~Razavinia, J.-T. Yu, and W.~Zhang.
\newblock Noncommutative {Bia\l{}ynicki-Birula.} theorem.
\newblock  arXiv preprint arXiv:1808.04903, 2018.

\bibitem{TA2} A.~Elishev, A.~Kanel-Belov, F.~Razavinia, J.-T. Yu, and W.~Zhang.
\newblock Torus actions on free associative algebras, lifting and
  {Bia\l{}ynicki-Birula.} type theorems.
\newblock arXiv preprint arXiv:1901.01385, 2019.

\bibitem{VDEssenImageA}
 A.~Van den Essen, \ {\it The Amazing Image Conjecture},
 arXiv:1006.5801.


\bibitem{VDEssenBondt} A.~Van den Essen and M.~de Bondt, {\it Recent progress on the Jacobian Conjecture},
    Proc. of the Int. Conf. Singularity Theory in honour of S. Lojawiewicz, Cracow, 22-26 March 2004, in Annales
    Polonici Mathematici, {\bf 87} (2005), 1--11.

\bibitem{VDEssenBondt1}
 A.~Van den Essen and M.~de Bondt, {\it  The Jacobian Conjecture for symmetric Dru\'zkowski mappings},
Annales Polonici Mathematici {\bf 86} No.~1 (2005), 43--46.


\bibitem{VDEssenImage}
 A.~Van den Essen, D.~Wright, and W.~Zhao,  {\it On the Image Conjecture}, Journal of Algebra
{\bf 340} (2011),  211---224.

\bibitem{F} R.H.~Fox, {\it Free differential calculus, I. Derivation in the free group ring},  Ann. of Math. (2) {\bf
    57} (1953), 547--560.


\bibitem{Gizatullin} M.H.~Gizatullin, {\it Automorphisms of affine surfaces, $I$, $II$,} Mathematics of the
    USSR-Izvestiya  {\bf 11(1)} (1977), 54--103.

\bibitem{GorniZampieri} G.~Gorni and  G.~Zampieri, {\it Yagzhev polynomial mappings: on the structure of the
    Taylor expansion of their local inverse.} Polon. Math.~{\bf 64} (1996), 285--290.

\bibitem{fedosov1994simple} B.~Fedosov.
\newblock A simple geometrical construction of deformation quantization.
\newblock  Journal of Differential Geometry, 40(2):213--238, 1994.

\bibitem{FMS} T.~Frayne, A.~C. Morel, and D.~S. Scott.
\newblock Reduced direct products.
\newblock Journal of Symbolic Logic, 31(3):506--507, 1966.

\bibitem{FulHar} W. Fulton and J. Harris.
\newblock Representation Theory. A First Course (2nd edition). Springer-Verlag.

\bibitem{FurKr} J.-P. Furter and H.~Kraft.
\newblock On the geometry of the automorphism groups of affine varieties.
\newblock arXiv preprint arXiv:1809.04175, 2018.

\bibitem{Gutwi} A~Gutwirth.
\newblock The action of an algebraic torus on the affine plane.
\newblock Transactions of the American Mathematical Society,
  105(3):407--414, 1962.

\bibitem{Jung} H. W. E.  Jung.
\newblock {{\"U}ber ganze birationale Transformationen der Ebene.}
\newblock { Journal f{\"u}r die reine und angewandte Mathematik},
  184:161--174, 1942.

\bibitem{KMMLR} S.~Kaliman, M.~Koras, L.~Makar-Limanov, and P.~Russell.
\newblock {$C^*$}-actions on {$C^3$} are linearizable.
\newblock { Electron. Res. Announc. Amer. Math. Soc}, 3:63--71, 1997.

\bibitem{KZ} S.~Kaliman and M.~Zaidenberg,  {\it Families of affine planes: the existence of a cylinder},
 Michigan Math.~J.  {\bf 49} (2001),
353--367.

\bibitem{Kuroda} S.~Kuroda, {\it Shestakov-Umirbaev reductions and Nagata's conjecture on a polynomial
    automorphism}, Tohoku Math. J.  {\bf 62} (2010), 75--115.

\bibitem{KuS} E.~Kuzmin and I.P.~Shestakov, Non-associative structures. (English), Algebra VI. Encycl. Math.
    Sci.  {\bf  57}, (1995) 197--280; translation from Itogi Nauki Tekh., Ser. Sovrem. Probl. Mat., Fundam.
    Napravleniya {\bf  57} (1990), 179--266.

\bibitem{Marek} M.~Karas',  {\em Multidegrees of tame automorphisms of $C^n$},  Dissertationes Math.  {\bf
    477} (2011), 55~pp.



\bibitem{PiontkovskiKhoroshkin}
 A.~Khoroshkin and D.~Piontkovski, {\it On generating series of finitely
presented operads}, preprint, 2012, arXiv:1202.5170.

\bibitem{Kam96} T.~Kambayashi.
\newblock Pro-affine algebras, {Ind}-affine groups and the {Jacobian} problem.
\newblock { Journal of Algebra}, 185(2):481--501, 1996.

\bibitem{Kam03} T.~Kambayashi.
\newblock Some basic results on pro-affine algebras and {Ind}-affine schemes.
\newblock { Osaka Journal of Mathematics}, 40(3):621--638, 2003.

\bibitem{KR} T.~Kambayashi and P.~Russell.
\newblock On linearizing algebraic torus actions.
\newblock { Journal of Pure and Applied Algebra}, 23(3):243--250, 1982.

\bibitem{BBL} A.~Kanel-Belov, V.~Borisenko, and V.~Latysev.
\newblock Monomial algebras.
\newblock { Journal of Mathematical Sciences}, 87:3463--3575, 1997.


\bibitem{K-BE2} A.~Kanel-Belov, A.~Elishev, and J.-T. Yu.
\newblock Independence of the {B-KK} isomorphism of infinite prime.
\newblock { arXiv preprint arXiv:1512.06533}, 2015.

\bibitem{K-BE4} A.~Kanel-Belov, A.~Elishev, and J.-T. Yu.
\newblock Augmented polynomial symplectomorphisms and quantization.
\newblock { arXiv preprint arXiv:1812.02859}, 2018.

\bibitem{KGE} A.~Kanel-Belov, S.~Grigoriev, A.~Elishev, J.-T. Yu, and W.~Zhang.
\newblock Lifting of polynomial symplectomorphisms and deformation
  quantization.
\newblock { Communications in Algebra}, 46(9):3926--3938, 2018.

\bibitem{K-BMR1} A. Kanel-Belov, S. Malev and L. Rowen.  The images of non-commutative
    polynomials evaluated on $2\times 2$ matrices. Proc. Amer. Math. Soc. 140 (2012), 465-478.

\bibitem{K-BMR2} A. Kanel-Belov, S. Malev and L. Rowen. The images of multilinear polynomials
    evaluated on $3\times 3$ matrices. Proc. Amer. Math. Soc. 144 (2016), 7-19 .

\bibitem{zhang2017} A.~Kanel-Belov, F.~Razavinia, and W.~Zhang.
\newblock Bergman's centralizer theorem and quantization.
\newblock { Communications in Algebra}, 46(5):2123--2129, 2018.

\bibitem{zhang2018berg} A.~Kanel-Belov, F.~Razavinia, and W.~Zhang.
\newblock Centralizers in free associative algebras and generic matrices.
\newblock  arXiv preprint arXiv:1812.03307, 2018.

\bibitem{BelovUzyRowenSerdicaBachtur}
 A. Kanel-Belov, L. H. Rowen  and U. Vishne.  Full
exposition of Specht's problem. Serdica Math. J. 38 (2012) 313-370.

\bibitem{KBYu} A.~Kanel-Belov, J.-T. Yu, and A.~Elishev.
\newblock On the augmentation topology of automorphism groups of affine spaces
  and algebras.
\newblock  International Journal of Algebra and Computation,
  28(08):1449--1485, 2018.

\bibitem{keller2003} B.~Keller.
\newblock Notes for an introduction to {Kontsevich's} quantization theorem,
  2003.

\bibitem{keller1939} O.-H. Keller.
\newblock {Ganze Cremona-Transformationen}.
\newblock  Monatshefte f{\"u}r Mathematik, 47(1):299--306, 1939.

\bibitem{kole2000} P.~S. Kolesnikov.
\newblock The {Makar-Limanov} algebraically closed skew field.
\newblock  Algebra and Logic, 39(6):378--395, 2000.

\bibitem{kole2001} P.~S. Kolesnikov.
\newblock Different definitions of algebraically closed skew fields.
\newblock Algebra and Logic, 40(4):219--230, 2001.

\bibitem{kontsevich2003} M.~Kontsevich.
\newblock Deformation quantization of {Poisson} manifolds.
\newblock  Letters in Mathematical Physics, 66(3):157--216, 2003.

\bibitem{Kon} M.~Kontsevich.
\newblock Holonomic {$D$}-modules and positive characteristic.
\newblock Japanese Journal of Mathematics, 4(1):1--25, 2009.

\bibitem{KoRu2} M.~Koras and P.~Russell.
\newblock {$C^*$-actions on $C^3$}: The smooth locus of the quotient is not of
  hyperbolic type.
\newblock Journal of Algebraic Geometry, 8(4):603--694, 1999.

\bibitem{KPZ} S.~Kovalenko, A.~Perepechko, M.~Zaidenberg, et~al.
\newblock On automorphism groups of affine surfaces.
\newblock In  {\em Algebraic Varieties and Automorphism Groups}, pages 207--286.
  Mathematical Society of Japan, 2017.

\bibitem{KrReg} H.~Kraft and A.~Regeta.
\newblock Automorphisms of the {Lie} algebra of vector fields.
\newblock  J. Eur. Math. Soc.(to appear), 2015.

\bibitem{KraftStampfli} H.~Kraft and I.~Stampfli.
\newblock On automorphisms of the affine {Cremona} group.
\newblock In  Annales de l'Institut Fourier, volume 63(3), pages
  1137--1148, 2013.

\bibitem{Kulik92} V.~S. Kulikov.
\newblock Generalized and local {Jacobian} problems.
\newblock  Russian Academy of Sciences. Izvestiya Mathematics, 41(2):351,
  1993.

\bibitem{Kulik01} V.~S. Kulikov.
\newblock The {Jacobian} conjecture and nilpotent maps.
\newblock  Journal of Mathematical Sciences, 106(5):3312--3319, 2001.

\bibitem{LLS} R.~Levy, P.~Loustaunau, and J.~Shapiro.
\newblock The prime spectrum of an infinite product of copies of {$Z$}.
\newblock  Fundamenta Mathematicae, 138:155--164, 1991.

\bibitem{YuYungChang} Y.-C. Li and J.-T. Yu. Degree estimate for subalgebras. J. Algebra 362
    (2012), 92-98.

\bibitem{lothaire1997} M.~Lothaire.
\newblock  Combinatorics on words, volume~17.
\newblock Cambridge university press, 1997.

\bibitem{M} L.~Makar-Limanov  {\it  A new proof of the Abhyankar-Moh-Suzuki
 Theorem}, arXiv:1212.0163,
    18 pages.

\bibitem{ML2} L.~Makar-Limanov.
\newblock Automorphisms of a free algebra with two generators.
\newblock  Functional Analysis and Its Applications, 4(3):262--264, 1970.

\bibitem{ML1} L.~Makar-Limanov.
\newblock On automorphisms of {Weyl} algebra.
\newblock  Bulletin de la Soci{\'e}t{\'e} Math{\'e}matique de France,
  112:359--363, 1984.

\bibitem{MLY} L. Makar-Limanov and J.-T. Yu. Degree estimate for subalgebras generated by two
    elements, J. Eur. Math. Soc. (JEMS) 10 (2008), 533-541.

\bibitem{makar85skew} L.~Makar-Limanov.
\newblock Algebraically closed skew fields.
\newblock  Journal of Algebra, 93(1):117--135, 1985.

\bibitem{ML3} L.~Makar-Limanov, U.~Turusbekova, and U.~Umirbaev.
\newblock Automorphisms and derivations of free {Poisson} algebras in two
  variables.
\newblock  Journal of Algebra, 322(9):3318--3330, 2009.

\bibitem{MSS} M.~Markl, S.~Shnider, J.~Stasheff, {\it Operads in algebra, topology and physics}, Mathematical
    Surveys and Monographs {\bf 96} (2002), AMS, Providence, RI.

\bibitem{MiySugie} M.~Miyanishi and  T.~Sugie, {\it   Affine surfaces containing cylinderlike open sets},  J. Math.
    Kyoto Univ. {\bf  20} (1980), 11--42.

\bibitem{Nag} M.~Nagata, {\it   On the automorphism group of $k[x,y]$}, Department of Mathematics, Kyoto
    University, Lectures in Mathematics, No. 5. Kinokuniya Book-Store Co., Ltd., Tokyo, 1972. v+53 pp.

\bibitem{Nie1} J.~Nielsen, {\it   Die Isomorphismen der allgemeinen, undendlichen Gruppen mit zwei
    Eerzeugenden}, Math. Ann.~{\bf  78} (1918), 385--397.

\bibitem{Nie2} J.~Nielsen, {\it   Die Isomorphismengruppe der freien Gruppen}, Math. Ann.~{\bf  91} (1924),
    169--209.

\bibitem{Ol} A.Yu.~Ol'shanskij, {\it Groups of bounded period with subgroups of prime order}, Algebra i Logika
    {\bf 21} (1982), 553--618, translation in Algebra and Logic {\bf  21} (1983), 369--418.

\bibitem{Peretz} R.~Peretz, {\it Constructing polynomial mappings using non-commutative algebras.} Aff. Algebr.
    geometry, 197--232, Cont.~Math.~{\bf  369} (2005), Amer.~Math.~Soc., Providence, RI.


\bibitem{Piontkovski} D.~Piontkovski, {\it Operads versus Varieties: a dictionary of universal algebra}, preprint,
    2011.

\bibitem{Piontkovski1} D.~Piontkovski, {\it On Kurosh problem in varieties of algebras}, Translated from
    Proceedings of Kurosh conference (Fund. Prikl. Matematika {\bf 14} (2008), 5, 171--184), Journal of
    Mathematical Sciences {\bf 163} No. 6  (2009), 743--750.



\bibitem{RazmyslovIANCapely} Yu.P.~Razmyslov,  {\it Algebras satisfying identity relations of Capelli type.}
    (Russian) Izv. Akad. Nauk SSSR Ser. Mat.  {\bf 45} (1981),  143--166, 240.

\bibitem{RazmyslovBook} Yu.P.~Razmyslov,  {\it Identities of algebras and their representations.} {\sl
    Sovremennaya Algebra}. [Modern Algebra] ``Nauka'', Moscow (1989), 432 pp. Translations of Mathematical
    Monographs  {\bf 138}, American Mathematical Society, Providence, RI, 1994, xiv+318 pp.

\bibitem{RazmZubrilin} Yu.P.~Razmyslov and K.A.~Zubrilin,  {\it Nilpotency of obstacles for the representability
    of algebras that satisfy Capelli identities, and representations of finite type.} (Russian) Uspekhi Mat. Nauk  {\bf
    48} (1993),   171--172; translation in Russian Math. Surveys  {\bf }48 (1993),  183--184.

\bibitem{Reu} C.~Reutenauer, {\it Applications of a noncommutative Jacobian matrix}, J.~Pure Appl.~Algebra
    {\bf 77} (1992), 634--638.


\bibitem{Row}  L.H.~Rowen, {\it {Graduate algebra: Noncommutative view}}, AMS Graduate Studies in
    Mathematics~\textbf{91} (2008).

\bibitem{Moh1} T.-T. Moh.
\newblock On the global {Jacobian} conjecture for polynomials of degree less
  than 100.
\newblock {\em preprint}, 1983.

\bibitem{Moh2} T.-T. Moh.
\newblock On the {Jacobian} conjecture and the configurations of roots.
\newblock  Journal f{\"u}r Mathematik. Band, 340:19, 1983.

\bibitem{moyal1949} J.~E. Moyal.
\newblock Quantum mechanics as a statistical theory.
\newblock In  Mathematical Proceedings of the Cambridge Philosophical
  Society, volume 45(1), pages 99--124. Cambridge University Press, 1949.

\bibitem{Orevkov2} S.~Y. Orevkov.
\newblock The commutant of the fundamental group of the complement of a plane
  algebraic curve.
\newblock  Russian Mathematical Surveys, 45(1):221, 1990.

\bibitem{Orevkov} S.~Y. Orevkov.
\newblock An example in connection with the {Jacobian} conjecture.
\newblock  Mathematical Notes, 47(1):82--88, 1990.

\bibitem{Orevkov3} S.~Y. Orevkov.
\newblock The fundamental group of the complement of a plane algebraic curve.
\newblock Mathematics of the USSR-Sbornik, 65(1):267, 1990.

\bibitem{BIP2} B~Plotkin.
\newblock Varieties of algebras and algebraic varieties.
\newblock  Israel Journal of Mathematics, 96(2):511--522, 1996.

\bibitem{BIP1} B.~Plotkin.
\newblock Algebras with the same (algebraic) geometry.
\newblock  arXiv preprint math/0210194, 2002.

\bibitem{Popov} V.~L. Popov.
\newblock Around the {Abhyankar-Sathaye} conjecture.
\newblock  arXiv preprint arXiv:1409.6330, 2014.

\bibitem{Procesi2} C.~Procesi.
\newblock  Rings with polynomial identities, volume~17.
\newblock M. Dekker, 1973.

\bibitem{Procesi} C.~Procesi.
\newblock The invariant theory of $n\times n$ matrices, 1976.

\bibitem{Razar} M.~Razar.
\newblock Polynomial maps with constant {Jacobian}.
\newblock  Israel Journal of Mathematics, 32(2-3):97--106, 1979.

\bibitem{robinson2016non} A.~Robinson.
\newblock  Non-standard analysis.
\newblock Princeton University Press, 2016.

\bibitem{rosset1976} S.~Rosset.
\newblock A new proof of the {Amitsur-Levitzki} identity.
\newblock  Israel Journal of Mathematics, 23(2):187--188, 1976.

\bibitem{rowen2008} L.~H. Rowen.
\newblock  Graduate Algebra: Noncommutative View, volume~9.
\newblock American Mathematical Society, Providence, RI, 2008.

\bibitem{Schof} A.H. Schofield.
\newblock Representations of rings over skew fields.
\newblock London Math. Soc. Lecture
    Note Series 92, Cambridge, Cambridge University Press (1985).

\bibitem{Sch} G.~Schwarz.
\newblock Exotic algebraic group actions.
\newblock  CR Acad. Sci. Paris, 309:89--94, 1989.

\bibitem{Shafarevich} I.~R. Shafarevich.
\newblock On some infinite-dimensional groups, II.
\newblock  Izvestiya Rossiiskoi Akademii Nauk. Seriya Matematicheskaya,
  45(1):214--226, 1981.

\bibitem{sha2013} Y.~Sharifi.
\newblock Centralizers in Associative Algebras.
\newblock PhD thesis, Science: Department of Mathematics, 2013.

\bibitem{ShestakovWedderburn} I.P.~Shestakov. Finite-dimensional algebras with a nil basis. ( in Russian)
    Algebra i Logika  {\bf 10} (1971),  87--99.

\bibitem{Shestakov1} I.P.~Shestakov. A quantization of Poisson superalgebras and a speciality of Jordan Poisson
    superalgebras. Algebra i Logika, 32, N 5 (1993), 572-585; English transl.: Algebra and Logic, 32, N 5 (1993), 309-317.

\bibitem{Shestakov2} I.P.~Shestakov. A quantization of Poisson algebras and a weak speciality of related Jordan
    superalgebras. Doklady RAN, 334, N 1 (1994), 29-31; English transl.: Russian Acad. Sci. Dokl. Math. 49 (1994), N 1, 34-37.

\bibitem{Shestakov3} I.P.~Shestakov. The speciality problem for Malcev algebras and deformations of Malcev
    Poisson algebras, in "Non-Associative Algebra and Its Applications". edited by R.Costa, A.Grishkov, H.Guzzo, and L.Peresi, Proceedings of the IV International Conference on Non-Associative Algebra and Its Applications, July 1998, S~ao Paulo, 365-371, Marcel Dekker, NY, 2000.

\bibitem{Shestakov4} I.P.~Shestakov. Speciality and Deformations of Algebras. in "Algebra: Proceedings of the
    International Algebraic Conference on Occasion of the 90th Birthday of A.G.Kurosh, Moscow, Russia, May 25-30, 1998" / Ed. Yu.Bahturin - Berlin; New York: de Gruyter, 2000; p. 345-356.




\bibitem {SU1} I. P. Shestakov and U. U. Umirbaev. Degree estimate and two-generated subalgebras of
    rings of polynomials. J. Amer. Math. Soc. 17 (2004), 181-196.

\bibitem{Shes2} I.~Shestakov and U.~Umirbaev.
\newblock The {Nagata} automorphism is wild.
\newblock  Proceedings of the National Academy of Sciences,
  100(22):12561--12563, 2003.

\bibitem{Shes1} I.~Shestakov and U.~Umirbaev.
\newblock {Poisson} brackets and two-generated subalgebras of rings of
  polynomials.
\newblock Journal of the American Mathematical Society, 17(1):181--196,
  2004.

\bibitem {SU2} I. P. Shestakov and U. U. Umirbaev. The tame and the wild automorphisms of
    polynomial rings in three variables. J. Amer. Math. Soc.  17 (2004), 197-220.

    \bibitem{Shp}
V.~Shpilrain, {\it On generators of $L/R^2$ Lie algebras},
 Proc. Amer. Math. Soc. {\bf 119} (1993), 1039--1043.





\bibitem{Singer} D.~Singer,  {\it On Catalan trees and the Jacobian conjecture,}
 Electron. J. Combin. {\bf 8} (2001), no.~1, Research Paper 2, 35 pp.~(electronic).

\bibitem{ShYu1} V.~Shpilrain and J.-T.~Yu, {\it Affine varieties with equivalent cylinders,} J. Algebra {\bf 251}
    (2002), no.~1, 295--307.

\bibitem{ShpYuFactor} V.~Shpilrain and J.-T.~Yu, {\it Factor algebras of free algebras: on a problem of
    G.Bergman,\/} Bull. London Math. Soc. {\bf 35} (2003), 706--710.

\bibitem{Su} M.~Suzuki, {\it  Propi\'et\'es topologiques des polyn\^omes de deux variables complexes, et
    automorphismes alg\'ebraique de l'espace $ C^2$}, J.~Math.~Soc.~Japan, {\bf 26} (1974), 241--257.



\bibitem{Tsu1} Y.~Tsuchimoto.
\newblock Preliminaries on {Dixmier} conjecture.
\newblock  Mem. Fac. Sci. Kochi Univ. Ser. A Math, 24:43--59, 2003.

\bibitem{Tsu2} Y.~Tsuchimoto.
\newblock Endomorphisms of {Weyl} algebra and $ p $-curvatures.
\newblock  Osaka Journal of Mathematics, 42(2):435--452, 2005.

\bibitem{Tsu3} Y. Tsuchimoto. Auslander regularity of norm based extensions of Weyl algebra.
    arXiv:1402.7153.

\bibitem{umirbaev1995ext} U.~Umirbaev.
\newblock On the extension of automorphisms of polynomial rings.
\newblock  Siberian Mathematical Journal, 36(4):787--791, 1995.

\bibitem{Umi} U.U.~Umirbaev, {\it On Jacobian matrices of Lie algebras}, 6th All-Union Conference on Varieties
    of Algebraic Systems, Magnitogorsk (1990), 32--33.

\bibitem{Umirbaev} U.U.~Umirbaev, {\it Shreer varieties of algebras}.(Russian). Algebra i Logika {\bf 33}(1994),
    no.~3,  317--340, 343; translation in Algebra and Logic {\bf 33} (1994),   180--193.

\bibitem{UmirbaevAnic} U.U.~Umirbaev, {\it Tame and wild automorphisms of polynomial algebras and free
    associative algebras},  Max-Planck-Institute f\"ur Mathematics, Bonn, Preprint MPIM 2004--108.

\bibitem{U} U.~Umirbaev.
\newblock The {Anick} automorphism of free associative algebras.
\newblock  Journal f{\"u}r die reine und angewandte Mathematik (Crelles
  Journal), 2007(605):165--178, 2007.

  \bibitem{U6} U.U.~Umirbaev, {\it Defining relations of the tame automorphism group of polynomial algebras in three variables}, J. Reine Angew. Math. {\bf 600} (2006),
203--235.

\bibitem{U5} U.U.~Umirbaev, {\it Defining relations for automorphism groups of free algebras}, J.~Algebra {\bf
    314} (2007), 209--225.

\bibitem{UY}  U. U. Umirbaev and J.-T. Yu. The strong Nagata conjecture. Proc. Natl. Acad. Sci. USA
 101 (2004), 4352-4355.

\bibitem{UrZ} C.~Urech and S.~Zimmermann.
\newblock Continuous automorphisms of {Cremona} groups.
\newblock  arXiv preprint arXiv:1909.11050, 2019.

\bibitem{vdE} A.~Van~den Essen.
\newblock Polynomial Automorphisms and the Jacobian Conjecture, volume
  190.
\newblock Birkh{\"a}user, 2012.

\bibitem{VdK} W.~Van~der Kulk.
\newblock On polynomial rings in two variables.
\newblock  Nieuw Arch. Wisk.(3), 1:33--41, 1953.

\bibitem{Vitushkin1} A.G.~Vitushkin, {\it A criterion for the representability of a chain of $\sigma$-processes by a
    composition of triangular chains} (Russian), Mat. Zametki {\bf 65} no.~5  (1999),  643--653; transl. in Math.
    Notes {\bf 65}  no.~5--6 (1999), 539--547.

\bibitem{Vit2} A.~G. Vitushkin.
\newblock On the homology of a ramified covering over {$C^2$}.
\newblock Mathematical Notes-New York, 64(5):726--731, 1998.

\bibitem{Vit1} A.~G. Vitushkin.
\newblock Evaluation of the {Jacobian} of a rational transformation of {$C^2$}
  and some applications.
\newblock Matematicheskie Zametki, 66(2):308--312, 1999.

\bibitem{Wedderburn} J.H.M.~Wedderburn,  {\it Note on algebras}, Annals of Mathematics {\bf 38}  (1937),
    854--856.

\bibitem{Wright_B} D.~Wright, {\it The Jacobian Conjecture as a problem in combinatorics},  in the monograph
    Affine Algebraic Geometry, in honor of Masayoshi Miyanishi,edited by Takayuki Hibi, published by Osaka
    University Press 2007,
 ArXiv: math.Co/0511214 v2, 22 Mar 2006.

\bibitem{Wright1} D.~Wright, {\it The Jacobian Conjecture: ideal membership questions and recent advances},
 Affine algebraic geometry, 261--276, Contemp. Math. {\bf 369} (2005).


\bibitem{yagzhevKoethe} A.V.~Yagzhev,  {\it On the Koethe problem} (Russian), Unpublished.

\bibitem{Yagzhev4} A.V.~Yagzhev, {\it Finiteness of the set of conservative polynomials of a given degree}
    (Russian), Mat. Zametki {\bf 41} (1987), no.~2, 148--151, 285.

\bibitem{Yagzhev5} A.V.~Yagzhev, {\it Nilpotency of extensions of an abelian group by an abelian group}
    (Russian), Mat. Zametki {\bf 43} (1988), no.~3, 424--427, 431; translation in Math. Notes {\bf 43} (1988),
    no.~3--4, 244--245.

\bibitem{Yagzhev6} A.V.~Yagzhev, {\it Locally nilpotent subgroups of the holomorph of an abelian group.}
    (Russian) Mat. Zametki {\bf 46} (1989), no.~6, 118.

\bibitem{Yagzhev7} A.V.~Yagzhev, {\it A sufficient condition for the algebraicity of an automorphism of a group.}
    (Russian) Algebra i Logika {\bf 28} (1989), no.~1, 117--119, 124; translation in Algebra and Logic {\bf 28}
    (1989), no.~1, 83--85.


\bibitem{Yagzhev} A.V.~Yagzhev, {\it The generators of the group of tame automorphisms of an algebra of
    polynomials} (Russian), Sibirsk. Mat. \v Z. {\bf 18} (1977), no.~1, 222--225, 240.

\bibitem{wang1980jacobian} S.~Wang.
\newblock A {Jacobian} criterion for separability.
\newblock  Journal of Algebra, 65(2):453--494, 1980.

\bibitem{Wright} D.~Wright.
\newblock On the {Jacobian} conjecture.
\newblock  Illinois Journal of Mathematics, 25(3):423--440, 1981.


\bibitem{Yag1} A.V. Yagzhev.
\newblock Invertibility of endomorphisms of free associative algebras (in Russian).
\newblock   Mat. Zametki 49 (1991), No. 4, 142–147, 160; translation in Math. Notes 49 (1991), No. 3–4,
    426–430.

\bibitem{Yag2} A.V. Yagzhev.
\newblock On endomorphisms of free algebras (in Russian).
\newblock Sibirsk. Mat. Zh. 21 (1980), No. 1, 181–192.

\bibitem{Yag3} A.V. Yagzhev.
\newblock Algorithmic problem of recognizing automorphisms among endomorphisms
    of free associative algebras of finite rank.
\newblock   Sib. Math. J. 21 (1980), 142--146.

\bibitem{Yag4} A.V. Yagzhev.
\newblock Keller's Problem.
\newblock Sib. Math. J., Vol. 21 (1980), 747--754.

\bibitem{Yagzhev10}
 A.V.~Yagzhev, {\it Engel algebras satisfying Capelli identities.}
(Russian) Proceedings of Shafarevich Seminar, Moscow, 2000; pages 83--88.

\bibitem{Yagzhev9}
 A.V.~Yagzhev, {\it Endomorphisms of polynomial rings and free
algebras of different varieties.} (Russian) Proceedings of Shafarevich Seminar, Moscow, 2000. pages 15--47.

\bibitem{YagzhevLast} A.V.~Yagzhev, {\it Invertibility criteria of a polynomial mapping.} (Russian, unpublished).




\bibitem{Zaks} A.~Zaks, {\it Dedekind subrings of $K[x_1,\dots,x_n]$ are rings of polynomials.} Israel Journal of
    math {\bf 9}, (1971), 285--289.

\bibitem{Zelmanov} E.~Zelmanov, {\it On the nilpotence of nilalgebras}, Lect. Notes Math.~{\bf 1352}  (1988),
    227--240.


\bibitem{VDEssenWenhua} W.~Zhao, {\it New Proofs for the Abhyankar-Gurjar Inversion Formula and the
    Equivalence of the Jacobian Conjecture and the Vanishing Conjecture},
   Proc. Amer. Math. Soc. {\bf 139} (2011), 3141--3154.

\bibitem{VDEssenWenhua1} W.~Zhao, {\it Mathieu Subspaces of Associative Algebras}, Journal of Algebra {\bf
    350} (2012),  245--272,
   arXiv:1005.4260

\bibitem{ZhevlakovKo} K.A.~Zhevlakov, A.M.~Slin'ko, I.P.~Shestakov, and A.I.~Shirshov, Nearly Associative
    Rings [in Russian], Nauka, Moscow (1978).


\bibitem{Zubrilin1} K.A.~Zubrilin,  {\it Algebras that satisfy the Capelli identities} (Russian), Mat. Sb. 186 (1995),
    no.~3, 53--64; translation in Sb. Math. {\bf 186} (1995), no.~3, 359--370.

\bibitem{Zubrilin2} K.A.~Zubrilin, {\it On the class of nilpotence of obstruction for the representability of algebras
    satisfying Capelli identities} (Russian), Fundam. Prikl. Mat. {\bf 1} (1995), no.~2, 409--430.



\bibitem{Zubrilin4} K.A.~Zubrilin, {\it On the Baer ideal in algebras that satisfy the Capelli identities} (Russian),
    Mat. Sb. {\bf 189} (1998), 73--82; translation in Sb. Math. {\bf 189} (1998),   1809--1818.

\bibitem{Za96} M.~G. Zaidenberg.
\newblock On exotic algebraic structures on affine spaces.
\newblock In {\em Geometric complex analysis}, pages 691--714. World
  Scientific, 1996.

\bibitem{Zhang-mas} W.~Zhang.
\newblock Alternative proof of {B}ergman's centralizer theorem by quantization.
\newblock Bar-Ilan University, Master Thesis, 2017.

\bibitem{Zhang-phd}  W.~Zhang.
\newblock Polynomial Automorphisms and Deformation Quantization.
\newblock Bar-Ilan University, PhDThesis, 2019.

\bibitem{Zubkov2} A.~N. Zubkov.
\newblock Matrix invariants over an infinite field of finite characteristic.
\newblock Siberian Mathematical Journal, 34(6):1059--1065, 1993.

\bibitem{Zubkov} A.~N. Zubkov.
\newblock A generalization of the {Razmyslov-Procesi} theorem.
\newblock Algebra and Logic, 35(4):241--254, 1996.



\end{thebibliography}
\end{document}